\documentclass{amsart}
\usepackage{amssymb, amsmath}
\usepackage{bbm}
\usepackage{calrsfs}
\usepackage{mathrsfs}

\def\largerightarrow{   -\negthinspace\negthinspace -\negthinspace\negthinspace
                                       \negthinspace\longrightarrow }

\begin{document}
\title [The Jordan Lattice completion ] {The Jordan lattice completion and a note on injective envelopes and von Neumann algebras}
\vbox{\hfil {\Large\bf  }\hfil}
\author{U. Haag}
\date{\today \\ \texttt{\hfil Contact:haag@zedat.fu-berlin.de}}
\maketitle
\par\noindent
\qquad\qquad\qquad\qquad\qquad\qquad $\> $ {\it In memory of my parents}
\par\noindent
\qquad\qquad\qquad\qquad\qquad\qquad\quad $\>\> $ {\it Kaethe and Rudolf}
\par\bigskip\noindent
\begin{abstract} The article exhibits, among other things, certain relations between the injective envelope $\, I ( A )\, $ of a 
$C^*$-algebra $\, A\, $ and the von Neumann algebra generated by a representation $\,\lambda\, $ of $\, A\, $ provided it is injective. In the first section two canonical lattice constructions associated with any function system $\, \mathfrak X\, $ (= underlying real ordered Banach space of an operator system) are defined and their general properties and interrelations are investigated. The first corresponds to the 
(positive) injective envelope of $\, \mathfrak X\, $ which is denoted 
$\, {\mathfrak L}_1 ( \mathfrak X )\, $ and is a monotone complete regular unital Banach lattice (= injective abelian $C^*$-algebra). It is the unique minimal complete Banach lattice containing $\, \mathfrak X\, $. The second is a somewhat maximal reasonable enveloping lattice construction which is called the {\it Jordan lattice associated with $\,\mathfrak X\, $} denoted 
$\, \mathfrak L ( \mathfrak X )\, $. Both constructions are used to prove the main theorem of the second section relating the injective envelope of a separable $C^*$-algebra with its enveloping von Neumann algebra in a given faithful $*$-representation and exhibiting the existence of weak injective envelopes for any separable $C^*$-algebra (for this notion see section 2). The last (third) section deals with lattices of projections of specific (noncommutative) $C^*$-algebras like von Neumann algebras and injective $C^*$-algebras or more generally (sequentially) monotone complete $C^*$-algebras. The notion of $\mathcal P$-map is defined and certain properties of these maps are
exhibited. Also the related notion of $\mu $-map (resp. $\mu $-functional) is established which includes and transcends the notion of a positive linear map (functional) as well as the notion of $\mathcal P$-measure which sort of generalizes the usual notion of measure (resp. $\mu $-measure) to a noncommutative context.
The last theorem of section 3 shows that any $C^*$-quotient of a sequentially monotone complete $C^*$-algebra is again sequentially monotone complete if the quotient admits an order isomorphic representation by bounded operators on a separable Hilbert space. \end{abstract}
\par\bigskip\bigskip\bigskip\noindent
{\bf 0.\quad Introduction.}
\par\bigskip\noindent
The following text is fairly selfcontained assuming only some general knowledge in the fields of 
operator theory, $C^*$-algebras and (real ordered) Banach spaces, so it should be accessible also to the interested reader who is not an outright expert in these subjects. Specific results which are used without proof may be found in general textbooks available as cited in the references. The author believes that many of the results and methods presented here are completely new opening perspectives for many different future applications beyond those given in the paper. This accounts in particular to the use of enveloping lattice constructions and consideration of general monotonous maps (as opposed to positive linear maps) which provide a powerful tool for many problems in a lattice environment. In a more limited sense this might also be the case for the (only partially monotonous) $\mathcal P$-maps and $\mathcal P$-measures described in section 3 although any computations in a noncommutative projection lattice are much more delicate than in the setting of a (linear) function lattice. The reader should be warned that the logical dependence of the presented results is not always in chronological order, but that the proofs of certain assertions might depend on some results which are proved afterwards so that the reader is requested to skip to the corresponding passages when needed. This is done in favour of grouping together statements which are closely related in content rather than obeying an accurate  choronlogy of logical dependence. We now proceed to give a brief account of the contents of this paper.
\par\noindent
In the first section the two fundamental lattice constructions are presented. This is done in the context of real (dually) ordered Banach spaces (see below for the precise definition). It is shown that such a space $\, \mathfrak V\, $ admits an order isomorphic embedding into a $C^*$-algebra (of continuous functions on a compact space) so that if $\,\mathfrak V\, $ is unital it is order isomorphic with a function system (in the sense of \cite{Ka1}). In the latter case the prototype for such an embedding is given by the lattice construction $\, {\mathfrak L}_1 ( \mathfrak V )\, $ which turns out to be a monotone complete abelian $C^*$-algebra and is equal to the (positive) injective envelope of $\, \mathfrak V\, $, i.e. for any positive linear embedding of unital real ordered Banach spaces  $\, \mathfrak W \subseteq \mathfrak Z\, $ there exists a positive linear map $\, \mathfrak Z \longrightarrow {\mathfrak L}_1 ( \mathfrak V )\, $ extending a given positive linear map $\, \mathfrak W \hookrightarrow {\mathfrak L}_1 ( \mathfrak V )\, $. Moreover $\, {\mathfrak L}_1 ( \mathfrak V )\, $ posseses the following {\it rigidity} property: any monotonous extension $\, \gamma : {\mathfrak L}_1 ( \mathfrak V ) \longrightarrow {\mathfrak L}_1 ( \mathfrak V )\, $ of the identity map of $\, \mathfrak V\, $ is the identity map. Using this construction it is shown that any unital complete real ordered Banach lattice is order isomorphic with a monotone complete (hence injective) abelian $C^*$-algebra. As the author learned in the course of writing this paper the construction of $\, {\mathfrak L}_1 \bigl( \mathfrak V \bigr)\, $ is originally due to Peressini, cf. \cite{Wr}, but the author would like to stress that the construction given here was done completely independent and without knowledge of the former results. If $\, {\mathfrak L}_1 ( \mathfrak V )\, $ is the minimal (complete) Banach lattice containing $\, \mathfrak V\, $ as a subspace, then the second lattice construction $\, \mathfrak L ( \mathfrak V )\, $ gives the somewhat maximal reasonable Banach lattice generated by $\, \mathfrak V\, $. Contrary to the prior case it is not generally monotone complete and in order to obtain a monotone complete lattice one defines the {\it Jordan lattice completion} as 
$$ {\mathfrak L}_q ( \mathfrak V )\> =\> {\mathfrak L}_1 \bigl( \mathfrak L ( \mathfrak V ) \bigr) \> . $$
It is then shown in Proposition 1 that there exists a unique surjective monotonous map 
$$ \pi :\> \mathfrak L ( \mathfrak V )\> \twoheadrightarrow\> {\mathfrak L}_1 ( \mathfrak V ) $$
extending the identity map of $\, \mathfrak V\, $. It turns out that $\, \pi\, $ is a lattice map so that it  extends to a (nonunique) surjective $*$-homomorphism from $\, {\mathfrak L}_q ( \mathfrak V )\, $ onto the latter by Theorem A of section 2. In case that $\, A = \mathfrak V\, $ is a (noncommutative) $C^*$-algebra there is a so called {\it Jordan squaring operation} defined on $\, {\mathfrak L}_q ( A )\, $ which is closely related to the $C^*$-square in $\, A\, $ and another partially inverse {\it Jordan squareroot operation} which extends the squareroot operation on $\, A\, $. These are investigated in detail and the Jordan squaring operation is then related to the (commutative) $C^*$-square in 
$\, {\mathfrak L}_q ( A )\, $ and certain of its quotients which results make up a major part of Theorem 1. In addition one obtains for any embedding of function systems $\, \mathfrak X \subseteq \mathfrak Y\, $ a canonical {\it normal} embedding of Jordan lattices $\, \Lambda : {\mathfrak L}_q ( \mathfrak X )\hookrightarrow {\mathfrak L}_q ( \mathfrak Y )\, $. The fact that $\, \Lambda\, $ is normal (and $*$-homomorphic) is quite surprising and is established in Theorem 1. In the course of its proof two canonical monotonous retractions 
$$ \overline r\, ,\, \underline r :\> {\mathfrak L}_q ( \mathfrak Y )\> \twoheadrightarrow\> {\mathfrak L}_q ( \mathfrak X ) $$
are employed which themselves are seen to posess very specific properties and these maps remain essential for large parts of the paper (for example in the proofs of Proposition 3 and Theorem 2). 
A somewhat isolated result in the first section is Proposition 2 which states that any completely positive complete quotient map 
$\, q : \mathfrak X \twoheadrightarrow M_n ( \mathbb C )\, $ of an operator system onto a matrix algebra which is an operator order epimorphism admits approximate completely positive splittings which are completely bounded by $\, n\, $. This result apart from being of interest for its own sake is used in the proof of Theorem  2 as well as Proposition 3 whose main outcome is to show that certain natural monotonous maps considered below are {\it boundedly basically increasing normal} (for an explanation of this technically rewarding but rather unintuitive notion we refer the reader to section 1 below).
\par\noindent
The second section is shorter than the first and contains in Theorem 2 one of the main results of this article as far as operator algebras are concerned. It exhibits a close relation between the injective envelope $\, I ( A )\, $ of a separable $C^*$-algebra $\, A\, $ and its enveloping von Neumann algebra $\, R\, $ in some given faithful $*$-representation on a separable Hilbert space provided $\, R\, $ is injective. In addition the notion of a weak injective envelope (weakly injective $C^*$-algebra) is defined and it is shown that any separable $C^*$-algebra $\, A\, $ posesses a weak injective envelope which can be identified with a certain monotone complete subspace of the positive injective envelope $\, {\mathfrak L}_1 ( A )\, $. It is also shown that any von Neumann algebra $\, R \subseteq \mathcal B ( \mathcal H )\, $ is weakly injective, i.e. there exists a positive linear projection $\, P : \mathcal B ( \mathcal H ) \rightarrow \mathcal B ( \mathcal H )\, $ with range $\, R\, $. Another fundamental result of section 2 is Theorem A which may be seen as a starting point of the whole paper as far as injective abelian $C^*$-algebras are concerned (whence the name). It reminds us that such an algebra is injective not only for positive linear maps, but also for $*$-homomorphisms (of abelian $C^*$-algebras) and linear lattice maps. This is used throughout the paper.
\par\noindent
Finally the last section is a bit different in flavour than the first two sections, but it is related in spirit since lattice constructions of different type are regarded. It can be considered as an autochtonous essay whose perception is in large parts independent of the prior article and is of interest for its own sake. The section investigates on conditions needed for a $C^*$-algebra in order that its subset of positive projections is a lattice (resp. complete lattice). The notion of $\mathcal P$-algebra and the notion of $\mathcal P$-map of $\mathcal P$-algebras are defined. These maps are essentially determined by their values on projections and they send projections to projections. Theorem P gives an account of the available results concerning these concepts. It seems that the state of this theory is a bit premature and the benefit of the included investigations lies more in posing questions than giving satisfactory answers.
In the last part of the section another class of monotonous maps called $\mu $-maps are considered together with the notions of $\mu $-measures and $\mu $-integrals as well as $\mathcal P$-measures of various types both which we expect to have widespread natural applications (for example in convex geometry). The paper proposes its own definition of a {\it convex metric space} tailored to suit state spaces and $\mu $-state spaces of 
(separable) $C^*$-algebras. Not all results are stated in the form of theorems or propositions but some are embedded in the main text. Examples include the decomposition of a $C$-contractive $\mathcal P$-measure over an ideal, the Jordan type decomposition of a (signed) $A$-linear $\mathcal P$-measure and the geometric center of a compact convex metric space. To give the reader a guideline with the rather confusing terminology: $\mu $-measures and $\mu $-integrals usually relate to a commutative context whereas the notion of $\mathcal P$-measure and $\mathcal P$-integral is a  generalization of this concept to the case where the domain of definition is a (noncommutative) $\mathcal P$-algebra. To the contrary $\mu $-functionals or $\mu $-maps are monotonous maps with certain additional features which may be defined on general function systems (resp. $\mu $-subsystems thereof). This is because only in the commutatve context a $\mu $-measure integrates to a fully monotonous $\mu $-integral thus defining a $\mu $-functional, whereas the canonical integral of a $\mathcal P$-measure is monotonous only on abelian $\mathcal P$-subalgebras. 
The final result of this treatise shows that quotients of sequentially monotone complete $C^*$-algebras are again sequentially monotone complete if they admit an order isomorphic representation on separable Hilbert space (the precise result of Theorem 3, (i) is that for any monotone increasing sequence in the quotient and given a countable set of majorants of this sequence there exists a majorant smaller than any of the elements in the countable set so that the quotient is in some sense "close" to being sequentially monotone complete).  Some other results are added which we do not spell out here in detail.
\par\bigskip\bigskip\bigskip\noindent
{\bf 1.\quad The Jordan lattice completion.}
\par\bigskip\noindent
In this paper we are primarily interested in $C^*$-algebras. Some of the following constructions however can be done in the more general setting of ordered real Banach spaces. By this we mean a real Banach space $\, \mathfrak V\, $ with a distinguished closed convex cone of positive elements 
$\, {\mathfrak V}_+\, $ such that $\, {\mathfrak V}_+ \cap - {\mathfrak V}_+ = \{ 0 \}\, $ and 
$\, \mathfrak V = {\mathfrak V}_+ - {\mathfrak V}_+\, $. The ordered Banach space is {\it unital} if the supremum of the unit ball,  denoted $\, {\bf 1}\, $, exists in $\,\mathfrak V\, $  and is of of norm one, i.e.  
$\, {\bf 1} - x\geq 0\, $ whenever $\, \Vert x\Vert\leq 1\, $.  The usual definition of a unital ordered Banach space includes a regularity condition, i.e. $\, 0\leq x\leq y\, $ implies $\, \Vert x\Vert\leq\Vert y\Vert\, $ or equivalently that every element dominated by $\, \bf 1\, $ is of norm less or equal to one, cf. 
\cite{A-E}. We will refer to such a space as a {\it regular (unital)} ordered Banach space.
A (linear) positive injection $\, \phi : \mathfrak V \hookrightarrow \mathfrak W\, $ of ordered real Banach spaces is termed an {\it order isomorphic embedding} if  $\, \phi ( x )\geq 0\, $ implies $\, x\geq 0\, $
and $\,\phi\, $ is bounded from below. The real ordered Banach space $\,\mathfrak V\, $ is called {\it dually ordered} if its dual is again an ordered real Banach space with positive cone 
$\, {\mathfrak V}^*_+\, $ given by the functionals attaining only nonnegative values on 
$\, {\mathfrak V}_+\, $ and if $\, {\mathfrak V}_+ = {\mathfrak V}^{**}_+ \cap \mathfrak V\, $. This is equivalent to the property that $\, \mathfrak V\, $ admits an order isomorphic embedding into a $C^*$-algebra.
If such an embedding exists one sees by restriction of the positive functionals of $\, A\, $ to the image of $\, \mathfrak V\, $ that any given element of $\, {\mathfrak V}^*\, $, which is bounded also on the image of $\,\mathfrak V\, $ in $\, A\, $ because the inverse map of the embedding is bounded, can be extended to a continuous functional on  $\, A^{sa}\, $ from the Hahn-Banach theorem, which then decomposes into the difference of two positive functionals. Since $\, A_+ = A \cap A^{**}_+\, $ also 
$\, {\mathfrak V}_+ = \mathfrak V \cap {\mathfrak V}^{**}_+\, $ follows.
On the other hand if the two conditions are satisfied then the first of the two lattice constructions below gives an order isomorphic embedding of $\,\mathfrak V\, $ into a regular complete Banach lattice 
$\, {\mathfrak L}_1 ( \mathfrak V )\, $  which, in turn may be embedded into a regular unital ordered Banach space itself, and thus admits an order isomorphic embedding into the commutative $C^*$-algebra of continuous functions on its state space (endowed with the relative $w^*$-topology, cf. \cite{A-E}, chap. 2). 
Let $\, \mathfrak V\, $ be a dually ordered real Banach space. To each subset 
$\, \mathcal C \subseteq {\mathfrak V}_+\, $ of positive elements one attributes an element 
$\, \underline{\mathcal C}\, $ in a real vector space 
$\, {\mathfrak L}_1 ( \mathfrak V )\, $ called the {\it infimum of $\,\mathcal C\, $}. There are several natural equivalence relations available in order to build a linear space, we consider the following:
$\, \mathcal C {\sim }_1  {\mathcal C}'\, $ if and only if for each element 
$\,\, a\in \mathfrak V\, $ such that 
$\, a\leq c\, $ for every $\, c\in\mathcal C\, $ this implies $\, a\leq c'\, $ for all $\, c'\in {\mathcal C}'\, $ and vice versa.  For each set $\, \mathcal C\, $ let 
$\, {\mathcal C}_c = \{ a\in \mathfrak V\,\vert\, a\leq\mathcal C \}\, $ denote the subset of elements all of which are smaller or equal than each element in $\,\mathcal C\, $. Then putting 
$\, ({\mathcal C}_c)^c = \{ c\in \mathfrak V_+\,\vert\, c\geq {\mathcal C}_c \}\, $ one obtains a maximal representative for the equivalence class of $\, \mathcal C\, $. 
Given two elements 
$\, \underline{\mathcal C} = \inf\, \{ c_{\lambda } {\}}_{\lambda }\, $ and 
$\, \underline{\mathcal D} = \inf\, \{ d_{\mu } {\}}_{\mu }\, $ as above their {\it sum} is defined to be 
$\, \underline{\mathcal C} + \underline{\mathcal D} = 
\inf\, \{ c_{\lambda } + d_{\mu }\,\vert\, c_{\lambda }\in {\mathcal C}\, ,\, d_{\mu }\in {\mathcal D} {\}}_{\lambda , \mu }\, $ and one checks that this definition is compatible with the equivalence relation hence well defined.  One also has  
$$ 2\, \mathcal C\> =\> \{ 2 c_{\lambda }\,\vert\, c_{\lambda }\in\mathcal C \}\> {\sim }_1\> 
\{ c_{\lambda } + c_{\mu }\,\vert\, c_{\lambda }\, ,\, c_{\mu }\in \mathcal C \}\> =\> \mathcal C\> +\> 
\mathcal C   $$
since if $\, a\in \mathfrak V\, $ is any element such that $\, a\leq 2 c_{\lambda }\, $ for all $\, \lambda\, $ then also $\, a \leq c_{\lambda } + c_{\mu }\, $ for all $\,\lambda\, ,\, \mu\, $. Therefore the scalar multiplication of elements $\, \underline{\mathcal C}\, $ with positive rational numbers is well defined. Approximating a real positive number by positive rational numbers from below one sees that scalar multiplication extends to the (positive) reals. 
The convex cone of such elements generates a 
real linear space $\, {\mathfrak L}_1 ( \mathfrak V )\, $ by taking arbitrary differences 
$\, \mathcal A = \underline{\mathcal C} - \underline{\mathcal D}\, $, subject to the equivalence relation
generated by $\, {\sim }_1\, $ and 
$$ ( \underline{\mathcal C} + \underline{\mathcal E} ) - ( \underline{\mathcal D} + \underline{\mathcal E} ) \> {\simeq }_1\> \underline{\mathcal C} - \underline{\mathcal D} $$
with $\, \mathcal E\subseteq {\mathfrak V}_+\, $ a subset of positive elements, and an order on this space is obtained by 
$$ \mathcal A \geq \mathcal B\> \iff\> \mathcal A - \mathcal B\> {\simeq }_1\> \underline{\mathcal C} \geq 0\> . $$
One finds that for positive elements the equivalence relation $\, {\simeq }_1\, $ which by definition is a stabilized version of $\, {\sim }_1\, $, i.e. 
$\, \mathcal C\> {\simeq }_1\> \mathcal D\, $ iff there exists a positive element
$\, \mathcal E\, $ such that $\, \mathcal C + \mathcal E {\sim }_1 \mathcal D + \mathcal E\, $, in fact reduces to $\, {\sim }_1\, $ so this equivalence relation has cancellation. Namely, suppose 
$\, \mathcal C + \mathcal E{\sim }_1 \mathcal D + \mathcal E\, $. 
Let $\, e_0\in \mathcal E\, $ be a fixed element and denoting 
$\, {\mathcal E}_c = \{ b\in \mathfrak V\,\vert\, b\leq \mathcal E \}\, $ the set of elements all of which are smaller or equal than every element in $\,\mathcal E\, $ consider the positive subset 
$\,  e_0 - {\mathcal E}_c = \{ e_0 - b\,\vert\, b\in {\mathcal E}_c \}\, $. 
Then since taking sums is compatible with $\, {\sim }_1\, $ one gets 
$\, \mathcal C + \mathcal E + ( e_0 - {\mathcal E}_c ) {\sim }_1
\mathcal D + \mathcal E + ( e_0 - {\mathcal E}_c )\, $. On the other hand 
$\, \mathcal E + ( e_0 - {\mathcal E}_c ) {\sim }_1 \{ e_0 \}\, $ since  
$\, a\leq e + e_0 - b\, $ for all $\, e\in\mathcal E\, ,\, b\in {\mathcal E}_c\, $ implies 
$\, b + ( a - e_0 )\in {\mathcal E}_c\, $ and by induction 
$\, b + n ( a - e_0 )\in {\mathcal E}_c\, $ for all $\, b\in {\mathcal E}_c\, ,\, n\in\mathbb N\, $. This then implies that 
$\, a - e_0 \leq 0\, $ since otherwise there would exist a positive functional $\,\phi\in {\mathfrak V}_+\, $ with $\, \phi ( a - e_0 ) > 0\, $ so that 
$\, \sup_{b\in {\mathcal E}_c} \{ \phi ( b ) \} = +\infty\, $ which is impossible since $\, {\mathcal E}_c\, $ is bounded above.  Therefore $\, a\leq e_0\, $ follows and the reverse implication is trivial. Thus 
$\, \mathcal C + \{ e_0 \} {\sim }_1 \mathcal D + \{ e_0 \}\, $ and since 
$\, \{ e_0 \}\, $ consists of a single element one concludes that $\,\mathcal C {\sim }_1 \mathcal D\, $.
One then sees by a similar argument that a general element $\, \mathcal A = \underline{\mathcal C} - 
\underline{\mathcal D}\, $ is equal to the infimum of the set 
$\, \{ c - a\,\vert\, c\in\mathcal C\, ,\, a\in {\mathcal D}_c \}\, $. Conversely any set of selfadjoint elements which is bounded below uniquely determines an element in 
$\, {\mathfrak L}_1 ( \mathfrak V )\, $ by taking its infimum, since each subset bounded from below is given as the difference of a positive subset and a single positive element of $\, \mathfrak V\, $ (viewed as a positive subset).
In particular $\, \underline{\mathcal C} \geq \underline{\mathcal D} \iff 
{\mathcal D}_c\subseteq {\mathcal C}_c\, $.
For an arbitrary subset $\, \{ \underline{\mathcal C}_{\lambda } {\}}_{\lambda }\, $ of positive elements its 
infimum $\, \inf_{\lambda }\, \underline{\mathcal C}_{\lambda }\, $ which is the largest element smaller or equal than each $\, \underline{\mathcal C}_{\lambda }\, $ is well defined and is given by the element corresponding to the union $\, \cup_{\lambda }\, {\mathcal C}_{\lambda }\, $. 
By shifting with a suitable positive element (resp. positive scalar if $\, \mathfrak V\, $ is unital) one may equally define an infimum for every set of general (selfadjoint) elements $\, \{ {\mathcal A}_{\lambda } {\}}_{\lambda }\, $ which is bounded below and correspondingly, by the symmetry $\, \mathcal A\mapsto - \mathcal A\, $ a supremum for every set of elements bounded above.  For example taking each $\, {\mathcal C}_{\lambda } = \{ c_{\lambda } \}\, $ to consist of a single positive element one gets 
$\, \inf_{\lambda }\, \{ {\mathcal C}_{\lambda } \} = \inf\, {\mathcal C} = \sup\, {\mathcal C}_c\, $ with 
$\, \mathcal C = \{ c_{\lambda } {\}}_{\lambda } = \cup\, {\mathcal C}_{\lambda }\, $. 
In particular one has the lattice operations
$$  \mathcal A\wedge \mathcal B\> ,\quad 
\mathcal A \vee \mathcal B $$
denoting the (unique !) maximal element smaller or equal to $\, \mathcal A\, $ and $\,\mathcal B\, $, resp. the unique minimal element larger or equal to both $\,\mathcal A\, $ and $\,\mathcal B\, $. 
If $\, \mathfrak V\, $ is not unital and $\, \mathcal A = \underline{\mathcal C} - \underline{\mathcal D}\, ,\, \mathcal B = \underline{\mathcal E} - \underline{\mathcal F}\, $ these operations may be defined as 
$$ \mathcal A \wedge \mathcal B\> =\> \inf\, \bigl\{\underline{\mathcal C} + \underline{\mathcal F}\, ,\, \underline{\mathcal E} + \underline{\mathcal D} \bigr\}\> -\> \bigl( \underline{\mathcal D} + \underline{\mathcal F} \bigr)\> , $$
$$ \mathcal A \vee \mathcal B\> =\> \sup\, \bigl\{ \underline{\mathcal C} + \underline{\mathcal F}\, ,\, \underline{\mathcal E} + \underline{\mathcal D} \bigr\}\> -\> \bigl( \underline{\mathcal D} + \underline{\mathcal F} \bigr)\> . $$
Taking 
$\, \mathcal B = 0\, $ the unique minimal positive decomposition of $\, \mathcal A\, $ is given by 
$$ \mathcal A\> =\> {\mathcal A}_+ - {\mathcal A}_-\> =\> ( \mathcal A \vee 0 ) + ( \mathcal A \wedge 0 )
\> .  $$ 
Indeed, suppose given two different positive decompositions 
$$ \mathcal A\> =\> \underline{\mathcal C}\, -\, \underline{\mathcal D}\> =\> \underline{\mathcal C}' \, -\, 
\underline{\mathcal D}'  $$
we claim that $\, \mathcal A = ( \underline{\mathcal C}\wedge \underline{\mathcal C}'  )\, -\, 
( \underline{\mathcal D}\wedge \underline{\mathcal D}' )\, $. To see this it is sufficient by the symmetry 
$\, \mathcal A\mapsto - \mathcal A\, $ to prove 
$$  \underline{\mathcal C}\, -\, \underline{\mathcal D}\> =\> \mathcal A\> \leq\> ( \underline{\mathcal C} \wedge \underline{\mathcal C}' )\, 
-\, ( \underline{\mathcal D} \wedge \underline{\mathcal D}' ) $$
$$\qquad \iff \underline{\mathcal C} + ( \underline{\mathcal D}\wedge \underline{\mathcal D}' )\> \leq\> 
( \underline{\mathcal C} \wedge \underline{\mathcal C}' )\, +\, \underline{\mathcal D}  $$
using the relation $\, \mathcal C + {\mathcal D}' {\sim }_1 {\mathcal C}' + \mathcal D\, $. Let 
$\, a\in A^{sa}\, $ be an element smaller or equal than each element in 
$\, \mathcal C + ( \mathcal D \cup {\mathcal D}' )\, $ and given an arbitrary element 
$\, b\in ( \mathcal C \cup {\mathcal C}' ) + \mathcal D\, $. Then either $\, b\in \mathcal C + \mathcal D\, $ in which case 
$\, a\leq b\, $ follows trivially, or $\, b\in {\mathcal C}' + \mathcal D\sim \mathcal C + {\mathcal D}' \, $ which again implies $\, a\leq b\, $. Using induction the minimal positive decomposition of 
$\, \mathcal A\, $ is given by 
$$ \mathcal A\> =\> {\mathcal A}_+\, -\, {\mathcal A}_-\> =\>  \inf_{\lambda }\, \{ \underline{\mathcal C}_{\lambda }\,\vert\, \mathcal A = 
\underline{\mathcal C}_{\lambda } - \underline{\mathcal D}_{\lambda } \}\, -\, 
\inf_{\lambda }\, \{ \underline{\mathcal D}_{\lambda }\,\vert\, 
\mathcal A = \underline{\mathcal C}_{\lambda }- \underline{\mathcal D}_{\lambda } \} \> . $$
Then $\, {\mathcal A}_+ \geq \mathcal A \vee 0\, $ with 
$\,\mathcal A = ( \mathcal A\vee 0 )\, -\, ( ( \mathcal A \vee 0 ) - \mathcal A )\, $ a positive decomposition of 
$\, \mathcal A\, $. 
From uniqueness one must have $\, {\mathcal A}_+ =  \mathcal A \vee 0\, $ and hence 
$\, {\mathcal A}_- = - ( \mathcal A \wedge 0 )\, $. 
For $\, \underline{\mathcal C} \geq 0\, $ a natural norm (satisfying the triangle inequality) is given by 
$$ \Vert \underline{\mathcal C} \Vert\> =\> \inf\, \{ \Vert d\Vert\,\vert\,  d\in \mathcal D\, ,\, \mathcal D {\sim }_1 \mathcal C \}\> =\> \inf\, \{ \Vert c\Vert\,\vert\, c\in ( {\mathcal C}_c )^c \}\> . $$
This norm is extended to general elements $\, \mathcal A = \underline{\mathcal C} - \underline{\mathcal D}\, $ by defining 
$$ \Vert \mathcal A \Vert = \max\, \bigl\{ \Vert {\mathcal A}_+\Vert\, ,\, 
\Vert {\mathcal A}_- \Vert \bigr\} \geq 0 \> . $$
The triangle inequality is readily checked since 
$\,  ( \mathcal A + \mathcal B )_{\pm } \leq  {\mathcal A }_{\pm } + {\mathcal B}_{\pm }\, $. Define 
$\, {\mathfrak L}_1 ( \mathfrak V )\, $ to be the resulting normed linear space which clearly is a complete vector lattice. Being monotone complete it is a Banach space, hence a regular unital Banach lattice. Namely if $\, \bigl\{ {\mathcal A}_n \bigr\} \subseteq {\mathfrak L}_1 ( \mathfrak V )\, $ is a Cauchy sequence then the sequence converges in norm towards the element 
$$ \mathcal A\> =\> \liminf_n\, {\mathcal A}_n\> =\> \sup_n\, \bigl\{ \inf_{m\geq n}\, \bigl\{ {\mathcal A}_m \bigr\} \bigr\}\> . $$
Thus $\, {\mathfrak L}_1 ( \mathfrak V )\, $ is a regular complete Banach lattice
 containing a canonical subspace order isomorphic with $\,\mathfrak V\, $. In many cases, e.g. if $\, \mathfrak V\, $ is a (unital) 
$C^*$-algebra (or for that matter a regular unital ordered Banach space = function system), the embedding $\, \mathfrak V \hookrightarrow {\mathfrak L}_1 ( \mathfrak V )\, $ is isometric with respect to the original norm of $
\, \mathfrak V\, $ (see the analogous argument below for 
$\, \mathfrak L ( A )\, $). In any case if $\,\mathfrak V\, $ is a dually ordered Banach space the embedding $\, \mathfrak V \hookrightarrow {\mathfrak L}_1 ( \mathfrak V )\, $ is order isomorphic. 
To see this consider a general selfadjoint element $\, x\in \mathfrak V\, $ and let 
$\,\mathfrak x = \inf\, \{ x \} = \sup\, \{ x \}\, $ be its image in $\, {\mathfrak L}_1 ( \mathfrak V )\, $. Then its minimal positive decomposition is given by $\, \mathfrak x = {\mathfrak x}_+ - {\mathfrak x}_-\, $ with 
$\, {\mathfrak x}_{\pm } = \inf\, \{\, c\, \vert\, 0\leq c\, ,\, \pm x\leq c \}\, $ so that $\, {\mathfrak x}_- = 0\, $ if and only if $\, x\geq 0\, $. Since the dual $\, {\mathfrak V}^*\, $ is positively generated, there exists 
$\, \alpha \geq 1\, $ such that $\, \mathfrak V\, $ is $\alpha $-normal which implies 
$$ \Vert \mathfrak x \Vert\> =\> \max\, \{ \Vert {\mathfrak x}_+\Vert\, ,\, \Vert {\mathfrak x}_-\Vert \}\> \geq\> {{\alpha }^{-1}\over 2 + {\alpha }^{-1}}\> \Vert x \Vert\>\geq\> {{\alpha }^{-1}\over 3} \>\Vert x\Vert $$
showing that the embedding $\, \mathfrak V \hookrightarrow {\mathfrak L}_1 ( \mathfrak V )\, $ is bounded below. Clearly, $\, {\mathfrak L}_1 ( \mathfrak V )\, $ is regular, i.e. 
$\, 0\leq x\leq y\, $ implies $\, \Vert x\Vert\leq \Vert y\Vert\, $. If $\, \mathfrak V\, $ is unital $\, {\mathfrak L}_1 ( \mathfrak V )\, $ is a unital complete regular Banach lattice, and as such isometrically order isomorphic with an (injective) commutative $C^*$-algebra by the Corollary below. 
In case of a general dually ordered Banach space one may adjoin a unit element by considering the linear space $\, {\widetilde{\mathfrak L}}_1 = {\mathfrak L}_1 ( \mathfrak V )\, +\, \mathbb R\, {\bf 1}\, $, where 
$\, {\bf 1} = \sup\, \{ x\in {\mathfrak L}_1 ( \mathfrak V )\,\vert\, \Vert x\Vert \leq 1 \} = \sup\, \{ x \in {\mathfrak L}_1 ( \mathfrak V )\,\vert\, x \geq 0\, ,\, \Vert x\Vert \leq 1 \}\, $ and ordering given by the positive cone 
$$ {\widetilde{\mathfrak L}}_1 ( \mathfrak V )_+\> =\> \bigl\{ p = \alpha\, a\, +\, \beta\, ( {\bf 1} - b )\,\bigm\vert\, \alpha\, ,\,\beta \geq 0\, ,\, 0\leq a\, ,\, b\, \leq {\bf 1} \bigr\} \> . $$
A norm is imposed on positive elements by the formula 
$$ \Vert p {\Vert}_1\> =\> \sup\, \bigl\{ \Vert c_+\Vert\,\bigm\vert\, c\leq p\, ,\, c\in {\mathfrak L}_1 ( \mathfrak V ) \bigr\} $$
which is well defined due to $\, {\mathfrak L}_1 ( \mathfrak V )\, $ being a regular lattice and satisfies the triangle inequality. Check that $\, p \geq c\, $ implies $\, p \geq c_+\, $ for $\, c\in {\mathfrak L}_1 ( \mathfrak V )\, $. This follows since $\, c_+\wedge c_- = 0\, $ implies $\, ( b + c_+ ) \wedge c_- = 
( b \wedge c_- ) + c_+\, $ and $\, c = d + e\, $ implies $\, c_+ \leq d_+ + e_+\, $. For $\, p\in {\mathfrak L}_1 ( \mathfrak V )\, $ one recovers the original norm. For general selfadjoint elements define 
$$ \Vert x {\Vert }_1\> =\> \inf\, \bigl\{ \max\, \{ \Vert p {\Vert }_1\, ,\, \Vert q {\Vert }_1 \}\,\bigm\vert\, x = p - q\, ,\, p , q \geq 0 \bigr\} $$
and check that the axioms for a norm are satisfied. In particular if 
$\, x = p - q\in {\mathfrak L}_1 ( \mathfrak V )\, $ with $\, p = a + \lambda ( {\bf 1} - b )\, ,\, q = c + \mu ( {\bf 1} - d )\, ,\, 
a\, ,\, b\, ,\, c\, ,\, d\, \in {\mathfrak L}_{1 , +}\, $ one has $\, \lambda = \mu\, $ and the triangle inequality gives $\, \Vert x {\Vert }_1 \geq \Vert x \Vert\, $ since 
$$ \Vert p {\Vert }_1\> \geq\> \Vert x_+ \Vert\> ,\> \Vert q {\Vert }_1\> \geq\> \Vert x_- \Vert \quad\Longrightarrow\quad \Vert x {\Vert }_1\>\geq\> \Vert x_+\Vert + \Vert x_- \Vert\> =\> \Vert x \Vert $$ 
so the embedding $\, {\mathfrak L}_1 ( \mathfrak V ) \hookrightarrow {\widetilde{\mathfrak L}}_1 ( \mathfrak V )\, $ is isometric and order isomorphic proving that 
$\,\mathfrak V\, $ admits an order isomorphic embedding into a regular unital ordered Banach space and a forteriori into a (commutative) $C^*$-algebra. The details of this argument are left to the reader.
In general the construction is not functorial with respect to arbitrary positive linear 
maps, meaning that in general there is no natural positive linear extension 
$\, {\mathfrak L}_1 ( \phi ) : {\mathfrak L}_1 ( \mathfrak V ) \rightarrow {\mathfrak L}_1 ( \mathfrak W )\, $ of a given positive linear map $\, \phi : \mathfrak V \rightarrow \mathfrak W\, $. If however $\, \pi :\,  \mathfrak V \twoheadrightarrow \mathfrak W\, $ is a complete order epimorphism so that the image set of the lower complement of any complemented positive subset $\, \mathcal C\subseteq {\mathfrak V}_+\, ,\, \mathcal C = \bigl( {\mathcal C}_c \bigr)^c $ is equal to the complement of the image set of $\, \mathcal C\, $ in $\,\mathfrak W\, $ 
there exists a unique monotonous extension $\, {\mathfrak L}_1 ( \phi ) :\, {\mathfrak L}_1 \bigl( \mathfrak V \bigr) \twoheadrightarrow {\mathfrak L}_1 \bigl( \mathfrak W \bigr)\, $ of 
$\, \phi\, $ which must be linear if $\, \phi\, $ is linear (see below).
One may wonder if there are some simpler notions ensuring functoriality. One notes the following: if 
$\, \pi :\,  \mathfrak V \twoheadrightarrow \mathfrak W\, $ is a simple order epimorphism, i.e. if $\, \pi ( {\mathfrak V}_+ ) = {\mathfrak W}_+\, $ and in addition the kernel of $\,\pi\, $ is positively generated, then for any two elements 
$\, x\, ,\, y\in {\mathfrak V}_+\, $ any element $\, b\in \mathfrak W\, $ with 
$\, b\leq \pi ( x )\, ,\, b\leq \pi ( y )\, $ lifts to an element $\, a\in\mathfrak V\, $ with $\, a\leq x\, ,\, a\leq y\, $. Since $\,\pi\, $ is an order epimorphism there exist lifts $\, a'\, ,\, a''\, $ of $\, b\, $ with $\, a'\leq x\, ,\, a''\leq y\, $. Then $\, a' - a'' = c - d\, $ with $\, c\, ,\, d\geq 0\,,\, c\, ,\, d\in \ker\,\pi\, $, so that putting 
$\, a = a' - c = a'' - d\, $ gives the result. This argument however only applies to finite subsets of 
$\, {\mathfrak V}_+\, $ and fails for arbitrary complemented positive subsets. A necessary and sufficient condition that $\, \phi\, $ admits a unique monotonous extension is the following: given any 
complemented positive subset $\, \mathcal C\subseteq {\mathfrak V}_+\, $ then $\, \inf\, \phi 
\bigl( \mathcal C \bigr)\, $ must be larger or equal to 
$\, \sup\, \phi \bigl( {\mathcal C}_c \bigr)\, $ which is encoded in the relation 
$\, \phi \bigl( {\mathcal C}_c \bigr) \subseteq \phi \bigl( \mathcal C \bigr)_c\, $. Then uniqueness of any extension can be encoded by the condition 
$\, \phi \bigl( {\mathcal C}_c \bigr)^c = \bigl( \phi \bigl( \mathcal C \bigr)_c \bigr)^c\, $ or equivalently 
$\, \sup\, \phi \bigl( \mathcal C \bigr)_c \leq \inf\, \phi \bigl( {\mathcal C}_c \bigr)^c\, $ since the reverse implication is trivial. Indeed if this relation holds for all complemented subsets then any monotonous extension is uniquely determined. Conversely suppose that the condition fails for some 
positive complemented subset $\, \mathcal C\, $. For any positive element $\, \mathfrak C\in {\mathfrak L}_1 ( \mathfrak V )\, $ write 
$$ \mathfrak C\> =\> \inf\, \mathcal C\> =\> \inf\, \bigl\{ c\in \mathfrak V\bigm\vert c \geq \mathfrak C  \bigr\}\> =\> \sup\, \bigl\{ b\in \mathfrak V\bigm\vert b \leq \mathfrak C \bigr\}\> =\> \sup\, {\mathcal C}_c  \> . $$
Putting $\, \underline\phi \bigl( \mathfrak C \bigr)\> =\> \inf\, \bigl\{ \phi ( \mathcal C ) \bigr\}\, $ and 
$\, \overline\phi \bigl( \mathfrak C \bigr) = \sup\, \bigl\{ \phi ( {\mathcal C}_c ) \bigr\}\, $ gives two different monotonous extensions of $\, \phi\, $ to $\, {\mathfrak L}_1 ( \mathfrak V )\, $.
Since the functoriality of 
$\, {\mathfrak L}_1\, $  for linear (positive) maps is so restricted one is forced or maybe seduced to consider a wider notion of functoriality.
\par\medskip\noindent
{\bf Definition.}\,  
(a)\, A positive map $\, m :\, \mathfrak V \rightarrow \mathfrak W\, $ of  two ordered Banach spaces is called {\it monotonous} iff
$$  x\>\leq\> y\>\Rightarrow\> m ( x )\>\leq\> m ( y ) \> . $$
It is called {\it (positively) homogenous} iff
$$ m ( \alpha\, x )\> =\> \alpha\, m ( x )\> ,\qquad \forall\> \alpha\in \mathbb R ( {\mathbb R}_+)\, ,\, x\in \mathfrak V\> .  $$
(b)\, A positive map $\, c :\, \mathfrak V \rightarrow \mathfrak W\, $ of two ordered Banach spaces is called {\it convex} if it is positively homogenous, monotonous and if the following condition holds
$$ c ( x + y )\> \leq\> c ( x )\> +\> c ( y )\> ,\qquad \forall\, x\, ,\, y\,\geq\, 0  \> . $$
Similar $\, c\, $ is called {\it concave} if it is monotonous, positively homogenous and if 
$$ c ( x + y )\> \geq\> c ( x )\> +\> c ( y )\> ,\qquad \forall\, x\, ,\, y\,\geq\, 0  \> . $$
(c)\, A monotonous map of unital dually ordered Banach spaces $\, \phi : \mathfrak V 
\twoheadrightarrow \mathfrak W\, $ is called {\it a bounded complete order epimorphism} iff 
for any bounded subset $\, \mathcal C\subseteq \mathfrak V\, $ the set 
$\, \phi \bigl( {\mathcal C}_c \bigr)\, $ is normdense in $\, \phi \bigl( \mathcal C \bigr)_c\, $.
\par\bigskip\noindent
Usually we will only consider monotonous maps which in addition are {\it positively homogenous}, i.e. $\, m ( \alpha\, x ) = \alpha\, m ( x )\, $ for $\,\alpha\in {\mathbb R}_+\, $. On the other hand there may be situations where it is convenient to drop this extra condition.
We also use the term {\it convex map} loosely for monotonous maps which are convex restricted to the positive cone and real homogenous (so that a homogenous convex map is in fact concave restricted to the cone of negative elements) as well as occasionally for maps which are convex only in the strict sense, i.e. $\, c\, \bigl( \lambda x + ( 1 - \lambda ) y \bigr) \leq \lambda\, c ( x ) + ( 1 - \lambda )\, c ( y )\, $.
An amiable aspect of a real homogenous convex map $\, c\, $ defined on a Banach lattice $\, \mathfrak L\, $ and assuming 
$\, c ( x ) = c ( x_+ ) - c ( x_- )\, $ with respect to the minimal positive decomposition of $\, x\, $ and further assuming that the kernel $\, \ker c = \{ x\in\mathfrak L\,\vert\, c ( x ) = 0 \}\, $ is {\it positively generated}, so that $\, x\in \ker\, c\, $ implies 
$\, x_{\pm }\in \ker\, c\, $, is that the kernel must be a linear subspace, and $\, c\, $ is linear for addition by elements from $\, \ker\, c\, $ since for every positive $\, z\in {\mathfrak L}_+\, ,\, d\in \ker\, c\, ,\, d\geq 0\, $ one has 
$$ c ( z )\>\leq\> c ( z + d )\>\leq\> c ( z )\> +\> c ( d )\> =\> c ( z )\> . $$
It is also true that a positive map ($\, c ( x ) \geq 0\, $ whenever $\, x\geq 0\, $) which is suplinear is automatically monotonous since $\, x \geq y\, $ implies $\, c ( x ) = c \bigl( y + ( x - y ) \bigr)\geq 
c ( y ) + c\bigl( x - y \bigr) \geq c ( y )\, $. By the symmetry $\, c ( x ) \mapsto - c ( - x )\, $ turning suplinear maps to sublinear maps a negative sublinear map is automatically monotonous.
We will later see some interesting examples of convex and concave maps. A real ordered Banach space $\, \mathfrak V\, $ whose positive cone $\, {\mathfrak V}_+\, $ admits an order unit 
$\, e\in {\mathfrak V}_+\, $ such that the relation $\, x\leq e\, $ holds for any 
$\, x\in B_1 ( \mathfrak V )\, $ in the unit ball is called a {\it preunital} ordered Banach space.  
\par\medskip\noindent
{\bf Corollary.}\quad A preunital dually ordered real Banach space $\, \mathfrak V\, $ is injective in the category of  preunital ordered Banach spaces and positive linear maps if and only if it is a complete Banach lattice. Also, any regular unital complete Banach lattice is isometrically order isomorphic to an injective commutative $C^*$-algebra and any preunital complete Banach lattice is order isomorphic to a regular unital complete Banach lattice, hence to a commutative $C^*$-algebra.
\par\bigskip\noindent
{\it Proof.}\quad Suppose that $\,\mathfrak V\, $ is a complete Banach lattice and given a (preunital) inclusion of preunital ordered real Banach spaces $\, \mathfrak X\subseteq \mathfrak Y\, $ together with a positive linear map $\, \phi :\,  \mathfrak X \rightarrow \mathfrak V\, $. Pick an arbitrary element $\, y_0\in {\mathfrak Y}_+ \backslash {\mathfrak X}_+\, $ and define 
$\, {\phi }_0 : \mathfrak X\, +\, \mathbb R y_0 \rightarrow \mathfrak V\, $ by linear extension of 
$\, {\phi }_0 {\vert }_{\mathfrak X} = \phi\, $ and 
$$ {\phi }_0 ( y_0 )\> =\ \sup\, \{ \phi ( x )\,\vert\, x\in \mathfrak X\, ,\, x\leq y_0 \} \> . $$
Since the inclusion $\, \mathfrak X \hookrightarrow \mathfrak Y\, $ is preunital the supremum is well defined.
Suppose that $\, y = x + \alpha y_0 \geq 0\, ,\, \alpha \in \mathbb R\, ,\, x\in \mathfrak X\, $ is a positive element in $\, \mathfrak X + \mathbb R\, y_0\, $. If $\,\alpha \leq 0\, $ then 
$$ {\phi }_0 ( y ) = \inf\, \{ \phi ( z )\, \vert\, z\in\mathfrak X\, ,\, z\geq y\geq 0 \}\> \geq \> 0 \> . $$
On the other hand if $\, \alpha > 0\, $ then 
$$ {\phi }_0 ( y )\> =\> \sup\, \{ \phi ( z )\,\vert\, z\in\mathfrak X\, ,\, z\leq y_0 \}\> \geq\> \phi ( 0 )\> =\> 0\> . 
$$
One proceeds by (transfinite) induction to extend $\,\phi\, $ to a positive map 
$\, \overline\phi : \mathfrak Y \rightarrow \mathfrak V\, $. This proves that any complete Banach lattice is injective for positive linear maps of preunital ordered Banach spaces. For the converse assume that the preunital dually ordered Banach space $\, \mathfrak V\, $ is injective for such maps. Then $\, \mathfrak V\, $ admits an order isomorphic preunital embedding into some preunital complete Banach lattice $\, \mathfrak L\, $ from the argument above  and the identity map of
$\, \mathfrak V\, $ extends to a positive map 
$\, \iota : \mathfrak L \rightarrow \mathfrak V\, $. 
Since any subset 
$\, \mathcal C\subseteq \mathfrak V\, $ which is bounded below has an infimum 
$\, x\in \mathfrak L\, $ the element 
$\, \iota ( x )\, $ is an infimum for $\,\mathcal C\, $ in $\, \mathfrak V\, $ by positivity of $\,\iota\, $. 
Therefore any subset $\,\mathcal C \subset \mathfrak V\, $ which is bounded below has an infimum and $\,\mathfrak V\, $ is a complete lattice. 
To prove the second statement one notes that if $\, \mathfrak X\, $ is a regular unital complete Banach lattice, it is in particular a regular unital Banach space, so it admits a unital isometric and order isomorphic embedding into a unital commutative (injective) $C^*$-algebra $\, A\, $. By positive injectivity of 
$\,\mathfrak X\, $ there exists a positive linear projection $\, \Phi : A \rightarrow A\, $ with 
$\, \Phi ( A ) = \mathfrak X\, $. Then the unital operator space $\,\mathfrak X\, $ (with respect to the operator space structure determined by the embedding into $\, A\, $) admits a unique (!) structure as a commutative $C^*$-algebra (compare with \cite{E-R}, Theorem  6.1.3).
Now suppose that $\, \mathfrak X\, $ is any preunital complete Banach lattice. As such it is certainly a dually ordered Banach space, and applying the functor $\, {\mathfrak L}_1\, $ gives an order isomorphism with a regular preunital complete Banach lattice, i.e. $\, 0\leq x\leq y\, $ implies $\, \Vert x\Vert\leq\Vert y\Vert\, $, and moreover the norm of a general element is given in terms of the minimal positive decomposition 
$\, x = x_+ - x_-\, $ as 
$$ \Vert x\Vert\> =\> \max\, \{ \Vert x_+\Vert\, ,\, \Vert x_-\Vert \} \> . $$
Let $\, {\bf 1}\, $ denote the supremum over the unit ball of the regular complete Banach lattice $\,{\mathfrak L}_1 ( \mathfrak X )\, $ and $\, \Vert {\bf 1}\Vert \geq 1\, $ its norm. Define a new norm on $\,\mathfrak X\, $ by 
$$  \Vert\vert x\vert\Vert\> =\> \inf\,\bigl\{ \alpha\geq 0\,\bigm\vert\, x \leq \alpha {\bf 1} \bigr\}\> ,\quad x\geq 0   $$
and $\, \Vert\vert x\vert\Vert = \max\, \{ \Vert\vert x_+ \vert\Vert\, ,\, \Vert\vert x_-\vert\Vert \}\, $ in the general case and check that the axioms of a norm are satisfied, defining a structure of a regular unital complete Banach lattice such that $\, \Vert\vert x\vert\Vert \leq\Vert x\Vert\leq \ \Vert {\bf 1}\Vert\, \Vert\vert
x\vert\Vert\, $\qed
\par\bigskip\noindent
If $\,\mathfrak V\, $ is a unital dually ordered Banach space we may call $\, {\mathfrak L}_1 ( \mathfrak V )\, $ the {\it positive injective envelope} of $\,\mathfrak V\, $ since it is the smallest positively injective Banach space containing 
$\,\mathfrak V\, $ (by an order isomorphic embedding if not by isometric embedding) and has a corresponding rigidity property: 
every monotonous map $\, m : {\mathfrak L}_1 ( \mathfrak V ) \rightarrow {\mathfrak L}_1 ( \mathfrak V )\, $ which restricts to the identity map of $\,\mathfrak V\, $ is equal to the identity map. 
\par\bigskip\noindent
{\it Remark.}\quad There are several generalizations of the ${\mathfrak L}_1$-construction, one is given by the following scheme: a regular unital real ordered Banach cone is a closed unital affine subset $\, \mathfrak C\subseteq \mathfrak V\, $ of a regular unital real ordered Banach space $\, \mathfrak V\, $, i.e. for all
$\, c\, ,\, d\in \mathfrak C\, ,\, \alpha\in {\mathbb R}_+\, $ and $\,\gamma\in\mathbb R\, $ one has 
$$ c\> +\> d\in \mathfrak C\> ,\quad \alpha\, c\,\in\, \mathfrak C\> ,\quad \gamma\, {\bf 1}\in \mathfrak C\> .  $$
The order relation on the set $\, \mathfrak C\, $ is induced by the order on 
$\,\mathfrak V\, $. Then also the subset $\, - \mathfrak C\subseteq \mathfrak V\, $ is a regular unital real ordered Banach cone, and one defines an equivalence relation for subsets $\, \mathcal C \subseteq \mathfrak C\, $ bounded below (resp. positive) by 
$$ \mathcal C\> \sim_1\> {\mathcal C}'\> \iff\> {\mathcal C}_c \> :=\> \bigl\{ b\in - \mathfrak C\bigm\vert b \leq c\, ,\, \forall c\in \mathcal C \bigr\}\> =\> \bigl\{ b\in - \mathfrak C\bigm\vert b\leq c'\, ,\, \forall c'\in {\mathcal C}' \bigr\}\> =\> {\mathcal C}_c' \> . $$
Dividing by the equivalence relation $\, \sim_1\, $ yields a set $\,  {\mathfrak L}_1 ( \mathfrak C )\, $. As above one checks that addition of two elements as well as multiplication with positive scalars is well defined and that $\, {\sim }_1\, $ has cancellation, i.e. 
$$  \mathcal C\> +\> \mathcal E\> {\sim }_1\> \mathcal D\> +\> \mathcal E\quad\Longrightarrow\quad 
\mathcal C\> {\sim }_1\> \mathcal D \> . $$
Then also subtraction of two elements is defined by 
$$ \bigl[ \mathcal C \bigr]\> - \bigl[ \mathcal D \bigr]\> =\> \bigl[ \mathcal C \bigr]\> +\> 
\bigl[ - {\mathcal D}_c \bigr] $$
and one obtains the structure of a real ordered normed space by imposing the corresponding norms  and ordering as above which in addition is a monotone complete lattice. 
For $\,\mathfrak C = \mathfrak V\, $ one obtains the original definition.
\par\noindent
Another generalization is to arbitrary partially ordered sets $\, \bigl( \mathcal S , \geq \bigr) \, $ by which we understand ordered sets such that for each pair of elements $\, x\, ,\, x\in \mathcal S\, $ there exist elements $\, w\, ,\, z\in \mathcal S\, $ with $\, w\leq x\, ,\, y\leq z\, $. Given a partially ordered set $\, \mathcal S\, $ let 
$\, C^c ( \mathcal S )\, $ denote the set of complemented subsets $\, C \subseteq \mathcal S\, $ which are bounded below so there exists an element $\, b\in \mathcal S\, $ with $\, b \leq c\, $ for all $\, c\in \mathcal C\, $ and $\, \mathcal C\, $ is the upper complement of its lower complement 
$\, {\mathcal C}_c = \{ b\in \mathcal S\,\vert\, b \leq c\, ,\, \forall c\in \mathcal C \}\, $. The subset $\, C^c ( \mathcal S )\, $ is naturally ordered by inclusion, i.e. $\, \mathcal C \geq \mathcal D \iff \mathcal C \subseteq \mathcal D\, $. Also let 
$\, C_c ( \mathcal S )\, $ denote the set of complemented subsets which are bounded above ordered by inclusion, so that there is a natural order isomorphism $\, C_c ( \mathcal S ) \simeq C^c ( \mathcal S )\, $ with $\, \mathcal C\, $ corresponding to its complement $\, {\mathcal C}_c\, $. 
Note that an arbitrary intersection of complemented subsets bounded below is a complemented subset bounded below, and an arbitrary intersection of complemented subsets bounded above is a complemented subset bounded above if one agrees to letting the void set be both upper and lower complement of the whole set $\, \mathcal S\, $. Then the partially ordered set $\, {\mathcal L}_1 ( \mathcal S ) = C^c ( \mathcal S ) \simeq C_c ( \mathcal S )\, $ is a monotone complete lattice with respect to the operations 
$$  \mathcal C \wedge \mathcal D\> =\> \bigl( {\mathcal C}_c \cap {\mathcal D}_c \bigr)^c\> ,\quad  \mathcal C \vee \mathcal D\> =\>  \mathcal C\cap \mathcal D \> , $$
i.e. any subset of elements which is bounded below has an infimum and any subset bounded above has a supremum, and there is a natural order isomorphic inclusion $\, \mathcal S \subseteq {\mathcal L}_1 ( \mathcal S )\, $ by $\, s \mapsto \{ s {\}}^c \equiv \{ s {\}}_c\, $. Then 
$\, \mathcal S\, $ is called {\it comparable} iff there exist monotonous functions
$$ r_+ :\> \mathcal S\> \longrightarrow\> \mathbb R\> ,\quad r_- :\> \mathcal S\>\longrightarrow\> \mathbb R  $$
such that any subset $\, \mathcal C\subseteq \mathcal S\, $ with $\, r_- ( \mathcal C ) := 
\inf\, \{ r_- ( c )\,\vert\, c\in \mathcal C \} > - \infty\, $ is bounded below, and any subset 
$\, \mathcal B\subseteq \mathcal S\, $ with $\, r_+ ( \mathcal B ) := \sup\, \{ r_+ ( b )\,\vert\, b\in \mathcal B \} < \infty\, $ is bounded above. Then if $\, \mathcal S\, $ is comparable  the monotone complete lattice $\, {\mathcal L}_1 ( \mathcal S )\, $ is injective in the category of comparable partially ordered sets and $r_+$-bounded (resp. $r_-$-bounded) monotonous maps. 
As an example of a comparable partially ordered set consider any unital subset $\, \mathcal S \subseteq \mathfrak V \, $ in a regular unital real ordered Banach space meaning that $\, \alpha\, {\bf 1}\in \mathcal S\, $  for any $\, \alpha\in \mathbb R\, $. Define 
$$ r_+ \bigl( x \bigr)\> =\> \inf\, \bigl\{ \alpha \in \mathbb R\bigm\vert \alpha\, {\bf 1}\geq x \bigr\}\> ,\quad r_- \bigl( x \bigr)\> =\> \sup\, \bigl\{ \beta\in \mathbb R\bigm\vert \beta\, {\bf 1} \leq x \bigr\}\> . $$ 
\par\bigskip\noindent
We now consider a second equivalence relation by which subsets of positive elements generate a linear lattice. In contrast with the previous construction this one is functorial for arbitrary monotonous (and positive linear) maps.
The construction can be done for an arbitrary ordered Banach space as above, but since we want to impose some additional structure we assume from the start that the base space is a (unital) $C^*$-algebra $\, A\, $. We only note that if $\,\mathfrak V\, $ is an arbitrary ordered real Banach space the positive linear map sending $\, \mathfrak V\, $ onto its canonical image in the regular Banach lattice $\,  \mathfrak L ( \mathfrak V )\, $ is order isomorphic if and only if $\,\mathfrak V\, $ is dually ordered by an argument as above. Thus the construction assigns to an arbitrary ordered real Banach space a canonically defined dually ordered Banach space together with a surjective positive linear map relating the two which is an order isomorphism only if $\,\mathfrak V\, $ is dually ordered itself.
Assume now that $\, \mathfrak V = A^{sa}\, $ is the real subspace of selfadjoint elements of a unital 
$C^*$-algebra $\, A\, $.
Much of the alternative construction parallels the preceding one. To each subset 
$\, \mathcal C \subseteq A_+\, $ of positive elements one attributes an element 
$\, \underline{\mathcal C}\, $ in a real vector space 
$\, \mathfrak L ( A )\, $ called the {\it infimum of $\,\mathcal C\, $}. 
There is a simple equivalence relation 
$\, \mathcal C \sim  {\mathcal C}'\, $ if and only if 
each element $\,\, c\in \mathcal C\, $ is larger or equal to some element in the closed convex hull of 
$\, {\mathcal C}'\, $ and vice versa.  
For each set $\, \mathcal C\, $ let $\, {\mathcal C}^c\, $ denote the set of positive elements each of which is larger than some element in the closed convex hull of $\, \mathcal C\, $.
Then $\, {\mathcal C}^c\, $ is the maximal representative in the simple equivalence class 
$\,  \underline{\mathcal C}\, $ of the set $\,\mathcal C\, $.
The equivalence class of $\, \mathcal C\, $ is then denoted 
$\, \underline{\mathcal C}\, $ and sometimes abusively identified with the maximal representative $\, {\mathcal C}^c\, $.
A preorder on the set of equivalence classes is obtained by the corresponding onesided relation, i.e. 
$\, \underline{\mathcal C} \geq  
\underline{\mathcal D}\, $ iff $\, {\mathcal C}^c \subseteq {\mathcal D}^c\, $. Then one may define the upper complement of the virtual positive element 
$\, \mathcal P = \underline{\mathcal C} - \underline{\mathcal D} \geq 0\, $ to consist of all elements 
$\, {\mathcal P}^c = \{ a\in A_+\,\vert\, (  a + \mathcal D )^c\subseteq {\mathcal C}^c \}\, $ which are larger or equal than  $\, \mathcal P\, $. This coincides for elements defined by single positive sets with the maximal representative set whence the notation but note that the positive element defined by the upper complement may be strictly larger then the difference $\, \underline{\mathcal C} - \underline{\mathcal D}\, $ in $\, \mathfrak L ( A )\, $ as defined below. 
Given two elements 
$\, \underline{\mathcal C} = \inf\, \{ c_{\lambda } {\}}_{\lambda }\, $ and 
$\, \underline{\mathcal D} = \inf\, \{ d_{\mu } {\}}_{\mu }\, $ as above their {\it sum} is defined to be 
$\, \underline{\mathcal C} + \underline{\mathcal D} = 
\inf\, \{ c_{\lambda } + d_{\mu }\,\vert\, c_{\lambda }\in {\mathcal C}\, ,\, d_{\mu }\in {\mathcal D} {\}}_{\lambda , \mu }\, $ and one checks that this definition is compatible with the equivalence relation hence well defined.  One also has  
$$ 2\, \mathcal C\> =\> \{ 2 c_{\lambda }\,\vert\, c_{\lambda }\in\mathcal C \}\> \sim\> 
\{ c_{\lambda } + c_{\mu }\,\vert\, c_{\lambda }\, ,\, c_{\mu }\in \mathcal C \}\> =\> 
\mathcal C\> +\> \mathcal C   $$
since each $\, c_{\lambda } + c_{\mu }\, $ is a convex combination of the elements $\, 2 c_{\lambda}\, ,\, 2 c_{\mu }\, $. Therefore as in the preceding construction the scalar multiplication of elements 
$\, \underline{\mathcal C}\, $ with positive real numbers is well defined. 
The convex cone of such elements generates a 
real linear space $\, \mathfrak L ( A )\, $ (the same notation will be used for the norm completion with respect to the norm defined below) by taking arbitrary differences 
$\, \mathcal A = \underline{\mathcal C} - \underline{\mathcal D}\, $ (in general the notation 
$\, \underline{\mathcal C}\, $ is reserved for elements defined by a specific representative positive set 
$\, \mathcal C\, $, it may happen in certain instances that $\,\mathcal C\, $ is also used to denote general selfadjoint elements in $\, \mathfrak L ( A )\, $ which meaning should be clear from the context), subject to the equivalence relation 
$$ ( \underline{\mathcal C} + \underline{\mathcal E} ) - ( \underline{\mathcal D} + \underline{\mathcal E} ) \>\simeq\> \underline{\mathcal C} - \underline{\mathcal D} $$
with $\, \mathcal E\subseteq A_+\, $ an arbitrary subset of positive elements, and an order on this space is obtained by 
$$ \mathcal A \geq \mathcal B\> \iff\> \mathcal A - \mathcal B \simeq 
\underline{\mathcal C} - \underline{\mathcal D} \geq 0  $$
with $\, {\mathcal D}^c \subseteq {\mathcal C}^c\, $.
For elements defined by positive sets the equivalence relation $\, \simeq\, $ coincides with 
$\,\sim\, $ so this equivalence relation has cancellation. To see this let 
$\, \mathcal D\, ,\, \mathcal E\, $ be given such that 
$\, \mathcal C + \mathcal E \sim \mathcal D + \mathcal E\, $. One can assume without loss of generality that $\, \mathcal C\, ,\,\mathcal D\, $ and $\, \mathcal E\, $ are maximal for $\,\sim\, $ so they are equal to their closed convex hull and the same is true for the sets
$\, \mathcal C + \mathcal E\, ,\, \mathcal D + \mathcal E\, $. We must show that 
$\, \mathcal D \subseteq \mathcal C\, $. Let $\, d\in \mathcal D\, ,\, e\in \mathcal E\, $ and 
$\,\epsilon > 0\, $  be given. 
There exist elements 
$\, c\in \mathcal C\, ,\, e'\in \mathcal E\, $ satisfying 
$$  c + e' \> \leq\> d + e + \epsilon 1\> . $$ 
Then again there exist $\, c'\in \mathcal C\, ,\, e''\in \mathcal E\, $ with 
$\, c' + e''\> \leq\> d + e' + \epsilon 1\, $ and inductively sequences
$\, ( c^{(k)} )_k\subseteq\mathcal C\, ,\, ( e^{(k)} )_k\subseteq\mathcal E\, $ with 
$$ c^{(k)} + e^{(k+1)} \> \leq\> d + e^{(k)} + \epsilon 1 \> . $$
For given $\, N\in\mathbb N\, $ one gets on adding these relations
$$ \sum_{k=1}^N\, {1\over N}\, ( c^{(k)} + e^{(k+1)} )\> \leq\> 
d\> +\> \sum_{k=1}^N\, {1\over N}\, e^{(k)} \> +\> \epsilon 1 $$
which implies 
$$ \sum_{k=1}^N\, {1\over N}\, c^{(k)}\> +\> {1\over N} e^{(N+1)}\>  \leq\> 
d\> +\> {1\over N}\, e\> +\> \epsilon 1\> . $$ 
Define the positive set $\, {\mathcal E}'\, $ to be the union of $\, \mathcal E\, $ and all elements 
of the form $\, \{ {1\over k} e^{(k+1)} \} \, $ with respect to all elements $\, d\, ,\, e\, $ and 
$\, \epsilon > 0\, $ as above.  
One concludes on letting $\, N\to \infty\, $  that 
$\, \mathcal D + \epsilon 1\subseteq \mathcal C + {\mathcal E}' \subseteq \mathcal C\, $ and by symmetry one has 
$\, \mathcal C + \epsilon 1\subseteq \mathcal D\, $. Since $\, \epsilon > 0\, $ was arbitrary one concludes $\, \mathcal C \sim \mathcal D\, $.
For an arbitrary set $\, \{ \underline{\mathcal C}_{\lambda } {\}}_{\lambda }\, $ of basic positive elements its infimum $\, \inf_{\lambda }\, \underline{\mathcal C}_{\lambda }\, $ which is the largest element smaller or equal than each $\, \underline{\mathcal C}_{\lambda }\, $ is well defined and is given by the element corresponding to the union $\, \cup_{\lambda }\, {\mathcal C}_{\lambda }\, $. 
By shifting with a fixed positive element $\, \underline{\mathcal C}\, $ one may then define an infimum for  subsets of general (selfadjoint) elements $\, \{ {\mathcal A}_{\lambda } {\}}_{\lambda }\, $ which are bounded below and are brought to the form above by the shift $\, {\mathcal A}_{\lambda } \mapsto {\mathcal A}_{\lambda } + \underline{\mathcal C}\, $. This applies in particular to finite subsets
$\, \{ {\mathcal A}_1 = \underline{\mathcal C}_1 - \underline{\mathcal D}_1\, ,\,\cdots\, ,\, {\mathcal A}_n = \underline{\mathcal C}_n - \underline{\mathcal D}_n \}\, $. By shifting each $\, {\mathcal A}_k\, $ with
$\, \sum_{k=1}^n\, \underline{\mathcal D}_k\, $, then taking the infimum of the corresponding basic positive elements and shifting back the result by $\, - \sum_{k=1}^n\, \underline{\mathcal D}_k\, $ one obtains the infimum of the set $\, \{ {\mathcal A}_k \}\, $.
In particular one has the lattice operations
$$  \mathcal A\wedge \mathcal B\> ,\quad 
\mathcal A \vee \mathcal B $$
denoting the (unique !) maximal element smaller or equal to $\, \mathcal A\, $ and $\,\mathcal B\, $, resp. the unique minimal element larger or equal to both $\,\mathcal A\, $ and $\,\mathcal B\, $. Taking 
$\, \mathcal B = 0\, $ the unique minimal positive decomposition of $\, \mathcal A\, $ is given by 
$$ \mathcal A\> =\> {\mathcal A}_+ - {\mathcal A}_-\> =\> ( \mathcal A \vee 0 ) + ( \mathcal A \wedge 0 )
\> .  $$ 
If $\, \mathcal A = \mathcal P - \mathcal Q\, $ is any positive decomposition of $\,\mathcal A\, $ then 
$\, \mathcal P \geq \mathcal A \vee 0\, $  with 
$\,\mathcal A = ( \mathcal A\vee 0 )\, -\, ( ( \mathcal A \vee 0 ) - \mathcal A )\, $ a positive decomposition of $\, \mathcal A\, $ which therefore is necessarily minimal. 
From uniqueness one gets $\, {\mathcal A}_+ =  \mathcal A \vee 0\, $ and hence 
$\, {\mathcal A}_- = - ( \mathcal A \wedge 0 )\, $. There also exists a minimal positive decomposition by {\it basic positive elements}, i.e. infima of positive sets.
Indeed, suppose given two different positive decompositions by basic positive elements
$$ \mathcal A\> =\> \underline{\mathcal C}\, -\, \underline{\mathcal D}\> =\> \underline{\mathcal C}' \, -\, 
\underline{\mathcal D}'  $$
we claim that $\, \mathcal A = ( \underline{\mathcal C}\wedge \underline{\mathcal C}'  )\, -\, 
( \underline{\mathcal D}\wedge \underline{\mathcal D}' )\, $. To see this it is sufficient by the symmetry 
$\, \mathcal A\mapsto - \mathcal A\, $ to prove 
$$  \underline{\mathcal C}\, -\, \underline{\mathcal D}\> =\> \mathcal A\> \leq\> ( \underline{\mathcal C} \wedge \underline{\mathcal C}' )\, 
-\, ( \underline{\mathcal D} \wedge \underline{\mathcal D}' ) $$
$$\qquad \iff \underline{\mathcal C} + ( \underline{\mathcal D}\wedge \underline{\mathcal D}' )\> \leq\> 
( \underline{\mathcal C} \wedge \underline{\mathcal C}' )\, +\, \underline{\mathcal D}  $$
using the relation $\, \mathcal C + {\mathcal D}' \sim {\mathcal C}' + \mathcal D\, $. Let 
$\, a = c + d \in ( \mathcal C \cup {\mathcal C}' ) + \mathcal D \, $ be given. If 
$\, c\in\mathcal C\, $ then $\, a\in \mathcal C + ( \mathcal D \cup {\mathcal D}' )\, $ follows trivially. If 
$\, c\in {\mathcal C}'\, $ there exists for each $\, \epsilon > 0\, $ a convex combination 
$\, b = \sum_k\, {\lambda }_k b_k\in \mathcal C + {\mathcal D}' \subseteq 
\mathcal C + ( \mathcal D \cup {\mathcal D}' )\, $ such that 
$\, a \geq b + \epsilon 1\, $ whence the result.
It follows that any selfadjoint element has a unique minimal decomposition as a difference of basic positive elements. In particular any positive element $\, \mathcal P = \underline{\mathcal C} - \underline{\mathcal D}\, $ with $\, \underline{\mathcal C} \geq \underline{\mathcal D}\, $ has a unique minimal representation by such elements. A (general) {\it basic} element of $\, \mathfrak L ( A )\, $ is an element which is the infimum of a given subset of selfadjoint elements of $\, A\, $ which is bounded below. Correspondingly, an element will be called {\it antibasic} iff it is the supremum of a given subset of selfadjoint elements bounded above. A positive antibasic element will sometimes be denoted 
$\, \overline{\mathcal P}\, $ to indicate that it is the supremum of its lower complement 
$\, {\mathcal P}_c = \{ a\in A^{sa}\,\vert\, a \leq \overline{\mathcal P} \}\, $. Then any element 
$$ \mathcal A\> =\> \underline{\mathcal C}\> -\> \underline{\mathcal D}\> =\> 
\bigl( R\> -\> \underline{\mathcal D} \bigr)\> -\> \bigl( R\> -\> \underline{\mathcal C} \bigr)\> =\> 
\overline{\mathcal P}\> -\> \overline{\mathcal Q} $$
can be decomposed as a difference of antibasic positive elements as well where $\, R\geq 0\, $ is a positive scalar exceeding both 
$\, \Vert \underline{\mathcal C}\Vert\, $ and $\, \Vert \underline{\mathcal D}\Vert\, $. 
For given $\,\mathcal P\geq 0\, $ consider the set of all positive elements $\, a\geq 0\, $ such that 
$\, \inf\, \{ a \} \geq \mathcal P\, $. This set $\, {\mathcal P}^c\, $ is a maximal set for $\,\sim\, $ called the upper complement of $\,\mathcal P\, $ which constitutes a basic positive element 
$\, {\mathcal P}^c \geq \mathcal P\, $ abusively identified with the corresponding maximal subset. By definition it is the smallest such element having this property and agrees with $\,\mathcal P\, $ if and only if $\, \mathcal P\, $ is basic itself.  By symmetry the supremum of the lower complement 
$\, {\mathcal P}_c = \{ a\in A^{sa}\,\vert\, a\leq \mathcal P \}\, $ constitutes an element $\, {\mathcal P}_c \leq \mathcal P\, $ of $\,\mathfrak L ( A )\, $ which is the largest antibasic element with this property. One may also consider the double lower complement $\, {\mathcal P}_{cc} = ( {\mathcal P}_c )^c\, $, which is basic, and the antibasic double upper complement respectively given by $\, {\mathcal P}^{cc} = ( {\mathcal P}^c )_c\, $. Then the triple upper and lower complement are defined recursively by 
$\, {\mathcal P}^{ccc} = ( {\mathcal P}^{cc} )^c = ( {\mathcal P}^c )_{cc}\, $ and 
$\, {\mathcal P}_{ccc} = ( {\mathcal P}_{cc} )_c = ( {\mathcal P}_c )^{cc}\, $ respectively. One has the relations 
$$ {\mathcal P}_c\>\leq\> \mathcal P\>\leq\> {\mathcal P}^c\> ,\quad {\mathcal P}_c\>\leq\> {\mathcal P}_{ccc}\> =\> {\mathcal P}^{cc}\>\leq\> {\mathcal P}^{ccc}\> =\> {\mathcal P}_{cc}\>
\leq\> {\mathcal P}^c  \> ,\quad {\mathcal P}^{cccc}\> =\> {\mathcal P}^{cc}\> ,\quad {\mathcal P}_{cccc}\> =\> {\mathcal P}_{cc} $$
but $\, \mathcal P\, $ need not to be related by order to neither of the elements $\, {\mathcal P}^{cc}\, ,\, {\mathcal P}_{cc}\, ,\, {\mathcal P}^{ccc}\, $ and $\, {\mathcal P}_{ccc}\, $. If $\,\mathcal P\, $ is basic 
then taking lower complements stabilizes already at the first step, i.e. $\, {\mathcal P}_{ccc} = {\mathcal P}_c\, $ and similarly taking upper complements stabilizes at the first step for antibasic elements.  For $\, \underline{\mathcal C} \geq 0\, $ a natural norm (satisfying the triangle inequality) is given by 
$$ \Vert \underline{\mathcal C} \Vert\> =\> \inf\, \{ \Vert d\Vert\,\vert\,  d\in \mathcal D\, ,\, \mathcal D \sim \mathcal C \}\> =\> \inf\, \{ \Vert c\Vert\,\vert\, c\in {\mathcal C}^c \}\> . $$
This norm is extended to general positive elements by 
$$ \Vert \mathcal P \Vert\> =\> \Vert {\mathcal P}^c \Vert $$
and to general selfadjoint elements  
$\, \mathcal A = \underline{\mathcal C} - \underline{\mathcal D}\, $ by defining 
$$ \Vert \mathcal A \Vert = \max\, \bigl\{ \Vert {\mathcal A}_+\Vert\, ,\, 
\Vert {\mathcal A}_- \Vert \bigr\} \geq 0 \> . $$
The triangle inequality is readily checked since 
$\,  ( \mathcal A + \mathcal B )_{\pm } \leq  {\mathcal A }_{\pm } + {\mathcal B}_{\pm }\, $.
If $\, \mathfrak a = \inf\, \{ a \}\, ,\, a\in A^{sa}\, $ is the image of the selfadjoint element $\, a = a_+ - a_-\, ,\, a_+a_- = 0\, $ in $\, \mathfrak L ( A )\, $ then $\, {\mathfrak a}_+ = \underline{\mathfrak a}_+ - \underline{\mathfrak a}_+ \wedge \underline{\mathfrak a}_-\, $ where 
$\, \underline{\mathfrak a}_+\, $ is the infimum of the set of all 
positive elements majorizing $\, a\, $, and $\, \underline{\mathfrak a}_-\, $ is the infimum of all positive elements majorizing $\, -a\, $, both of which are obviously maximal for the equivalence relation 
$\, \sim\, $. 
One easily checks that $\, b\geq a\, ,\, b\geq 0\, $ implies $\, \Vert b \Vert \geq 
\Vert a_+ \Vert\, $ and $\, b\geq -a\, ,\, b\geq 0\, $ implies $\, \Vert b \Vert \geq \Vert a_-\Vert\, $ 
(for example by extending a positive linear functional $\,\rho\, $ of $\, C^* ( a , 1 )\, $ satisfying 
$\, \rho ( a ) = \rho ( a_+ ) = \Vert a_+ \Vert\, $ to a positive linear functional of $\, A\, $ by the Hahn-Banach Theorem, 
then $\, \rho ( b ) \geq \Vert a_+ \Vert\, $ for each $\, b\in \underline{\mathfrak a}_+\, $). Thus the norm of 
$\, \mathfrak a\, $ is the same as the usual $C^*$-norm of $\, a\, $. Let us just remark that if instead considering a general ordered (say $\alpha $-normal) real Banach space the construction still yields a positive homomorphism $\, \mathfrak V\rightarrow \mathfrak L ( \mathfrak V )\, $ which is bounded from below on positive elements, but not necessarily bounded from below  or even injective on general elements. Thus if $\,\mathfrak V\, $ is not dually ordered there must exist for each $\,\epsilon > 0\, $ an element $\, x\in\mathfrak V\, $ such that every positive element which is larger than $\, x\, $ is also larger than $\, - x \, $ up to addition of a positive element of norm less or equal than $\,\epsilon\, $ since the minimal positive decomposition of the image $\,\mathfrak x\, $ in $\, \mathfrak L ( \mathfrak V )\, $ is given by 
$$ \mathfrak x\> =\> {\mathfrak x}_+\> -\> {\mathfrak x}_-\> ,\quad {\mathfrak x}_{\pm }\> =\> \underline{\mathfrak x}_{\pm }\> -\> \underline{\mathfrak x}_+\wedge \underline{\mathfrak x}_-  $$
where $\, \underline{\mathfrak x}_{\pm } = \inf\, \{ c\in {\mathfrak V}_+\,\vert\, c\geq \pm x \}\, $. 
We want to extend both the squaring operation as well as taking square roots of positive elements to 
$\, \mathfrak L ( A )\, $. It turns out that defining a squaring operation which is literally an extension of the square for single elements and has the properties to be expected is a rather strenuous adventure and appears to be unnatural. Instead we will give a stratified extension of the squaring operation which 
compensates for not being a proper extension by some other favourable aspects, e.g. being 
monotone.  For an element $\, \underline{\mathcal C}\, $ represented by a positive set define
a basic squaring operation by the formula
$$ \underline{\mathcal C}^2 \> =\> \inf\,\bigl\{ c^2\,\vert\, c\in {\mathcal C}^c \bigr\} \leqno{( 1 )} $$
to be the element represented by the set of squares of the maximal representative of 
$\,\underline{\mathcal C}\, $ and similarly 
$$ \sqrt{\underline{\mathcal C}}\> =\> \inf\, \bigl\{ \sqrt c\,\vert\, c\in {\mathcal C}^c \bigr\} \> . \leqno{( 2 )}  
$$
If $\, \mathfrak c = \inf\, \{ c \}\, ,\, c\geq 0\, $ is the image of a positive element of $\, A\, $ then 
$\,  {\mathfrak c}^2 \leq \inf\, \{ c^2 \}\, $ whereas 
$\, \sqrt{\mathfrak c} = \inf\, \{ \sqrt c \}\, $ by operator monotonicity of the square root. Of course for 
commutative $\  A\, $ both operations are proper extensions since in this case also the square of operators is monotone. Note that in general these operations need not be strictly inverse to each other, i.e. $\, ( \sqrt{ \underline{\mathcal C}} )^2\, $ may be strictly smaller than $\, \underline{\mathcal C}\, $. Yet one does have the onesided inversion formula
$$ \sqrt{ \underline{\mathcal C}^2}\> =\> \underline{\mathcal C} \> . \leqno{ ( 3 )} $$ 
This follows since the set of squares of the maximal representative of $\, \underline{\mathcal C}\, $ is again maximal for any convex combination of elements in this set is larger or equal to some element in the set. Indeed, for a given convex combination of two positive elements
$\, \lambda a + ( 1 - \lambda ) b\, ,\, 0\leq \lambda\leq 1\, $  one has the inequality
$$ ( \lambda a + ( 1 - \lambda ) b )^2 \leq \lambda a^2 + ( 1 - \lambda ) b^2 \> , \leqno{( 4 )} $$
and taking square roots of this set will give back the original element $\, \underline{\mathcal C}\, $ by operator monotonicity of the square root.
To extend both definitions to general positive elements one is guided by consideration of certain convexity (concavity) properties in case of single operators. Any reasonable definition of the square root for general positive elements should be monotone, in particular should satisfy 
$\, \sqrt{\mathcal P} \leq \sqrt{ {\mathcal P}^c }\, $. 
The inverse of inequality $\, ( 4 )\, $ for convex combinations of squareroots for two positive elements 
$\, a\, ,\, b \geq 0\, $ reads 
$$  \lambda\> \sqrt a\> +\> ( 1 - \lambda )\> \sqrt b\> \leq\> \sqrt{ \lambda\, a\, +\, ( 1 - \lambda )\, b} \> . \leqno{( 5 )} $$
Both inequalites $\, ( 4 )\, ,\, ( 5 )\, $ generalize to  $\, \mathfrak L ( A )\, $ on replacing $\, a\, ,\, b\, $ by 
{\it basic} positive elements $\, \underline{\mathcal C}\, ,\, \underline{\mathcal D}\, $, 
$$ ( \lambda \underline{\mathcal C}\> +\> ( 1 - \lambda ) \underline{\mathcal D} )^2\> \leq\> 
\lambda\, \underline{\mathcal C}^2\> +\> ( 1 - \lambda )\, \underline{\mathcal D}^2\> . \leqno{( 4' )} $$ 
$$ \lambda\, \sqrt{\underline{\mathcal C}}\> +\> ( 1 - \lambda )\, \sqrt{\underline{\mathcal D}}\> \leq\> 
\sqrt{\lambda\, \underline{\mathcal C}\> +\> ( 1 - \lambda )\, \underline{\mathcal D} }\> .\leqno{( 5' )} $$
Indeed, if $\, p\, $ is an element of the upper complement on the right side of $\, ( 5' )\, $ it can be represented in the form $\, \sqrt{ \lambda\, c + ( 1 - \lambda )\, d}\, $ with $\, c\in {\mathcal C}^c\, ,\, d\in {\mathcal D}^c\, $ up to considering larger elements or convex combinations. $\, ( 4' )\, $ follows much in the same way.
Assuming $\, a + d = c\, $ the formula $\, ( 5 )\, $ applied to $\, a\, $ and $\, d\, $ replacing $\, b\, $ renders 
$$ \sqrt a\> \leq\> {1\over \sqrt{\lambda }}\> \bigl( \sqrt c\>
 -\> \sqrt{ 1 - \lambda } \sqrt d \> \bigr) \leqno{( 6 )} $$
for every $\, 0 < \lambda \leq 1\, $. For $\, \mathcal P \geq 0\, $ one defines 
$$ \sqrt{\mathcal P}\> =\> \inf_{\lambda\, ,\, \underline{\mathcal C}\, ,\, \underline{\mathcal D}}\> 
\left\{\, {1\over\sqrt{\lambda }}\> \bigl(\, \sqrt{\underline{\mathcal C}}\> -\> \sqrt{ 1 - \lambda }\, 
\sqrt{\underline{\mathcal D} } \bigr)\>\bigm\vert\> 0 < \lambda \leq 1\> ,\> \mathcal P\, \leq\, \underline{\mathcal C}\, -\, \underline{\mathcal D}\> \right\} \> . \leqno{ ( 7 )} $$
From this definition it is immediate that the square root is monotone, i.e. $\, \mathcal P \leq \mathcal Q\, $ implies $\, \sqrt{\mathcal P} \leq \sqrt{\mathcal Q}\, $. There is only one problem that the infimum may not be well defined in $\, \mathfrak L ( A )\, $ which contrary to $\, {\mathfrak L}_1 ( A )\, $ is not a complete lattice. The solution of this problem is to consider the following iterated lattice construction 
$\, {\mathfrak L}_q ( A ) = {\mathfrak L}_1 ( \mathfrak L ( A ) )\, $. The canonical isometric embedding 
$\, \mathfrak L ( A ) \hookrightarrow {\mathfrak L}_q ( A )\, $ is supplemented by a (nonunique) extension $\, {\pi }_q : {\mathfrak L}_q ( A ) \twoheadrightarrow {\mathfrak L}_1 ( A )\, $ of the canonical surjection $\, \pi : \mathfrak L ( A ) \twoheadrightarrow {\mathfrak L}_1 ( A )\, $ which is a $*$-homomorphism of the underlying commutative $C^*$-algebras (see below). The Banach completion of $\, \mathfrak L ( A )\, $ is naturally a $C^*$-subalgebra of $\, {\mathfrak L}_q ( A )\, $ since $\, \mathfrak L ( A )\subseteq {\mathfrak L}_q ( A )\, $ is a sublattice (compare the Remark after the Corollary of Theorem 2) and we will usually not distinguish between $\, \mathfrak L ( A )\, $ as defined above and its completion since practically all interesting properties shared by $\, \mathfrak L ( A )\, $ extend to its completion from continuity.  
Then the squareroot of $\, \mathcal P\, $ is well defined as an element of $\, {\mathfrak L}_q ( A )\, $ and we may extend the definition to arbitrary positive elements of $\, {\mathfrak L}_q ( A )\, $ represented by a subset 
$\,  \mathfrak P = \inf\, \{ {\mathcal P}_{\lambda } {\}}_{\lambda\in\Lambda }\, $ by 
$$ \sqrt{\mathfrak P} = \sup\, \bigl\{ \sqrt{\mathcal Q}\,\bigm\vert\,  0\leq \mathcal Q\leq {\mathcal P}_{\lambda }\, ,\, \forall \lambda\in\Lambda  \bigr\} \> . \leqno{( 8 )} $$ 
From monotonicity of the square root and the fact that $\, \mathfrak L ( A )\, $ is a lattice already (so that the lower complement of a positive subset can be replaced by the positive part of the lower complement in the notion of equivalence) this is well defined irrespective of the representative set and monotone.
Check that the definition $\, ( 7 )\, $ coincides with the previous definition for basic positive elements. It is then sufficient to verify 
$$ \sqrt{\underline{\mathcal P}}\> \leq\> {1\over \sqrt{\lambda }}\, \bigl( \sqrt{\underline{\mathcal C}}\> -\> 
\sqrt{ 1 - \lambda }\, \sqrt{\underline{\mathcal D}} \bigr) $$
whenever $\, \underline{\mathcal P} \leq \underline{\mathcal C} - \underline{\mathcal D}\, $ or 
$$ \sqrt{\underline{\mathcal P}}\> +\> \sqrt{ 1 - \lambda\over \lambda }\, \sqrt{\underline{\mathcal D}}\> \leq\> {1 \over \sqrt{\lambda }}\, \sqrt{\underline{\mathcal C}} \> . $$
Rewriting the left side as
$$ \lambda\, {\sqrt{\underline{\mathcal P}}\over\lambda }\> +\> ( 1 - \lambda )\, 
{\sqrt{\underline{\mathcal D}}\over \sqrt{\lambda }\, \sqrt{1 - \lambda } } $$ 
and applying $\, ( 5' )\, $ with $\, \underline{\mathcal C}\, $ replaced by $\, {\lambda }^{-1}\, \underline{\mathcal P}\, $ and $\, \underline{\mathcal D}\, $ replaced by 
$\, ( 1 - \lambda )^{-1}\, \underline{\mathcal D}\, $ gives the result. 
In complete analogy with the case of the square root one defines the square in general by the formula 
$$ {\mathcal P}^2\> =\> \sup_{\lambda\, ,\, \underline{\mathcal C}\, ,\, \underline{\mathcal D}}\> 
\left\{\, \lambda\, \bigl( \underline{\mathcal C}^2\> -\> {1\over 1 - \lambda }\, \underline{\mathcal D}^2\, \bigr)\,\bigm\vert\, 0\leq\lambda < 1\, ,\, \mathcal P \geq \underline{\mathcal C} - \underline{\mathcal D} \right\} \> . \leqno{ ( 9 )}$$
From property $\, ( 4' )\, $ one checks that this coincides with the previous definition for basic positive elements. Moreover the formula makes sense for arbitrary selfadjoint elements. For 
$\, \mathcal A\in \mathfrak L ( A )\, $ define 
$$ sqr ( \mathcal A )\> =\> \sup_{\lambda\, ,\, \underline{\mathcal C}\, ,\, \underline{\mathcal D}}\> 
\left\{\, \lambda\, \bigl( \underline{\mathcal C}^2\> -\> {1\over 1 - \lambda }\, \underline{\mathcal D}^2\, \bigr)\,\bigm\vert\, 0\leq\lambda < 1\, ,\, \mathcal A \geq \underline{\mathcal C} - \underline{\mathcal D} \right\} \> . \leqno{ ( 10 )}  $$
It is immediate that $\, sqr ( \mathcal A )\, $ is monotone and positive, and that 
$\, sqr ( \mathcal A ) = 0\, $ for $\, \mathcal A \leq 0\, $. Now put 
$$ {\mathcal A}^2\> =\> sqr ( \mathcal A )\> +\> sqr ( - \mathcal A )\> \leqno{ ( 11 )}$$
to extend the definition of the square to arbitrary selfadjoint elements.  
The definition of $\, sqr\, $ (and a forteriori the squaring operation for selfadjoint elements) extends to 
$\, {\mathfrak L}_q ( A )\, $ in a similar fashion representing a general element
$\, \mathfrak A\, $ as the supremum of its lower complement in $\, \mathfrak L ( A )\, $ and putting 
$\, sqr\,\mathfrak A = \sup\, \{ sqr\, \mathcal A\,\vert\, \mathcal A \leq \mathfrak A \}\, $.
The complete lattice $\, {\mathfrak L}_q ( A )\, $ endowed with the two operations of taking squares and squareroots (of positive elements) will be called 
the {\it Jordan lattice} associated with $\, A\, $. Of course the Jordan squaring operations $\, (1),\, (9),\, (10)\, $ can also be defined if $\, A\, $ is only a Jordan-$*$-algebra (JC-algebra), and the Jordan squareroot operations $\, (2),\, (7),\, (8)\, $ can even be defined for an  {\it inverse Jordan algebra}, i.e. for an operator subsystem $\, \mathfrak X\subseteq \mathcal B ( \mathcal H )\, $ whose positive part is closed under taking squareroots. 
For each element in $\, {\mathfrak L}_q ( A )\, $ there are defined upper (basic) and lower (antibasic) complements defined analogously as above such that one has the relations
$$ {\mathfrak A}_c\>\leq\> {\mathfrak A}_{cc}\> ,\quad {\mathfrak A}^{cc}\>\leq\> {\mathfrak A}^c\> ,\quad 
{\mathfrak A}_c\>\leq\> \mathfrak A\>\leq\> {\mathfrak A}^c  $$
but note that one no longer has the relation $\, {\mathfrak A}^{cc} \leq {\mathfrak A}_{cc}\, $ nor either of the relations $\, {\mathfrak A}_{cc} = {\mathfrak A}^{ccc}\, ,\, {\mathfrak A}^{cc} = {\mathfrak A}_{ccc}\, $ in this case. Examples for a convex positive map are given by the projections 
$\, uc\, ,\, uccc :\, {\mathfrak L}_q ( A ) \rightarrow \mathfrak L ( A )\, $ onto the subset of basic (resp.complemented basic) elements defined by
$$ \mathfrak A\>\mapsto\> {\mathfrak A}^c\> ,\quad \mathfrak A\>\mapsto\> {\mathfrak A}^{ccc} $$
and natural examples of concave maps are given  by the projections onto the subset of antibasic (resp. complemented antibasic) elements  
$\, lc\, ,\, lccc :\, {\mathfrak L}_q ( A ) \rightarrow \mathfrak L ( A )\, $  
$$ \mathfrak A\> \mapsto\> {\mathfrak A}_c\> ,\quad \mathfrak A\>\mapsto\> {\mathfrak A}_{ccc} \> . $$
The fact that $\, uccc\, $ is convex follows from the fact that it factors over the convex surjection 
$\, {\pi }^c : {\mathfrak L}_q ( A ) \twoheadrightarrow {\mathfrak L}_1 ( A )\, $ by composition with the (maximal) convex lift $\, cv : {\mathfrak L}_1 ( A ) \rightarrow \mathfrak L ( A )\, $ sending the equivalence class of a basic element to the basic element of its maximal representative, i.e. $\, cv \bigl( [ \underline{\mathcal C} ] \bigr) = \underline{\mathcal C}_{cc}\, $. The maps $\, ucc\, $ and $\, lcc\, $ cannot be expected to be either convex nor concave in general.
From the squaring operation one derives a Jordan type binary operation by the formula 
$$ \underline{\langle}\> \mathcal A\, ,\, \mathcal B \underline{\rangle }\> =\> {1\over 2}\>
\bigl[ ( \mathcal A + \mathcal B )^2 \> -\> {\mathcal A}^2\>  -\> 
{\mathcal B}^2 \bigr]\> .  \leqno{( 12 )} $$
In order that this operation should be distributive with respect to addition in either variable one must have
$$ ( \underline{\mathcal C} + \underline{\mathcal D} + \underline{\mathcal E} )^2\> -\> 
( \underline{\mathcal C} + \underline{\mathcal D} )^2\> -\> 
( \underline{\mathcal C} + \underline{\mathcal E} )^2\> -\> 
( \underline{\mathcal D} + \underline{\mathcal E} )^2\> +\> 
\underline{\mathcal C}^2\>  +\> \underline{\mathcal D}^2\>  +\> \underline{\mathcal E}^2\> \bumpeq\> 0  
\leqno{( * )} $$
for every triple $\, ( \underline{\mathcal C}\, ,\, \underline{\mathcal D}\, ,\, \underline{\mathcal E} )\, $ of basic positive elements.  The expression on the left side of $\, ( * )\, $ is a {\it quadratic cloud} associated with the given triple and a main goal in the following is to find suitable conditions and relations making the quadratic cloud disappear. An important special case when this happens is if  $\, A\, $ is commutative. One then checks that the product $\, (12)\, $ gives a distributive, commutative  and associative Banach algebra product on $\, \mathfrak L ( A )\, $ which extends the product of $\, A\, $ and can be shown to satisfy the $C^*$-condition $\, \Vert \underline\langle \mathcal A\, ,\, \mathcal A \underline\rangle \Vert\, =\, \Vert \mathcal A {\Vert }^2\, $. The details are left to the reader.
Consider the complexification 
$\, {\mathfrak L}_q ( A )_c = {\mathfrak L}_q ( A ) {\otimes }_{\mathbb R} \mathbb C\, $. 
Although we have not yet imposed a norm on the complexification we may consider the original real subspace $\, {\mathfrak L}_q ( A )\, $ as consisting of selfadjoint elements with respect to the antilinear involution 
$\, \mathcal X = \mathcal A + i\, \mathcal B \mapsto \mathcal A - i\, \mathcal B = {\mathcal X}^*\, $ for
$\, \mathcal A\, ,\,\mathcal B\in \mathfrak L ( A )\, $. It is then convenient to define
$$ \Vert \mathcal X \Vert\> =\> \sup\, \bigl\{ \Vert \omega \mathcal X\, +\, \overline{\omega } {\mathcal X}^*\Vert\,\bigm\vert\, \omega\in\mathbb C\, ,\, \vert\omega\vert = 1 \bigr\}\> . $$
As a preamble to the construction of an associative product (on a certain quotient space of 
$\, {\mathfrak L}_q ( A )_c\, $) based on the Jordan bracket consider the following extension of the basic squaring operation on the subset of elements of the form 
$\,\mathcal X = \underline{\mathcal C} + i \underline{\mathcal D}\, $ with 
$\, \underline{\mathcal C}\, ,\, \underline{\mathcal D} \geq 0\, $. One defines  
$$ ( \underline{\mathcal C} + i \underline{\mathcal D} )\> \underline{\cdot }\> 
( \underline{\mathcal C} - i \underline{\mathcal D} )
\> =\> \inf\, \bigl\{ c_{\mu , \nu } c_{\mu , \nu }^*\,\vert\, c_{\mu , \nu } = c_{\mu } + i\, d_{\nu }\, ,\, c_{\mu }\in {\mathcal C}^c\, ,\, d_{\nu }\in {\mathcal D}^c \bigr\} \leqno{( 13 )} $$
which operation is well defined by uniqueness of the maximal representatives 
$\, {\mathcal C}^c\, ,\, {\mathcal D}^c\, $ and coincides with the basic squaring operation in case of a positive element. Note that the operation $\, ( 13 )\, $ is monotone in the sense that 
$\, ( \underline{\mathcal C} + i \underline{\mathcal D} )\,\underline\cdot\, ( \underline{\mathcal C} - i \underline{\mathcal D} )\,\geq\, ( \underline{\mathcal C}' + i \underline{\mathcal D}' )\,\underline\cdot\, ( \underline{\mathcal C}' - i \underline{\mathcal D}' )\, $ whenever 
$\, \underline{\mathcal C} \geq \underline{\mathcal C}'\, $ and $\, \underline{\mathcal D}\geq \underline{\mathcal D}'\, $.
From $\, ( 13 )\, $ one may derive another type of binary operation (of two positive elements say) which corresponds to the Lie bracket of two operators. 
The Lie type operation is given by 
$$ \lfloor\> \underline{\mathcal C}\, ,\, \underline{\mathcal D}\> \rfloor\> =\> {i\over 2}\>\bigl[\,
\bigl( \underline{\mathcal C} + i \underline{\mathcal D} \bigr)\, \underline{\cdot }\, 
\bigl( \underline{\mathcal C} - i \underline{\mathcal D} \bigr)\, -\, \underline{\mathcal C}\> \underline{\cdot }\> \underline{\mathcal C}\, -\, \underline{\mathcal D}\> \underline{\cdot }\> \underline{\mathcal D}\,\bigr] \> . \leqno{( 14 )} $$
In order that this operation should be distributive with respect to addition in the second variable one must have
$$ \bigl(\> \underline{\mathcal C} + i\, ( \underline{\mathcal D} + \underline{\mathcal E} )\> \bigr)\> \underline{\cdot }\> 
\bigl(\> \underline{\mathcal C} - i ( \underline{\mathcal D} + \underline{\mathcal E} )\> \bigr)\> -\> 
(\> \underline{\mathcal C} + i \underline{\mathcal D}\> )\> \underline{\cdot }\> 
(\> \underline{\mathcal C} - i \underline{\mathcal D}\> )\quad \leqno{(**)} $$
$$\qquad -\> (\> \underline{\mathcal C} + i \underline{\mathcal E}\> )\> \underline{\cdot }\> (\> \underline{\mathcal C} - i \underline{\mathcal E }\> )\> -\> 
(\> \underline{\mathcal D} + \underline{\mathcal E}\> )^2
\> +\> \underline{\mathcal C}^2\> +\> \underline{\mathcal D}^2\> +\> 
\underline{\mathcal E}^2\qquad \bumpeq\quad 0  $$
and similarly for the first variable (which case can be reduced to the one considered above by multiplying with the complex number $i$). The expression on the left side of $\, (**)\, $ is a 
{\it quadratic cloud of second type} associated with the given triple of positive elements. 
Assuming that all quadratic clouds of first and second type disappear under suitable restrictions and relations one may define a general (bilinear) product by linear extension of 
$$ \underline{\mathcal C}\> \underline\cdot\> \underline{\mathcal D}\> =\> 
{1\over 2}\, \underline{\langle}\, \underline{\mathcal C}\, ,\,  \underline{\mathcal D}\, \underline{\rangle }\> +\> {1\over 2}\, \lfloor\, \underline{\mathcal C}\, ,\, \underline{\mathcal D}\, \rfloor \> .  \leqno{( 15 )} $$
The relations needed to get associativity will be referred to as a {\it quadratic cloud of third type} associated with a given triple of elements and determining their precise form is left to the reader.
\par\bigskip\noindent
{\it Remark.}\quad As with the ${\mathfrak L}_1$-functor the $\mathfrak L$-construction admits various generalizations. One is to 
unital affine subcones $\, \mathfrak C\subseteq\mathfrak V\, $ in a regular unital real ordered Banach space $\, \mathfrak V\, $. Define $\, \mathfrak L \bigl( \mathfrak C \bigr) \subseteq \mathfrak L \bigl( \mathfrak V \bigr)\, $ to be the subspace which is linearly generated by (positive) basic elements 
$\, \{ \underline{\mathcal C} \}\, $ with $\, \mathcal C\subseteq \mathfrak C\, $ so that an arbitrary element of $\, \mathfrak L \bigl( \mathfrak C \bigr)\, $ can be written as a difference of two such elements or equivalently as a difference of two positive antibasic elements which are suprema of some given subsets in $\, - \mathfrak C\, $ which are bounded above. The Proposition below then also applies to this more general situation.
\par\bigskip\noindent
{\bf Proposition 1.}\quad For any unital ordered Banach space $\, \mathfrak V\, $ the canonical surjection 
$$ \pi :\> \mathfrak L ( \mathfrak V )\> \twoheadrightarrow {\mathfrak L}_1 ( \mathfrak V ) $$
obtained by assigning with each positive basic element in $\, \mathfrak L ( \mathfrak V )\, $ its corresponding equivalence class in $\, {\mathfrak L}_1 ( \mathfrak V )\, $ is a $*$-homomorphism (lattice map) and is the unique monotonous extension of the identity map of $\, \mathfrak V\, $. It extends (nonuniquely) to a $*$-homomorphic surjection 
$$ {\pi }_q : {\mathfrak L}_q ( \mathfrak V )\> \twoheadrightarrow {\mathfrak L}_1 ( \mathfrak V ) $$
of the underlying (commutative) $C^*$-algebras. Moreover two canonical monotonous extensions of $\, \pi\, $ to $\, {\mathfrak L}_q ( \mathfrak V )\, $ are given by the maps 
$$ {\pi }_c :\> {\mathfrak L}_q ( \mathfrak V )\> \largerightarrow\> {\mathfrak L}_1 ( \mathfrak V )\> ,\qquad {\pi }_c ( \mathfrak A )\> =\> \pi \bigl( {\mathfrak A}_c \bigr)\> , $$
$$ {\pi }^c :\> {\mathfrak L}_q ( \mathfrak V )\> \largerightarrow\> {\mathfrak L}_1 ( \mathfrak V )\> ,\qquad {\pi }^c ( \mathfrak A )\> =\> \pi \bigl( {\mathfrak A}^c \bigr) \> . $$
Of these $\, {\pi }_c\, $ is the minimal monotonous extension of $\, \pi\, $, whereas $\, {\pi }^c\, $ is the maximal monotonous extension of $\, \pi\, $.
If $\, \phi : \mathfrak V \rightarrow \mathfrak W\, $ is any monotonous map of unital ordered Banach spaces there exists a canonical maximal, and a canonical minimal extension of $\, \phi\, $ to monotonous maps 
$$ cv ( \phi )\, ,\, cc ( \phi ) :\> {\mathfrak L}_1 \bigl( \mathfrak V \bigr)\> \largerightarrow\> {\mathfrak L}_1 \bigl( \mathfrak W \bigr) \> . $$
Of these $\, cv ( \phi )\, $ is convex and $\, cc ( \phi )\, $ is concave (for positive linear $\,\phi\, $). 
\par\bigskip\noindent
{\it Proof.}\quad There are two canonical monotonous  extensions of any monotonous map $\, \phi : \mathfrak V \rightarrow \mathfrak W\, $ to 
$$  cv ( \phi ) : {\mathfrak L}_1 ( \mathfrak V )\> \longrightarrow\> {\mathfrak L}_1 ( \mathfrak W ) \> , $$
$$  cc ( \phi ) : {\mathfrak L}_1 ( \mathfrak V )\>\longrightarrow\> {\mathfrak L}_1 ( \mathfrak W )\> .  $$
These are obtained by composing the functorial composition 
$$ \mathfrak L ( \mathfrak V ) \> \buildrel \mathfrak L ( \phi )\over \largerightarrow\> 
\mathfrak L ( \mathfrak W )\> \buildrel \pi\over\largerightarrow {\mathfrak L}_1 ( \mathfrak W ) $$
on the left with the two lifts 
$$ cv\, ,\, cc :\> {\mathfrak L}_1 \bigl( \mathfrak V \bigr)\> \largerightarrow\> \mathfrak L \bigl( \mathfrak V \bigr) $$
induced for elements $\, x\in {\mathfrak L}_1 ( \mathfrak V )\, $ by the assignment 
$$ x\> =\> [ \underline{\mathcal C} ]\quad\mapsto\quad\underline{\mathcal C}_{cc}\> =\> cv ( x )\> ,  $$
 $$ x\> =\> [ \underline{\mathcal C} ]\quad\mapsto\quad \underline{\mathcal C}_c\> =\> cc ( x ) $$
respectively where $\, \underline{\mathcal C}\, $ is any basic element representing $\, x\, $. 
Note that always $\, cc ( x )\>\leq\> cv ( x )\, $.
Obviously $\, cv ( \phi ) : {\mathfrak L}_1 ( \mathfrak V ) \twoheadrightarrow {\mathfrak L}_1 ( \mathfrak W )\, $ is the maximal monotonous extension of $\, \phi\, $ and is checked to be convex if $\, \phi\, $ is positive linear (or convex) whereas $\, cc ( \phi )\, $ is the minimal monotonous extension of $\,\phi\, $ and is checked to be concave whenever $\,\phi\, $ is positive linear (or concave). Writing $\, {\mathfrak L}_q ( \mathfrak V ) = 
{\mathfrak L}_1 ( \mathfrak L ( \mathfrak V ) )\, $ and $\, {\mathfrak L}_1 ( \mathfrak V ) \simeq 
{\mathfrak L}_1 ( {\mathfrak L}_1 ( \mathfrak V ) )\, $ one obtains $\, cv ( \pi ) = {\pi }^c\, $ since 
$\, cv ( \pi )\, $ is obviously the maximal monotonous extension of the identity of $\,\mathfrak V\, $ as any such extension is uniquely determined on $\, \mathfrak L ( \mathfrak V )\, $ and equal to $\,\pi\, $ from the argument below. On the other hand 
$\, cv ( \pi ) \leq {\pi }^c\, $ whence both maps must be equal. Similarly $\, cc ( \pi ) = {\pi }_c\, $.
The kernel of $\, \pi\, $ is positively generated by elements of the form 
$\, \bigl\{ {\mathcal A}^c - \mathcal A\, ,\, \mathcal A - {\mathcal A}_c \bigr\}\, $. This is because 
$$ {\mathcal A}_c \leq {\mathcal A}^{cc}\leq {\mathcal A}_{cc}\leq {\mathcal A}^c $$ 
for $\, \mathcal A \in \mathfrak L ( \mathfrak V )\, $ with $\, \pi ( {\mathcal A}_c ) = \pi ( {\mathcal A}_{cc} )\, ,\, \pi ( {\mathcal A}^{cc} ) = \pi ( {\mathcal A}^c )\, $. Hence the kernel must be a $*$-ideal and $\, \pi\, $ a surjective $*$-homomorphism which admits a $*$-homomorphic extension $\, {\pi }_q : {\mathfrak L}_q ( \mathfrak V ) \twoheadrightarrow {\mathfrak L}_1 ( \mathfrak V )\, $ from Theorem A of section 2.  It is also clear that any monotonous map extending the identity of $\, \mathfrak V\, $ must be smaller than $\, {\pi }^c\, $ and larger than $\, {\pi }_c\, $ while $\, {\pi }^c = \pi = {\pi }_c\, $ restricted to $\,\mathfrak L ( \mathfrak V )\, $ proving that $\,\pi\, $ is uniquely determined\qed
\par\bigskip\noindent
{\it Remark.} (i)\quad For any order isomorphic unital embedding of unital ordered Banach spaces 
$\, \lambda : \mathfrak V \hookrightarrow \mathfrak W\, $ there exists a canonical convex, monotonous and positively homogenous retraction 
$\, \underline r : {\mathfrak L}_q ( \mathfrak W ) \twoheadrightarrow {\mathfrak L}_q ( \mathfrak V )\, $ for the functorial map $\, \Lambda : \mathfrak L ( \mathfrak V ) \hookrightarrow \mathfrak L ( \mathfrak W )\, $ induced by restriction of a positive basic element $\, \underline{\mathcal C} \in \mathfrak L ( \mathfrak W )_+\, $ to the subset of elements $\, {\mathcal C}^r = \{ c\in {\mathfrak V}_+\,\vert\, 
\lambda ( c )\in  \mathcal C \}\, $, i.e. $\, \underline r ( \underline{\mathcal C} ) = \underline{\mathcal C}^r\, $. Clearly if 
$\, \underline{\mathcal C}^{\lambda } = \mathfrak L ( \lambda ) ( \underline{\mathcal C} )\, $ is the image of a positive basic element in $\, \mathfrak L ( \mathfrak V )\, $ then 
$\, \underline r ( \underline{\mathcal C}^{\lambda } ) = \underline{\mathcal C}\, $. Also it is easy to see that 
$\, ( \underline{\mathcal C} + \underline{\mathcal D} )^r \leq \underline{\mathcal C}^r + \underline{\mathcal D}^r\, $ and $\, \underline{\mathcal C}^r \leq \underline{\mathcal D}^r\, $ whenever 
$\, \underline{\mathcal C} \leq \underline{\mathcal D}\, $, so that $\, r\, $ is convex and monotonous on the convex subcone of positive basic elements. To extend the definition to more general elements define 
$$ \underline r ( \mathcal A )\> =\> \sup\, \bigl\{ \underline{\mathcal C}^r\, -\, \underline{\mathcal D}^r\,\bigm\vert\, \underline{\mathcal C}\, -\, \underline{\mathcal D} \leq \mathcal A \bigr\}\>\in {\mathfrak L}_q ( \mathfrak V ) \> . $$
From definition $\, \underline r\, $ is monotonous on $\,\mathfrak L ( \mathfrak W )\, $ (if well defined).
If $\, \underline{\mathcal C}^{\lambda }\, $ is an image of a basic element in 
$\, \mathfrak L ( \mathfrak V )\, $ and $\, \underline{\mathcal D}\, $ is an arbitrary basic element then 
$\, ( \underline{\mathcal C}^{\lambda } + \underline{\mathcal D} )^r =   
{\underline{\mathcal C}}^{\lambda , r} + \underline{\mathcal D}^r\, $ whence one has the relation 
$$ \underline{\mathcal C}^{\lambda , r}\> -\> \underline{\mathcal D}^{\lambda , r}\>\geq\> 
\underline{\mathcal C}^r\> -\> \underline{\mathcal D}^r $$
whenever $\, \underline{\mathcal C}^{\lambda } - \underline{\mathcal D}^{\lambda } \geq \underline{\mathcal C} - \underline{\mathcal D}\, $ showing that $\, \underline r\, $ is a retraction for elements in the image of $\, \mathfrak L ( \mathfrak V )\, $. Also for a positive basic element 
$\, \underline{\mathcal P}\, $ one gets 
$\, \underline{\mathcal P}^r \geq \underline{\mathcal C}^r - \underline{\mathcal D}^r\, $ whenever 
$\, \underline{\mathcal P} \geq \underline{\mathcal C} - \underline{\mathcal D}\, $ showing for one thing that $\, \underline r\, $ is well defined since 
$\, \underline r ( \mathcal P ) \leq \underline r ( {\mathcal P}^c ) = ( {\mathcal P}^c )^r\, $ and that one recovers the original definition in case of a positive basic element. In fact it is easily seen that for any (even nonpositive) 
basic element one has $\, \underline r ( \underline{\mathcal A} ) = \underline{\mathcal A}^r\, $ where 
$\, {\mathcal A}^r = \mathcal A \cap \mathfrak V\, $. Then the condition of positivity of basic elements in the definition of the extension of $\,\underline r\, $ is seen to be superfluous and can be replaced by considering differences of arbitrary basic elements. 
Positive homogeneity of $\,\underline r\, $ is more than obvious, then we may define  
$$ \underline r ( \mathfrak A )\> =\> \sup\, \bigl\{ \underline r ( \mathcal A )\,\bigm\vert\, \mathcal A \leq \mathfrak A \bigr\} $$
in order to extend the definition of $\, \underline r\, $ to $\, {\mathfrak L}_q ( \mathfrak W )\, $. Combining this with the convex lift $\, cv : {\mathfrak L}_1 ( \mathfrak W ) \rightarrow \mathfrak L ( \mathfrak W )\, $ and the linear surjection $\, \pi : {\mathfrak L}_q ( \mathfrak V ) \twoheadrightarrow {\mathfrak L}_1 ( \mathfrak V )\, $ of Lemma 1 one obtains a convex monotonous retraction
$$  {\underline r}_1 : {\mathfrak L}_1 ( \mathfrak W )\buildrel cv\over\longrightarrow \mathfrak L ( \mathfrak W ) \buildrel \underline r\over\longrightarrow {\mathfrak L}_q ( \mathfrak V ) \buildrel\pi\over\longrightarrow {\mathfrak L}_1 ( \mathfrak V ) $$
for the functorial maps $\, cv ( \lambda )\, ,\, cc ( \lambda )\, $. Note also that $\,\underline r\, $ is generally magnifying, i.e. $\, ( \Lambda\circ \underline r ) ( \mathfrak A ) \geq \mathfrak A\, $. Namely any element 
can be represented as a supremum of basic elements $\, \mathfrak A = \sup_{\mu }\, \underline{\mathcal A}_{\mu }\, $. Then obviously $\, ( \Lambda \circ \underline r ) ( \underline{\mathcal A}_{\mu } ) \geq 
\underline{\mathcal A}_{\mu }\, $ for each index $\,\mu\, $ whence 
$$ ( \Lambda\circ \underline r ) ( \mathfrak A )\> =\> ( \Lambda\circ \underline r ) ( \sup\, \underline{\mathcal A}_{\mu } )\>\geq\> \sup\, \underline{\mathcal A}_{\mu }^{r , \lambda }\>\geq\> \mathfrak A\> . $$
This in turn shows that $\, \underline r\, $ is {\it increasing normal}, because 
$$ \underline r ( \sup\, {\mathfrak A}_{\mu } )\>\geq\> \sup\, \underline r ( {\mathfrak A}_{\mu } )\> =\> \underline r ( \sup\, ( \Lambda\circ \underline r ) ( {\mathfrak A}_{\mu } ) )\>\geq\> 
\underline r ( \sup\, {\mathfrak A}_{\mu } ) $$
hence equality in each instance.
\par\noindent
(ii)\quad Since any element $\, \mathcal A\in \mathfrak L ( \mathfrak W )\, $ can also be written as a difference of (positive or arbitrary) antibasic elements $\, \mathcal A = \overline{\mathcal P} - \overline{\mathcal Q}\, $ we may consider another functorial monotonous retraction $\, \overline r : {\mathfrak L}_q ( \mathfrak W ) \twoheadrightarrow {\mathfrak L}_q ( \mathfrak V )\, $ induced by restriction of antibasic elements, i.e. 
$$ \overline r ( \overline{\mathcal P} )\> =\> {\overline{\mathcal P}}^r\> =\> 
\sup\, \{ a\in \mathfrak V\,\vert\, \lambda ( a )\leq \overline{\mathcal P} \} \> . $$ 
By symmetry with the above argument the general assignment
$$ \overline r ( \mathfrak A )\> =\> \inf\, \bigl\{ {\overline{\mathcal P}}^r\> -\> {\overline{\mathcal Q}}^r\,\vert\, \overline{\mathcal P}\> -\> \overline{\mathcal Q} \geq \mathfrak A \bigr\} $$
defines a monotonous retraction $\, \overline r : {\mathfrak L}_q ( \mathfrak W ) \twoheadrightarrow {\mathfrak L}_q ( \mathfrak V )\, $ for $\,\Lambda\, $ which extends the above definition for antibasic positive elements. Moreover $\, \overline r\, $ is generally reducing, i.e., 
$\, ( \Lambda\circ \overline r ) ( \mathfrak A ) \leq \mathfrak A\, $ and {\it decreasing normal}, i.e. 
$\, \overline r ( \inf\, {\mathfrak A}_{\mu } )\, =\, \inf\, \overline r ( {\mathfrak A}_{\mu } )\, $. One obtains the general relation 
$$ \overline r ( \mathfrak A\, +\, \mathfrak B )\> =\> \inf\, \bigl\{ \overline{\mathcal P}^r\,\bigm\vert\, \overline{\mathcal P}\geq \mathfrak A\, +\, \mathfrak B \bigr\}\>\leq\> \inf\, \bigl\{ ( \overline{\mathcal P}\, -\, \underline{\mathcal A} )^r\, +\, \underline{\mathcal A}^r\,\bigm\vert\, \underline{\mathcal A}\,\leq\, \mathfrak A\, ,\, \overline{\mathcal P}\,\geq\, \mathfrak A\, +\, \mathfrak B \bigr\} $$
$$\quad \leq\> 
\sup\, \bigl\{ \underline{\mathcal A}^r\,\bigm\vert\, \underline{\mathcal A}\,\leq\, \mathfrak A \bigr\}\> +\> 
\inf\, \bigl\{ \overline{\mathcal Q}^r\,\bigm\vert\, \overline{\mathcal Q}\,\geq\,\mathfrak B \bigr\}
\> =\> \underline r ( \mathfrak A )\> +\> \overline r ( \mathfrak B ) \> .  $$
\par\bigskip\noindent 
In the following $\, A\, $ denotes a $C^*$-algebra as usual, but note that the same concepts can be applied to a regular unital ordered Banach space $\, \mathfrak V\, $ (isometrically order isomorphic with an operator system $\, \mathfrak X \subseteq \mathcal B ( \mathcal H )\, $).
Let $\, \underline{\mathfrak L} ( A )\subseteq \mathfrak L ( A )\, $ denote the subcone of basic elements 
which may be identified with their maximal representative set $\, {\mathcal A}^c\subseteq A^{sa}\, $, and is closed under addition and multiplication by positive scalars. Also let 
$\, \underline C ( \mathcal S ( A ) )\, $ be the subcone of upper semicontinuous functions on the state space of $\, A\, $ and $\, \widehat C ( \mathcal S ( A ) )\, $ its linear envelope. Define a positive linear map $\, s : \mathfrak L ( A )\, \rightarrow \widehat C ( \mathcal S ( A ) )\, $ by linear extension of the affine map 
$$ s\, :\> \underline{\mathfrak L} ( A )\> \longrightarrow\> \underline C ( \mathcal S ( A ) )\> ,\quad 
s ( {\mathcal A}^c ) ( \rho )\> =\> \inf\, \{ \rho ( a )\,\vert\, a\in {\mathcal A}^c \} \leqno{( 16 )} $$
where $\, \rho\in \mathcal S ( A )\, $ and check that it is well defined and {\it strictly monotone} in the sense that $\, \mathcal A \lneqq \mathcal B\, $ implies 
$\, s ( \mathcal A ) \lneqq s ( \mathcal B )\, $. This implies that $\, s\, $ is an order isomorphic injection. 
Indeed, suppose that $\, s ( \mathcal A ) \geq 0\, $ for some $\, \mathcal A = \underline{\mathcal C} - \underline{\mathcal D}\, $. Then $\, s ( \underline{\mathcal C} ) \geq s ( \underline{\mathcal D} )\, $. Choosing  an arbitrary $\, c\in {\mathcal C}^c\, $ and any state $\,\sigma\, $ one gets 
$\, s_{\sigma } (   \underline{{\mathcal D}^c \cup \{ c \} }  ) = s_{\sigma } ( \underline{\mathcal D} )\, $ which implies $\, c\in {\mathcal D}^c\, $ by strict monotonicity of $\, s\, $. Therefore 
$\, \mathcal C\geq \mathcal D\, $ and $\,\mathcal A \geq 0\, $.
To see that $\, s\, $ is strictly monotone it is sufficient to prove that whenever 
$\, \mathcal C \subseteq A_+\, $ is a maximal closed convex subset of positive elements (i.e. $\, x\in\mathcal C\, $ implies $\, y\in\mathcal C\, $ for $\, y\geq x\, $) and 
$\, d\geq 0\, $ is  a single positive element not contained in $\, \mathcal C\, $ then there exists a positive functional $\, \phi\in ( A^* )_+\, $ separating $\, \mathcal C\, $ and $\, \{ d \}\, $. To prove this choose 
a small convex open neighbourhood $\, {\mathcal C}_0 \supset \mathcal C\, $ such that 
$\, d\, $ is not contained in the closure of $\, {\mathcal C}_0\, $ and a point $\, c_0 \in\mathcal C\, $. Let 
$\, \mathfrak C\, =\, {\mathcal C}_0\, -\ c_0 = 
\{ c - c_0\,\vert\, c\in {\mathcal C}_0 \}\, $. Then $\,\mathfrak C\, $ is a convex open neighbourhood of 
$\, 0\, $ with associated Minkowski functional 
$$ m ( x ) = \inf\, \{ s > 0\,\vert\, s^{-1} x\in \mathfrak C \}  \> . $$
One checks that $\, m ( x ) = 0\, $ whenever $\, x\geq 0\, $ since $\, A_+ \subseteq \mathfrak C\, $. This implies that any functional dominated by $\, m\, $ must be negative. Put 
$\, z = d - c_0\notin\mathfrak C\, $ so that $\, m ( z ) > 1\, $. Define $\, {\psi }_0 ( s z ) = s\, $ for 
$\, s\in\mathbb R\, $ and extend $\, {\psi }_0\, $ to a (necessarily negative) functional $\,\psi\, $ on $\, A\, $ dominated by $\, m\, $ from the Hahn-Banach theorem. Then $\, \psi\, $ separates $\, \mathcal C\, $ and $\, \{ d \}\, $ and the same holds for the positive functional $\,  \phi = - \psi\, $ proving the statement. For a given state $\, \rho\in \mathcal S ( A )\, $ let 
$\, s_{\rho } : \mathfrak L ( A ) \twoheadrightarrow \mathbb C\, $ be the induced state on 
$\, \mathfrak L ( A )\, $, i.e. $\, s_{\rho } ( \mathcal A ) = s ( \mathcal A ) ( \rho )\, $. Then the kernel of 
$\, s_{\rho }\, $ is positively generated, since if $\, \underline{\mathcal C}\, ,\, \underline{\mathcal D}\, $ are given with $\, s_{\rho } ( \underline{\mathcal C} ) = s_{\rho } ( \underline{\mathcal D} ) = r\, $ then also 
$\, s_{\rho } ( \underline{\mathcal C} \wedge \underline{\mathcal D} ) = r\, $. 
Also consider the restriction of $\, s\, $ to the {\it pure} states 
$$ r\, : \mathfrak L ( A )\>\longrightarrow\> \widehat C ( \mathcal P ( A ) ) $$
where $\, \mathcal P ( A )\, $ denotes the set of pure states of $\, A\, $, and in case that 
$\, M\, $ is a von Neumann algebra consider the restriction to the {\it normal} states
$$ s^{\nu } : \mathfrak L ( M )\> \longrightarrow\> \widehat C ( {\mathcal S}^{\nu } ( M ) )\> , $$
resp. if $\, M = \mathcal B ( \mathcal H )\, $ the restriction to the {\it vector states} (= pure normal states) 
$$ r^{\nu } : \mathfrak L ( \mathcal B ( \mathcal H ) )\> \longrightarrow\> 
\widehat C ( {\mathcal P}^{\nu } ( \mathcal B ( \mathcal H ) ) )\>  . $$
The maps $\, r\, ,\, s^{\nu }\, $ and $\, r^{\nu }\, $ fail to be injective in general. 
We will sometimes also consider the subset 
$\, {\mathcal S}^{\nu }_f\, $ of (normal) states which are finite convex combinations of vector states. 
\par\bigskip\noindent
For the following constructions we consider concrete operator subsystems
$\, \mathfrak X \subseteq \mathcal B ( \mathcal H )\, $, i.e. a specific unital completely isometric representation is taken into account and we therefore consider the subset of states on 
$\, \mathfrak X\, $ which are restrictions of states in $\, {\mathcal S}^{\nu }_f\, $.  For a fixed normal state 
$\,\rho \in {\mathcal S}^{\nu }_f\, $ let 
$\, \underline{\mathfrak C}_{0 , \rho } \subseteq \underline{\mathfrak L} ( \mathfrak X )_+\, $ be the subcone generated by basic positive  elements which are restrictions $\, \underline{\mathcal E}^r = \underline r ( \underline{\mathcal E} )\, $ of basic positive elements 
$\, \underline{\mathcal E}\in \mathfrak L ( \mathcal B ( \mathcal H ) )\, $ such that
$\, s^{\nu }_{\rho } ( \underline{\mathcal E} ) = s^{\nu } ( \underline{\mathcal E} ) ( \rho ) = 0\, $ and 
$\, \overline{\mathfrak C}_{0 , \rho }\subseteq {\mathfrak L}_+ ( \mathfrak X )\, $ the positive subcone generated by images of basic positive elements $\, \underline{\mathcal E}\, $ as above under the {\it antibasic} restriction map $\, \overline r\, $ as in part (ii) of the Remark after Proposition 1. Then 
$\, \underline{\mathfrak C}_{0 , \rho }\, $ (resp. $\, \overline{\mathfrak C}_{0 , \rho }\, $) generates an order ideal $\, \underline{\mathfrak J}_{0 , \rho }\, $ (resp. $\, \overline{\mathfrak J}_{0 , \rho }\, $) in  
$\, \mathfrak L ( \mathfrak X )\, $ (by restriction of the corresponding order ideal in $\, {\mathfrak L}_q ( \mathfrak X )\, $) whose positive part consists of all positive elements 
$\, \mathfrak Q\in \mathfrak L ( \mathfrak X )\, $ which are dominated by some element 
$\, \underline{\mathcal E}\in \underline{\mathfrak C}_{0 , \rho }\, $ (resp. $\, \mathcal P\in \overline{\mathfrak C}_{0 , \rho }\, $). Since always $\, \overline r \leq \underline r\, $ one obtains 
$\, \overline{\mathfrak J}_{0 , \rho }\subseteq \underline{\mathfrak J}_{0 , \rho }\, $. The corresponding (unrestricted) order ideal of $\, {\mathfrak L}_q ( \mathcal B ( \mathcal H ) )\, $ is denoted $\, {\mathfrak J}_{0 , \rho}\, $. Put $\, {\mathfrak J}_0 = \bigcap_{\rho\in {\mathcal P}^{\nu } } {\mathfrak J}_{0 , \rho}\, $ and let 
$\, \underline{\mathfrak J}_0 = \underline r \bigl( {\mathfrak J}_0  \bigr) \subseteq \mathfrak L ( \mathfrak X )\, $ be the order ideal which is generated by the basic restriction of the intersection of all kernels 
$$ \bigcap_{\rho\in {\mathcal P}^{\nu }}  \ker\, r^{\nu }_{\rho } $$
which is the same as the order ideal generated by 
$\, \bigcap_{\rho }\, \underline{\mathfrak C}_{0 , \rho }\, $ from the fact that $\, {\mathfrak J}_{0 , \rho } = \ker\, r^{\nu }_{\rho }\, $ whenever $\, \rho\in {\mathcal P}^{\nu }\, $ to be proved in Theorem 1 below.
If $\, a\in {\mathfrak X}^{sa}\, $ is any selfadjoint element 
put $\, \mathfrak a = \inf\, \{ a \}\, $, and $\, \underline{\mathfrak a}_{\pm } = {\mathfrak a}_{\pm }^c\, $ where $\, \mathfrak a = {\mathfrak a}_+ - {\mathfrak a}_-\, $ denotes the minimal positive decomposition of $\,\mathfrak a\, $ in $\, \mathfrak L ( A )\, $.
Note that each element $\, \underline{\mathfrak a}_+\wedge \underline{\mathfrak a}_-\, $ is the (basic) restriction of the corresponding element in $\, \mathfrak L ( \mathcal B ( \mathcal H ) )\, $ viewing $\, a\, $ as an element in $\,\mathcal B ( \mathcal H )\, $. Then the same holds for the squares of such elements (in case that $\,\mathfrak X = A\, $ is a sub-$C^*$-algebra) which are all contained in $\, \underline{\mathfrak J}_0\, $ by Lemma 2 below.
For $\, \rho\in {\mathcal S}^{\nu }_f\, $ put $\, \underline{\mathfrak J}_{\rho }\, =\, \underline{\mathfrak J}_{0 , \rho } +\underline{\mathfrak J}_0\, $ (resp. $\, \overline{\mathfrak J}_{\rho }\, =\,\overline{\mathfrak J}_{0 , \rho } +\underline{\mathfrak J}_0\, $).
For any collection of normal states
$\, \mathcal F \subseteq {\mathcal S}^{\nu }_f\, $ let 
$\, \underline{\mathfrak J}_{0 , \mathcal F}\, $ and $\, \overline{\mathfrak J}_{0 , \mathcal F}\, $ be the order ideal generated by $\, \underline r\, \bigl( 
\bigcap_{\rho \in \mathcal F }\, {\mathfrak J}_{0 , \rho } \bigr)\, $ and $\, \overline r\, \bigl( 
\bigcap_{\rho \in \mathcal F }\, {\mathfrak J}_{0 , \rho } \bigr)\, $ respectively, and  
$\, \underline{\mathfrak J}_{\mathcal F } = \underline{\mathfrak J}_{0 , \mathcal F}\, +\, \underline{\mathfrak J}_0\, ,\,  \overline{\mathfrak J}_{\mathcal F } = \overline{\mathfrak J}_{0 , \mathcal F}\, +\, \underline{\mathfrak J}_0\, $. Put 
$\, \underline{\pi }_{0 , \mathcal F} : \mathfrak L ( \mathfrak X ) \twoheadrightarrow \mathfrak L ( \mathfrak X ) / 
\underline{\mathfrak J}_{0 , \mathcal F}\, $, resp. 
$\, \underline{\pi }^{\mathcal F} : \mathfrak L ( \mathfrak X ) \twoheadrightarrow L^{\mathfrak F} ( \mathfrak X ) = 
\mathfrak L ( \mathfrak X ) / \underline{\mathfrak J}_{\mathcal F}\, $ and similarly 
$\, \overline{\pi }^{0 , \mathcal F}\, ,\, \overline{\pi }^{\mathcal F}\, $ denote the induced surjections all of which are (restrictions of) homomorphisms due to the fact that the kernel is positively generated hence an ideal in $\, {\mathfrak L}_q ( \mathfrak X )\, $.
It will follow from Theorem 1 that for a pure normal state of $\, \mathcal B ( \mathcal H )\, $ one has 
$\, {\mathfrak J}_{\rho } = {\mathfrak J}_{0 , \rho }\,=\,\ker\, r^{\nu }_{\rho }\, $, and 
$\, {\mathfrak J}_{0 , \rho }\subseteq \ker\,  s^{\nu }_{\rho }\, $, resp. 
$\, {\mathfrak J}_{\rho } \subseteq \bigcap_{i=1}^n\, \ker\,  r^{\nu }_{{\rho }_i}\, ,\, \rho = 
\sum_{i=1}^n\, {\lambda }_i {\rho }_i\, $ for a state 
$\, \rho\in {\mathcal S}^{\nu }_f\, $, so the definition "makes sense" (at least for 
$\, \mathfrak X = \mathcal B ( \mathcal H )\, $) in that these order ideals are not unreasonably large. 
Since the antibasic restriction map is in general reducing, i.e. $\, ( \lambda\circ \overline r ) ( \mathfrak A )\leq \mathfrak A\, $ with $\, \lambda\, $ as above one easily finds that 
$\, \overline{\mathfrak J}_{0 , \rho }\subseteq \ker\, s^{\nu }_{\rho }\, $.
However the basic restriction of 
$\, {\mathfrak J}_{0 , \rho }\, $ to $\, {\mathfrak L}_q ( \mathfrak X )\, $  may not lie in the kernel of 
$\, s^{\nu }_{\rho }\, $. We let 
$\, {\mathcal S}^{\nu }_0\, $ denote the collection of vector states which is the same as 
$\, {\mathcal P}^{\nu } ( \mathcal B ( \mathcal H ) )\, $. Restricted to a general operator subsystem however these may no longer be pure states so that in this setting the notation 
$\, {\mathcal S}^{\nu }_0\, $ seems more appropriate. 
\par\noindent
Suppose that $\, M\subseteq \mathcal B ( \mathcal H )\, $ is a von Neumann algebra.
Consider the subspace of $\, \mathfrak L ( M )\, $ generated by positive basic elements which are infima of $w^*$-closed subsets $\, \mathcal C\subseteq M_+\, $. Since the sum of two $w^*$-compact convex sets is again $w^*$-compact and convex it follows from the Krein-Smulian Theorem (c.f. Theorem 2.5.9 of \cite{Pe2})  that the sum of two maximal $w^*$-closed positive subsets is again $w^*$-closed. 
Taking $w^*$-closures gives a natural projection from positive (closed) convex subsets to $w^*$-closed positive convex subsets. The subspace of differences of basic elements corresponding to $w^*$-closed positive convex sets forms a sublattice 
$\, {\mathfrak L}^{\nu } ( M )\subseteq \mathfrak L ( M )\, $ since the sum of two  
$w^*$-compact sets as well as the convex hull of the union of finitely many $w^*$-compact sets is again $w^*$-compact. Then the assignment 
$\, \underline{\mathcal C} \mapsto \underline{\mathcal C}^{\nu }\, $ extends to a well defined linear and positive projection $\, p^{\nu } :\, \mathfrak L ( M )\, \twoheadrightarrow\, {\mathfrak L}^{\nu } ( M ) \subseteq \mathfrak L ( M )\, $ onto the sublattice $\, {\mathfrak L}^{\nu } ( M )\, $. Also the subspace of differences of boundedly generated positive basic elements which are infima of $w^*$-closed positive subsets constitute a sublattice 
$\, {\mathfrak L}^{\nu }_0 ( M ) \subseteq {\mathfrak L}^{\nu } ( M )\, $, and it will be shown that this sublattice has an order isomorphic representation by $\, s^{\nu }_f\, $, i.e. the restriction of $\, s\, $ to the subset $\, {\mathcal S}^{\nu }_f\, $. Note first that for a bounded set 
the notions $w^*$-closed (= $\sigma $-weakly closed) and weak operator closed agree, so the same is true for positive boundedly generated maximal (= closed under norm and addition of positive elements) convex subsets $\, \mathcal C\, $ which are $w^*$-closed, resp. weak operator closed if and only if there exists a $w^*$-closed (= weak operator closed) bounded subset $\, \mathcal B \subset\mathcal C\, $ having the property that for each $\, c\in \mathcal C\, $ there is a $\, b\in \mathcal B\, $ with $\, b\leq c\, $. For suppose that $\, \mathcal C\, $ is not weak operator closed. Then there exists a net 
$\, ( c_{\lambda } ) \subseteq \mathcal C\, $ converging weakly to $\, c\notin\mathcal C\, $. For each 
$\, \lambda\, $ choose an element $\, b_{\lambda } \leq c_{\lambda }\, ,\, b_{\lambda }\in \mathcal B\, $. Then the bounded net $\, ( b_{\lambda } {)}_{\lambda }\, $ admits a convergent subnet 
$\, ( b_{{\lambda }_{\mu }} )\to b\, $ converging weakly to an element $\, b\in \mathcal B\, $. 
Since $\, b_{\lambda } \leq c_{\lambda }\, $ for each index $\,\lambda\, $ this implies 
$\, \rho ( b ) \leq \rho ( c )\, $ for each vector state $\,\rho\, $ whence $\, b\leq c\, $ and 
$\, c\in\mathcal C\, $ follows. Considering the quotient lattices of $\, \mathfrak L ( M )\, $ modulo the kernels of the map $\, s^{\nu }\, $ and $\, s^{\nu }_f\, $ both of which are order ideals one gets
$\, {\mathfrak L}^{\nu} ( M ) \simeq \mathfrak L ( M ) / \ker\, s^{\nu } \simeq s^{\nu } \bigl( \mathfrak L ( M ) \bigr)\, $ 
and may correspondingly define the lattice 
$\, {\mathfrak L}^w ( M ) = \mathfrak L ( M ) / \ker\, s^{\nu }_f \simeq s^{\nu }_f \bigl( \mathfrak L ( M ) \bigr)\, $  
employing the weak operator topology instead of the $w^*$-topology. It is obvious that the image of a positive basic element $\, \underline{\mathcal C}\, $ in $\, {\mathfrak L}^w ( M )\, $ is the same as the image of the element corresponding to the weak operator closure of $\, {\mathcal C}^c\, $.  Although this is not a sublattice of $\, \mathfrak L ( M )\, $ the sublattice 
$\, {\mathfrak L}^{\nu }_0 ( M ) \subseteq \mathfrak L ( M )\, $ also embeds naturally, isometrically and order isomorphic into $\, {\mathfrak L}^w ( M )\, $ by the obvious identifications. 
To see this it is sufficient to prove that for any weak operator closed maximal convex subset 
$\, \mathcal C\subseteq \mathcal B ( \mathcal H )_+\, $ and 
$\, d\geq 0\, ,\, d\notin \mathcal C\, $ there exists a state in $\, {\mathcal S}^{\nu }_f\, $ separating 
$\, \mathcal C\, $ and $\, \{ d \}\, $. Since $\, \mathcal C\, $ is weak operator closed there exists a weak operator open convex neighbourhood of $\, \{ d \}\, $ not intersecting $\,\mathcal C\, $. The Hahn-Banach separation theorem (\cite{Pe2}, Theorem 2.4.7) gives the existence of a normal selfadjoint functional 
$\, \rho\, $ in the linear span of  functionals 
$\, x\mapsto \langle x\, \xi\, ,\, \eta \rangle\, ,\, \xi\, ,\,\eta\in\mathcal H\, $ separating $\,\mathcal C\, $ and 
$\, \{ d \}\, $. Since the set $\,\mathcal C\, $ is closed under addition of positive elements this functional must be definite, for otherwise if $\,\rho = {\rho }_+ - {\rho }_-\, $ denotes the Jordan decomposition of 
$\,\rho \, $ there exists for each $\,\epsilon > 0\, $ a positive element 
$\, k\in M_+\, $ satisfying $\, {\rho }_+ ( k )\leq \epsilon\, ,\, {\rho }_- ( 1 - k )\leq \epsilon\, $ so that 
$\, \rho ( \mathcal C ) = \mathbb R\, $ is the real line which means that the functional cannot separate 
$\, \mathcal C\, $ and $\, \{ d\}\, $. Therefore $\, \rho\, $ may be assumed to be a normal state in 
$\, {\mathcal S}^{\nu }_f\, $ proving (order isomorphic) injectivity of $\, s^{\nu }_f\, $ restricted to $\, {\mathfrak L}^{\nu }_0 ( M )\, $ (compare with the argument above). A similar separation argument by arbitrary normal states works for maximal $w^*$-closed positive sets.
In case of 
$\, {\mathfrak L}^{\nu }_0 ( M )\, $ we may also substitute the dense subset 
$\, {\mathcal S}^{\nu }_{\Lambda }\subseteq {\mathcal S}^{\nu }_f\, $ of normal states which are finite convex combinations of vector states whose support projection 
$\, P_F\, $ is subordinate to some given orthonormal basis 
$\, \{ {\xi }_{\lambda } {\}}_{\lambda\in \Lambda }\, $ of $\, \mathcal H\, $, where  
$\, P_F = \sum_{i=1}^n\, p_{{\lambda }_i}\, ,\, F = \{ {\lambda }_1 ,\cdots , {\lambda }_n \}\subseteq \Lambda\, $ and $\, p_{\lambda }\, $ is the minimal projection corresponding to the basis vector 
$\, {\xi }_{\lambda }\in \mathcal H\, $,
since this subset is normdense in $\, {\mathcal S}^{\nu }\, $ with respect to the norm in $\, M^*\, $  hence $w^*$-dense in $\, {\mathcal S}^{\nu }_f\subseteq {\mathcal S}^{\nu }\, $ viewed as a subset of 
$\, {\mathfrak L}^{\nu }_0 ( M )^*\, $.  
\par\bigskip\noindent
{\bf Lemma 1.} (i)\quad For an arbitrary collection 
 $\, \mathcal F \subseteq {\mathcal S}^{\nu }_f\, $ and any $C^*$-subalgebra 
 $\, A \subseteq \mathcal B ( \mathcal H )\, $ the basic squaring operation $\, ( 1 )\, $ as well as the operation 
 $\, ( 13 )\, $ for basic positive elements drops to the quotient
$\, {\underline L}^{0 , \mathcal F} ( A )\, =\, \mathfrak L ( A )\, /\, 
\underline{\mathfrak J}_{0 ,\mathcal F}\, $, and 
to the quotient $\, {\overline L}^{\mathcal F} ( A )\, =\, 
\mathfrak L ( A )\, /\, \overline{\mathfrak J}_{\mathcal F}\, $. 
In particular if $\, \underline{\mathcal C}\, ,\, \underline{\mathcal D} \geq 0\, $ with 
$\,\underline{\mathcal C} - \underline{\mathcal D}\in \overline{\mathfrak J}_{\mathcal F}, $ then 
$$ \underline{\pi }_{\mathcal F} ( \underline{\mathcal C}^2 )\> =\> \underline{\pi }_{\mathcal F} ( \underline{\mathcal D}^2 )\> ,\quad 
\overline{\pi }_{\mathcal F} ( \underline{\mathcal C}^2 )\> =\> \overline{\pi }_{\mathcal F} ( \underline{\mathcal D}^2 )\>  . $$
\par\smallskip\noindent
(ii)\quad For any orthogonal projection $\, P\in \mathcal B ( \mathcal H )\, $ the map 
$$ Ad\, P :\> \mathcal B ( \mathcal H )\> \twoheadrightarrow\> \mathcal B ( P\, \mathcal H )\> ,\quad x\>\mapsto\> Ad\, P ( x )\> =\> P\, x\, P $$
is a bounded complete order epimorphism.
\par\bigskip\noindent
{\it Proof.}\quad To prove the first statement of the Lemma we begin by considering the special case 
$\, A = \mathcal B ( \mathcal H )\, $. Let $\, \rho\, $ be any normal state. Then there exists a 
sequence of pairwise orthogonal vector states $\, \{ {\rho }_i {\}}_{i = 1}^{\infty }\, $ and a decreasing sequence of positive numbers $\, \{ {\lambda }_i {\}}_{i = 1}^{\infty }\, ,\, 0 \leq {\lambda }_i \leq 1\, ,\, \sum_{i = 1}^\infty\, {\lambda }_i = 1\, $ such that $\, \rho ( x ) = \sum_{i = 1}^{\infty }\, {\lambda }_i {\rho }_i ( x )\, $. To each $\, {\rho }_i\, $ corresponds a minimal projection $\, p_i\in \mathcal B ( \mathcal H )\, $ such that $\, {\rho }_i ( x ) = tr ( p_i\, x )\, $. Assume first that $\, \rho\, $ is a pure (vector) state with corresponding minimal projection $\, p\, $. Let 
$\, \underline{\mathcal E} \geq 0\, $ be a basic positive element such that 
$\, r^{\nu }_{\rho } ( \underline{\mathcal E} ) = 0\, $. Then for any given $\,\epsilon > 0\, $ there exists 
$\, e\in\mathcal E\, $ with $\, \rho ( e ) \leq \epsilon\, $. We claim that this implies 
$\, r^{\nu }_{\rho } ( {\mathfrak e}^2 ) \leq 4 {\epsilon }^2\, $ where as usual $\, \mathfrak e = \inf\, \{ e \}\, $. If $\, \rho ( e ) = 0\, $ then $\, \rho ( e^2 ) = 0\, $ and we are done. Otherwise assume $\, 0 < \rho ( e ) \leq \epsilon\, $. One has $\, p e p = \rho ( e ) p\, $ so that we may find $\, e' \geq e\, $ with 
$\, e' p = p e' =  (\rho ( e ) + \epsilon ) p\, $  
by adding an element $\, \epsilon p\, +\,  R\, ( 1 - p )\, $ to $\, e\, $ with 
$\, \epsilon > 0\, $ as above and 
$\, R\gg 0\, $ a sufficiently large scalar, and simultanously erasing the matrix coefficients 
$\, p e ( 1 - p )\, $ and $\, ( 1 - p ) e p\, $. Then 
$\, r^{\nu }_{\rho } ( {\mathfrak e}^2 ) \leq \rho ( ( e' )^2 )\> =\> ( \rho ( e ) + \epsilon )^2\leq 4 
{\epsilon }^2\, $. Then also $\, r^{\nu }_{\rho } ( \underline{\mathcal E}^2 ) \leq 4 {\epsilon }^2\, $ and since $\, \epsilon\, $ can be chosen arbitrarily small this settles the case of a vector state. In fact, basically the same argument gives $\, r^{\nu }_{\rho } ( \underline{\mathcal C}^2 ) \leq r^{\nu }_{\rho } ( \underline{\mathcal C} )^2\, $, and the reverse implication is an easy consequence of the Schwarz inequality for 
$\, \rho\, $. Therefore 
$\, r^{\nu }_{\rho } ( \underline{\mathcal C}^2 ) = r^{\nu }_{\rho } ( \underline{\mathcal C} )^2\, $ for any positive basic element $\, \underline{\mathcal C}\in \mathfrak L ( \mathcal B ( \mathcal H ) )\, $.
The method of proof immediately generalizes to the case of a finite convex combinations of pairwise orthogonal vector states provided that the matrix coefficients $\, p_i e p_j = 0\, $ for $\, i\neq j\, $ are all zero. In the general case there exists for any given $\, \epsilon > 0\, $ an element $\, e\in\mathcal E\, $ 
with $\, {\rho }_i ( e ) \leq \epsilon\, $ for each $\, i = 1\, ,\cdots\, ,\, n\, $. Upon adding a positive element of the form $\, P c P\, $ with $\, P = \sum_{i=1}^n\, p_i\, $ one can substitute $\, e\, $ by an element 
$\, e' \geq e\, $ with $\, p_i e' p_j = 0\, $ whenever $\, i\neq j\, $ and $\, {\rho }_i ( e' ) \leq 2^{i}\, \epsilon\, $ so that if $\, \epsilon\, $ is chosen sufficiently small these values become arbitrarily small. Then one may apply the argument above to show that $\, s^{\nu }_{\rho } ( \underline{\mathcal E}^2 ) = 0\, $. It is rather surprising that the corresponding result for an arbitrary normal state seems to be false (or extremely hard to prove), comparing with the fact that any normal state is in the norm closure of 
$\, {\mathcal S}^{\nu }_f ( \mathcal B ( \mathcal H ) )\, $.
By monotonicity of the square it follows a forteriori that 
$\, {\mathcal Q}^2\in {\mathfrak J}_{0 , \rho }\, $ for $\, 0\leq \mathcal Q \leq \underline{\mathcal E}\, $.  Let $\, A \subseteq \mathcal B ( \mathcal H )\, $ be a $C^*$-algebra. If  
$\, \underline{\mathcal E}\in \underline{\mathfrak J}_{0 , \rho } ( A )\, $ is the restriction of an element 
$\, \underline{\mathcal E}'\in {\mathfrak J}_{0 , \rho } ( \mathcal B ( \mathcal H ) )\, $ then 
$\, \underline{\mathcal E}^2\, $ is the restriction of $\, ( \underline{\mathcal E}' )^2\, $ since the squareroot will give back the original elements. Therefore also 
$\, \underline{\mathcal E}^2\in \underline{\mathfrak J}_{0 , \rho }\, $. Then 
$\, \underline{\pi }_{0 ,\rho } ( ( \underline{\mathcal C} + \underline{\mathcal E} )^2 )\geq \underline{\pi }_{0 ,\rho } ( \underline{\mathcal C}^2 )\, $ by monotonicity of $sqr$ and positivity of 
$\, \underline{\pi }_{0 ,\rho }\, $. 
On the other hand from $\, ( 2' )\, $ one gets 
$$ \underline{\pi }_{0 ,\rho } ( ( \underline{\mathcal C} + \underline{\mathcal E} )^2 ) \leq 
{1\over\lambda }\, \underline{\pi }_{0 ,\rho } ( \underline{\mathcal C}^2 ) + {1\over 1 - \lambda }\, 
\underline{\pi }_{0 ,\rho } ( \underline{\mathcal E}^2 )\> =\> {1\over\lambda }\, \underline{\pi }_{0 ,\rho } ( \underline{\mathcal C}^2 ) $$
for each $\, 0 < \lambda < 1\, $. This proves $\, \underline{\pi }_{0 ,\rho } ( ( \underline{\mathcal C} + \underline{\mathcal E} )^2 ) = \underline{\pi }_{0 ,\rho } ( \underline{\mathcal C}^2 )\, $. Similarly 
$\, \underline{\pi }_{0 ,\rho } ( ( \underline{\mathcal C} - \underline{\mathcal E} )^2 ) \leq \underline{\pi }_{0 ,\rho } ( \underline{\mathcal C}^2 )\, $ whenever 
$\, \underline{\mathcal C} \geq \underline{\mathcal E}\, $, while on the other hand 
$$ \underline{\pi }_{0 ,\rho } ( ( \underline{\mathcal C} - \underline{\mathcal E} )^2 )\> \geq\> 
\lambda\, \bigl( \underline{\pi }_{0 ,\rho } ( \underline{\mathcal C}^2 )\> -\> {1\over 1 - \lambda }\, 
\underline{\pi }_{0 ,\rho } ( \underline{\mathcal E}^2 ) \bigr)\> =\> \lambda\, 
\underline{\pi }_{0 ,\rho } ( \underline{\mathcal C}^2 ) $$
for $\, 0\leq \lambda < 1\, $ proving the reverse inequality. 
If $\, \mathcal X = \underline{\mathcal C} + i\, \underline{\mathcal D} =  \mathcal Y + ( \mathcal P + i\, \mathcal Q )\, $ with $\, {\mathcal P}^c\, ,\,
{\mathcal Q}^c\in \underline{\mathfrak J}_{0 , \rho }^{\nu }\, $ then putting $\, \overline{\mathcal X} = \mathcal Y + 
( {\mathcal P}^c + i\, {\mathcal Q}^c )\, $ it is easy to see from monotonicity of the basic operation 
$\, ( 13 )\, $ that 
$$ \underline{\pi }_{0 ,\rho } ( \mathcal X\underline{\cdot } {\mathcal X}^* )\> \geq\> 
\underline{\pi }_{0 ,\rho } ( \mathcal Y\underline{\cdot } {\mathcal Y}^* ) \> . $$
On the other hand one gets for each $\, 0 < \lambda < 1\, $
$$ \underline{\pi }_{0 ,\rho } ( \mathcal X\underline{\cdot } {\mathcal X}^* )\> \leq\> 
\underline{\pi }_{0 ,\rho } ( \overline{\mathcal X}\underline{\cdot } {\overline{\mathcal X}}^* )\>\leq\> 
 {1\over\lambda }\, \underline{\pi }_{0 ,\rho } ( \mathcal Y\underline{\cdot } {\mathcal Y}^* )\> +\> {1\over 1 - \lambda }\, \underline{\pi }^{0 ,\rho } ( ( {\mathcal P}^c + i\, {\mathcal Q}^c )\underline{\cdot } ( {\mathcal P}^c - i\, {\mathcal Q}^c ) ) $$
$$\quad =\> {1\over \lambda }\, \underline{\pi }_{0 ,\rho } ( \mathcal Y\underline{\cdot } {\mathcal Y}^* ) $$
by operator convexity of $\, x\mapsto x x^*\, $ proving the reverse inequality. Note however that the identity $\, s^{\nu }_{\rho } ( \underline{\mathcal C}^2 ) = s^{\nu }_{\rho } ( \underline{\mathcal C} )^2\, $ is no longer valid if $\, \rho\, $ is not a pure state even in the case $\, A = \mathcal B ( \mathcal H )\, $.  
In order to prove the second statement we have to presuppose certain results to be proved in Theorem 1 below, namely that the inclusion $\, \Lambda : \mathfrak L ( A ) \hookrightarrow \mathfrak L ( \mathcal B ( \mathcal H ) )\, $ extends canonically to a normal multiplicative inclusion 
$\, \Lambda : {\mathfrak L}_q ( A ) \hookrightarrow {\mathfrak L}_q ( \mathcal B ( \mathcal H ) )\, $ which then drops to an injective $*$-homomorphism  
$$ \tilde\Lambda :\> {\mathfrak L}_q ( A )\, /\, \overline{\mathfrak J}_{\mathfrak F}\> \longrightarrow\> 
{\mathfrak L}_q ( \mathcal B ( \mathcal H ) )\, /\, \bigl( {\mathfrak J}_0\, +\, \lambda ( \overline{\mathfrak J}_{\mathcal F} ) \bigr) \> . $$
Also the (homogeneized) basic restriction map $\, \underline r\, $ as well as the antibasic restriction map  drop to  monotonous retractions 
$$ \tilde{\underline r}\> ,\> \tilde{\overline r} :\> {\mathfrak L}_q ( \mathcal B ( \mathcal H ) )\, /\, \bigl( {\mathfrak J}_0\, +\, \lambda ( \overline{\mathfrak J}_{\mathcal F} ) \bigr)\> \twoheadrightarrow\> {\mathfrak L}_q ( A )\, /\, \overline{\mathfrak J}_{\mathcal F}  $$
for $\, \tilde\lambda\, $ since $\, \underline r\, $ is convex on positive elements while 
$$ \overline r ( \mathfrak P )\leq \overline r ( \mathfrak P + \mathfrak Q )\leq \overline r ( \mathfrak P ) + \underline r ( \mathfrak Q ) \equiv \overline r ( \mathfrak P )  $$
in case that $\, 0\leq \mathfrak Q\in {\mathfrak J}_0\, +\, \lambda ( \overline{\mathfrak J}_{\mathcal F} )\, $.
Then one easily sees that $\, \tilde{\underline r}\, $ is generally magnifying whereas 
$\, \tilde{\overline r}\, $ is generally reducing. Then modulo $\, \underline{\mathfrak J}_0\, $ one has the identity
$$  \bigl[ \underline{\mathcal C}^2 \bigr]\> =\>  \bigl[ \underline{\mathcal C}\cdot \underline{\mathcal C} \bigr]\> =\>  \bigl[ \underline{\mathcal C}\bigr]\cdot  \bigl[ \underline{\mathcal C} \bigr]   $$
for any positive basic element $\, \underline{\mathcal C}\, $ as 
$$ [ \underline{\mathcal C}^2 ]\> =\> \tilde{\underline r} \bigl[ \Lambda ( \underline{\mathcal C}^2 ) \bigr]\> =\> \tilde{\underline r} \bigl[ \Lambda ( \underline{\mathcal C} )^2 \bigr]\> =\> 
\tilde{\underline r} \bigl[ \Lambda ( \underline{\mathcal C}\cdot\underline{\mathcal C} ) \bigr]\> =\> 
\bigl[ \underline{\mathcal C} \bigr]\cdot \bigl[ \underline{\mathcal C} \bigr] \> . $$
Here we have used the relation $\, \Lambda ( \underline{\mathcal C} )^2 \equiv \Lambda ( \underline{\mathcal C}\cdot \underline{\mathcal C} )\, $ modulo $\, {\mathfrak J}_0\, $ which is a special case of a general result to be proved below in Theorem 1 but in principle follows from the result 
$\, r^{\nu }_{\rho } ( \underline{\mathcal C}^2 ) = r^{\nu }_{\rho } ( \underline{\mathcal C} )^2\, $ derived above. Then if 
$\, \underline{\mathcal C} = \underline{\mathcal D}\, +\, \mathcal Q\, $ with $\,  \mathcal Q \in \overline{\mathfrak J}_{\mathcal F }\, $ we can use the $C^*$-square to see that
$$ \overline{\pi }_{\mathcal F } \bigl(  \underline{\mathcal C}^2 \bigr)\> =\>  \overline{\pi }_{\mathcal F } \bigl( \underline{\mathcal C} \bigr)\cdot \overline{\pi }_{\mathcal F } \bigl( \underline{\mathcal C} \bigr)\> =\>
\overline{\pi }_{\mathcal F } \bigl ( \underline{\mathcal D} \bigr)\cdot \overline{\pi }_{\mathcal F } \bigl( \underline{\mathcal D} \bigr)\> =\>
\overline{\pi }_{ \mathcal F } \bigl( \underline{\mathcal D}^2 \bigr)\>  . $$
The argument in the case of $\, ( 13 )\, $ is quite similar and we leave it to the reader. This proves (i).
\par\noindent
To prove (ii) let $\, P\in\mathcal B ( \mathcal H )\, $ be an orthogonal projection so that $\, Ad\, P \bigl( \mathcal B ( \mathcal H ) \bigr) \simeq \mathcal B \bigl( P\, \mathcal H \bigr)\, $. If $\, \mathcal C\subseteq \mathcal B ( \mathcal H )\, $ is any bounded positive convex subset then its image 
$\, {\mathcal C}_P = P\, \mathcal C\, P\subseteq \mathcal B \bigl( P\, \mathcal H \bigr)\, $ is again bounded positive and convex. Therefore if $\, b_P\in \bigl( {\mathcal C}_P \bigr)_c\subseteq \mathcal B \bigl( P\, \mathcal H \bigr)\, $ is any element in the lower complement of $\, {\mathcal C}_P\, $ there exists $\, r\in {\mathbb R}_-\, $ such that $\, b_P\oplus r\, \bigl( {\bf 1} - P \bigr) \leq c\, $ for any $\, c\in \mathcal C\, $ whence the result follows\qed
\par\bigskip\noindent
Let $\, a = a_+ - a_-\, $ be a selfadjoint element of $\, A\, $ with $\, a_{\pm}\geq 0\, ,\, a_+a_- = 0\, $ and $\, C\subseteq A\, $ an abelian $C^*$-subalgebra containing $\, a\, $. Define 
$$ sqr ( a )\> =\> \sup_{\lambda , C , c , d }\, \left\{ \lambda\, \bigl( c^2\> -\> {1\over 1 - \lambda }\, d^2 \bigr)\,\Bigm\vert\, c\, ,\, d\in C^+\, ,\, c - d\leq  a \right\} \> . \leqno{( 18 )} $$
If $\, a\geq 0\, $ then simply $\, sqr ( a ) = a^2\, $ by monotonicity of the square in a commutative $C^*$-algebra and $\, ( 4 )\, $, but note that $\, ( 18 )\, $ is in  general not the same as 
$\, sqr ( \mathfrak a )\, ,\, \mathfrak a = \inf\, \{ a \}\, $ as defined above. Let  
$\, {\mathfrak a}_{\pm } = ( \pm \mathfrak a ) \vee 0\, $ and 
$\, {\underline{\mathfrak a }}_{\pm } = \inf\, \{\, b\geq 0\,\vert\, b\geq \pm a \} = {\mathfrak a}_{\pm }^c\, $. Then 
$\, {\mathfrak a}_{\pm } = {\underline{\mathfrak a}}_{\pm } - ( {\underline{\mathfrak a}}_+ \wedge 
{\underline{\mathfrak a}}_- )\, $ so that 
$\, \mathfrak a = {\underline{\mathfrak a}}_+ - {\underline{\mathfrak a}}_-\, $ is the minimal basic positive decomposition of $\, \mathfrak a\, $.
One has the following result
\par\bigskip\noindent
{\bf Lemma 2.}\quad With notation as above 
$$ sqr ( a )\> =\> sqr ( a_+ )\> =\> a_+^2\> , \quad sqr ( \mathfrak a )\> \leq\> {\mathfrak a}_+^2\> . $$
If $\, A = \mathcal B ( \mathcal H )\, $ then
$$  sqr ( \mathfrak a )^{cc}\> = \>  ( {\underline{\mathfrak a}}_+^2 )_c\>  ,\quad \pi ( sqr ( \mathfrak a ) )\> =\> \pi (  {\mathfrak a}_+^2 )\> =\> \pi (  {\underline{\mathfrak a}}_+^2 )\> .  $$
Moreover $\, ( \underline{\mathfrak a}_+\wedge \underline{\mathfrak a}_- )_c = 0\, $ for every 
$\, a\in A^{sa}\, $, resp. 
$\, r^{\nu }_{\rho } ( \underline{\mathfrak a}_+\wedge \underline{\mathfrak a}_- ) = 0\, $ for every pure normal state $\, \rho\in {\mathcal P}^{\nu } ( \mathcal B ( \mathcal H ) )\, $.
\par\bigskip\noindent
{\it Proof.}\quad  $\, sqr ( a ) \leq sqr ( a_+ ) = a_+^2\, $ is obvious since any element which commutes with $\, a\, $ also commutes with $\, a_+\, $ so we only need to prove the converse on fixing $\, C = C^* ( a )\, $. 
Let $\, R\geq 0\, $ be a positive scalar. One has
$$ sqr ( a )\>\geq\> \sup_{\lambda , R }\, \bigl\{ \lambda\, \bigl( ( R + 1 )^2 a_+^2\> -\> {1\over 1 - \lambda }\, ( R a_+ + a_- )^2 \bigr)\, \bigr\} $$
and hence putting $\, {\lambda }_R = ( R + 1 )^{-1}\, $
$$ sqr ( a )\> \geq\> \sup_R\, \bigl\{ {\lambda }_R\, \bigl( ( R + 1 )^2 a_+^2\> -\> 
{1\over 1 - {\lambda }_R }\, ( R^2 a_+^2\> +\> a_-^2 ) \bigr)\, \bigr\}\> \geq\> 
a_+^2\> -\> {1\over R}\, a_-^2 $$
which converges uniformly to $\, a_+^2\, $ as $\, R\to\infty\, $.
This proves the first result.
\par\noindent
By monotonicity of $\, sqr\, $ one has 
$$ sqr ( \mathfrak a )\> \leq\>  sqr ( {\mathfrak a}_+ )\> =\> {\mathfrak a}_+^2 $$
which implies $\, sqr ( \mathfrak a )^{cc} \leq  ( {\underline{\mathfrak a}}_+^2 )_c\, $
so we only need to prove the reverse inequality, i.e. given an arbitrary state $\, \rho\in \mathcal S ( A )\, $ we must prove the inequality $\, s_{\rho } ( sqr ( \mathfrak a )^{cc} ) \geq 
s_{\rho } ( ( {\underline{\mathfrak a}}_+^2 )_c )\, $ provided that 
$\, A = \mathcal B ( \mathcal H )\, $. 
We first do the case where $\, A = M_n ( \mathbb C )\, $ is a matrix algebra. In this case it is sufficient to consider only vector states since both elements are in the image of the canonical concave lift 
$$ {\mathfrak L}_1 \bigl( M_n ( \mathbb C ) \bigr) \buildrel cc\over\longrightarrow\> \mathfrak L \bigl( M_n ( \mathbb C ) \bigr) $$
to the quotient map which factors over $\, r^{\nu }\, $ as shown in the proof of Theorem 1. Given 
$\, a\in M_n ( \mathbb C )^{sa}\, $ one can assume that $\, a_{ij } = 0\, $ for $\, i\neq j\, $ and there exist two projections $\, p_+\, ,\, p_-\, $ commuting with $\, a\, $ such that $\, p_+ + p_- = 1\, $ and 
$\, a_+ = p_+ a\, ,\, a_- = p_- a\, $. Then $\, \rho\, $ decomposes as 
$\, \rho = {\rho }_+ + {\rho }_- + {\rho }_{+/-}\, $ where 
$$ {\rho }_+ ( x )\> =\> \rho ( p_+\, x\, p_+ )\> ,\quad 
{\rho }_- ( x )\> =\> \rho ( p_-\, x\, p_- )\> ,\quad {\rho }_{+ / -} ( x )\> =\> 
\rho ( p_+\, x\, p_- + p_-\, x\, p_+ ) \> . $$
One may further decompose $\, p_+\, $ as a sum of minimal projections 
$\, p_+ = \sum_{k=1}^m\, e_k\, $ which correspond to eigenvectors of $\, a_+\, $ and $\, p_-\, $ as a sum of minimal projections $\, e_{m+1}\, ,\,\cdots\, ,\, e_n\, $ corresponding to eigenvectors of $\, a_-\, $. 
Assume now that there exists an index $\, 1\leq k\leq m\, $ and an index $\, m+1 \leq l\leq n\, $ such that $\, t_{r s} = t_{s r} = 0\, $ for all $\, 1\leq r\, ,\, s \leq n\, $ unless $,\, r\, ,\, s\in \{ k\, ,\, l \}\, $. This means that $\, \rho\, $ factors canonically over the completely positive projection $\, Q\, $ onto the corresponding copy of $\, M_2 ( \mathbb C )\, $ with $\, Q\, sqr ( \mathfrak a )\, Q = sqr ( Q\, \mathfrak a\, Q )\, ,\, Q\, {\mathfrak a}_+\, Q = (Q\, \mathfrak a\, Q)_+\, $ and we may as well assume $\, n = 2\, $. Then $\, \rho ( x ) = t r\, ( t\, x )\, $ with 
$$ t\> =\> \begin{pmatrix} t_+ & - \tau \\ - \overline\tau & t_- \end{pmatrix}\> ,\quad t_1\> +\> t_2 = 1\> ,\quad 
\tau\overline\tau\> =\> t_+\, t_- \> . $$
Without loss of generality $\, \tau\in {\mathbb R}_+\, $.
Suppose first that $\, t_+\, a_+ < t_-\, a_-\, $. Then the element 
$$ c\> =\> \lambda\> \begin{pmatrix} t_- & \tau \\ \tau & t_+ \end{pmatrix} $$
is contained in $\, {\underline{\mathfrak a}}_+\, $ if $\, \lambda > 0\, $ is chosen large enough, and 
$\, \rho ( c^2 ) = 0\, $ so that $\, s_{\rho } ( sqr ( \mathfrak a ) ) \geq
 s_{\rho } ( {\underline{\mathfrak a}}_+^2 )\, $ follows trivially. Approximating $\, a\, $ in norm by  elements with slightly smaller positive part one finds that the same result holds in case that $\, t_+\, a_+ = t_- a_-\, $.
On the other hand if $\, t_+\, a_+ > t_-\, a_-\, $ the element 
$$ d\> =\> \lambda\, \begin{pmatrix} t_- & \tau \\ \tau & t_+ \end{pmatrix} $$
is contained in $\, {\underline{\mathfrak a}}_-\, $ if $\, \lambda > 0\, $ is chosen large enough with 
$\, \rho ( x\, d )\, =\, \rho ( d\, x )\, =\, 0\, $ for every $\, x\in A\, $. In particular 
$$ \rho ( c\, d\, +\, d\, c )\> =\> 0 \leqno{( 19 )} $$
with $\, c\in {\underline{\mathfrak a}}_+\, $ any element minimizing the value of $\, \rho ( c^2 )\, $. As we will see this implies the same result. 
Assume for the moment that we are given a functional $\, \rho\, $ such that for the minimal value 
$\, \rho ( c^2 )\, ,\, c\in {\underline{\mathfrak a}}_+\, $ there exists $\, d\in {\underline{\mathfrak a}}_-\, $ with $\, \rho ( c\, d\, +\, d\, c )\, =\, 0\, $. We will see that this implies 
$\, s_{\rho } ( sqr ( \mathfrak a ) )\> \geq\> s_{\rho } ( {\underline{\mathfrak a}}_+^2 )\, $. 
For each $\, 0 \leq \lambda < 1\, $  and each 
$\, R\geq 0\, $ one has the estimate 
$$ sqr ( \mathfrak a )\> \geq\> 
\lambda\, \left( \bigl( ( R + 1 )\, {\underline{\mathfrak a}}_+ \bigr)^2\> -\> {1\over 1 - \lambda }\, 
\bigl( R\, {\underline{\mathfrak a}}_+\> +\>  {\underline{\mathfrak a}}_- \bigr)^2 \right) $$
by definition of $\, sqr\, $. Then, given $\, c\in {\underline{\mathfrak a}}_+\, ,\, 
d\in {\underline{\mathfrak a}}_-\, $ as above
$$ s_{\rho } ( sqr ( \mathfrak a ) )\> \geq\> \lambda\, \bigl( ( R + 1 )^2 \rho ( c^2 )\> -\> {1\over 1 - \lambda }\, ( R^2 \rho ( c^2 )\, +\, \rho ( d^2 ) ) \bigr) $$
since $\, \rho ( c\, d\, +\, d\, c )\, =\, 0\, $. At $\, {\lambda }_R = ( R + 1 )^{-1} \, $ the expression to the right attains the value 
$\, \rho ( c^2 )\, +\, R^{-1}\, \rho ( d^2 )\, $ which converges to the desired estimate 
$\, s_{\rho } ( sqr ( \mathfrak a ) )\, \geq\, \rho ( c^2 ) = s_{\rho } ( {\underline{\mathfrak a}}_+^2 )\, $ as 
$\, R\to\infty\, $. In particular one gets 
$$ s_{\rho } ( sqr ( \mathfrak a )^c )\>\geq\> s_{\rho } ( ( {\underline{\mathfrak a}}_+^2 )_c ) \> . $$
The whole argument given above extends to the slightly more general setting if one only assumes $\, a_{ij} = 0\, $ whenever $\, 1\leq i \leq m\, ,\, m+1\leq j\leq n\, $ or 
$\, m+1\leq i\leq n\, ,\,1\leq j\leq m\, $. If $\, \overline\rho\, $ is a given state of $M_2$-type as above so that $\, \overline\rho = \rho\circ \varphi\, $  with 
$\, \varphi : M_n ( \mathbb C ) \twoheadrightarrow M_2 ( \mathbb C )\, $ the canonical retraction corresponding to the specified pair of indices $\, ( k , l )\, $ with $\, 1\leq k\leq m\, ,\, m+1\leq l\leq n\, $
inducing a surjection of $\, {\underline{\mathfrak a}}_+\, $ onto 
$\, ({\underline{\mathfrak a}}_{\varphi })_+\, $ where $\, {\mathfrak a}_{\varphi } = \inf\, \{ \varphi ( a ) \}\, $, 
and if $\, \varphi ( c )\, $ minimizes the value of $\, \rho ( \varphi ( c )^2 )\, $ one needs to check that it is possible to extend the $2\times 2$-matrix $\, \varphi ( c )^2\, $ (or some positive element which takes the same value at $\,\rho\, $) to an $n\times n$-matrix 
$\, {\overline c}^2\, $ satisfying $\, \overline c\in {\underline{\mathfrak a}}_+\, $. Since we may add arbitrarily large elements on the diagonal entries different from $\, ( k , k )\, $ and $\, ( l , l )\, $ and scalar multiples of the complementary onedimensional projection with respect to $\, t\, $ (assuming that 
$\, \rho\, $ is a pure state) this is readily achieved. This means that the whole argument is invariant under unitary transformations of the form $\, Ad\, ( U \oplus V )\, $ with $\, U\, $ a unitary matrix in 
$\, M_m ( \mathbb C )\, $, and $\, V\, $ a unitary in the complementary copy of 
$\, M_{n-m} ( \mathbb C )\, $. If $\, \rho\, $ is any pure state on $\, M_n ( \mathbb C )\, $ it corresponds to a minimal projection $\, t\in M_n ( \mathbb C )\, ,\, t^2 = t = t^*\, $ via $\, \rho ( x ) = t r ( t\, x )\, $. Put 
$\, t_1 = p_+ t p_+\, ,\, t_{12} = p_+ t p_-\, ,\, t_2 = p_- t p_-\, $. Applying a unitary transformation of type 
$\, Ad\, ( U \oplus V )\, $ it is possible to diagonalize $\, t_1\, $ and $\, t_2\, $ both of which are rank one operators so that only two nonzero diagonal entries remain, i.e. $\, \rho\circ Ad\, ( U \oplus V )\, $ is of the type considered above and 
$\, s_{\rho } ( sqr ( \mathfrak a ) ) \geq s_{\rho } ( {\underline{\mathfrak a}}_+^2 )\, $ follows. 
Since for each fixed $\,\mathcal A\, $ the map $\, s ( \mathcal A )\, $ (resp. $\, s ( \mathcal B )\, $) viewed as a map from $\, \mathcal S ( A )\, $ to $\,\mathbb R\, $ extends to a ${\mathbb R}_+$-homogenous map on the cone of positive functionals which is an infimum of affine maps if $\, \mathcal A\, $ is basic, resp. a supremum of affine maps if $\, \mathcal B\, $ is antibasic one concludes that if given positive functionals 
$\, \rho\, ,\, {\rho }'\, ,\, {\rho }''\, $ such that $\, \rho = {\rho }'\, +\, {\rho }''\, $ and 
$$ s_{{\rho }'} ( \mathcal A )\>\geq\> s_{{\rho }'} ( \mathcal B )\> ,\quad 
s_{{\rho }''} ( \mathcal A )\> \geq \> s_{{\rho }''} ( \mathcal B ) $$
then $\, s_{\rho } ( \mathcal A )\,\geq\, s_{\rho } ( \mathcal B )\, $ follows. Therefore 
$\, s_{\rho } ( sqr ( \mathfrak a )^c )\> \geq\> s_{\rho } ( ( {\underline{\mathfrak a}}_+^2 )_c )\, $ for every 
state $\, \rho\, $ (being a convex combination of irreducible states).
Also, since 
$\, \mathcal Q = sqr ( \mathfrak a )^c\, -\, ( {\underline{\mathfrak a}}_+^2 )_c\, $ is a basic element with 
$\, s ( \mathcal Q ) \geq 0\, $ one concludes that $\, \mathcal Q \geq 0\, $ whence 
$\, sqr ( \mathfrak a )^c\,\geq\, ( {\underline{\mathfrak a}}_+^2 )_c\, $. This again implies 
$$ sqr ( \mathfrak a )^{cc}\> \geq\> ( {\underline{\mathfrak a}}_+^2 )_c \> . $$
As we have seen the reverse inequality is trivial. This accounts for the case $\, A = M_n ( \mathbb C )\, $.
The argument given above can be generalized to infinite dimensions, i.e. the case 
$\, A = \mathcal B ( \mathcal H )\, $ where $\, \mathcal H\, $ is a Hilbert space of arbitrary dimension, the most important ingredient being the existence of a matrix decomposition $\, p_+ + p_- = 1\, $ corresponding to the positive and negative part of $\, a\, $ plus pure normal states corresponding to minimal projections via the trace functional. The rest of the argument easily adapts to cover infinite dimensions since the normal states suffice to determine positivity of basic elements in 
$\, \mathfrak L ( \mathcal B ( \mathcal H ) )\, $. This then also proves the assertion 
$\, r^{\nu }_{\rho } ( \underline{\mathfrak a}_+\wedge \underline{\mathfrak a}_- ) = 0\, $ for any pure normal state $\, \rho\in {\mathcal P}^{\nu } ( \mathcal B ( \mathcal H ) )\, $ and any selfadjoint element $\, a\in \mathcal B ( \mathcal H )\, $\qed
\par\bigskip\noindent
{\it Remark.}\quad The formula $\, \pi ( sqr ( \mathfrak a ) ) = \pi ( \underline{\mathfrak a}_+^2 )\, $ is in fact valid for arbitrary $C^*$-algebras as follows from Theorem 1. The Lemma also implies 
$\, ( {\mathfrak a}_+^2 )^{cc}\> =\> 
( {\underline{\mathfrak a}}_+^2 )_c\, $, as $\, sqr ( \mathfrak a )^{cc} \leq ( {\mathfrak a}_+^2 )^{cc}\leq 
( \underline{\mathfrak a}_+^2 )_c\, $. Moreover the proof of the Lemma shows that for 
$\, A = \mathcal B ( \mathcal H )\, $ and given pure normal state $\,\rho\, $ one either has 
$\, r^{\nu }_{\rho } ( \underline{\mathfrak a}_+ ) = 0\, $ or else 
$\, r^{\nu }_{\rho } ( \underline{\mathfrak a}_- ) = 0\, $ so that $\, \underline{\mathfrak a}_+ \wedge \underline{\mathfrak a}_- \in {\mathfrak J}_{0 , \rho }\, $. In particular 
$\, \underline{\mathfrak a}_+ \wedge \underline{\mathfrak a}_- \in ker\,\pi\, $ for arbitrary $C^*$-algebras $\, A\, $ by the following scheme. Although there is in general no natural linear map $\, {\mathfrak L}_1 ( A ) \rightarrow {\mathfrak L}_1 ( \mathcal B ( \mathcal H ) )\, $ associated with a given $*$-representation $\, \lambda : A \hookrightarrow \mathcal B ( \mathcal H )\, $
there is some connection, namely a selfadjoint element $\, \mathcal A\in \mathfrak L ( A )\, $ whose image in $\, \mathfrak L ( \mathcal B ( \mathcal H ) )\, $ is in the kernel of $\,\pi\, $ (with respect to 
$\, {\mathfrak L}_1 ( \mathcal B ( \mathcal H ) )\, $) is in the kernel of $\,\pi\, $ for $\, A\, $. To see this it is sufficient consider positive elements $\,\mathcal P\, $ which are identified with their lower complements in $\, {\mathfrak L}_1 ( A )\, $. If $\,\mathcal P\, $ maps to $\, 0\, $ in 
$\, {\mathfrak L}_1 ( \mathcal B ( \mathcal H ) )\, $ then any element in the lower complement with respect to $\, \mathcal B ( \mathcal H )\, $ is negative or zero, which implies that the same holds for the lower complement with respect to $\, A\, $ whence $\, \pi ( \mathcal P ) = 0\, $. In particular 
$\, \underline{\mathfrak a}_+ \wedge \underline{\mathfrak a}_- \in ker\,\pi\, $ holds for arbitrary 
$C^*$-algebras $\, A\, $ since if $\, \Lambda : \mathfrak L ( A ) \hookrightarrow \mathfrak L ( \mathcal B ( \mathcal H ) )\, $ denotes the functorial linear map associated with some faithful $*$-representation $\, \lambda : A \hookrightarrow \mathcal B ( \mathcal H )\, $ then $\, \Lambda\, $ is (the restriction of) a $*$-homomorphism by Theorem 1 hence a lattice map and $\, \Lambda ( \underline{\mathfrak a}_{\pm } ) = \underline{\mathfrak a}^{\lambda }_{\pm }\, $ where $\, {\mathfrak a}^{\lambda } = \inf\, \{ \lambda ( a ) \}\, $. In case of 
injective $\, A\, $ one also has a positive linear retraction $\, \Upsilon : \mathfrak L ( \mathcal B ( \mathcal H ) ) \twoheadrightarrow \mathfrak L ( A )\, $ functorially extending a given completely positive linear retraction $\, \upsilon : \mathcal B ( \mathcal H ) \twoheadrightarrow A\, $ for $\,\lambda\, $ and there is a commutative diagram 
$$ \vbox{\halign{ #&#&#&#&#\cr 
&\hfil ${\mathfrak L} ( A )$\hfil &\hfil $\buildrel\pi\over\largerightarrow $\hfil &\hfil 
${\mathfrak L}_1 ( A )$\hfil & \cr
\hfil $ ( 20 )\qquad\qquad\qquad\qquad\qquad\qquad$\hfil &\hfil $\Lambda\Bigm\downarrow $\hfil &&\hfil 
$ cv ( \upsilon )\Bigm\uparrow cc ( \upsilon ) $\hfil &\hfil $\qquad\qquad\qquad\qquad\qquad\qquad $\hfil \cr
&\hfil ${\mathfrak L} ( \mathcal B ( \mathcal H ) )$\hfil &\hfil $\buildrel\pi\over\largerightarrow $\hfil &\hfil ${\mathfrak L}_1 ( \mathcal B ( \mathcal H ) )$\hfil & \cr }} \> . $$
To see this consider the composition of the upper horizontal map of the diagram with either map $\, cc ( \lambda )\, $ and $\, cv ( \lambda )\, $ respectively and check that the image of a positive element $\, \mathcal P\, $ by the former composition is smaller or equal than its image under the linear composition 
$\, \mathfrak L ( A ) \hookrightarrow \mathfrak L ( \mathcal B ( \mathcal H ) ) \twoheadrightarrow {\mathfrak L}_1 ( \mathcal B ( \mathcal H ) )\, $ whereas the image of the same element by the latter (convex) composition is larger or equal than $\, \pi\circ \Lambda ( \mathcal P )\, $. On the other hand any of the compositions 
$$ cc ( \lambda )\circ cc ( \upsilon ) = cc ( \lambda )\circ cv ( \upsilon ) = cv ( \lambda )\circ cc ( \upsilon ) = cv ( \lambda )\circ cv ( \upsilon ) = id $$
equals the identity map from rigidity. Therefore by the sandwich principle both compositions 
$\, cv ( \upsilon )\circ\pi\circ \Lambda = cc ( \upsilon )\circ\pi\circ\Lambda\, $ must equal $\,\pi\, $.
\par\bigskip\noindent
{\bf Lemma 3.}\quad Let $\, \mathcal H\, $ be a Hilbert space of arbitrary dimension and 
$\, \mathcal B ( \mathcal H )\, $ the algebra of bounded operators on $\,\mathcal H\, $. 
Then for any pure normal state $\, \rho \in {\mathcal P}^{\nu } ( \mathcal B ( \mathcal H ) )\, $ the induced functional 
$$ r^{\nu }_{\rho } : \mathfrak L ( \mathcal B ( \mathcal H ) )\> \rightarrow\> \mathbb R $$ 
satifies $\, r^{\nu }_{\rho } ( \underline{\mathcal C}^2 ) = r^{\nu }_{\rho } ( \underline{\mathcal C} )^2 \, $ for every basic positive element $\, \underline{\mathcal C}\geq 0\, $ and every element $\, \mathcal A\, $ with $\, r^{\nu } ( \mathcal A ) = 0\, $ is contained in $\, {\mathfrak J}_{0 , \rho }\, $. This applies in particular to the quadratic clouds of first, second and third type with respect to any given triple 
$\, ( \underline{\mathcal C}\, ,\, \underline{\mathcal D}\, ,\, \underline{\mathcal E} )\, $ of basic positive elements which are in the kernel of 
$\, {\pi }_0 = {\pi }_{0 , {\mathcal P}^{\nu } ( \mathcal B ( \mathcal H ) )}\, $. 
\par\bigskip\noindent
{\it Proof. }\quad Let $\, A = \mathcal B ( \mathcal H )\, $ and $\,\rho\, $ be a given pure normal state.
The first statement is a consequence of the proof of Lemma 1. For $\, \underline{\mathcal E}\in \mathfrak L ( \mathcal B ( \mathcal H ) )\, $ a (positive) basic element put $\, \underline{\mathcal E}_{\rho } = \underline{\mathcal E} - r^{\nu }_{\rho } \bigl( \underline{\mathcal E} \bigr)\, {\bf 1}\in \ker\, r^{\nu }_{\rho }\, $. It is sufficient to show that $\, \underline{\mathcal E}_{\rho }\in {\mathfrak J}_{0 , \rho }\, $. Since $\, \underline{\mathcal E}_{\rho }\, $ is basic its minimal basic positive decomposition is given by 
$$ \underline{\mathcal E}_{\rho }\> =\> \underline{\mathcal E}_{\rho }^+\> -\> \underline{\mathcal E}_{\rho }^-\> =\> ( \underline{\mathcal E}_{\rho , +} )^c\> -\> ( \underline{\mathcal E}_{\rho , -} )^c $$
where $\, \underline{\mathcal E}_{\rho } = \underline{\mathcal E}_{\rho , +} 
- \underline{\mathcal E}_{\rho , -}\, $ denotes the minimal positive decomposition. 
Given $\, \epsilon > 0\, $ choose $\, e\in {\mathcal E}^c\, $ with $\, \rho ( e )\leq r^{\nu }_{\rho } ( \underline{\mathcal E} ) + \epsilon\, $. Then 
$\, 0 = r^{\nu }_{\rho } ( \underline{\mathcal E}_{\rho } ) \leq r^{\nu }_{\rho } ( {\mathfrak e}_{\rho } ) \leq \epsilon\,,$ where 
$\, {\mathfrak e}_{\rho } = \inf\, \{ e - r^{\nu }_{\rho } ( \underline{\mathcal E} ) {\bf 1} \}\, $. By Lemma 2 one always has either 
$\, r^{\nu }_{\rho } ( ( \underline{\mathfrak e}_{\rho } )_+ ) = 0\, $ or else 
$\, r^{\nu }_{\rho } ( ( \underline{\mathfrak e}_{\rho } )_- ) = 0\, $ which implies 
$$ 0\>\leq\> r^{\nu }_{\rho } ( \underline{\mathcal E}^+_{\rho } )\> \leq\> r^{\nu }_{\rho } ( ( \underline{\mathfrak e}_{\rho } )_+ )\> \leq \epsilon\> . $$ 
Then $\, r^{\nu }_{\rho } ( \underline{\mathcal E}_{\rho }^-  ) \leq \epsilon\, $ by linearity of 
$\, r^{\nu }_{\rho }\, $ and therefore since $\,\epsilon\, $ was arbitrary 
$$  \underline{\mathcal E}\> \equiv\> r^{\nu }_{\rho } ( \underline{\mathcal E} )\, {\bf 1}\quad \mod {\mathfrak J}_{0 , \rho }\> .  $$
The argument above shows that 
$$ \mathcal A\> =\> \underline{\mathcal C}\> -\> \underline{\mathcal D}\>\in \ker\, r^{\nu }_{\rho }\quad\Longrightarrow\quad \mathcal A\> =\> \underline{\mathcal C}_{\rho }\> -\> \underline{\mathcal D}_{\rho }\>\in\> {\mathfrak J}_{0 , \rho }\> , $$ 
i.e. the kernel of $\, r^{\nu }_{\rho }\, $ is equal to $\, {\mathfrak J}_{0 , \rho }\, $\qed
\par\bigskip\noindent
The Lemma shows in particular that the induced map 
$$ r^{\nu } : {\pi }_{00} ( \mathfrak L ( \mathcal B ( \mathcal H ) ) )\> \rightarrow\> \widehat C ( {\mathcal P}_{\nu } ( \mathcal B ( \mathcal H ) ) ) $$
is injective where $\, {\pi }_{00} = {\pi }_{0 , {\mathcal S}^{\nu }_0}\, $ in the terminology introduced before Lemma 1.  
\par\bigskip\noindent
We now consider some specific notions defined for injective operator systems.
Although injective operator systems carry a uniquely determined structure as a $C^*$-algebra it is sometimes convenient to "forget" the multiplicative structure. 
In the following we will occasionally consider unital completely positive linear representations of (mostly injective) operator systems $\, \lambda : \mathfrak X \rightarrow \mathcal B ( \mathcal H )\, $. Such a representation will be called {\it transitive} iff the image contains an irreducible subalgebra of $\, \mathcal B ( \mathcal H )\, $, in particular this implies that
given any finitedimensional orthogonal projection 
$\, P\in \mathcal B ( \mathcal H )\, $ one has $\, P\, \lambda ( {\mathfrak X}_+ )\, P\, =\, \mathcal B ( P\, \mathcal H )_+\, $ and $\, P\, \lambda ( {\mathfrak X}_1 )\, P\, =\, \mathcal B ( P\, \mathcal H )_1\, $ where $\, {\mathfrak X}_+\, $ denotes the positive cone and $\, {\mathfrak X}_1\, $ denotes the unit ball of $\,\mathfrak X\, $. By Kadison's transitivity theorem (c.f. \cite{Pe1}, Theorem 2.7.5)  these properties hold for any irreducible $*$-representation of a $C^*$-algebra $\, A\, $. 
The representation $\,\lambda\, $ will be called {\it relatively transitive} iff the image contains a subalgebra which is strongly dense in the strong closure $\, \lambda ( \mathfrak X )''\, $.
An {\it injective representation} of an injective operator system $\, I\, $ is supposed to mean a unital completely positive linear map
$\, \lambda : I \rightarrow \mathcal B ( \mathcal H )\, $ such that $\, \lambda\, $ factors as a product of a complete quotient map $\, \pi : I \twoheadrightarrow J\, $ with injective image $\, J = \lambda ( I )\, $ and a completely isometric unital embedding $\, \iota : J \subseteq \mathcal B ( \mathcal H )\, $.  The representation is called {\it multiplicative} iff $\, \pi\, $ is a surjective $*$-homomorphism. Note however that a multiplicative injective representation need not be a $*$-homomorphism since 
$\, \iota\, $ is only assumed to be completely isometric. The representation is called {\it split injective} iff the quotient map $\,\pi\, $ admits a completely positive unital linear cross section $\, s : J \rightarrow I\, $.
Two injective representations 
$\, \lambda : I \rightarrow \mathcal B ( {\mathcal H}_1 )\, $ and 
$\, \mu : I \rightarrow \mathcal B ( {\mathcal H}_2 )\, $
are (spatially) equivalent iff there exists a unitary 
$\, U\in \mathcal B ( {\mathcal H}_1\, ,\, {\mathcal H}_2 )\, $ such that $\, \mu ( x )\, =\, U\, \lambda ( x )\, U^*\, $. Note however that in general a subrepresentation of an injective representation need not be injective if not completely isometric. Similarly the direct sum of two injective representations need not be injective itself. 
If $\, \rho\, ,\, \sigma\in \mathcal S ( \mathfrak X )\, $ are two states of an operator system 
$\, \mathfrak X\, $ then $\, \rho\, $ will be called {\it separated from $\,\sigma\, $} iff given any convex decomposition $\, \sigma = \lambda\, \rho\, +\, ( 1 - \lambda )\, {\sigma }'\, $ with 
$\, {\sigma }'\in \mathcal S ( \mathfrak X )\, $ this implies 
$\, \lambda = 0\, $. The pair of states $\, ( \rho\, ,\, \sigma )\, $ will be called {\it separated} iff 
$\, \rho\, $ is separated from $\,\sigma\, $ and vice versa $\,\sigma\, $ is separated from $\,\rho\, $.
An easy example of separated states is given by a pair of different pure states $\, \rho\, ,\, \sigma\in \mathcal P ( \mathfrak X )\, $.
To get a better feeling for the notion of separated states consider the case of a matrix algebra 
$\, \mathfrak X = M_n ( \mathbb C )\, $ and check that in this case $\, \rho\, $ is separated from 
$\,\sigma\, $ if and only if the contraction of $\,\rho\, $ to the orthogonal complement 
of the support of $\,\sigma\, $ is nontrivial. More generally in case of a finite dimensional operator system $\,\mathfrak X\, $ one has that $\,\rho\, $ is separated from $\,\sigma\, $ if and only if the restriction of $\,\rho\, $ to the positive kernel of $\,\sigma\, $ (= intersection of the kernel of $\,\sigma\, $ with the positive cone $\, {\mathfrak X}_+\, $) is nontrivial. If $\, \Sigma\subseteq \mathcal S ( \mathfrak X )\, $ is a subset of states of $\,\mathfrak X\, $ then $\,\Sigma\, $ will be called 
{\it positively separating for $\,\mathfrak X\, $} iff given $\, x\in\mathfrak X\, $ the condition $\, \rho ( x )\geq 0\, $ for all 
$\, \rho\in\Sigma\, $ implies $\, x\geq 0\, $. The pair
$\, ( \mathfrak X\, ,\, \Sigma )\, $ will be called {\it a separating pair} iff $\, \Sigma\, $ consists of mutually separated states and is positively separating for $\,\mathfrak X\, $.  
Given a function subsystem $\, \mathfrak X \subseteq \mathcal B ( \mathcal H )\, $ and a subset $\, \mathcal F \subseteq {\mathcal S}^{\nu }_0\, $ of vector states $\, \mathfrak X\, $ will be called {\it separating} for $\, \mathcal F\, $ if the set of restrictions of $\, \mathcal F\, $ to $\, \mathfrak X\, $ consists of mutually separated states. An injective representation
$\, \lambda : \mathfrak X \rightarrow \mathcal B ( \mathcal H )\, $
will be called {\it separating} iff there exists a 
subset $\, \Sigma\subseteq {\mathcal S}_0^{\nu }\, $ of vector states such that 
$\, ( R , \Sigma )\, $ is a separated pair if $\, R = \lambda ( \mathfrak X )''\, $ denotes the strong closure of $\,\mathfrak X\, $, and {\it supertransitive} iff it is the direct sum of separable transitive representations (in the following the term {\it supertransitive representation} will usually be understood to include the notion {\it split injective representation} when speaking of injective operator systems whereas a (relatively) transitive representation of an injective $C^*$-algebra is not necessarily required to be injective itself unless stated explicitely). In case a transitive decomposition exists it is necessarily unique up to equivalence. In case of a supertransitive representation we assume a fixed decomposition into  relatively transitive separable split injective factors and in case of a separating injective representation given a fixed subset of vector states $\, \Sigma\, $ as above as being part of the representation data.
\par\bigskip\noindent
{\it Example.}\quad Suppose given an injective $C^*$-algebra $\, I\, $. Let $\, A\subseteq I\, $ be a unital separable 
$C^*$-subalgebra and $\, {\lambda }_{\rho }\, $ an irreducible $*$-representation of $\, A\, $ associated with some pure state $\, \rho\in \mathcal P ( A )\, $.
Then there exists a completely isometric 
embedding $\, {\iota }_{\rho } : I ( A_{\rho } ) \rightarrow \mathcal B ( {\mathcal H }_{\rho } )\, $ extending $\, {\lambda }_{\rho }\, $ where $\, I ( A_{\rho } )\, $ denotes the injective envelope of 
$\, A_{\rho } = {\lambda }_{\rho } ( A )\, $. By injectivity of $\, I\, $ and rigidity of $\, I ( A )\, $ there exists a completely positive projection $\, {\Phi }_A : I \rightarrow I\, $ with range completely isometric to $\, I ( A )\, $ which factors as a product of a completely positive retraction 
$\, {\rho }_A : I \twoheadrightarrow I ( A )\, $ and a completely isometric unital embedding 
$\, {\iota }_A : I ( A ) \hookrightarrow I\, $ extending the embedding $\, A \subseteq I\, $. Then there is a completely positive map 
$\, I \twoheadrightarrow I ( A ) \rightarrow {\iota }_{\rho }\, ( I ( A_{\rho } ) )\, $ factoring over $\, {\rho }_A\, $ which extends $\, {\lambda }_{\rho }\, $ (and for simplicity is denoted by the same letter). In case that $\, A\twoheadrightarrow 
A_{\rho }\, $ admits a completely positive cross section (i.e. if the associated extension is semisplit) one also obtains a completely positive cross section $\, I ( A_{\rho } )\rightarrow I ( A ) \hookrightarrow I\, $ whose existence is due to injectivity of $\, I ( A )\, $ plus rigidity of $\, I ( A_{\rho } )\, $ showing that the representation $\, {\lambda }_{\rho }\, $ is injective. In any case it is transitive from the fact that it extends an irreducible $*$-representation of $\, A\, $.
This shows that any injective $C^*$-algebra admits a faithful supertransitive representation by summing up over all such transitive building blocks with respect to arbitrary separable subalgebras of $\, I\, $. Clearly this representation is completely isometric restricted to an arbitrary separable subalgebra 
$\, A\subseteq I\, $ hence completely isometric for $\, I\, $ which is the inductive limit of its separable subalgebras. An injective $C^*$-algebra will be called {\it separably determined} iff it identifies with the injective envelope of some separable $C^*$-subalgebra $\, A\, $ and {\it separably representable} iff it admits a faithful injective representation on separable Hilbert space. In the first case $\, I = I ( A )\, $ admits a faithful separable and separating supertransitive representation. 
With notation as above let  $\,{\rho }_I\, $ denote the (possibly mixed) state of $\, I = I ( A )\, $ represented by the vector state corresponding to $\, \rho\, $ with respect to the representation $\, {\lambda }_{\rho }\, $. 
Since $\, {\lambda }_{\rho } ( A )\subseteq \mathcal B ( {\mathcal H}_{\rho } )\, $ is an irreducible subalgebra there exist for any two different unit vectors $\, \xi\, ,\, \eta\in {\mathcal H}_{\rho }\, $ a unitary 
$\, u\in A\, $ with $\, \overline u\, \xi = \eta\, $ where $\, \overline u = {\lambda }_{\rho } ( u )\, $ (compare \cite{Pe1}, Theorem 2.7.5). Then we may regard 
$\, u\, $ as a unitary of $\, I ( A )\, $ and since $\, A\, $ is represented homomorphically one gets 
$\, {\lambda }_{\rho } ( u\, x\, u^* )\, =\, \overline u\, {\lambda }_{\rho } ( x )\, {\overline u}^*\, $ which can be seen from the Stinespring dilation of $\, {\lambda }_{\rho }\, $. Therefore any two vector states are unitarily equivalent (and separated being extensions of different pure states of $\, A\, $) for $\, I ( A )\, $. Then one may consider another inequivalent pure state 
$\, {\rho }'\, $ of $\, A\, $ and repeat the whole construction with respect to $\, {\rho }'\, $. If a state of $\, I\, $ in the new set should coincide with some state already constructed one must have $\, A_{\rho } = A_{{\rho }'}\, $ and $\, \rho  \sim {\rho }'\, $ contradicting the assumption that $\, \rho\, $ and $\, {\rho }'\, $ are inequivalent. Otherwise one obtains a subset of states having trivial intersection with the previous ones. Suppose there exists a vector state $\, {\rho }_{\xi }\, $ in the first transitive building block and a vector state $\, {\rho }_{{\xi }'}\, $ in the second block such that 
$\, {\rho }_{\xi }\, $ and $\, {\rho }_{{\xi }'}\, $ are not separated for $\, I\, $. Then their restrictions to $\, A\, $ are also not separated giving a contradiction since these restrictions represent different (even inequivalent) pure states of $\, A\, $.
Proceeding in this manner one arrives at a separating  supertransitive representation 
$\, {\lambda }_A\, $ which is completely isometric restricted to $\, A\, $ (since the atomic representation of $\, A\, $ is faithful) hence completely isometric for $\, I ( A )\, $. In particular the representation is injective. By separability of $\, A\, $ it contains a faithful separable supertransitive subrepresentation. Now if $\, I\, $ is the injective envelope of the separable subalgebra $\, A\, $ then it is also the injective envelope of any larger separable subalgebra 
$\, A \subseteq B\subseteq I\, $ so we may do the whole construction above with respect to 
$\, B\, $. 
For any pure state $\, \rho\, $ of $\, A\, $ there exists a pure state $\, {\rho }'\, $ of $\, B\, $ restricting to 
$\, \rho\, $. Then $\, A_{\rho } \subseteq B_{{\rho }'}\, $. On the other hand there may exist pure states of 
$\, B\, $ which restrict to mixed states of $\, A\, $. For each $\, \rho\, $ as above choose a pure state 
$\, {\rho }'\, $ lying above $\, \rho\, $ and consider the corresponding subrepresentation $\, {\lambda }_B^{A}\, $ of $\, {\lambda }_B\, $ which still is faithful hence injective because the GNS-representation space $\, {\mathcal H}_{{\rho }'}\, $ for $\, {\rho }'\, $ contains the representation space $\, {\mathcal H}_{\rho }\, $ for $\, \rho\, $. Moreover choose 
$\, {\lambda }_{{\rho }'}\, $ so that its contraction to $\, {\mathcal H}_{\rho }\, $ equals 
$\, {\lambda }_{\rho }\, $. Let $\, V_{\rho } : {\mathcal H}_{\rho }\hookrightarrow 
{\mathcal H}_{{\rho }'}\, $ denote the corresponding isometry. Then one sees that 
$\, {\lambda }_A\, $ is a 
quotient of $\, {\lambda }_B^{A}\, $ by the contraction $\, \prod_{\rho }\, V_{\rho }^*\, $. In this way one obtains a projective system $\, \{ {\lambda }_{A_{\mu }} {\}}_{\mu } \, $ of separating supertransitive injective representations of $\, I\, $ with $\, \mu \leq \nu\, $ iff $\, A_{\nu } \subseteq A_{\mu }\, $. 
\par\bigskip\noindent
{\bf Lemma 4.}\quad Let $\, \mathfrak X\, $ be an operator system. If given $\, \rho\, ,\, \sigma\in \mathcal S ( \mathfrak X )\, $ such that $\,\rho\, $ is separated from $\,\sigma\, $ then for each 
$\,\epsilon > 0\, $ there exists a positive element $\, c\in {\mathfrak X}_+\, $ with $\,\rho ( c ) = 1\, $ and $\, \sigma ( c )\leq \epsilon\, $. 
\par\bigskip\noindent
{\it Proof.}\quad Consider the set of twodimensional operator subsystems $\, \{ {\mathfrak X}_{\kappa } \}\, ,\, {\mathfrak X}_{\kappa }\subseteq\mathfrak X $ generated by a single positive element $\, c_{\kappa }\, $ of norm one and the order unit $\, \bf 1\, $. Then for each 
$\,\epsilon > 0\, $ there exists an index $\, \kappa = \kappa ( \epsilon )\, $ such that for the restrictions $\, {\rho }_{\kappa }\, ,\, {\sigma }_{\kappa }\in \mathcal S ( {\mathfrak X}_{\kappa } )\, $
the maximal number $\, 0\leq {\lambda }_{\kappa }\leq 1\, $ with 
$\, {\sigma }_{\kappa }\, =\, {\lambda }_{\kappa }\, {\rho }_{\kappa }\, +\, ( 1 - {\lambda }_{\kappa } )\, {\sigma }_{\kappa }'\, $ satisfies $\, {\lambda }_{\kappa }< \epsilon\, $. For suppose this is not the case so there exists $\, {\lambda }_0 > 0\, $ with $\, {\lambda }_{\kappa } \geq {\lambda }_0\, $ for all $\, \kappa\, $. Since any positive element is contained in some twodimensional operator subsystem this implies $\, \sigma \geq {\lambda }_0\, \rho\, $ whence there exists a convex decomposition $\, \sigma\, =\, {\lambda }_0\, \rho\, +\, ( 1 - {\lambda }_0 )\, {\sigma }'\, $ giving a contradiction. Given $\,\epsilon > 0\, $ choose $\, \kappa\, $ with $\, {\lambda }_{\kappa } \leq\epsilon\, $. By maximality of $\, {\lambda }_{\kappa }\, $ and finite dimensionality of $\, {\mathfrak X}_{\kappa }\, $ one gets that $\, {\rho }_{\kappa }\, $ is separated from $\, {\sigma }_{\kappa }'\, $. Then the positive kernel of $\, {\sigma }_{\kappa }'\, $ must be nontrivial and there exists a positive element $\, c\geq 0\,$ in this kernel with 
$\, \rho ( c ) = 1\, $ and $\, \sigma ( c ) = {\lambda }_{\kappa }\, \rho ( c ) \leq\epsilon\, $\qed
\par\bigskip\noindent
The following Proposition is used in the proof of Theorem 2 but is of interest in its own right. The remark is that given an irreducible representation of a $C^*$-algebra on Hilbert space any contraction of the image by a finitedimensional projection gives a complete quotient map onto $\, M_n ( \mathbb C )\, $ which is an operator order epimorphism whence by Proposition 2 there exist approximate completely positive cross sections for each such map. One would like to improve the statement to an exact result but for most purposes in analysis 'approximately' is sufficient. 
If $\, V\, ,\, W\, $ are matrix ordered spaces admitting an order unit $\, {\bf 1}_V\, $ and 
$\, {\bf 1}_W\, $ respectively (are order isomorphic to operator systems, see \cite{C-E} for the definition) then given $\, \epsilon > 0\, $ a linear map $\phi : V \rightarrow W\, $ is {\it completely $\epsilon $-positive} with respect to the order units iff  for all positive elements 
$\, x\in M_n ( V )_+\, ,\,  x \leq {\bf 1}_{V , n}\, $ one has 
$\, \phi ( x )\geq - \epsilon {\bf 1}_{W , n}\, $. A map $\, \phi : V \rightarrow W\, $ of (ordered) Banach spaces is {\it (positively) $\epsilon $-contractive} if $\, \Vert \phi ( x ) \Vert \leq 1 + \epsilon\, $ for $\, \Vert x \Vert \leq 1\, $ (and $\, x\geq 0\, $ respectively). 
\par\bigskip\noindent
{\bf Proposition 2.} \quad Let $\, \mathfrak X\, $ be an operator system. If 
$$ \psi :\> \mathfrak X\> \twoheadrightarrow\> M_n ( \mathbb C ) $$
is a unital linear completely positive complete quotient map which is an operator order epimorphism,
i.e. $\, {\psi }_k \bigl( M_k ( \mathfrak X )_+ \bigr) = M_k \bigl( M_n ( \mathbb C ) \bigr)_+\, $ for each $\, k\geq 1\, $ then for any given 
$\, \epsilon > 0\, $ the quotient map $\,\psi\, $ admits a completely positive linear map $\, {\phi }_{\epsilon } : M_n ( \mathbb C ) \hookrightarrow \mathfrak X\, $ with $\, \Vert {\phi }_{\epsilon } {\Vert }_{cb} \leq n\, $ such that $\, {\phi }_{\epsilon }\, $ is $\epsilon $-close to a cross section with respect to the completely bounded norm, i.e. $\, \Vert ( \psi\circ {\phi }_{\epsilon } )\, -\, id_{M_n ( \mathbb C )} {\Vert }_{cb} \leq \epsilon\, $.
\par\bigskip\noindent
{\it Proof.}\quad 
Let $\, \{ {\mathcal F}_{\lambda } {\}}_{\lambda\in\Lambda }\, $ denote the collection of finitedimensional operator subsystems of $\, \mathfrak X\, $ ordered by inclusion  with corresponding inclusion map $\, {\iota }_{\lambda } : {\mathcal F}_{\lambda } \hookrightarrow \mathfrak X\, $. Then given 
$\, \epsilon > 0\, $ there exists an index $\, {\lambda }_{\epsilon }\, $ such that the restriction of $\, \psi\, $ to $\, {\mathcal F}_{\lambda }\, $ is approximately an operator order epimorphism, i.e. 
$\, {\psi }_k \bigl( M_k ( {\mathcal F}_{\lambda } )_{+ , 1} \bigr)\, $ is $\epsilon $-normdense in $\, M_k \bigl( M_n ( \mathbb C ) \bigr)_{+ , 1}\, $ for each $\, k = 1 ,\cdots , n\, $ and the completely bounded norm of the inverse of the induced completely contractive map 
$$ \overline{\psi }_{\lambda }\, :\> {\mathcal G}_{\lambda }\> =\> {\mathcal F}_{\lambda }\, \bigm/\, \bigl( \ker\, \psi \cap {\mathcal F}_{\lambda } \bigr)\>\buildrel\sim\over\longrightarrow\> M_n ( \mathbb C )\> $$
is bounded by $\, 1 + \epsilon\, $ for all $\, \lambda \geq {\lambda }_{\epsilon }\, $, i.e. the restriction $\, {\psi }_{\lambda }\, $ of $\,\psi\, $ to $\, {\mathcal F}_{\lambda }\, $ is nearly a complete quotient map. Let $\, {\mathcal G}_{\lambda }\, $ be matrix ordered by the image matrix order cones of $\, {\mathcal F}_{\lambda }\, $ under the quotient map, i.e. the positive cone of 
$\, M_k ( {\mathcal G}_{\lambda } )\, $ identifies with some $\epsilon $-normdense subcone of 
$\, M_k ( M_n ( \mathbb C ) )_+\, $ by the linear completely contractive isomorphism 
$\, {\mathcal G}_{\lambda } \buildrel\sim\over\longrightarrow M_n ( \mathbb C )\, $. 
The complete quotient map $\, {\mathcal F}_{\lambda } \twoheadrightarrow {\mathcal G}_{\lambda }\, $ induces a completely isometric  complete order isomorphic injection of preduals 
$\, {\mathcal G}_{\lambda *} \hookrightarrow {\mathcal F}_{\lambda *}\, $. As we will see 
the completely $\epsilon $-positive identification $\, ( {\overline\psi }_{\lambda }^{-1} )_* : {\mathcal G}_{\lambda *} \buildrel\sim\over\longrightarrow M_n ( \mathbb C )_*\, $ has a completely positive perturbation $\, {\overline\phi }_{\lambda *} : {\mathcal G}_{\lambda *} \buildrel\sim\over\longrightarrow M_n ( \mathbb C )_*\, $ which is $\epsilon $-close to 
$\, ( {\overline\psi }_{\lambda }^{-1} )_*\, $ with respect to the completely bounded norm.
Since the (pre)dual  of $\, M_n ( \mathbb C )\, $ being selfdual with respect to the matrix order structure is injective, see Theorem 5.1 of \cite{C-E}, the completely 
positive map $\, {\overline\phi }_{\lambda *}\, $ admits a completely positive extension 
$$ {\phi }_{\lambda *} : {\mathcal F}_{\lambda *}\>\twoheadrightarrow\> M_n ( \mathbb C )_* \> , $$ 
which can be assumed completely bounded by $\, n\, $ as will turn out. 
One has the following general scheme: given a (finite dimensional) operator system $\, \mathcal F\, $ its  dual $\, {\mathcal F}^*\, $ is a dually ordered Banach space as defined above (resp. matrix ordered space in the sense of \cite{C-E}) and for each selfadjoint bounded functional $\, \phi\in {\mathcal F}^*\, $ there exists a decomposition as a difference of positive functionals $\, \phi = {\phi }_+\, -\, {\phi }_-\, $ with $\, \Vert\phi\Vert\, =\, \Vert {\phi }_+\Vert\, +\,\Vert {\phi }_-\Vert\, $. The case of a $C^*$-algebra is given by the Jordan decomposition. In the general case a unital completely positive representation $\, \iota : \mathcal F \hookrightarrow \mathcal B ( \mathcal H )\, $ dualizes to a positive complete quotient map $\, {\iota }^* : \mathcal B ( \mathcal H )^* \twoheadrightarrow {\mathcal F}^*\, $ so we may lift a selfadjoint element 
$\, \phi\in  {\mathcal F}^* \, $ to a selfadjoint element 
$\, \overline{\phi }\in \mathcal B ( \mathcal H )^*\, $ with $\, \Vert \phi\Vert\, = \Vert \overline{\phi }\Vert\, $ and consider the image 
$\, \phi = {\phi }_+\, -\, {\phi }_-\, $ of the Jordan decomposition of $\,\overline\phi\, $ to get the result. Then there exists a canonical unitization $\, V^+\, $ of any real ordered Banach space $\, V\, $ having the property that given any element $\, x\in V\, $ there is a positive decomposition $\, x = x_+\, -\, x_-\, ,\, x_{\pm }\in V_+\, $ such that $\, \Vert x\Vert\geq \max\, \{ \Vert x_+\Vert\, ,\, \Vert x_-\Vert \}\, $. Namely consider the space $\, V\, +\, \mathbb R\, {\bf 1}\, $ with order structure determined by the closed convex cone generated by $\, V_+\, $ and all positive scalar multiples of elements $\, \{ {\bf 1}\, -\, x\,\vert\, x\in V_+\, ,\, \Vert x\Vert\leq 1 \}\, $. Define the norm of an element $\, y = x\, +\, \lambda\, {\bf 1}\, $ by 
$\, \Vert y\Vert\, =\, \Vert x\Vert\, +\, \vert\lambda\vert\, $ and check that this makes $\, V^+\, $ a unital ordered Banach space in the sense given at the beginning of the section, i.e. 
$\, 0\leq y\, ,\, \Vert y\Vert\leq 1\, $ implies $\, {\bf 1}\, -\, y\geq 0\, $. The space $\, V^+\, $ admits a canonical decomposition as a direct $1$-sum $\, V^+ \simeq V {\oplus }_1\, \mathbb C\, {\bf 1}\, $ since there exists a natural contractive retraction $\, p : V^+ \twoheadrightarrow V\, $ which sends $\, x\, +\, \lambda {\bf 1}\, $ to $\, x\, $, and another natural contractive retraction 
$\, V^+ \twoheadrightarrow  \mathbb C\, {\bf 1}\, $ sending $\, x + \lambda\, {\bf 1}\, $ to 
$\, \lambda\, {\bf 1}\, $. The second retraction is clearly positive while the first is not. However in case of $\, M_n ( \mathbb C )\, $ which is order isomorphic with its dual by the map 
which sends a positive functional $\, \rho\,$ to $\, t\in M_n ( \mathbb C )_+\, $ if $\, \rho ( x ) = tr ( tx )\, $ one obtains a positive retraction $\, p_+ : M_n ( \mathbb C )_*^+ \twoheadrightarrow M_n ( \mathbb C )_*\, $ by sending the unit element $\, {\bf 1}\, $ to the preimage of $\, {\bf 1}_n\in M_n ( \mathbb C )\, $ under this identification. Since the order isomorphic identification $\, \iota : M_n ( \mathbb C )_* \buildrel\sim\over \longrightarrow M_n ( \mathbb C )\, $ is contractive one finds that the resulting retraction is positive. Define a matrix order on 
$\, M_n ( \mathbb C )_*^+\, $ by taking the closure of the affine sum of all cones of the form 
$\, \alpha\, \bigl( M_r ( M_n ( \mathbb C )_* )^+\bigr)_+ {\alpha }^* \subseteq M_k ( M_n ( \mathbb C )_*^+ )\, $ with $\, \alpha\in M_{k\, r} ( \mathbb C )\, $ a scalar matrix and 
$\, \bigl( M_r ( M_n ( \mathbb C )_* )^+\bigr)_+\, $ the positive cone of the unitization of 
$\, M_k ( M_n ( \mathbb C )_* )\, $ as defined above as positive cone for 
$\, M_k ( M_n ( \mathbb C )_*^+ )\, $. One checks that this defines a matrix order restricting to the previously defined order on each subspace $\, M_k ( M_n ( \mathbb C )_* )^+\, $ from the relation 
$\, \Vert \sum_i\, {\alpha }_i\, x {\alpha }_i^* \Vert \leq \Vert \sum_i\, {\alpha }_i\, {\alpha }_i^* {\Vert }\, \Vert x \Vert\, $ and 
$\, \alpha\, {\bf 1}_r\, {\alpha }^* = \alpha\, {\alpha }^*\, {\bf 1}_k\, $ for $\, \alpha\in M_{k\, r} ( \mathbb C )\, $ a scalar matrix and $\, x\in M_r ( M_n ( \mathbb C )_* )\, $. In the same manner one extends the given matrix orders on $\, {\mathcal F}_{\lambda *}\, $ and $\, {\mathcal G}_{\lambda *}\, $ to matrix orders on $\, {\mathcal F}_{\lambda *}^+\, $ and $\, {\mathcal G}_{\lambda *}^+\, $ respectively and checks that the canonical unital extensions of the maps considered above remain completely positive or completely $\epsilon $-positive   as the case may be.
Letting $\, r \bigl( M_n ( \mathbb C )_*^+ \bigr)\, ,\, r \bigl( {\mathcal F}_{\lambda *}^+ \bigr)\, $ and 
$\, r \bigl( {\mathcal G}_{\lambda *}^+ \bigr)\, $ denote the operator systems associated with the unital matrix ordered spaces $\, M_n ( \mathbb C )_*^+\, ,\, {\mathcal F}_{\lambda *}^+\, $ and $\, {\mathcal G}_{\lambda *}^+\, $ respectively, 
consider the composition of contractive resp. completely order isomorphic maps
$$ \rho  :\> r ( M_n ( \mathbb C )_*^+ )\> \buildrel r^{-1}\over\longrightarrow\> M_n ( \mathbb C )_*^+\> \buildrel p\over\longrightarrow\> M_n ( \mathbb C )_*\>\buildrel\iota\over\longrightarrow\> M_n ( \mathbb C )  $$
and
$$ {\rho }_+  :\> r ( M_n ( \mathbb C )_*^+ )\> \buildrel r^{-1}\over\longrightarrow\> M_n ( \mathbb C )_*^+\> \buildrel p_+\over\longrightarrow\> M_n ( \mathbb C )_*\>\buildrel\iota\over\longrightarrow\> M_n ( \mathbb C )\> ,  $$
where  the first map is an expansive unital complete order isomorphism, the last map is a completely contractive order isomorphism and the middle map is  contractive in the first instance and completely positive unital in the second instance. Note that in each case the maps
$$ M_n ( \mathbb C )_*\> \buildrel r\over\largerightarrow\> r \bigl( M_n ( \mathbb C )_*^+ \bigr)\> , \quad {\mathcal F}_{\lambda *}\buildrel r\over\largerightarrow\>
r \bigl( {\mathcal F}_{\lambda *}^+ \bigr)\> ,\quad {\mathcal G}_{\lambda *}\>\buildrel r\over\largerightarrow\> r ( {\mathcal G}_{\lambda *}^+ \bigr) $$ 
are positively isometric. This follows since for a positive element $\, x\geq 0\, $ in $\, {\mathcal F}_{\lambda *}\, $ say, its norm in $\, r \bigl( {\mathcal F}_{\lambda *}^+ \bigr)\, $ is given by the minimal value $\, \alpha \geq 0\, $ such that $\, x \leq \alpha {\bf 1}\, $ which is equal to $\, \Vert x\Vert\, $ by construction of $\, {\mathcal F}_{\lambda *}^+\, $.
Also we may write $\, p_+ = p + {\pi }_0\, $ with $\, {\pi }_0 : M_n ( \mathbb C )_*^+ \rightarrow M_n ( \mathbb C )_*\, $ the completely positive map with kernel $\, M_n ( \mathbb C )_*\, $ and sending the unit element to the (preimage of the) unit matrix $\, {\bf 1}_n\, $. The unital completely positive map $\, {\rho }_+\, $ corresponds to a unital positive map (state)
$$ s_{{\rho }_+} : M_n ( M_n ( \mathbb C )_*^+ )\> \longrightarrow\> \mathbb C $$ 
defined by linear extension of 
$$ s_{{\rho }_+} \bigl( E_{k l} \otimes x \bigr)\> =\> {1\over n}\, ( {\rho }_+ ( x ) )_{k l} $$
where $\, E_{k l}\in M_n ( \mathbb C )\, $ denotes the matrix unit with $( i , j )$-th matrix coefficient equal to $\, {\delta}_{i k}\, {\delta }_{j l}\, $ (Kronecker $\delta $-function, see \cite{Pa}, chap. 6). Similarly $\, \rho \, $ corresponds to a linear functional 
$\, s_{\rho } : M_n ( M_n ( \mathbb C )_*^+ ) = M_n ( \mathbb C ) \otimes M_n ( \mathbb C)_*^+ \rightarrow \mathbb C\, $ such that 
$\, s_{{\rho }_+} = s_{\rho } + s_{{\pi }_0}\, $ using the corresponding formulas. Via the completely $\epsilon $-positive unital identifications $\, {\mathcal G}_{\lambda *}^+ \buildrel\sim\over \longrightarrow M_n ( \mathbb C )_*^+\, $ one obtains corresponding linear (resp. linear $\epsilon $-positive) maps
$$ s^{\lambda }_{{\rho }_+}\, ,\, s^{\lambda }_{\rho } :\> 
M_n ( {\mathcal G}_{\lambda *}^+ )\> \longrightarrow\> \mathbb C  $$
corresponding to $\, {\rho }_+\circ {\overline\phi }_{\lambda *}^+\, $ and $\, \rho\circ {\overline\phi }_{\lambda *}^+\, $ and transforming correspondingly under the adjoint action of the unitary group of 
$\, M_n ( \mathbb C )\, $.  
One now defines a positive suplinear positively homogenous map 
$$ {\sigma }^{\lambda } :\> M_n ( {\mathcal F}_{\lambda *}^+ )\> \longrightarrow\> 
\mathbb C $$
using the natural unital order isomorphic injection $\, {\iota }_{\lambda } :\, M_n ( {\mathcal G}_{\lambda *}^+ ) \hookrightarrow M_n ( {\mathcal F}_{\lambda *}^+ )\, $ by the formula 
$$ {\sigma }^{\lambda } ( c )\> =\> \sup\, \bigl\{ s_{\rho }^{\lambda } ( a )\,\bigm\vert\,  
{\iota }_{\lambda } ( a ) \leq c\, ,\, a\in M_n ( {\mathcal G}_{\lambda *}^+ ) \bigr\}  $$
and easily checks that the restriction of $\, {\sigma }^{\lambda }\, $ to $\, {\mathcal G}_{\lambda *}^+\, $ is $\epsilon $-close to $\, s^{\lambda }_{\rho }\, $. Define 
$\, {\sigma }^{\lambda }_+ ( x ) = {\sigma }^{\lambda } ( x )\, +\, s_{{\pi }_0} ( x )\, $ where the latter is the obvious (linear positive) extension to $\, M_n ( {\mathcal F}_{\lambda *}^+ )\, $ of the map as above. Then there exists from the Hahn-Banach Theorem a linear extension $\, {\widetilde s}^{\lambda }_{\rho }\, $ of $\, s^{\lambda }_{\rho }\, $ satisfying $\, {\widetilde s}^{\lambda }_{\rho } ( x ) \geq {\sigma }^{\lambda } ( x )\, $ for each $\, x\in M_n ( {\mathcal F}_{\lambda *}^+ )\, $ and similarly $\, {\widetilde s}^{\lambda }_{{\rho }_+} ( x ) = {\widetilde s}^{\lambda }_{\rho } ( x )\, +\, s_{{\pi }_0} ( x )\,\geq\, {\sigma }^{\lambda }_+ ( x )\, $. This necessitates that $\, {\widetilde s}^{\lambda }_{{\rho }_+}\, $ is unital and positive. Then the associated linear map 
$$ \widetilde{\phi }_{\lambda * }^+ : {\mathcal F}_{\lambda *}^+\>\longrightarrow\> M_n ( \mathbb C )  $$
which extends $\, {\rho }_+\, $ via the identification $\, {\mathcal G}_{\lambda *}^+ \buildrel\sim\over\rightarrow M_n ( \mathbb C )_*^+\, $ is completely positive (compare with Theorem 6.1 of \cite{Pa}). But then the restriction of $\, \widetilde{\phi }_{\lambda *}\, $ associated with $\, {\widetilde s}^{\lambda }_{\rho }\, $ to $\, {\mathcal F}_{\lambda *}\, $ which extends $\, \rho \, $  is also completely positive since $\, {\pi }_0\, $ is trivial for 
$\, {\mathcal F}_{\lambda *}\, $.  We still need to show that the linear map 
$$ {\phi }_{\lambda *} :\> {\mathcal F}_{\lambda *}\>\buildrel {\widetilde \phi }_{\lambda *}\over\largerightarrow\> M_n ( \mathbb C )\> \buildrel\sim\over\largerightarrow\> M_n ( \mathbb C )_* $$
is completely bounded by $\, n\, $ which will be the case if $\, \widetilde{\phi }_{\lambda *}^+\, $ is positively contractive implying that $\, {\phi }_{\lambda *}\, $ restricted to $\, {\mathcal F}_{\lambda *}\, $ is positively bounded by $\, n\, $ whence also the dual map 
$$ {\phi }_{\lambda } :\> M_n ( \mathbb C )\> \largerightarrow\> {\mathcal F}_{\lambda } $$ 
is positively bounded by $\, n\, $ and being completely positive is completely bounded by $\, n\, $. 
However the map 
$$ r \bigl( {\mathcal F}_{\lambda *}^+ \bigr)\>\buildrel r^{-1}\over\largerightarrow\> {\mathcal F}_{\lambda *}^+\> \buildrel {\widetilde\phi }_{\lambda }^+\over\largerightarrow\> M_n ( \mathbb C ) $$
is completely contractive since unital and completely positive with $\, r^{-1}\, $ positively isometric restricted to $\, {\mathcal F}_{\lambda *}\, $ which proves the assertion. It is then obvious that 
$\, {\phi }_{\lambda }\, $ is $\epsilon $-close to a cross section for the quotient map as
$$ \Vert {\psi }_{\lambda }\circ {\phi }_{\lambda }\> -\> id_{M_n ( \mathbb C )} {\Vert }_{cb}\> =\> 
\Vert \overline{\psi }_{\lambda }\circ \overline{\phi }_{\lambda }\> -\> id_{M_n ( \mathbb C )} {\Vert }_{cb}\> =\> \Vert \overline{\psi }_{\lambda }\circ \bigl( \overline{\phi }_{\lambda }\, -\, {\overline\psi }_{\lambda }^{-1} \bigr) {\Vert }_{cb}\> \leq\> \epsilon $$
proving the Proposition\qed
\par\bigskip\noindent
Let $\, \mathfrak L\, $ denote a complete Banach lattice. A monotonous map 
$\, \phi : {\mathfrak L}_q ( \mathfrak X ) \rightarrow \mathfrak L\, $ is called {\it boundedly basically increasing normal} iff given a monotone increasing net of uniformly boundedly generated basic elements $\, \bigl( \underline{\mathcal A}_{\omega } {\bigr)}_{\omega }\, $ with supremum $\, \mathfrak A = \sup_{\omega } \underline{\mathcal A}_{\omega }\in {\mathfrak L}_q ( \mathfrak X )\, $ where each element $\, \underline{\mathcal A}_{\omega }\, $ is the infimum over a bounded subset $\, {\mathcal A}_{\omega }\subseteq B_r ( \mathfrak X ) = \bigl\{ x\in \mathfrak X\bigm\vert \Vert x\Vert \leq r \bigr\}\, $ with $\, r > 0\, $ fixed independent of $\, \omega\, $, then one has 
$$ \phi \bigl( \mathfrak A \bigr)\> =\> \sup_{\omega\in\Omega }\, \phi \bigl( \underline{\mathcal A}_{\omega } \bigr)\> . $$
Similarly, if $\, M\, $ is a von Neumann algebra a monotonous map $\, \psi : {\mathfrak L}^{\nu }_q ( M ) \rightarrow \mathfrak L\, $ is called boundedly basically increasing normal iff given a monotone increasing net $\, \bigl\{ \underline{\mathcal A}_{\omega } {\bigr\}}_{\omega\in\Omega }\, $ of uniformly boundedly generated $w^*$-closed basic elements with supremum $\, \mathfrak A\in {\mathfrak L}_q^{\nu } ( M )\, $ one has 
$$ \psi \bigl( \mathfrak A \bigr)\> =\> \sup_{\omega\in\Omega }\, \psi \bigl( \underline{\mathcal A}_{\omega } \bigr)\> .  $$ 
Caution: note that the definition requires the index set $\, \Omega\, $ to be directed, and that the last statements of the following Proposition may fail if one considers arbitrary sets of uniformly bounded basic elements since the supremum of two basic elements is no longer a basic element.
\par\bigskip\noindent
{\bf Proposition 3.}\quad Let $\, A \subseteq \mathfrak X \subseteq \mathcal B ( \mathcal H )\, $ be a relatively transitive injective operator subsystem with strong closure $\, R\subseteq \mathcal B ( \mathcal H )\, $ acting on the separable Hilbert space $\,\mathcal H\, $ and containing the strongly dense subalgebra $\, A \subseteq R\, $. Let $\, \mathcal J\subseteq R\, $ be the operator subsystem generated by the affine cone $\, {\mathcal J}^-\, $ of elements which are infima of monotone decreasing nets $\, ( x_{\lambda } {)}_{\lambda } \searrow z\, $ with $\, \{ x_{\lambda } \} \subseteq {\mathfrak X}^{sa}\, $ and $\, {\mathcal J}_A \subseteq \mathcal J\, $ the operator subsystem generated by the affine cone $\, {\mathcal J}_A^-\, $ of elements which are infima of monotone decreasing sequences $\, ( a_n {)}_n \searrow z\, $ with $\, \{ a_n \} \subseteq A^{sa}\, $. Let $\, \mathfrak L ( {\mathcal J}^- )\subseteq \mathfrak L ( \mathcal J )\, $ (resp. 
$\, \mathfrak L ( {\mathcal J}_A^- )\subseteq \mathfrak L ( {\mathcal J}_A )\, $) denote the sublattices generated by basic elements which are infima of elements in $\, {\mathcal J}^-\, $ (resp. $\, {\mathcal J}_A^-\, $).
There exist unique monotonous extensions
$$ {\pi }_{\mathcal J}^R :\> \mathfrak L \bigl( R \bigr)\> \largerightarrow\> {\mathfrak L}_1 \bigl( \mathcal J \bigr)\> ,\quad {\pi }_{{\mathcal J}_A}^R :\> \mathfrak L ( R )\> \largerightarrow\> {\mathfrak L}_1 ( {\mathcal J}_A ) $$
of the natural identification of the subspaces $\, \mathcal J\subseteq R\subseteq \mathfrak L ( R )\, $ and $\, \mathcal J\subseteq {\mathfrak L}_1 ( \mathcal J )\, $ (over the identity of $\, {\mathcal J}_A\, $ respectively). Similarly there exist unique monotonous extensions 
$$ {\pi }_{\mathfrak X}^{{\mathcal J}^-} :\> \mathfrak L ( {\mathcal J}^- )\> \largerightarrow\> {\mathfrak L}_1 ( \mathfrak X )\> ,\quad {\pi }_A^{{\mathcal J}_A^-} :\> \mathfrak L ( {\mathcal J}_A^- )\> \largerightarrow\> {\mathfrak L}_1 ( A ) $$
of the natural identification of the subspaces $\, \mathfrak X \subseteq {\mathcal J}_- \subseteq \mathfrak L ( {\mathcal J}_- )\, $ and $\, \mathfrak X \subseteq {\mathfrak L}_1 ( \mathfrak X )\, $ (over the identity of $\, A\, $ respectively).
\par\noindent
If $\, R = \prod_{\lambda }\, \mathcal B ( {\mathcal H}_{\lambda } )\, $ is a direct product of type I factors  and letting $\, s_{\Psi }\, $ denote the restriction of the function representation $\, s : \mathfrak L ( R ) \rightarrow l_{\infty } \bigl( \mathcal S ( R ) \bigr)\, $ to the subset of vector states corresponding to vectors in $\, \cup_{\lambda }\, {\mathcal H}_{\lambda }\, $ the uniquely determined surjective $*$-homomorphism 
$$ {\pi }_{\Psi } :\>  \mathfrak L \bigl( R \bigr)\,\bigm/\, \ker\, s_{\Psi }\> \largerightarrow\> {\mathfrak L}_1 \bigl( R \bigr) $$ 
admits a $*$-homomorphic cross section
$$ s_R :\> {\mathfrak L}_1 \bigl( R \bigr)\> \largerightarrow\> \mathfrak L \bigl( R \bigr)\,\bigm/\, \ker\, s_{\Psi } $$
extending the identity map on $\, R\, $.
\par\noindent
If $\, \mathfrak X \subseteq R \subseteq B\subseteq \mathcal B ( \mathcal H )\, $ is a relatively transitive operator subsystem with strong closure $\, R\, $  contained in the von Neumann algebra 
$\, B\, $ which is a direct product of type I factors
and $\, \bigl\{ \underline{\mathcal A}_{\omega } \bigr\}\subseteq {\mathfrak L}^{\nu } ( R )\subseteq {\mathfrak L}^{\nu } \bigl( B \bigr)\, $ is a directed monotone increasing net of $w^*$-closed basic elements with supremum $\, \mathfrak A\in {\mathfrak L}^{\nu }_q \bigl( B \bigr)\, $ then the minimal monotonous extension
$\, \underline\pi : {\mathfrak L}^{\nu }_q \bigl(  B \bigr) \rightarrow {\mathfrak L}_1 ( \mathfrak X )\, $ of the identity map of $\, \mathfrak X\, $ is generally concave and satisfies $\, \underline\pi ( \mathfrak A ) = \sup_{\omega }\ \underline\pi \bigl( \underline{\mathcal A}_{\omega } \bigr)\, $ (is boundedly basically increasing normal for monotone increasing nets of uniformly boundedly generated basic elements in the image of $\, {\mathfrak L}^{\nu } ( R )\, $).
In particular if $\, R\, $ is a direct product of type I factors the canonical (minimal) concave extension 
$$ {\pi }^{\nu }_c :\> {\mathfrak L}^{\nu }_q \bigl( R \bigr)\> \largerightarrow\> {\mathfrak L}_1 \bigl( R \bigr)\> ,\quad {\pi }^{\nu }_c \bigl( \mathfrak A \bigr)\> =\> \sup\, \bigl\{ \pi ( b )\bigm\vert b\in B ( X )\, ,\,  b \leq \mathfrak A \bigr\} $$ 
of $\, {\pi }^{\nu }\, $ is boundedly basically increasing normal (note that this applies in particular to the abelian $W^*$-algebra $\, R = B ( X )\, $ of bounded functions on the set $\, X\, $).
\par\noindent
If $\, \bigl\{ \underline{\mathcal A}_{\lambda }^{ccc} {\bigr\}}_{\lambda } \nearrow \mathfrak A\, $ is a monotone increasing net of complemented basic elements in $\, \mathfrak L ( \mathfrak X )\, $ then 
$\, {\pi }_c ( \mathfrak A )\, =\, {\pi }^c ( \mathfrak A )\, =\, \sup_{\lambda }\, \pi \bigl( \underline{\mathcal A}_{\lambda } \bigr)\, $ and both $\, {\pi }^c\, ,\, {\pi }_c\, $ are linear with respect to addition/subtraction of $\, \mathfrak A\, $, i.e. $\, {\pi }^c \bigl( \pm\mathfrak A + \mathfrak B ) = \pm {\pi }^c ( \mathfrak A ) + {\pi }^c ( \mathfrak B )\, ,\, {\pi }_c \bigl( \pm\mathfrak A + \mathfrak B \bigr) = \pm {\pi }_c ( \mathfrak A ) + {\pi }_c ( \mathfrak B )\, $ for any $\, \mathfrak B\in {\mathfrak L}_q ( \mathfrak X )\, $.
\par\bigskip\noindent
{\it Proof.}\quad Assume given the relatively transitive injective operator subsystem containing a strongly dense sub-$C^*$-algebra $\, A\subseteq \mathfrak X\subseteq R \subseteq\mathcal B ( \mathcal H )\, $ acting on the separable Hilbert space $\,\mathcal H\, $ and $\, \mathcal J , {\mathcal J}_A\, $ be defined as above. From the Up-Down-Theorem (Theorem 2.4.3 of 
\cite{Pe1}) any selfadjoint element $\, z\in R^{sa}\, $ is the infimum of a sequence $\, ( y_n )_n \searrow z\, $ with each $\, y_n\in {\mathcal J}_A\, $ the supremum of a sequence 
$\, ( a_{nm} )_m \nearrow y_n\, $ with $\, a_{nm}\in A^{sa}\subseteq {\mathfrak X}^{sa}\, $.
The (abstract) operator system $\,\mathfrak X\, $ is monotone complete by injectivity.
Given any monotonous retraction 
$\, r : {\mathfrak L}_q ( \mathcal J ) \twoheadrightarrow {\mathfrak L}_q ( \mathfrak X )\, $ for the (unique by Theorem 1 below) extension of the natural embedding $\, \mathfrak L ( \mathfrak X ) \subseteq \mathfrak L ( \mathcal J )\, $ consider the image of any element $\,y\in {\mathcal J}^+\, $ in $\, {\mathfrak L}_1 ( \mathfrak X )\, $ under the composition 
$$ {\mathfrak L}_q ( \mathcal J )\> \buildrel r\over\largerightarrow\> {\mathfrak L}_q ( \mathfrak X )\> \buildrel {\pi }^c , {\pi }_c\over\largerightarrow\> {\mathfrak L}_1 ( \mathfrak X )\> . $$
If $\, y = \sup_{\lambda }\, z_{\lambda }\, $ let $\, z\in \mathfrak X\, $ be the supremum of the same monotone increasing net in $\,\mathfrak X\, $, i.e. $\, y\leq z\, $, so that 
$\, ( {\pi }^c\circ r ) ( \inf\, \{ y \} )\leq ( {\pi }^c\circ r ) ( \inf\, \{ z \} ) = \pi ( \inf\, \{ z \} )\, $. On the other hand 
$\, ( {\pi }_c\circ r ) ( \inf\, \{ y \} ) \geq ( \pi\circ r ) ( \inf\, \{ z_{\lambda } \} ) = \inf\, \{ z_{\lambda } \}\, $ for each $\, \lambda\, $. Then since $\, \inf\, \{ z \} = \sup_{\lambda } \bigl\{ z_{\lambda } \bigr\}\, $ in 
$\, {\mathfrak L}_1 ( \mathfrak X )\, $ one gets that the value 
$\, ( {\pi }^c\circ r ) ( \inf\, \{ y \} ) = \inf\, \{ z \} = ( {\pi }_c\circ r ) ( \inf \{ y \} )\, $ is unique irrespective of the chosen retraction $\, r\, $. Next consider a general element
$\, y = y^+ - y^-\in {\mathcal J}^{sa}\, ,\, y^{\pm }\in {\mathcal J}^+$. Put $\, z = z^+ - z^-\, $ with $\, z^{\pm } = r ( y^{\pm } )\, $ which does not depend on $\, r\, $. Then 
$$ \pi ( z )\>=\>  ( {\pi }_c\circ r ) ( y^+ - z^- )\> \leq\> ( {\pi }_c\circ r ) ( y )\>\leq\> ( {\pi }^c\circ r ) ( y )\>\leq\> ( {\pi }^c\circ r ) ( z^+ - y_- )\> =\> \pi ( z ) $$
so that $\, {\pi }^c\circ r = {\pi }_c\circ r\, $ is uniquely determined on $\, \mathcal J\, $.  Since 
$\, \mathfrak X\, $ is injective there exists a completely positive linear retraction 
$\, \upsilon : \mathcal J \twoheadrightarrow \mathfrak X\, $ which as seen above is uniquely determined and induces a functorial linear quotient lattice map $\, \Upsilon = \mathfrak L ( \upsilon ) : \mathfrak L ( \mathcal J ) \twoheadrightarrow \mathfrak L ( \mathfrak X )\, $ extending (possibly nonuniquely) to a surjective $*$-homomorphism $\, \Upsilon : {\mathfrak L}_q ( \mathcal J ) \twoheadrightarrow {\mathfrak L}_q ( \mathfrak X )\, $ denoted by the same letter for simplicity (compare with Theorem 1).
Clearly $\, \Upsilon\, $ sends the $C^*$-subalgebra $\, C^* ( \mathcal J )\subseteq {\mathfrak L}_q ( \mathcal J )\, $ onto the $C^*$-subalgebra $\, C^* ( \mathfrak X ) \subseteq {\mathfrak L}_q ( \mathfrak X )\, $ generated by $\, \mathfrak X\, $. Since $\, {\pi }_c\circ \overline r\, $ agrees with 
$\, {\pi }^c\circ \underline r\, $ for elements in $\, C^* ( \mathcal J )\, $ by the argument below the homomorphism $\, \Upsilon\, $ sends the intersection of the kernel of $\, \pi\, $ with $\, C^* ( \mathcal J )\, $ into 
$\, \ker \pi\cap C^* ( \mathfrak X )\, $ so that $\, \Upsilon\, $ drops to a surjective $*$-homomorphism on the corresponding quotient $C^*$-algebras 
$$ \overline\upsilon :\> \pi \bigl( C^* ( \mathcal J ) \bigr)\> \twoheadrightarrow\> \pi \bigl( C^* ( \mathfrak X ) \bigr) $$
which by the Theorem A of section 2 has a $*$-homomorphic extension 
$$ \overline\Upsilon :\> {\mathfrak L}_1 ( \mathcal J )\> \twoheadrightarrow\> {\mathfrak L}_1 ( \mathfrak X ) \> . $$
There is a commutative diagram
$$ \vbox{\halign{ #&#&#\cr
\hfil $\mathfrak L ( {\mathcal J}^- )$\hfil &\hfil $\buildrel \overline r\, ,\, \underline r\over\largerightarrow $\hfil &\hfil 
${\mathfrak L}_q ( \mathfrak X )$\hfil \cr
\hfil $\pi\>\Bigm\downarrow$\hfil &&\hfil ${\pi }^c\> \Bigm\downarrow {\pi }_c $\hfil \cr
\hfil ${\mathfrak L}_1 ( \mathcal J )$\hfil &\hfil $\buildrel \overline\Upsilon\over\largerightarrow $\hfil &\hfil ${\mathfrak L}_1 ( \mathfrak X )$\hfil \cr }} \> . $$
This follows from Proposition 1 since the map 
$$ {\pi }^- :\> \mathfrak L \bigl( {\mathcal J}^- \bigr)\> \largerightarrow\> {\mathfrak L}_1 \bigl( {\mathcal J}^- \bigr)\> \simeq\> {\mathfrak L}_1 \bigl( \mathfrak X \bigr) $$
is uniquely determined with the property of being monotonous and extending the canonical map 
$\, {\upsilon }^- : {\mathcal J}^-\> \rightarrow\> \mathfrak X\, $. 
\par\noindent
A similar argument applies to show that there is a unique monotonous extension 
$$ {\upsilon }_R :\> R\> \largerightarrow\> {\mathfrak L}_1 ( \mathcal J ) $$
of the canonical identification of the subspaces $\, \mathcal J \subseteq R\, $ and $\, \mathcal J \subseteq {\mathfrak L}_1 ( \mathcal J )\, $.
We wish to show the existence of a corresponding commutative diagram of the form
$$ \vbox{\halign{ #&#&#\cr
\hfil $\mathfrak L ( R )$\hfil &\hfil $\buildrel \overline r\, ,\, \underline r\over\largerightarrow $\hfil &\hfil 
${\mathfrak L}_q ( \mathcal J )$\hfil \cr
\hfil $\pi\>\Bigm\downarrow$\hfil &&\hfil ${\pi }^c\> \Bigm\downarrow\> {\pi }_c $\hfil \cr
\hfil ${\mathfrak L}_1 ( R )$\hfil &\hfil $\buildrel \overline{\Upsilon }_R\over\largerightarrow $\hfil &\hfil ${\mathfrak L}_1 ( \mathcal J )$\hfil \cr }} \> . $$
The subset of elements 
$\, \bigl\{ \mathfrak A \bigm\vert ( {\pi }_c\circ\overline r ) \bigl( \mathfrak A \bigr) = ( {\pi }^c\circ\underline r ) \bigl( \mathfrak A \bigr) \bigr\}\subseteq {\mathfrak L}_q ( R )\, $ is a subalgebra. 
To see this suppose given two elements $\, \mathfrak A\, ,\, \mathfrak B\in {\mathfrak L}_q ( R )\, $ with $\, ( {\pi }_c\circ \overline r ) ( \mathfrak A ) = ( {\pi }^c\circ \underline r ) ( \mathfrak A )\, $ and $\, ( {\pi }_c\circ \overline r ) ( \mathfrak B ) = ( {\pi }^c\circ \underline r ) ( \mathfrak B )\, $. 
Then by concavity of $\, {\pi }_c\circ \overline r\, $ and convexity of $\, {\pi }^c\circ \underline r\, $ one gets 
$$ ( {\pi }_c\circ \overline r ) \bigl( \mathfrak A + \mathfrak B \bigr)\>\geq\> ( {\pi }_c\circ \overline r ) \bigl( \mathfrak A \bigr)\> +\> ( {\pi }_c\circ \overline r ) \bigl( \mathfrak B \bigr)\> =\> 
( {\pi }^c\circ \underline r ) \bigl( \mathfrak A \bigr)\> +\> ( {\pi }^c\circ \underline r ) \bigl( \mathfrak B \bigr)\> \geq\> ( {\pi }^c\circ \underline r ) \bigl( \mathfrak A + \mathfrak B \bigr) $$
hence equality since generally $\, {\pi }_c\circ \overline r \leq {\pi }^c\circ\underline r\, $. Therefore the subset of such elements is a linear space being closed under addition and change of sign. Also 
$$ ( {\pi }_c\circ \overline r ) \bigl( \mathfrak A \wedge \mathfrak B \bigr)\> =\> ( {\pi }_c\circ \overline r ) \bigl( \mathfrak A ) \wedge ( {\pi }_c\circ \overline r ) \bigl( \mathfrak B \bigr)\> =\> ( {\pi }^c\circ \underline r ) \bigl( \mathfrak A \bigr)\wedge ( {\pi }^c\circ \underline r ) \bigl( \mathfrak B \bigr)\> \geq\> ( {\pi }^c\circ \underline r ) \bigl( \mathfrak A \wedge \mathfrak B \bigr) $$
which again implies equality by $\, {\pi }_c\circ\overline r\leq {\pi }^c\circ\underline r\, $. Here we have used the fact that generally $\, {\pi }_c\, $ and $\,\overline r\, $ commute with the wedge-operation. Therefore also
$$ ( {\pi }_c\circ \overline r ) \bigl( \mathfrak A \vee \mathfrak B \bigr)\> =\> ( {\pi }^c\circ \underline r ) \bigl( \mathfrak A \vee \mathfrak B \bigr) $$
so that the subset in question is a sublattice hence a $C^*$-algebra being norm closed. Note that for a positive element $\, \mathfrak P \geq 0\, $ one in fact has the identities
$$ \overline r \bigl( \mathfrak P\cdot \mathfrak P \bigr)\> =\> \overline r \bigl( \mathfrak P \bigr)\cdot \overline r \bigl( \mathfrak P \bigr)\> ,\quad 
\underline r \bigl( \mathfrak P\cdot \mathfrak P \bigr)\> =\> \underline r \bigl( \mathfrak P \bigr)\cdot \underline r \bigl( \mathfrak P \bigr) $$
since by monotonicity of the square root also 
$$ \overline r \bigl( {\mathfrak P}^{1\over 2} \bigr)\> \geq\> \overline r \bigl( \mathfrak P \bigr)^{1\over 2}\> ,\quad \underline r \bigl( {\mathfrak P}^{1\over 2} \bigr)\> \leq\> \underline r \bigl( \mathfrak P \bigr)^{1\over 2} \> , $$
i.e. $\, \overline r\, $ and $\, \underline r\, $ are $\mathcal P$-maps in the sense of section 3 below. Hence the subset in question is a sub-$C^*$-algebra of $\, {\mathfrak L}_q ( R )\, $ which in particular contains the subalgebra $\, C^* ( R )\subseteq {\mathfrak L}_q ( R )\, $ generated by $\, R\, $ and is multiplicative restricted to $\, C^* ( R )\, $ so that as above $\, {\upsilon }_R\, $ drops to a $*$-homomorphism 
$$ \overline{\upsilon }_R : \pi \bigl( C^* ( R ) \bigr)\>\largerightarrow\> \pi \bigl( C^* ( \mathcal J ) \bigr) $$
which map extends to a $*$-homomorphism 
$$ \overline{\Upsilon }_R :\> {\mathfrak L}_1 ( R )\> \largerightarrow\> {\mathfrak L}_1 ( \mathcal J ) 
$$
giving the bottom line of the diagram above. We will see that $\, \overline{\Upsilon }_R\, $ is in fact the unique monotonous extension of 
$\, {\upsilon }_R : R \twoheadrightarrow {\mathfrak L}_1 ( \mathcal J )\, $. To this end construct two order isomorphic monotonous embeddings 
$\, {\rho }^{\pm } : R \hookrightarrow \mathfrak L ( \mathcal J )\,\bigm/\, \ker\, s_{\Psi }\, $ in the following way: put 
$$ {\rho }^+ ( z )\> \equiv\> \inf\, \bigl\{ y\in\mathcal J \bigm\vert y \geq z  \bigr\}\> \equiv\> \underline r ( \{ z \} )\quad \mod\> \ker\, s_{\Psi } $$
and
$$ {\rho }^- ( z )\> \equiv\> \sup\,\bigl\{ y\in \mathcal J\bigm\vert y \leq z \bigr\}\>\equiv\> \overline r ( \{ z \} )\quad \mod\> \ker\, s_{\Psi } $$
and check that this definition yields two monotonous maps on $\, R\, $. 
To see that $\, {\rho }^{\pm }\, $ is a cross section for the quotient map to $\, {\mathfrak L}_1 ( \mathcal J )\, $ note that $\, \mathfrak L ( \mathcal J )\cap\ker\, s_{\Psi } \subseteq \ker\,\pi\, $ and that the image of $\, {\rho }^{\pm } ( z )\, $ under the function representation $\, s_{\Psi }\, $ agrees with the image of 
$\, z\, $ modulo $\, \ker\, s_{\Psi }\, $. By monotonicity of $\, \rho\, $ one gets 
$$  \overline r ( z )\> =\> \sup\, \bigl\{ y\in\mathcal J \bigm\vert y \leq z \bigr\}\> \equiv\> 
\sup\, \bigl\{ \rho ( y )\bigm\vert y \leq z\, ,\, y\in \mathcal J \bigr\}\> =\> {\rho }^- ( z )\>\leq\> {\rho }^+ ( z )\qquad\qquad $$
$$\quad =\> \inf\, \bigl\{ \rho ( w )\bigm\vert w \geq z\, ,\, w\in \mathcal J \bigr\}\> \equiv\>  \inf\, \bigl\{ w\in\mathcal J\bigm\vert w\geq z \bigr\}\> =\> \underline r ( z )\>\>\quad\mod\> \ker\, s_{\Psi }\> .\qquad  $$
Then $\, {\rho }^{\pm }\, $ extends to monotonous maps 
$$ r^{\pm } :\> {\mathfrak L}_q ( R )\> \largerightarrow\> {\mathfrak L}_1 \left( \mathfrak L ( \mathcal J )\, \bigm/\, \ker\, s_{\Psi } \right)  $$
larger than $\, \overline r\, $ and smaller than $\, \underline r\, $ with respect to some ($*$-homomorphic) extension of the quotient map modulo $\, \ker\, s_{\Psi }\, $ to 
$$ {\mathfrak L}_q ( \mathcal J )\>\largerightarrow\> {\mathfrak L}_1 \left( \mathfrak L ( \mathcal J )\,\bigm/\, \ker\, s_{\Psi } \right) $$ 
(compare the proof of Theorem 1), e.g. putting 
$$ r^- ( \mathfrak A )\> \equiv\>\left( \sup\, \bigl\{ {\rho }^- ( z )\bigm\vert z\in R\, ,\, z \leq \mathfrak A \bigr\} \vee \bigl[ \overline r ( \mathfrak A ) \bigr] \right) \wedge \bigl[ \underline r ( \mathfrak A ) \bigr]\> , $$
$$ r^+ ( \mathfrak A )\> \equiv\>\left( \inf\, \bigl\{ {\rho }^+ ( z )\bigm\vert z\in R\, ,\, z \geq \mathfrak A \bigr\} \vee \bigl[ \overline r ( \mathfrak A ) \bigr] \right) \wedge \bigl[ \underline r ( \mathfrak A ) \bigr]\> . $$
Similarly the (inverse) order isomorphic identification $\, \rho ( R ) \simeq R\, $ extends to a monotonous map 
$$ q :\>  {\mathfrak L}_1 \left( \mathfrak L ( \mathcal J )\, \bigm/\, \ker\, s_{\Psi } \right)\> \twoheadrightarrow\> {\mathfrak L}_1 ( R ) $$
and one has $\,  q\circ r = \pi\, $ restricted to $\, \mathfrak L ( R )\, $ by uniqueness of $\, \pi\, $ showing on one hand that $\, q\, $ is  surjective, and secondly uniquely determined (on the image of $\, \mathfrak L ( R )\, $) , since two different choices for $\, q\, $ would lead to two different choices for $\, \pi\, $ which is impossible (compare with Proposition 1). Restricted to the subspace $\, \mathcal J\subseteq {\mathfrak L}_q ( \mathcal J )\, $ the composition of $\, q\circ s_{\Psi }\, $ with any monotonous extension $\, {\overline\Upsilon }_R : {\mathfrak L}_1 ( R ) \twoheadrightarrow {\mathfrak L}_1 ( \mathcal J )\, $ of $\, \upsilon\, $ is equal to $\, \pi\, $ so that again by uniqueness of $\, \pi\, $ on $\, \mathfrak L ( \mathcal J )\, $ extending the identity of $\, \mathcal J\, $ one finds that $\, {\overline\Upsilon }_R\, $ is uniquely determined on the image of $\, \mathfrak L ( \mathcal J )\, $ in $\, {\mathfrak L}_1 ( R )\, $ and is a $*$-homomorphism (a lattice map) since as we have seen there exists a $*$-homomorphic extension of $\, \overline{\upsilon }_R\, $. However the map $\, q\circ s_{\Psi } : \mathfrak L ( \mathcal J ) \rightarrow {\mathfrak L}_1 ( R )\, $ is clearly surjective since $\, s_{\Psi } \bigl( \mathfrak L ( \mathcal J ) \bigr) \simeq s_{\Psi } \bigl( \mathfrak L ( R ) \bigr)\, $ so that $\, {\overline\Upsilon }_R\, $ is the uniquely determined
monotonous extension of $\, {\upsilon }_R\, $. Then $\, q\, $ must be equal to the composition of the inclusion $\, \mathfrak L ( \mathcal J ) \subseteq \mathfrak L ( R )\, $ with $\, \pi\, $ modulo $\, \ker\, s_{\Psi }\, $ (on the image of $\, \mathfrak L ( \mathcal J )\, $). It is now obvious that the diagram 
$$ \vbox{\halign{ #&#&#\cr
\hfil $\mathfrak L ( R )$\hfil &\hfil $\quad\buildrel r^{\pm }\over\largerightarrow\>\>\> $\hfil &\hfil 
${\mathfrak L}_1 \left( \mathfrak L ( \mathcal J )\bigm/ \ker\, s_{\Psi } \right)$\hfil \cr
\hfil $\pi\>\Bigm\downarrow$\hfil &&\hfil $\pi\> \Bigm\downarrow\quad\quad $\hfil \cr
\hfil ${\mathfrak L}_1 ( R )$\hfil &\hfil $\quad\buildrel \overline{\Upsilon }_R\over\largerightarrow $\hfil &\hfil ${\mathfrak L}_1 ( \mathcal J )\quad$\hfil \cr }}  $$
is commutative (at least on the subalgebra of elements for which $\, {\pi }_c\circ \overline r = {\pi }^c\circ \underline r\, $ from $( 20 )$ since $\,\underline r\, $ is generally magnifying and $\, \overline r\, $ is generally reducing). 
If $\, \underline{\mathcal A}\in \underline{\mathfrak L} ( R )\, $ is a basic element one gets
$$ ( {\pi }_c\circ \overline r ) \bigl( \underline{\mathcal A} \bigr)\> \geq\> ( {\overline\Upsilon }_R\circ q ) \bigl( r^- \bigl( \underline{\mathcal A}_c \bigr) \bigr)\> =\> ( {\overline\Upsilon }_R\circ \pi ) \bigl( \underline{\mathcal A}_c \bigr)\> =\> ( {\overline\Upsilon }_R\circ \pi ) \bigl( \underline{\mathcal A} \bigr)  $$
and 
$$ ( {\pi }^c\circ \overline r ) \bigl( \underline{\mathcal A} \bigr)\>\leq\> ( {\overline\Upsilon }_R\circ q\circ r^+ ) \bigl( \underline{\mathcal A} \bigr)\> =\> ( {\overline\Upsilon }_R\circ\pi ) \bigl( \underline{\mathcal A} \bigr)  $$
whence equality follows in both instances. Now let 
$\, \overline{\mathcal B}\in \overline{\mathfrak L} ( R )\, $ be an antibasic element. Then 
$$ ( \pi\circ \overline r ) \bigl( \overline{\mathcal B} \bigr)\> =\> ( \overline{\Upsilon }_R\circ q\circ s_{\Psi } ) \bigl( \overline{\mathcal B} \bigr)\> =\>  ( \overline{\Upsilon }_R\circ\pi ) \bigl( \overline B \bigr)\> =\> ( \overline{\Upsilon }_R\circ q\circ r^- ) \bigl( \overline{\mathcal B} \bigr)\>\geq\> ( \pi \circ \overline r ) \bigl( \overline{\mathcal B} \bigr) \> . $$
For a general element $\, \underline{\mathcal A} + \overline{\mathcal B}\in \mathfrak L ( R )\, $ one gets by concavity of $\, \overline r\, $
$$ ( {\pi }_c\circ \overline r ) \bigl( \underline{\mathcal A} + \overline{\mathcal B} \bigr)\> \geq\> 
( {\pi }_c\circ \overline r ) \bigl( \underline{\mathcal A} \bigr)\> +\> ( \pi\circ \overline r ) \bigl( \overline{\mathcal B} \bigr)\> =\> ( {\overline\Upsilon }_R\circ \pi ) \bigl( \underline{\mathcal A} \bigr)\> +\> ( {\overline\Upsilon }_R\circ\pi ) \bigl( \overline{\mathcal B} \bigr)\> =\> ( {\overline\Upsilon }_R\circ\pi ) \bigl( \underline{\mathcal A} + \overline{\mathcal B} \bigr) $$
and applying the symmetry $\, \mathfrak A \mapsto - \mathfrak A\, $ gives 
$$ ( \overline{\Upsilon }_R\circ \pi ) \bigl( \underline{\mathcal A} + \overline{\mathcal B } \bigr)\> \geq\> ( {\pi }^c\circ \underline r ) \bigl( \underline{\mathcal A} + \overline{\mathcal B} \bigr) $$
so that $\, {\pi }_c\circ \overline r = {\pi }^c\circ \underline r\, $ is linear (and multiplicative) on $\, \mathfrak L ( R )\, $. The argument in case of $\, {\mathcal J}_A\, $ is much the same proving the first assertions of the Proposition.
\par\noindent
If $\, R\, $ is a direct product of type I factors and letting $\, \Psi = {\mathcal P}^{\nu } \bigl( R \bigr)\, $ denote the subset of pure normal states one has an isomorphism
$$ {\mathfrak L}_1 \left( \mathfrak L \bigl( R \bigr) \,\bigm/\, \ker\, s_{\Psi } \right)\> \simeq\> l_{\infty } \bigl( \Psi \bigr) $$
from Theorem 1. Consider the image of 
$\, \mathfrak L \bigl( R \bigr)\, $ modulo $\, \ker\, s_{\Psi }\, $. Any boundedly generated element of $\, \mathfrak L \bigl( R \bigr)\, $ (the subspace of differences of boundedly generated basic elements, see above) defines a (uniformly) continuous function in $\, C ( \Psi )\, $ with respect to the metric induced by the norm in $\, R_*\, $ since generally 
the image of a basic element in $\, l_{\infty } ( \Psi )\, $ is an upper semicontinuous function, but in order to get a discontinuity one needs elements of arbitrary large norm which are not dominated by any element in a given bounded subset, i.e. if $\, \underline{\mathcal A}\in \mathfrak L ( R )\, $ is boundedly generated, then there exists for every $\, \epsilon > 0\, $ a $\, \delta > 0\, $ such that $\, \Vert s_{\rho} ( \underline{\mathcal A} ) - s_{{\rho }'} ( \underline{\mathcal A} ) \Vert < \epsilon\, $ whenever $\, \Vert \rho - {\rho }' \Vert < \delta\, $.
On the other hand each continuous function in $\, C ( \Psi )\, $ is equal to the image of a complemented basic (and antibasic) element in $\, \mathfrak L \bigl( R \bigr)\, $. To see this choose for the given continuous function $\, f\in C ( \Psi )\, $, given $\, \epsilon > 0\, $ and for any given point $\, \rho\in \Psi \, $ an element in $\, x_{\rho , f}\in R\, $ with $\, \rho ( x_{\rho , f } ) = f ( \rho )\, $ and 
$\, {\rho }' ( x_{\rho , f} ) \geq f ( {\rho }' ) - \epsilon\, $ which is easily achieved on letting $\, \rho\, $  correspond to the eigenvector for the minimal eigenvalue of $\, x_{\rho , f }\, $ and  choosing 
$\, x_{\rho , f}\, $ sufficiently large on the orthogonal complement of $\, \rho\, $. Then the image of the basic element 
$\, \underline{\mathcal A}_{f , \epsilon } = \inf\, \{ x_{\rho , f} \}\, $ in $\, l_{\infty } ( \Psi )\, $ satisfies 
$$ s_{\Psi } \bigl( \underline{\mathcal A}_{f , \epsilon } \bigr)\> \leq\>  f\>\leq s_{\Psi } \bigl( \underline{\mathcal A}_{f , \epsilon } \bigr) + \epsilon\, {\bf 1} \> . $$ 
Similarly one can pointwise approximate the function $\, f\, $ from below by an antibasic element 
$\, \overline{\mathcal B}_{f , \epsilon }\, $ satisfying
$$ s_{\Psi } \bigl( \overline{\mathcal B}_{f , \epsilon } \bigr)\> -\> \epsilon\, {\bf 1}\> \leq\> f\>\leq 
s_{\Psi } \bigl( \overline{\mathcal B}_{f , \epsilon } \bigr) $$
by the corresponding procedure which also entails 
$$ \overline{\mathcal B}_{f , \epsilon }\> -\> \epsilon\, {\bf 1}\>\leq\> \underline{\mathcal A}_{f , \epsilon }\> +\> \epsilon\, {\bf 1} $$
since any element $\, b\leq \overline{\mathcal B}_{f , \epsilon }\, $ satisfies the relation 
$\, b \leq a + 2\, \epsilon\, {\bf 1}\, $ for every element $\, a\geq \underline{\mathcal A}_{f , \epsilon }\, $ by construction. Therefore one also has the relations 
$$ \overline{\mathcal B}_{f , \epsilon }\>\leq\> \overline{\mathcal B}_{f , \epsilon }^{cc}\>\leq\> \overline{\mathcal B}_{f , \epsilon }^c\> \leq\> \underline{\mathcal A}_{f , \epsilon }\> +\> 2\, \epsilon\, {\bf 1}\> ,\quad \overline{\mathcal B}_{f , \epsilon }\> -\> 2\,\epsilon\, {\bf 1}\> \leq\> \bigl( \underline{\mathcal A}_{f , \epsilon } \bigr)_c\>\leq\> \bigl( \underline{\mathcal A}_{f , \epsilon } \bigr)_{cc}\>\leq\> \underline{\mathcal A}_{f , \epsilon } $$
$$\quad\Longrightarrow\quad s_{\Psi } \bigl( \bigl( \underline{\mathcal A}_{f , \epsilon } \bigr)_c \bigr)\> ,\> s_{\Psi } \bigl( \bigl( \underline{\mathcal A}_{f , \epsilon } \bigr)_{cc} \bigr)\> \leq\> f\>\leq\> s_{\Psi } \bigl( \bigl( \underline{\mathcal A}_{f , \epsilon } \bigr)_c \bigr)\> +\> \epsilon\, {\bf 1}\> ,\> s_{{\Psi }_R } \bigl( \bigl( \underline{\mathcal A}_{f , \epsilon } \bigr)_{cc} \bigr)\> +\> \epsilon\, {\bf 1} $$ 
so that $\, f\, $ can be uniformly approximated by images of complemented basic (antibasic) elements which moreover approximately agree with the images of their complements. Thus if 
$\, f = s_{\Psi } \bigl( \mathcal A \bigr)\, $ for some $\,\mathcal A\in \mathfrak L \bigl( R \bigr)\, $ then the images of $\, \pi ( \mathcal A )\, $ in $\, l_{\infty } ( \Psi )\, $ using both the convex and the concave lift $\, cv , cc :\, {\mathfrak L}_1 \bigl( R \bigr) \rightarrow  \mathfrak L \bigl( R \bigr)\, $ composed with $\, s_{\Psi }\, $ must agree with each other and with 
$\, f = s_{\Psi } ( \mathcal A )\, $. Then the algebra of (bounded) continuous functions $\, C \bigl( \Psi \bigr)\, $ is contained in 
$\, s_{\Psi } \bigl( \mathfrak L \bigl( R \bigr) \bigr)\, $ and can also be identified with a subalgebra of 
$\, {\mathfrak L}_1 \bigl( R \bigr)\, $. This identification extends from the Corollary of Theorem 2 to a $*$-homomorphic cross section $\, s_R :\, {\mathfrak L}_1 ( R ) \rightarrow s_{\Psi } \bigl( \mathfrak L ( R ) \bigr)\, $ as stated. 
\par\noindent
Let  $\, \mathfrak X \subseteq R \subseteq \mathcal B ( \mathcal H )\, $ be a relatively transitive operator subsystem with strong closure $\, R\, $ acting on the Hilbert space $\, \mathcal H\, $ and let $\, \Psi = {\mathcal P}^{\nu } ( \mathcal B ( \mathcal H ) )\, $ denote  
the vector states of $\, \mathcal B ( \mathcal H )\, $. If $\, \Lambda : {\mathfrak L}_q \bigl( R \bigr) \rightarrow {\mathfrak L}_q \bigl( \mathcal B ( \mathcal H ) \bigr)\, $ denotes the functorial and normal extension of the inclusion $\, R \subseteq \mathcal B ( \mathcal H )\, $ as in Theorem 1 below 
and $\, \bigl\{ \Lambda \bigl( \underline{\mathcal A}_{\omega } \bigr) \bigr\}\subseteq {\mathfrak L}^{\nu } \bigl( \mathcal B ( \mathcal H ) \bigr)\, $ is a directed monotone increasing net of uniformly boundedly generated $w^*$-closed basic elements in the image of $\, \mathfrak L ( R )\, $ with supremum $\, \mathfrak A\in {\mathfrak L}^{\nu }_q \bigl( \mathcal B ( \mathcal H ) \bigr)\, $ we want to show that 
the minimal monotonous extension $\, \underline\pi : {\mathfrak L}_q^{\nu } \bigl( \mathcal B ( \mathcal H ) \bigr) \twoheadrightarrow {\mathfrak L}_1 ( \mathfrak X )\, $ of the identity of $\, \mathfrak X\, $ satisfies
$\, \underline\pi  \bigl( \mathfrak A \bigr) = \sup_{\omega }\, \underline\pi \bigl( \underline{\mathcal A}_{\omega } \bigr)\, $. The map $\, \underline\pi\, $ is obtained by taking the supremum of all elements  in the relative lower complement 
$$ {\mathfrak A}_{c , \mathfrak X}\> =\> \bigl\{ b\in\mathfrak X\bigm\vert b \leq \mathfrak A \bigr\} $$
in $\, {\mathfrak L}_1 ( \mathfrak X )\, $, i.e. 
$$ \underline\pi \bigl( \mathfrak A \bigr)\> =\> \sup\, \bigl\{ \pi ( b )\bigm\vert b\in\mathfrak X\, ,\, b \leq \mathfrak A \bigr\} $$
with respect to the natural embedding $\, \mathfrak X \subseteq \mathcal B ( \mathcal H ) \subseteq {\mathfrak L}^{\nu }_q \bigl( \mathcal B ( \mathcal H ) \bigr)\, $ which in particular shows that $\,\underline\pi\, $ is generally concave. Thus it is sufficient to consider the case where $\, \mathfrak A = \inf \{ b \}\, $ for some $\, b\in\mathfrak X\, $ on replacing $\, \mathfrak A\, $ by $\, \mathfrak B = \mathfrak A \wedge \mathfrak b\, $ and $\, \bigl\{ \underline{\mathcal A}_{\omega } \bigr\}\, $ by the elements $\, \bigl\{ \underline{\mathcal B}_{\omega }\bigm\vert \underline{\mathcal B}_{\omega } =  \underline{\mathcal A}_{\omega } \wedge \mathfrak b\bigr\}\, $ for $\, \mathfrak b = \inf \{ b \}\, ,\, b \leq \mathfrak A\, $ which again are $w^*$-closed uniformly boundedly generated basic elements in the image of $\, \mathfrak L ( R )\, $ and determine uniformly bounded Lipschitz normcontinuous monotone increasing concave functions $\, \{ {\widetilde f}_{\omega } : {\mathcal S}^{\nu } ( \mathcal B ( \mathcal H ) ) \rightarrow \mathbb R \}\, $ on the normal states of $\, \mathcal B ( \mathcal H )\, $
converging pointwise to a normcontinuous concave function $\, \widetilde f \leq s^{\nu } ( \mathfrak b )\, $ whose restriction $\, f\, $ to $\, \Psi\, $ being normcontinuous determines a complemented basic element $\, \underline{\mathcal A}_f \leq \mathfrak b\, $ from the argument above. Then $\, \widetilde f\, $ can be pointwise approximated from above by (normcontinuous) realvalued affine functions = elements of $\, \mathcal B ( \mathcal H )^{sa}\, $ thus defining a $w^*$-closed basic element $\, \underline{\mathcal A}_f \leq \underline{\mathcal A}_{\widetilde f} \leq \Lambda ( \mathfrak b )\, $ of $\, \mathfrak L ( \mathcal B ( \mathcal H ) )\, $. Since any complemented basic element is $w^*$-closed it is sufficient to check this relation modulo $\, \ker\, s^{\nu }\, $ (modulo $\, \ker\, s_{\Psi }\, $ for complemented basic elements and modulo $\, \ker\, s^{\nu }\, $ for $w^*$-closed basic elements). Let $\, \widetilde{\mathcal K} \subseteq \mathcal B ( \mathcal H )\, $ denote the unitization of the compact operators whence 
$\, \mathcal S \bigl( \widetilde{\mathcal K} \bigr) = co \bigl( {\mathcal S}^{\nu } ( \mathcal B ( \mathcal H ) ) \cup \{ * \} \bigr)\, $ where $\, *\, $ denotes the multiplicative functional $\, * : \widetilde{\mathcal K} \twoheadrightarrow \mathbb C\, $ with kernel $\, \mathcal K = \mathcal K ( \mathcal H )\, $. Let $\, \mathcal I : \mathfrak L \bigl( \widetilde{\mathcal K} \bigr) \rightarrow \mathfrak L ( \mathcal B ( \mathcal H ) )\, $ denote the functorial normal embedding as in Theorem 1 below. We may then consider the restricted basic elements $\, \bigl\{ \underline r \bigl( \Lambda \bigl( \underline{\mathcal B}_{\omega } \bigr) \bigr) \bigr\} \subseteq \mathfrak L \bigl( \widetilde{\mathcal K} \bigr)\, $ whose function representation on $\, {\mathcal S}^{\nu } ( \mathcal B ( \mathcal H ) )\, $ is given by $\, \{ {\widetilde f}_{\omega } \}\, $ due to the fact that $\, \widetilde{\mathcal K}\, $ is order dense in $\, \mathcal B ( \mathcal H )\, $. Then it is easy to see that the function representation of the supremum of these elements restricted to $\, {\mathcal S}^{\nu } ( \mathcal B ( \mathcal H ) )\, $ is given by $\, \widetilde f\, $ since the value of the pointwise supremum $\, {\widetilde f}^*\, $ of the extended functions 
$\, \bigl\{ {\widetilde f}_{\omega }^* = s \bigl( \underline r \bigl( \Lambda \bigl( \underline{\mathcal B}_{\omega } \bigr) \bigr) \bigr) \bigr\}\, $ on $\, \mathcal S \bigl( \widetilde{\mathcal K} \bigr)\, $ is concave and continuous restricted to $\, {\mathcal S}^{\nu } ( \mathcal B ( \mathcal H ) )\, $ hence it can be completed to a concave upper semicontinuous function defining a basic element of $\, \underline{\mathfrak L} \bigl( \widetilde{\mathcal K} \bigr)\, $ on raising the value at the point $\, *\, $ to $\, {\widetilde f}^* ( * ) :=  \limsup_{\rho\to *}\, \widetilde f ( \rho )\, $. This element is clearly larger than 
the supremum of the restricted sequence hence larger than $\, \Lambda \bigl( \mathfrak b \bigr)\, $ (viewing it as an element of $\, \mathfrak L ( \mathcal B ( \mathcal H ) )\, $ via $\, \mathcal I\, $). Therefore 
$\, \widetilde f \geq s^{\nu } \bigl( \Lambda \bigl( \mathfrak b \bigr) \bigr)\, $ while on the other hand $\, {\widetilde f}_{\omega } \leq s^{\nu } \bigl( \Lambda \bigl( \mathfrak b \bigr) \bigr)\, $ implies 
$\, \widetilde f = s^{\nu } \bigl( \Lambda \bigl( \mathfrak b \bigr) \bigr)\, $. Then $\, f = s_{\Psi } \bigl( \Lambda ( \mathfrak b ) \bigr)\, $ and
$\, \underline{\mathcal A}_{\widetilde f} = \Lambda \bigl( \mathfrak b \bigr)\, $ follows both elements being $w^*$-closed.
Put $\, {\widetilde g}_{\omega } = s^{\nu } \bigl( \Lambda \bigl( \underline{\mathcal B}_{\omega } {\bigr)}_{cc} \bigr)\, ,\, {\widetilde h}_{\omega } = s^{\nu } \bigl( \Lambda \bigl( \underline{\mathcal B}_{\omega } \bigr)_c \bigr)\, $ with restrictions $\, {\widetilde g}_{\omega } {\vert }_{\Psi } = f_{\omega } = {\widetilde h}_{\omega } {\vert }_{\Psi }\, $ from the argument above. Then the pointwise limit functions $\, \widetilde g = \sup_{\omega } {\widetilde g}_{\omega }\, ,\, \widetilde h = \sup_{\omega }\, {\widetilde h}_{\omega }\, $ satisfy $\, \widetilde h \leq\widetilde g\leq \widetilde f\, $. Moreover for each finitedimensional projection $\, P\in R\, $ and given $\,\epsilon > 0\, $ there exists 
$\, {\omega }^P_{\epsilon }\in\Omega\, $ such that the restriction of $\, {\widetilde g}_{\omega }\, $ to $\, \mathcal S ( P\, R\, P )\, $ is $\epsilon $-affine for $\, \omega\geq {\omega }^P_{\epsilon }\, $, i.e. given any finite collection of vector states $\, \{ {\rho }_1\, ,\cdots , {\rho }_n \}\, $ subordinate to $\, P\, $ and a convex combination $\, \sigma = \sum_{k = 1}^n\, {\lambda }_k\, {\rho }_k\, $ one has $\, {\widetilde g}_{\omega } ( \sigma ) \leq \sum_{k = 1}^n\, {\lambda }_k\, {\widetilde g}_{\omega } ( {\rho }_k ) + \epsilon\, $. It is then easy to see by an argument as used above that the corresponding holds for the functions $\, \{ {\widetilde h}_{\omega } \}\, $ which become approximately affine restricted to finite subspaces as $\, \omega \to\infty\, $ implying in particular that the limit functions $\, \widetilde h = \widetilde g = \widetilde f\, $ all coincide. If $\, {\mathcal S}^{\nu } ( \mathfrak X ) \subseteq \mathcal S ( \mathfrak X )\, $ denotes the restriction of $\, {\mathcal S}^{\nu } \bigl( \mathcal B ( \mathcal H ) \bigr)\, $ to $\,\mathfrak X\, $ define a map 
$$ {\Gamma }^{\nu } :\> l_{\infty } \bigl( {\mathcal S}^{\nu } \bigl( \mathcal B ( \mathcal H ) \bigr) \bigr)\> \largerightarrow\> l_{\infty } \bigl( {\mathcal S}^{\nu } ( \mathfrak X ) \bigr)\> ,\quad {\Gamma }^{\nu } ( f ) ( \sigma )\> =\> \inf\, \bigl\{ f ( \rho )\bigm\vert {\iota }^* ( \rho ) = \sigma \bigr\} $$
with respect to the functorial restriction map $\, {\iota }^* : \mathcal S \bigl( \mathcal B ( \mathcal H ) \bigr) \rightarrow \mathcal S ( \mathfrak X )\, $. One checks that the image of a lower semicontinuous convex function corresponding to an antibasic element under $\, {\Gamma }^{\nu }\, $ is again lower semicontinuous and convex. Moreover it is easy to see that $\, {\Gamma }^{\nu }\, $ is increasing normal restricted to functions $\, \{ f \}\, $ satisfying $\, {\iota }^* ( \rho ) = {\iota }^* ( {\rho }' ) \Longrightarrow f ( \rho ) = f ( {\rho }' )\, $ which is the case for the function representations of the elements $\, \bigl\{ \underline{\mathcal B}_{\omega } \bigr\}\, $. Extending the assignment to all of $\, \mathcal S ( \mathfrak X )\, $ by 
$$ \Gamma :\> l_{\infty } \bigl( {\mathcal S}^{\nu } \bigl( \mathcal B ( \mathcal H ) \bigr) \bigr)\> \largerightarrow\> \mathcal S ( \mathfrak X )\> ,\quad \Gamma ( f ) ( \sigma )\> =\> \liminf_{\tau\to\sigma }\, {\Gamma }^{\nu } ( f ) ( \tau ) $$
which is well defined due to the fact that $\, {\mathcal S}^{\nu } ( \mathfrak X ) \subseteq \mathcal S ( \mathfrak X )\, $ is $w^*$-dense one finds that $\, \Gamma ( f )\, $ is again lower semicontinuous and convex if $\, f\, $ has these properties so that $\, \Gamma ( f )\, $ defines an antibasic element of 
$\, \mathfrak L ( \mathfrak X )\, $. We let 
$$ {\Gamma }^{\nu }_0 :\> l_{\infty } \bigl( {\mathcal S}^{\nu } \bigl( \mathcal B ( \mathcal H ) \bigr) \bigr)\> \largerightarrow\> l_{\infty } \bigl( {\mathcal S}^{\nu }_0 ( \mathfrak X ) \bigr) $$
denote the composition of $\, {\Gamma }^{\nu }\, $ with the evaluation map on the vector states of $\, \mathfrak X\, $. 
Finally $\, {\Gamma }^{\nu }_0\, $ may be composed with the minimal monotonous extension 
$$ {\underline\pi }_{\mathfrak X} : l_{\infty } \bigl( {\mathcal S}^{\nu }_0 ( \mathfrak X ) \bigr)\> \largerightarrow\> {\mathfrak L}_1 ( \mathfrak X ) $$
extending the identity map of $\, \mathfrak X\, $ on the left and the function representation $\, s^{\nu }\, $ on the right to obtain a map 
$\, {\mathfrak L}^{\nu } \bigl( \mathcal B ( \mathcal H ) \bigr) \twoheadrightarrow {\mathfrak L}_1 ( \mathfrak X )\, $ which we claim to be equal to $\, \underline\pi\, $ at least restricted to antibasic elements. This follows since the antibasic element in $\, \mathfrak L ( \mathfrak X )\, $ corresponding to $\, \Gamma ( f )\, $ with 
$\, f = s^{\nu } \bigl( \overline{\mathcal B} \bigr)\, $ is clearly smaller than $\, \overline{\mathcal B}\in {\mathfrak L}^{\nu } \bigl( \mathcal B ( \mathcal H ) \bigr) \subseteq \mathfrak L \bigl( \mathcal B ( \mathcal H ) \bigr)\, $. One now checks that for $\, \widetilde f = s^{\nu } ( \mathfrak b )\, ,\, {\widetilde h}_{\omega } = s^{\nu } \bigl( \Lambda ( \underline{\mathcal B}_{\omega } )_c \bigr)\, $ one has 
$\, \bigl( ( \underline\pi\circ {\Gamma }^{\nu }_0 ) ( {\widetilde h}_{\omega } ) {\bigr)}_{\omega } \nearrow ( \underline\pi\circ {\Gamma }^{\nu } ) ( \widetilde f )\, $ which follows from the argument above since $\, s_0^{\nu } \bigl( \Lambda ( \underline{\mathcal B}_{\omega } )_c \bigr) = s^{\nu }_0 \bigl( \Lambda ( \underline{\mathcal B}_{\omega } ) \bigr)\, $ and $\, {\Gamma }^{\nu }\, $ is increasing normal for images of elements in $\, \mathfrak L ( R )\, $ and $\, \underline\pi\, $ agrees with $\,\pi\, $ for images of antibasic elements the latter which factors over the function representation $\, s^{\nu }_0\, $ from Theorem 1 below. The argument in the general case where $\, B\, $ is a direct product of type I factors is just the same on noting that there exists a normal ($w^*$-continuous) retraction 
$\, \upsilon : \mathcal B ( \mathcal H ) \twoheadrightarrow B\, $ inducing an affine $w^*$-continuous embedding $\, {\upsilon }_* : {\mathcal P}^{\nu } ( B ) \hookrightarrow {\mathcal P}^{\nu } \bigl( \mathcal B ( \mathcal H ) \bigr)\, $ of the pure normal states. If $\, \mathfrak X = R \subseteq \mathcal B ( \mathcal H )\, $ is a direct product of type I factors the result above specializes to show that the minimal monotonous extension 
$\, {\pi }^{\nu }_c :\, {\mathfrak L}_q^{\nu } \bigl( R \bigr) \twoheadrightarrow\, {\mathfrak L}_1 \bigl( R \bigr)\, $ of the identity map of $\, R\, $ is boundedly basically increasing normal. 
\par\noindent
Assume in particular given a monotone increasing net $\, \bigl\{ \underline{\mathcal A}_{\lambda }^{ccc} {\bigr\}}_{\lambda } \subseteq \mathfrak L ( \mathfrak X )\, $ of complemented basic elements 
with supremum $\, \mathfrak A\in {\mathfrak L}_q ( \mathfrak X )\, $. Then the maximal representative of the basic element $\, {\mathfrak A}^c\, $ is given by the intersection of the maximal representatives $\, \bigl\{ {\mathcal A}_{\lambda }^{ccc} \bigr\}\, $ of the complemented elements $\, \bigl\{ \underline{\mathcal A}_{\lambda }^{ccc} \bigr\}\, $ which again is complemented and gives the supremum of the corresponding image net $\, \bigl\{ \pi \bigl( \underline{\mathcal A}_{\lambda } \bigr) \bigr\}\, $ in $\, {\mathfrak L}_1 \bigl( \mathfrak X \bigr)\, $, i.e. 
$$ {\pi }^c \bigl( \mathfrak A \bigr)\> =\> \pi \bigl( {\mathfrak A}^c \bigr)\> =\> \sup_{\lambda }\, \pi \bigl( \underline{\mathcal A}_{\lambda } \bigr)\> =\> {\pi }_c \bigl( \mathfrak A \bigr) \> . $$
Since $\, - \mathfrak A\, $ is the infimum of a monotone decreasing net of complemented antibasic elements one deduces by symmetry $\, {\pi }^c ( \mathfrak A ) = - {\pi }_c \bigl( - \mathfrak A \bigr) = - {\pi }^c \bigl( - \mathfrak A \bigr)\, $.
Then 
$$ \pm {\pi }^c \bigl( \mathfrak A \bigr)\> +\> {\pi }^c \bigl( \mathfrak B \bigr)\> =\> - {\pi }^c \bigl( \mp \mathfrak A \bigr)\> +\> {\pi }^c \bigl( \mathfrak B \bigr)\> \leq\> {\pi }^c \bigl( \pm \mathfrak A\, +\, \mathfrak B \bigr)\> \leq \pm {\pi }^c \bigl( \mathfrak A \bigr)\> +\> {\pi }^c \bigl( \mathfrak B \bigr) $$
from convexity of $\, {\pi }^c\, $ and the case of $\, {\pi }_c\, $ follows by symmetry\qed
\par\bigskip\noindent
{\it Remark and Definitions.} (i) 
Let $\, {\lambda }_u : A \rightarrow \mathcal B ( {\mathcal H}_u )\, $ denote the universal representation of the abstract $C^*$-algebra $\, A\, $ with associated (canonical) 
order isomorphic embedding $\, {\Lambda }_u : {\mathfrak L}_q ( A ) \hookrightarrow 
{\mathfrak L}_q ( \mathcal B ( {\mathcal H}_u ) )\, $ (see Theorem 1). Writing $\, L_0^{\nu } ( \mathcal B ( {\mathcal H}_u ) )\, =\, {\mathfrak L}_1 \bigl( \mathfrak L ( \mathcal B ( {\mathcal H}_u ) )\bigm/ \ker\, r^{\nu } \bigr)\, $ there is a positive linear map 
$\, {\pi }^{\nu }_0 : {\mathfrak L}_q ( \mathcal B ( {\mathcal H}_u ) ) \twoheadrightarrow L^{\nu }_0 ( \mathcal B ( {\mathcal H}_u ) )\, $ extending the natural quotient map modulo the kernel of 
$\, r^{\nu }\, $. Consider the composite map
$$ {\Gamma }_u :\> {\mathfrak L}_q ( A ) \buildrel {\Lambda }_u\over\largerightarrow\> {\mathfrak L}_q ( \mathcal B ( {\mathcal H}_u ) )\>\buildrel {\pi }^{\nu }_0\over\largerightarrow\> 
L^{\nu }_0 ( \mathcal B ( {\mathcal H}_u ) $$  
which is injective and order isomorphic (any state of $\, A\, $ is represented by a vector state of $\, \mathcal B ( {\mathcal H}_u )\, $). Then again $\, {\Gamma }_u\, $ admits a positive linear retraction $\, \Psi : L_0^{\nu } ( \mathcal B ( {\mathcal H}_u ) ) \twoheadrightarrow {\mathfrak L}_q ( A )\, $
and the (commutative) $C^*$-product in $\, {\mathfrak L}_q ( A )\, $ may be defined recurring to the product in $\, L^{\nu }_0 ( \mathcal B ( {\mathcal H}_u ) )\, $ by the formula 
$\, \mathfrak P\cdot \mathfrak P = \Psi \bigl( {\Gamma }_u ( \mathfrak P )\cdot {\Gamma }_u ( \mathfrak P ) \bigr)\, $. For a basic positive element $\, \underline{\mathcal D}\in \mathfrak L ( \mathcal B ( {\mathcal H}_u ) )\, $ Theorem 1 below gives $\, {\pi }^{\nu }_0 ( \underline{\mathcal D} )\cdot {\pi }^{\nu }_0 ( \underline{\mathcal D} ) = {\pi }^{\nu }_0 ( \underline{\mathcal D}^2 )\, $. 
Therefore if $\, \underline{\mathcal C}\in {\mathfrak L} ( A )\, $ is a basic positive element then 
$\, \underline{\mathcal C}\cdot \underline{\mathcal C} \leq \underline{\mathcal C}^2\, $ since 
$\, \underline{\mathcal C}\cdot \underline{\mathcal C} = \Psi ( 
( \underline{\mathcal C}^{\gamma } )^2 )\, $ where 
$\, \underline{\mathcal C}^{\gamma } = {\Gamma }_u ( \underline{\mathcal C} ) \, $
which is represented by the basic positive element 
$\, \inf\, \{ {\lambda }_u ( c )\,\vert\, c\in \underline{\mathcal C} \}\, $. Then the element 
$\, ( \underline{\mathcal C}^{\gamma } )^2\, $ is smaller or equal to the element represented by 
the positive set $\, \{ {\lambda }_u ( c )^2 = {\lambda }_u ( c^2 )\,\vert\, c\in \underline{\mathcal C} \}\, $
which by the property of $\, \Psi\, $ being a retraction for $\, {\Gamma }_u\, $ is mapped back to $\, \underline{\mathcal C}^2\, $. 
\par\smallskip\noindent
(ii) For a $C^*$-algebra $\, A\, $ let $\, {\mathfrak I}_A\subseteq \mathfrak L ( A )\, $ denote the order ideal which is the restriction of the order ideal in $\, {\mathfrak L}_q ( A )\, $ generated by all positive elements $\, \{ \underline{\mathcal C}^2\, -\, \underline{\mathcal C}\cdot \underline{\mathcal C} \}\, $ with $\, \underline{\mathcal C}\in \mathfrak L ( A )\, $ a positive basic element. 
For simplicity $\, {\mathfrak I}_{\mathcal B ( \mathcal H )}\, $ will usually be abbreviated to 
$\, \mathfrak I\, $. Define $\, {\mathfrak L}_{\mathfrak I} ( A )\, $ to be the ${\mathfrak L}_1$-completion of the quotient lattice $\, \mathfrak L ( A )\bigm / {\mathfrak I}_A\, $ and
$$ {\pi }_{\mathfrak I} :\> {\mathfrak L}_q ( A )\> \longrightarrow\> {\mathfrak L}_{\mathfrak I} ( A )\> =\> 
{\mathfrak L}_1 \left( \mathfrak L ( A )\bigm / {\mathfrak I}_A \right) $$
any chosen positive linear extension of the natural quotient map. Then the basic squaring operation $\, (1)\, $ drops to 
$\, {\mathfrak L}_{\mathfrak I} ( A )\, $ where it coincides with the $C^*$-square. At this point it is not clear whether this quotient lattice is nontrivial or what it looks like for a general $C^*$-algebra but it will be shown in Theorem 1 that the quotient map to $\, {\mathfrak L}_{\mathfrak I} ( A )\, $ factors the quotient map to $\, {\mathfrak L}_1 ( A )\, $ so the definition makes good sense.
\par\smallskip\noindent
(iii) Let $\, \mathfrak X\subseteq \mathcal B ( \mathcal H )\, $ be an operator subsystem and given any collection of vector states $\, \mathcal F \subseteq {\mathcal S}^{\nu }_0\, $ the lattice 
$\, {\mathfrak L}_1 \Bigl( \mathfrak L ( \mathfrak X ) / {\mathfrak J}_{0 , \mathcal F} \Bigr)\, $ embeds into $\, {\mathfrak L}_1 \Bigl( \mathfrak L \bigl( \mathcal B ( \mathcal H ) \bigr) / {\mathfrak J}_{0 , \mathcal F} \Bigr)\, $ which is isomorphic to $\, l_{\infty } ( \mathcal F )\, $ by Theorem 1 below. Let $\, {\pi }^{\mathfrak X}_{0 , \mathcal F} : \mathfrak L ( \mathfrak X ) \twoheadrightarrow \mathfrak L ( \mathfrak X ) \bigm / {\mathfrak J}_{0 , \mathcal F}\, $ denote the quotient map. For each $\, \rho\in\mathcal F\, $ consider the basic element 
$\, \underline{\mathcal E}_{\rho , \mathfrak X}\in \mathfrak L ( \mathfrak X )\, $ which is the infimum of the set 
$\, \{ c\in \mathfrak X\, \vert\, c\geq 0\, ,\, \rho ( c ) \geq 1 \}\, $. Define a monotonous increasing normal map 
$$ j_{\mathfrak X} :\> l_{\infty } ( \mathcal F )\>\longrightarrow\> {\mathfrak L}_q ( \mathfrak X ) $$
by the assignment
$$  \alpha\quad\mapsto\quad j_{\mathfrak X} ( \alpha ) = \sup_{\rho\in\mathcal F}\, \alpha ( \rho )\, \underline{\mathcal E}_{\rho , \mathfrak X} $$
if $\, \alpha\in l_{\infty } ( \mathcal F )\, $ is a realvalued function. 
Restricted to images of positive basic elements $\, j_{\mathfrak X}\, $ is a lift for 
$\, {\pi }_{0 , \mathcal F}^{\mathfrak X}\, $. To see this note that for a positive basic element 
$\, \underline{\mathcal C}\in \mathfrak L ( \mathfrak X )_+\, $ one has $\, \sup_{\rho }\, \bigl\{ r^{\nu }_{\rho } ( \underline{\mathcal C} )\, \underline{\mathcal E}_{\rho , \mathfrak X}  \bigr\} \leq \underline{\mathcal C}\, $ by construction of the basic elements 
$\, \{ \underline{\mathcal E}_{\rho , \mathfrak X} \}\, $ so the result follows from monotonicity.
\par\smallskip\noindent
(iv) One often encounters the following situation: the quotient map $\, \mathfrak L ( \mathfrak V ) \twoheadrightarrow \mathfrak L ( \mathfrak V )\, /\, \mathfrak J\, $ modulo an order ideal $\, \mathfrak J\subseteq \mathfrak L ( \mathfrak V )\, $ extends (mostly nonuniquely) to a $*$-homomorphism of ${\mathfrak L}_1$-completions 
$$ {\mathfrak L}_q \bigl( \mathfrak V \bigr)\> \largerightarrow\> {\mathfrak L}_1 \left( \mathfrak L \bigl( \mathfrak V \bigr)\,\bigm/\, \mathfrak J \right) $$
which however need not be surjective unless the original quotient map admits a monotonous splitting. If this is not the case we say that the extended map is an {\it approximate quotient map}. 
\par\bigskip\noindent
{\bf Theorem 1.}\quad 
For any order isomorphic unital embedding of function systems $\, \mathfrak V \subseteq \mathfrak W \, $  the induced functorial positive linear map
$\, \mathfrak L ( \mathfrak V ) \hookrightarrow \mathfrak L ( \mathfrak W )\, $ extends uniquely to an injective normal $*$-homomorphism 
$$ {\mathfrak L}_q ( \mathfrak V )\> \longrightarrow\> {\mathfrak L}_q ( \mathfrak W ) \> . $$
Let $\, \widetilde{\mathfrak L} ( \mathfrak V )\subseteq {\mathfrak L}_q ( \mathfrak V )\, $ denote the linear subspace generated by $\, \mathfrak L ( \mathfrak V )\, $ plus arbitrary suprema of sets of uniformly boundedly generated basic elements.
Any monotonous extension 
$$ {\pi }_q :\> {\mathfrak L}_q \bigl( \mathfrak V \bigr)\> \largerightarrow\> {\mathfrak L}_1 \bigl( \mathfrak V \bigr) $$
of $\, \pi\, $ is uniquely determined and linear restricted to $\, \widetilde{\mathfrak L} \bigl( \mathfrak V \bigr)\, $ and is boundedly basically increasing normal (resp. boundedly antibasically decreasing normal). 
For any state $\, \sigma\in \mathcal S ( \mathfrak V )\, $ the induced functional 
$$ s_{\sigma } :\> \mathfrak L ( \mathfrak V )\> \longrightarrow\> \mathbb R $$
extends (nonuniquely) to a multiplicative functional on $\, {\mathfrak L}_q ( \mathfrak V )\, $.  
\par\noindent
Assume that $\,  \mathfrak X \subseteq \mathcal B ( \mathcal H )\, $ is an operator subsystem. Let $\, \mathcal F\subseteq {\mathcal S}_0^{\nu }\, $ be any collection of vector states and $\, {\mathfrak J}_{0 , \mathcal F} \subseteq \mathfrak L ( \mathcal B ( \mathcal H ) )\, $ the order ideal as defined above.
Put $\, L_0^{\mathcal F } ( \mathfrak X )\, =\, {\mathfrak L}_1 \bigl( \mathfrak L ( \mathfrak X )\, /\, {\mathfrak J}_{0 , \mathcal F , \mathfrak X} \bigr)\, $ where 
$\, {\mathfrak J}_{0 , \mathcal F , \mathfrak X} = {\mathfrak J}_{0 , \mathcal F} \cap \mathfrak L ( \mathfrak X )\, $. In case that $\, \mathcal F = {\mathcal S}_0^{\nu }\, $ we simply write 
$\, L^{\nu }_0 ( \mathfrak X )\, $ etc.. Then $\, {\mathfrak J}_{0 , \mathcal F}\, $ is equal to the intersection of all kernels $\, \bigcap_{\rho\in\mathcal F}\, \ker\,r^{\nu }_{\rho }\, $.  For any inclusion of operator subsystems $\, \mathfrak X \subseteq \mathfrak Y\subseteq \mathcal B ( \mathcal H )\, $ there exist $*$-homomorphic extensions of the canonical embeddings/quotients of $\, \mathfrak L ( \mathfrak X )\, ,\, \mathfrak L ( \mathfrak Y )\, ,\, \mathfrak L ( \mathcal B ( \mathcal H ) )\, $ making a commutative diagram
$$ \vbox{\halign{ #&#&#\cr
\hfil ${\mathfrak L}_q ( \mathfrak X )$\hfil &\hfil $\buildrel \mathcal I\over\largerightarrow$\hfil &\hfil ${\mathfrak L}_q ( \mathfrak Y )$\hfil \cr
\hfil ${\pi }_{0 , \mathfrak X}^{\mathcal F}\Biggm\downarrow\quad $\hfil &&\hfil ${\pi }_{0 , \mathfrak Y}^{\mathcal F}\Biggm\downarrow\quad $\hfil \cr
\hfil $L_0^{\mathcal F} ( \mathfrak X )$\hfil &\hfil $\buildrel {\mathcal I}_0^{\mathcal F}\over\largerightarrow $\hfil &\hfil $L_0^{\mathcal F} ( \mathfrak Y )$\hfil \cr }}  $$
with $\,\mathcal I\, $ the functorial normal $*$-homomorphism induced by the inclusion 
$\, \mathfrak X \subseteq \mathfrak Y\, $. The horizontal maps are injective and the vertical maps are approximate quotient maps (surjective $*$-homomorphisms if $\,\mathfrak X\, $ is separating for $\,\mathcal F\, $). One has
$\, L_0^{\mathcal F} ( \mathcal B ( \mathcal H ) ) \simeq l_{\infty } ( \mathcal F )\, $ and in case that $\,\mathfrak X\, $ is separating for $\,\mathcal F\, $ the map $\, {\mathcal I}_0^{\mathcal F}\, $ is an isomorphism $\, L_0^{\mathcal F } ( \mathfrak X ) \simeq L_0^{\mathcal F} ( \mathcal B ( \mathcal H ) )\, $.  The surjection $\,\pi : {\mathfrak L}_q ( \mathfrak X ) \twoheadrightarrow {\mathfrak L}_1 ( \mathfrak X )\, $ factors over 
$\, {\pi }_{0 , \mathfrak X}^{\nu }\, $. 
\par\noindent
Moreover if $\, A\subseteq \mathcal B ( \mathcal H )\, $ is an irreducible 
$C^*$-subalgebra the image of the basic squaring operations $\, ( 1 )\, ,\, ( 13 )\, $ in 
$\, L^{\nu }_0 ( A )\, $ agrees with the $C^*$-square and the $C^*$-product $\, \mathcal X\cdot {\mathcal X}^*\, $ for any element of the form $\, \mathcal X = \underline{\mathcal C}\, +\, i\, \underline{\mathcal D}\, $ respectively.
More generally for any (abstract) separable $C^*$-algebra $\, A\, $ the $\, {\mathfrak L}_1$-completion
$\, {\mathfrak L}^r_0 ( A ) = {\mathfrak L}_1 \bigl( r \bigl( \mathfrak L ( A ) \bigr) \bigr)\, $ of the quotient of 
$\, \mathfrak L ( A )\, $ modulo the kernel of 
$$ r : \mathfrak L ( A )\> \longrightarrow\> l_{\infty } \bigl( \mathcal P ( A ) \bigr) $$
is a commutative $W^*$-algebra isomorphic to $\, l_{\infty } ( \mathcal P ( A ) )\, $ such that the basic squaring operations 
$\, ( 1 )\, ,\, ( 13 )\, $ drop to this algebra where they coincide with the corresponding $C^*$-squares and -products. The canonical surjection $\, \mathfrak L ( A ) \twoheadrightarrow r ( \mathfrak L ( A ) )\, $ extends (nonuniquely) to a surjective $*$-homomorphism 
$$ {\pi }_0^r : {\mathfrak L}_q ( A )\> \largerightarrow {\mathfrak L}^r_0 ( A )\>\simeq l_{\infty } ( \mathcal P ( A ) ) $$
factoring $\, \pi\, $ and factoring over $\, {\pi }_{\mathfrak I}\, $.
\par\noindent
Suppose that $\, I\, $ is an injective $C^*$-algebra. For any (equivalence class of) multiplicative injective representation 
$\, \mu : I \longrightarrow \mathcal B ( \mathcal H )\, $ there exists a canonical order ideal 
$\, {\mathfrak J}^{\mu }_0 \subseteq  \mathfrak L ( I )\, $ such that the basic squaring operations 
$\, ( 1 )\, ,\, ( 13 )\, $ drop to the quotient modulo the ideal $\, {\mathfrak J}^{\mu }_0 \, $ which is a $C^*$-algebra where they coincide with the image of the corresponding $C^*$-squares and -products. The quotient map may be extended (nonuniquely) to a $*$-homomorphism
$$ {\pi }_0^{\mu } :\> {\mathfrak L}_q ( I )\> \largerightarrow\> {\mathfrak L}^{\mu }_0 ( I )\> =\> {\mathfrak L}_1 \bigl( \mathfrak L ( I )\, /\, {\mathfrak J}^{\mu }_0 \bigr)\> .   $$
In particular let $\, {\mathfrak J}_0^m\, =\, \bigcap_{\mu }\, {\mathfrak I}^{\mu }_0\, $ where 
$\, \{ \mu \}\, $ ranges over all multiplicative injective representations of $\, I\, $ up to equivalence 
and put 
$$ {\mathfrak L}_0^m ( I )\> =\> {\mathfrak L}_1 \bigl( \mathfrak L ( I )\, /\, {\mathfrak J}_0^m \bigr)   \> . $$
The canonical quotient maps modulo $\, {\mathfrak J}_0^m\, $ may be extended (nonuniquely) to a  $*$-homomorphism 
$$ {\pi }_0^m :\> {\mathfrak L}_q ( I )\> \largerightarrow\> {\mathfrak L}_0^m ( I ) $$
 factoring $\, \pi\, $ and factoring over $\, {\pi }_{\mathfrak I}\, $. 
\par\noindent
For any $C^*$-algebra $\, A\, $ the Jordan squaring operations 
$\, ( 1 )\, ,\, ( 9 )\, ,\, ( 10 )\, ,\, ( 11 )\, ,\, ( 13 )\, $ drop to $\, {\mathfrak L}_1 ( A )\, $ where they agree with the corresponding $C^*$-squares and -products, i.e. for any $\, \mathfrak A\in {\mathfrak L}_q ( A )\, $ one has $\, \pi   ( \mathfrak A )^2\, =\, \pi  ( {\mathfrak A}^2 )\, $. 
Moreover if $\, {\mathfrak P}^{1/2}\, $ denotes the $C^*$-squareroot of the positive element 
$\, \mathfrak P\geq 0\, $ and $\, \sqrt{\mathfrak P}\, $ is the Jordan lattice squareroot of $\, \mathfrak P\, $ as above then 
$$  \pi ( {\mathfrak P}^2 )\> =\> \pi ( \mathfrak P )^2\> ,\quad \pi ( \sqrt{\mathfrak P} )\> \leq \> \pi  ( {\mathfrak P}^{1/2} )\> .  \leqno{( 21 )} $$
\par\noindent
Also for any selfadjoint element $\, a\in A\, $
and $\, \mathfrak a = \{ a \} = \underline{\mathfrak a}_+ - \underline{\mathfrak a}_-\, $ the minimal basic positive decomposition of $\, a\, $ in $\, {\mathfrak L} ( A )\, $ the image of $\, \underline{\mathfrak a}_+ \wedge \underline{\mathfrak a}_-\, $ in $\, {\mathfrak L}_1 ( A )\, $ is trivial (so $\, \pi   ( \mathfrak a ) = 
\pi   ( \underline{\mathfrak a}_+ ) - \pi   ( \underline{\mathfrak a}_- )\, $ is the minimal positive decomposition of $\, a\, $ in $\, {\mathfrak L}_1 ( A )\, $).
\par\noindent
\par\bigskip\noindent
{\it Proof.}\quad 
Let $\, \mathfrak V \subseteq \mathfrak W\, $ be an order isomorphic unital embedding of function systems and $\, {\iota }_q : {\mathfrak L}_q ( \mathfrak V ) \hookrightarrow 
{\mathfrak L}_q ( \mathfrak W )\, $ any monotonous extension of the functorial linear map 
$\, \iota : \mathfrak L ( \mathfrak V ) \hookrightarrow \mathfrak L ( \mathfrak W )\, $. Let
$\, \mathfrak A\, =\, \sup_{\lambda }\, \{ \underline{\mathcal A}_{\lambda } \}\in {\mathfrak L}_q ( \mathfrak V )\, $ be the supremum of certain basic elements and note that since any element of 
$\, {\mathfrak L}_q ( \mathfrak V )\, $ is the supremum of elements in $\, \mathfrak L ( \mathfrak V )\, $ and every element 
$$ \mathcal A\> =\> \underline{\mathcal C}\, - \underline{\mathcal D}\> =\> 
\sup\, \bigl\{ \underline{\mathcal C}\> -\> \sup\, \{ d \}\,\bigm\vert\, d\in \mathcal D \bigr\} $$ 
is the supremum of basic elements, the same is true for an arbitrary element 
$\, \mathfrak A\in {\mathfrak L}_q ( \mathfrak V )\, $. If $\, \underline{\mathcal A}_{\lambda }^{\iota }\, =\, 
\iota ( \underline{\mathcal A}_{\lambda } )\, $ consider the supremum 
$\, \overline{\mathfrak A}\, =\, \sup_{\lambda }\, \{ \underline{\mathcal A}_{\lambda }^{\iota } \} \in 
{\mathfrak L}_q ( \mathfrak W )\, $. Then $\, \overline{\mathfrak A}\, $ is the infimum of all {\it antibasic} elements 
$$ \overline{\mathfrak A}\> =\> \inf\, \bigl\{ \overline{\mathcal B}_{\mu }\,\bigm\vert\, \overline{\mathcal B}_{\mu }\in \mathfrak L ( \mathfrak W )\, ,\, \overline{\mathfrak A}\,\leq\, \overline{\mathfrak B}_{\mu } \bigr\} $$
larger or equal to $\, \overline{\mathfrak A}\, $. From monotonicity one gets 
$\, \overline{\mathfrak A}\,\leq\, {\iota }_q ( \mathfrak A )\, $ so that any monotonous retraction 
$\,  r_q : {\mathfrak L}_q ( \mathfrak W ) \twoheadrightarrow {\mathfrak L}_q ( \mathfrak V )\, $ for 
$\, {\iota }_q\, $ must satisfy $\, r_q ( \overline{\mathfrak A} )\, =\, \mathfrak A\, $. Consider the functorial retraction $\, \overline r\, $ induced by restriction of antibasic elements as in part (ii) of the Remark after Proposition 1. Then $\, ( \iota\circ \overline r ) ( \overline{\mathcal B}_{\mu } )\, \leq\, \overline{\mathcal B}_{\mu }\, $ for each index $\, \mu\, $ which by monotonicity implies 
$$ {\iota }_q ( \mathfrak A )\> =\>  ( {\iota }_q\circ \overline r ) ( \overline{\mathfrak A} )\> \leq\> 
\overline{\mathfrak A} $$
hence equality. Thus the extension $\, {\iota }_q\, $ is necessarily unique and (increasing) normal since the image of the supremum of an arbitrary set $\, \{ {\mathfrak A}_{\lambda } \}\, $ is seen to agree with the image of the supremum of a corresponding subset of basic elements, each of which is dominated by some 
$\, {\mathfrak A}_{\lambda }\, $ and this image in turn agrees with the supremum of the images. 
Similarly one finds that $\, {\iota }_q\, $ is decreasing normal, hence normal. Now let a state $\, \sigma\in \mathcal S ( \mathfrak V )\, $ be given. It is easy to see that the kernel of the map 
$$ s_{\sigma } :\> \mathfrak L ( \mathfrak V )\> \longrightarrow\> \mathbb R $$
is positively generated. Indeed, if $\, s_{\sigma } ( \underline{\mathcal C} )\, =\, s_{\sigma } ( \underline{\mathcal D} )\, $ then this common value is equal to $\, s_{\sigma } ( \underline{\mathcal C}\wedge \underline{\mathcal D} )\, $ so that the positive elements 
$\, \underline{\mathcal C}\, -\, \underline{\mathcal C}\wedge \underline{\mathcal D}\, $ and 
$\, \underline{\mathcal D}\, -\, \underline{\mathcal C}\wedge \underline{\mathcal D}\, $ are contained in 
the kernel of $\, s_{\sigma }\, $. Thus $\, s_{\sigma }\, $ is a lattice map which from Theorem A of section 2 can be extended to a corresponding lattice map $\, {\overline s}_{\sigma } : {\mathfrak L}_q ( \mathfrak V ) \twoheadrightarrow \mathbb R\, $ the kernel of which is again positively generated hence must be a $*$-ideal of the underlying commutative $C^*$-algebra (since the kernel of a positive linear map of $C^*$-algebras is a $*$-ideal whenever positively generated). Then 
$\, {\overline s}_{\sigma }\, $ is uniquely determined on the sub-$C^*$-algebra generated by 
$\, \mathfrak L ( \mathfrak V )\, $. From this one easily sees that the normal linear injection $\, {\iota }_q :  {\mathfrak L}_q ( \mathfrak V ) \hookrightarrow {\mathfrak L}_q ( \mathfrak W )\, $ as above which certainly is a lattice map is in fact a $*$-homomorphism since every state of 
$\, \mathfrak W\, $ determines a state of $\, \mathfrak V\, $ so that the composition 
$$ \mathfrak L \bigl( \mathfrak V \bigr)\> \hookrightarrow\> \mathfrak L \bigl( \mathfrak W \bigr)\> \buildrel s\over\largerightarrow\> l_{\infty } \bigl( \mathcal S ( \mathfrak W ) \bigr) $$
extends to a multiplicative embedding
$$ {\mathfrak L}_q \bigl( \mathfrak V \bigr)\> \largerightarrow\> l_{\infty } \bigl( \mathcal S ( \mathfrak W ) \bigr)   $$
together with a compatible multiplicative (lattice map) extension 
$$ {\mathfrak L}_q \bigl( \mathfrak W \bigr)\> \largerightarrow\> l_{\infty } \bigl( \mathcal S ( \mathfrak W ) \bigr)  $$
from Theorem A of section 2.
\par\smallskip\noindent
Now let $\, \mathfrak X \subseteq \mathcal B ( \mathcal H )\, $ be an operator subsystem. From Lemmas 1 and 3 one gets that 
$\, L^{\nu }_{00} ( \mathcal B ( \mathcal H ) ) = \mathfrak L ( \mathcal B ( \mathcal H ) ) / {\mathfrak J}_0\, $ admits the structure of an associative algebra with corresponding Jordan product induced by the squaring operation of basic positive elements, i.e. 
$$ {\pi }_{00} ( \underline{\mathcal C} )\cdot {\pi }_{00} ( \underline{\mathcal D} )\> +\> 
{\pi }_{00} ( \underline{\mathcal D} )\cdot {\pi }_{00} ( \underline{\mathcal C} )\> =\> 
{\pi }_{00} ( ( \underline{\mathcal C} + \underline{\mathcal D} )^2\> -\> \underline{\mathcal C}^2\> -\> \underline{\mathcal D}^2 ) \> , $$
and that the induced map 
$\, l : L^{\nu }_{00} ( \mathcal B ( \mathcal H ) ) \rightarrow l_{\infty } ( {\mathcal P}^{\nu } ( \mathcal B ( \mathcal H ) ) )\, $ is an algebra homomorphism. Since the range is a commutative $C^*$-algebra and $\, l\, $ is injective one finds that $\, L^{\nu }_{00} ( \mathcal B ( \mathcal H ) )\, $ is commutative. If 
$\, r^{\nu } ( \underline{\mathcal C} )\geq r^{\nu } ( \underline{\mathcal D} )\, $ then 
$\, {\pi }_{00} ( \underline{\mathcal D} ) = {\pi }_{00} ( \underline{\mathcal C} \wedge \underline{\mathcal D} )\, $ which shows that $\, l\, $ is order isomorphic, so as far as the order and algebra structure are concerned $\, L^{\nu }_{00} ( \mathcal B ( \mathcal H ) )\, $ may be identified with a dense subalgebra of a selfadjoint real commutative $C^*$-algebra, i.e. realvalued continuous functions on a compact space 
$\, X\, $. Since the norm is monotone for positive elements and $\, l\, $ is contractive as well as order isomorphic one deduces the inequalities 
$$ - \Vert \overline{\mathcal A} \Vert \cdot \Vert \overline{\mathcal B} \Vert\>\leq\> \overline{\mathcal A}\cdot \overline{\mathcal B}\> \leq\> \Vert \overline{\mathcal A} \Vert\cdot \Vert \overline{\mathcal B} \Vert 
$$
from the corresponding relations in $\, C ( X )\, $ where $\, \overline{\mathcal A}\,,\, \overline{\mathcal B} \in L^{\nu }_{00} ( \mathcal B ( \mathcal H ) )\, $ are arbitrary selfadjoint elements. This implies 
$\, \Vert \overline{\mathcal A} \cdot \overline{\mathcal B} \Vert \leq \Vert \overline{\mathcal A} \Vert\cdot
\Vert \overline{\mathcal B} \Vert\, $, i.e. (the norm completion of) $\, L^{\nu }_{00} ( \mathcal B ( \mathcal H ) )\, $ is a Banach algebra. The inequality
$$ - \Vert \overline{\mathcal P} + \overline{\mathcal Q} \Vert\>\leq\>  \overline{\mathcal P} - \overline{\mathcal Q}\> \leq\> \Vert \overline{\mathcal P} + \overline{\mathcal Q} \Vert $$
valid for any two positive elements $\, \overline{\mathcal P}\, ,\, \overline{\mathcal Q} \geq 0\, $ implies 
$\, \Vert \overline{\mathcal P} - \overline{\mathcal Q} \Vert\leq \Vert \overline{\mathcal P} + \overline{\mathcal Q} \Vert\, $ so that from Theorem 4.2.5 of \cite{A-K} $\, L^{\nu }_{00} ( \mathcal B ( \mathcal H ) )\, $ is isometrically isomorphic to a commutative $C^*$-algebra. Then also 
$\, l\, $ being injective must be an isometry, so 
$\, L^{\nu }_{00} ( \mathcal B ( \mathcal H ) ) \simeq C_{\mathbb R} ( X ) \subseteq l_{\infty } ( {\mathcal S}^{\nu }_0 )\, $. From injectivity of $\, l_{\infty } ( {\mathcal S}^{\nu }_0 )\, $ there exists a monotonous (of course even linear positive) extension 
$\, \overline l : L_0^{\nu } \bigl( \mathcal B ( \mathcal H ) \bigr) = {\mathfrak L}_1 \bigl( L_{00}^{\nu } \bigl( \mathcal B ( \mathcal H ) \bigr) \bigr) \longrightarrow l_{\infty } ( {\mathcal S}^{\nu }_0 )\, $ of $\, l\, $. Every pure normal state 
$\, \rho\in {\mathcal S}^{\nu }_0\, $ defines a corresponding pure state 
$\, {\overline r}_{\rho }\, $ of $\, L^{\nu }_{00} ( \mathcal B ( \mathcal H ) )\, $. Considering the 
dual space $\, V^* = l_1 \bigl( {\mathcal S}^{\nu }_0 \bigr)\, $ of $\, V = c_0 \bigl( {\mathcal S}^{\nu }_0 \bigr)\, $ one has $\, V^{**} \simeq  
l_{\infty } \bigl( {\mathcal S}^{\nu }_0 \bigr)\, $.
Then the extension 
$\, \overline l :  L_0^{\nu } ( \mathcal B ( \mathcal H ) )\> \longrightarrow\> V^{**}\, $ is necessarily surjective.
Namely if $\, x\in V^{**}_+\, $ is any positive element then $\, x\, $ is the supremum of the elements 
$\, \{ x ( r^{\nu }_{\rho } ){\eta }_{\rho }  \}\, $ where 
$\, {\eta }_{\rho } ( \sum_\kappa\, {\lambda }_{\kappa }\, r^{\nu }_{\kappa } ) = {\lambda }_{\rho }\, $.
Now each $\, {\eta }_{\rho }\, $ is in the image of $\, L^{\nu }_{00} ( \mathcal B ( \mathcal H ) )\, $ and, if $\, \mathfrak X\subseteq \mathcal B ( \mathcal H )\, $ is transitive also in the image of 
$\, L_{00}^{\nu } ( \mathfrak X ) = \mathfrak L ( \mathfrak X )\, /\, {\mathfrak J}_{0 , \mathfrak X}\, $.
This follows from the fact that for every pure normal state $\, \kappa \neq \rho\, $ and $\, \epsilon > 0\, $ there exists an element $\, x_{\kappa , \epsilon }^{\rho }\in {\mathfrak X}_+\, $ with $\, \kappa ( x_{\kappa , \epsilon}^{\rho } ) \leq \epsilon\, ,\, \rho ( x_{\kappa , \epsilon }^{\rho } ) \geq 1\, $. Putting $\, \underline{\mathcal E}_{\rho , \mathfrak X}\, =\, \inf\, \{ c\in {\mathfrak X}_+\,\vert\, \rho ( c ) \geq 1 \}\, $ the image of the basic positive element $\, \underline{\mathcal E}_{\rho , \mathfrak X}\, $ is clearly equal to $\, {\eta }_{\rho }\, $ (transitive case). Now it is easy to see that the positive injective envelope of the commutative $C^*$-algebra $\, c_0 \bigl( {\mathcal S}_0^{\nu } \bigr)\, $ of functions vanishing at infinity is equal to $\, l_{\infty } \bigl( {\mathcal S}_0^{\nu } \bigr)\, $ whence the latter is also the injective envelope of $\, c_0 \bigl( {\mathcal S}_0^{\nu } \bigr)\subseteq L_{00}^{\nu } \bigl( \mathfrak X \bigr) \subseteq l_{\infty } \bigl( {\mathcal S}^{\nu }_0 \bigr) \, $ (transitive case).
Since any positive linear extension of $\, r^{\nu }\, $ to 
$\, C^* ( \mathfrak L ( \mathcal B ( \mathcal H ) ) )\subseteq {\mathfrak L}_q ( \mathcal B ( \mathcal H ) )\, $ must be a $C^*$-homomorphism from the Schwarz inequality and the fact that 
$\, \underline{\mathcal C}\cdot \underline{\mathcal C} \leq \underline{\mathcal C}^2\, $ for a positive basic element there always exists a multiplicative positive linear extension to 
$\, {\mathfrak L}_q ( \mathcal B ( \mathcal H ) )\, $ by the Corollary of Theorem 2. We denote any chosen 
(nonunique) such extension by 
$$ {\pi }^{\nu }_0 :\> {\mathfrak L}_q ( \mathcal B ( \mathcal H ) )\> \largerightarrow\> l_{\infty } ( {\mathcal S}^{\nu }_0 )\> =\> L_0^{\nu } ( \mathcal B ( \mathcal H ) ). $$
and note that it restricts canonically to a $*$-homomorphism $\, {\mathfrak L}_q ( \mathfrak X ) \rightarrow L_0^{\nu } ( \mathcal B ( \mathcal H ) )\, $ via the unique normal multiplicative embedding $\, {\mathfrak L}_q ( \mathfrak X ) \hookrightarrow {\mathfrak L}_q ( \mathcal B ( \mathcal H ) )\, $ yielding a $*$-homomorphic injection $\, L_{00 }^{\nu } ( \mathfrak X ) \hookrightarrow L_{00}^{\nu } ( \mathcal B ( \mathcal H ) )\, $ where $\, L_{00}^{\nu } ( \mathfrak X )\, $ denotes the image of 
$\, C^* \bigl( \mathfrak L ( \mathfrak X ) \bigr)\subseteq {\mathfrak L}_q \bigl( \mathfrak X \bigr)\, $. 
Choose any $*$-homomorphic extension 
$\, L_0^{\nu } ( \mathfrak X ) \hookrightarrow L_0^{\nu } ( \mathcal B ( \mathcal H ) )\, $ of this map. Then it admits an extremal retraction $\, p_{\mathfrak X} : L_0^{\nu } ( \mathcal B ( \mathcal H ) ) \twoheadrightarrow L_0^{\nu } ( \mathfrak X )\, $ which is a $*$-homomorphism from the Corollary of Theorem 2. Define 
$\, {\pi }_{0 , \mathfrak X}^{\nu }\, $ to be the composition 
$\, p_{\mathfrak X}\circ {\pi }_0^{\nu }\circ \mathcal I\, $. One may if necessary change the original homomorphism $\, {\pi }_0^{\nu }\, $ in order that $\, \mathcal I ( \ker\,{\pi }_{0 , \mathfrak X}^{\nu } ) \subseteq \ker\, {\pi }_0^{\nu }\, $ using the methods as above starting from the canonical $*$-homomorphism 
$$ C^* \left( \mathfrak L \bigl( \mathcal B ( \mathcal H ) \bigr) + \mathcal I \bigl( \ker\, {\pi }_{0 , \mathfrak X}^{\nu } \bigr) \right)\> \largerightarrow\> L_{00} \bigl( \mathcal B ( \mathcal H ) \bigr) $$ sending
$\, \mathcal I \bigl( \ker\, {\pi }_{0 , \mathfrak X}^{\nu } \bigr) + {\mathfrak J}_0\, $ to zero. Then assume this compatibility condition. The image of $\, {\mathfrak L}_q \bigl( \mathfrak X \bigr)\, $ under the composition $\, {\pi }_0^{\nu }\circ \mathcal I\, $ is $*$-isomorphic to 
$\, {\pi }^{\nu }_{0 , \mathfrak X} \bigl( {\mathfrak L}_q ( \mathfrak X ) \bigr)\subseteq L_0^{\nu } \bigl( \mathfrak X \bigr)\, $ so again from the Corollary of Theorem 2 one may find a $*$-homomorphic extension $\, {\mathcal I}_0 : L_0^{\nu } ( \mathfrak X ) \longrightarrow L_0^{\nu } ( \mathcal B ( \mathcal H ) )\, $ making a commutative diagram with the other maps as in the Theorem. The same line of arguments extends to the case 
of $\, L_0^{\nu , \mathcal F} ( \mathfrak X ) \subseteq L_0^{\nu  , \mathcal F} ( \mathcal B ( \mathcal H ) )\, $ for an arbitrary subset $\, \mathcal F\subseteq {\mathcal P}^{\nu } ( \mathcal B ( \mathcal H ) )\, $. In particular $\, L_0^{\nu , \mathcal F} ( \mathcal B ( \mathcal H ) ) \simeq l_{\infty } ( \mathcal F )\, $ and the same holds for $\, L_0^{\nu , \mathcal F} ( \mathfrak X )\, $ in case that $\,\mathfrak X\, $ is separating for the subset $\,\mathcal F \subseteq {\mathcal S}_0^{\nu }\, $ by an argument as before because $\, L_0^{\nu , \mathcal F} ( \mathfrak X )\, $ contains every minimal projection of $\, l_{\infty } \bigl( \mathcal F \bigr)\, $ in this case. In fact the positive part of the subspace $\, c_0 ( \mathcal F )\subseteq l_{\infty } ( \mathcal F )\, $ is contained in the image of the cone of positive basic elements $\, \underline{\mathfrak L} ( \mathfrak X )\subseteq \mathfrak L ( \mathfrak X )\, $ in this case so that  the monotonous map $\, j_{\mathfrak X} :\, l_{\infty } ( \mathcal F )_+ \rightarrow {\mathfrak L}_q ( \mathfrak X )_+\, $ as in part (iii) of the Remark above which extends canonically to a monotonous map on the linear envelopes must be a lift for $\, {\pi }^{\mathcal F}_{0 , \mathfrak X}\, $ from rigidity of the embedding $\, c_0 ( \mathcal F ) \subseteq l_{\infty } ( \mathcal F )\, $ proving that $\, {\pi }^{\mathcal F}_{0 , \mathfrak X}\, $ is surjective. 
The statement that $\, \pi\, $ factors over $\, {\pi }^{\nu }_{0 , \mathfrak X}\, $ is obvious from Proposition 1. This proves the existence of a commutative diagram as above in the case $\, \mathfrak Y = \mathcal B ( \mathcal H )\, $. However the argument in case of an arbitrary inclusion $\, \mathfrak X \subseteq \mathfrak Y\, $ is much the same. 
\par\noindent
Now let $\, A\subseteq \mathcal B ( \mathcal H )\, $ be an irreducible $C^*$-subalgebra. 
Then the natural inclusion 
$\, {\iota }_q : {\mathfrak L}_q ( A ) \hookrightarrow {\mathfrak L}_q ( \mathcal B ( \mathcal H ) )\, $ drops to an isomorphism 
$$ L_0^{\nu } ( A )\, =\, {\mathfrak L}_1 \bigl( \mathfrak L ( A )\,/\, {\mathfrak J}_{0 , A} \bigr)\> \buildrel\sim\over\longrightarrow\> 
L_0^{\nu } ( \mathcal B ( \mathcal H ) )\>\simeq\> l_{\infty } ( {\mathcal S}^{\nu }_0 ) $$ 
and as noted in the proof of Lemma 1 the basic squaring operations $\, ( 1 )\, ,\, ( 13 )\, $ agree with the corresponding $C^*$-squares and -products for the image of any positive basic elements in $\, \mathfrak L ( \mathcal B ( \mathcal H ) )\, $ implying in particular that for any positive basic element $\, \underline{\mathcal C}\in \mathfrak L ( A )\, $
the image of $\, {\iota }_q ( \underline{\mathcal C} )^2\, $ agrees with the image of 
$\, \underline{\mathcal C}\cdot \underline{\mathcal C}\, $ modulo $\, {\mathfrak J}_0\, $.  
Then we only need to show the congruence of $\, {\iota }_q ( \underline{\mathcal C}^2 )\, $ and $\, {\iota }_q ( \underline{\mathcal C} )^2\, $, and similarly 
the congruence of $\, {\iota }_q ( \mathcal X \underline\cdot {\mathcal X}^* )\, $ and 
$\, {\iota }_q ( \mathcal X ) \underline\cdot {\iota }_q ( {\mathcal X} )^*\, $ modulo $\, {\mathfrak J}_0\, $ where $\, \mathcal X = \underline{\mathcal C}\, +\, i\, \underline{\mathcal D}\,\in \mathfrak L ( A )_c\, $. In the first case write 
$$ 0\>\leq\> {\iota }_q ( \underline{\mathcal C}^2 )\> -\> {\iota }_q ( \underline{\mathcal C} )^2\> =\> 
\sup_a\,\left\{ {\iota }_q ( \underline{\mathcal C}^2 )\, -\, \inf\, \{ a^2 \}\,\bigm\vert\, a\in {\iota }_q ( \underline{\mathcal C} )^c \right\} $$
$$\quad =\> \sup_a\, \left\{\, \inf_c\, \left\{\, \inf\, \{  c^2 - a^2 \}\,\bigm\vert\, 
c\in\mathcal C\, \right\}\, \Bigm\vert\, a\in {\iota }_q ( \underline{\mathcal C} )^c\, \right\}\> . $$
If we can show the estimate  
$$ \inf_c\, \left\{\, \inf\, \{ c^2 - a^2 \}\,\bigm\vert\, c\in \mathcal C\,\right\}\>\leq\> \epsilon\, {\bf 1} $$
modulo $\, {\mathfrak J}_0\, $ for each given $\, a\in {\iota }_q ( \underline{\mathcal C} )^c\, $ and 
$\,\epsilon > 0\, $ we get the result since $\,\epsilon \, $ may be chosen arbitrary small and the image function of $\, {\iota }_q ( \underline{\mathcal C}^2 ) - {\iota }_q ( \underline{\mathcal C} )^2\, $ in $\, l_{\infty } ( {\mathcal S}^{\nu }_0 )\, $ is the pointwise supremum of the functions corresponding to basic elements as above fixing $\, a\, $.
Then it is sufficient to check that 
$$ \inf_c\, \left\{\, \inf\, \{ P\, ( c^2 - a^2 )\, P \}\,\bigm\vert\, c\in \overline{\mathcal C}^{\nu }\,\right\}\>\leq\> \inf_c\, \left\{\, \inf\, \{ P\, ( c^2 - a^2 )\, P \}\,\bigm\vert\, c\in \mathcal C\,\right\}\> \leq\>
\epsilon\, P $$
for any finitedimensional orthogonal projection $\, P\in\mathcal B ( \mathcal H )\, $. From the Schwarz inequality one gets $\, P\, a^2\, P \geq ( P a P )^2 = a_P^2\, $ and it is sufficient to 
check 
$$ \inf_c\, \left\{\, \inf\, \{ P c^2 P - a_P^2 \}\,\bigm\vert\, c\in \mathcal C\,\right\}\>\leq\> \epsilon\, P \> . $$
Since $\, a\in {\iota }_q ( \underline{\mathcal C} )^c\, $ there exists $\, c\in\mathcal C\, $ with 
$\, c \leq a\, $ which implies $\,  c_P \leq a_P\, $ and invoking Kadison's  transitivity theorem revisited, Theorem 2.7.5 of \cite{Pe1}, there exists a positive element 
$\, d\in A_+\, $ with $\, c_P + d_P = a_P\, $. Then for given $\,\epsilon > 0\, $ there exists a finitedimensional projection $\, Q \geq P\, $ with 
$\, \Vert P\, ( c + d )\, ( 1 - Q )\Vert \leq \epsilon\, $ and we may choose a positive element 
$\, h_{\epsilon }\in \mathcal B ( \mathcal H )_+\, $ with $\, \Vert P\, h_{\epsilon }\, P\Vert \leq 
{\epsilon\over 3\, \Vert a\Vert }\, $ and $\,  P\, ( c + d + h_{\epsilon } )\, ( Q - P )\, =\, ( Q - P )\, ( c + d + h_{\epsilon } )\, P = 0\, $.  
Once more invoking Theorem 2.7.5 of \cite{Pe1} there exists a positive element 
$\, d_{\epsilon }\in A_+\, $ of norm $\, \Vert h_{\epsilon } \Vert\, $ with $\, Q\, d_{\epsilon }\, Q = Q\, h_{\epsilon }\, Q\, $ and $\, Q\, d_{\epsilon }\, ( 1 - Q ) = 0 = ( 1 - Q )\, d_{\epsilon }\, Q\, $. Putting 
$\, c_{\epsilon } = c + d + d_{\epsilon }\in\mathcal C\, $ one finds that 
$\, \Vert P\, c_{\epsilon }^2 P - ( P\, c_{\epsilon }\, P )^2\Vert\leq {\epsilon }^2\, $ and 
$\, \Vert ( P\, c_{\epsilon }\, P )^2\, -\, a_P^2\Vert \leq {2\over 3}\,\epsilon + {\epsilon }^2\, $
so that altogether $\, \Vert P\, c_{\epsilon }^2\, P\, -\, a_P^2 \Vert \leq\epsilon\, $ if 
$\, \epsilon > 0\, $ is chosen small enough. The argument in case of $\, \mathcal X\cdot {\mathcal X}^*\, $ is very similar and we leave it to the reader. For an irreducible representation $\, {\lambda }_{\rho }\, $ obtained from a pure state $\, \rho\in \mathcal P ( A )\, $ the above result gives an isomorphism
$\, L_0^{\nu } \bigl( {\lambda }_{\rho } ( A ) \bigr)\> \simeq\> l_{\infty } \bigl( {\mathcal P}_{{\lambda }_{\rho }} ( A ) \bigr)\, $
where $\, {\mathcal P}_{\lambda } ( A )\, =\, \bigl\{ \rho\in \mathcal P ( A )\,\vert\,  \rho  = {\rho }_{\xi }\, ,\, \xi\in {\mathcal H}_{\lambda } \bigr\}\, $. Then $\, {\mathfrak L}^r_0 ( A )\, $ is given by a diagonal subalgebra of the direct product of such quotients where the product ranges over the different folia of pure states of $\, A\, $. Since each factor is a $C^*$-algebra such that the Jordan square coincides with the $C^*$-square the same is true for the direct product map. By an argument as above there exists a multiplicative extension 
$$ {\pi }_0^r : {\mathfrak L}_q ( A )\> \largerightarrow {\mathfrak L}_0^r ( A )\> =\> {\mathfrak L}_1 \bigl( r ( \mathfrak L ( A ) ) \bigr) $$
extending the natural map $\, r : \mathfrak L ( A ) \rightarrow l_{\infty } ( \mathcal P ( A ) )\, $. To see that $\, {\mathfrak L}^r_0 ( A ) \simeq l_{\infty } ( \mathcal P ( A ) )\, $ it is enough to note that the image of $\, \mathfrak L ( A )\, $ in $\, l_{\infty } ( \mathcal P ( A ) )\, $ contains each function $\, {\delta }_{\rho }\, ,\, {\delta }_{\rho } ( \rho ) = 1\, ,\, {\delta }_{\rho } ( \sigma ) = 0\, $ for $\,\sigma\neq\rho\, $ which is the case since any two states in $\,\mathcal P ( A )\, $ are separated. By an argument as used in the proof above the ${\mathfrak L}_1$-completion of the image of 
$\,\mathfrak L ( A )\, $ must agree with the whole algebra $\, l_{\infty } ( \mathcal P ( A ) )\, $ and by rigidity of the embedding $\, c_0 \bigl( \mathcal P ( A ) \bigr) \subseteq l_{\infty } \bigl( \mathcal P ( A ) \bigr)\, $ the map $\, {\pi }_0^r\, $ is necessarily surjective.
That $\, {\pi }_0^r\, $ factors $\, \pi\, $ is again obvious from Proposition 1.
This implies that for an arbitrary $C^*$-algebra $\, A\, $ the image of the Jordan square of a basic positive element in $\, {\mathfrak L}_1 ( A )\, $ coincides 
with the $C^*$-square since the induced map $\, {\mathfrak L}^r_0 ( A ) \twoheadrightarrow 
{\mathfrak L}_1 ( A )\, $ is necessarily a $*$-homomorphism. From monotonicity of the Jordan square one then gets $\, \pi ( {\mathfrak P}^2 )\leq \pi ( ( {\mathfrak P}^c )^2 ) = 
\pi ( \mathfrak P )\cdot \pi ( \mathfrak P )\, $ while on the other hand 
$\, {\mathfrak P}^2 \geq \mathfrak P\cdot \mathfrak P\, $ implying the reverse inequality. 
Similarly one gets for a general selfadjoint element $\, \mathfrak A\in {\mathfrak L}_q ( A )\, $ the estimate 
$$ sqr ( \mathfrak A )\>\geq\> \sup_{\lambda , \underline{\mathcal A}}\, \left\{ 
\lambda\,\left( ( \underline{\mathcal A} + r\, {\bf 1} ) ^2\> -\> {1\over 1 - \lambda }\, r^2\, {\bf 1} \right)\,\Bigm\vert\, \underline{\mathcal A}\leq\mathfrak A\> ,\> \underline{\mathcal A} + r\, {\bf 1} \geq 0\> .\> 0 \leq \lambda < 1 \right\} $$
$$\>\geq \sup_{\lambda , \underline{\mathcal A}}\,\left\{ 
\lambda \left(\, ( \underline{\mathcal A} + r\, {\bf 1} ) \cdot ( \underline{\mathcal A}\, +\, r\, {\bf 1} )\> -\> {1\over 1 - \lambda }\,r^2\, {\bf 1}\,\right) \Bigm\vert\,\underline{\mathcal A}\leq\mathfrak A ,\>  
\underline{\mathcal A} + r\, {\bf 1}\geq 0 ,\> 0 \leq \lambda < 1 \right\} $$
$$\quad =\> \sup_{\underline{\mathcal A}\leq\mathfrak A}\,\left\{\, \underline{\mathcal A}_+\cdot \underline{\mathcal A}_+\,\right\}\> =\> {\mathfrak A}_+\cdot {\mathfrak A}_+ $$ 
where the first identity in the third line can be checked pointwise on the spectrum of 
$\, {\mathfrak L}_q ( A )\, $ since for each homomorphism into $\,\mathbb R\, $ the supremum is taken for $\, r\to\infty\, $ and the second identity follows by normality of the commutative $C^*$-square. On the other hand 
$$ \pi ( sqr ( \mathfrak A ) )\>\leq\> \pi ( {\mathfrak A}_+^2 )\> =\> \pi ( {\mathfrak A}_+ )\cdot \pi ( {\mathfrak A}_+ ) $$
so that one gets an equality. Therefore also 
$$ \pi ( {\mathfrak A}^2 )\> =\> \pi ( sqr ( \mathfrak A ) )\> +\> \pi ( sqr ( - \mathfrak A ) )\> =\>
\pi ( {\mathfrak A}_+ )\cdot \pi ( {\mathfrak A}_+ )\> +\> \pi ( {\mathfrak A}_- )\cdot \pi ( {\mathfrak A}_- )\> =\> \pi ( \mathfrak A )\cdot \pi ( \mathfrak A )\> . $$
The second inequality of $\, (21)\, $ follows since 
$$ \pi ( \sqrt{\mathfrak P} )^2\> =\>\pi ( \sqrt{\mathfrak P}^2 )\>\leq\> \pi ( \sqrt{{\mathfrak P}^c}^2 )\> \leq \pi ( \mathfrak P ) $$ 
by taking squareroots. The statement that the image of $\, \underline{\mathfrak a}_+ \wedge \underline{\mathfrak a}_-\, $ is trivial in $\, {\mathfrak L}_1 ( A )\, $ follows from the Remark after Lemma 2. 
\par\noindent
Now suppose that $\, I\, $ is an injective $C^*$-algebra. Assume first that $\, \mu : I \rightarrow \mathcal B ( \mathcal H )\, $ is a faithful injective representation and let $\,\upsilon : \mathcal B ( \mathcal H ) \twoheadrightarrow I\, $ be a completely positive retraction (left inverse) for $\, \mu \, $. One gets induced maps 
$$ \mathcal M :\quad {\mathfrak L}_q ( I )\> \longrightarrow {\mathfrak L}_q ( \mathcal B ( \mathcal H ) )\> ,  $$
$$ \Upsilon :\quad  \mathfrak L ( \mathcal B ( \mathcal H ) )\> \longrightarrow\> \mathfrak L ( I ) 
\> , $$ 
the former which is a normal $*$-homomorphism. $\, \Upsilon\, $ sends the order ideal 
$\, {\mathfrak J}_0^r\subseteq C^* \bigl( \mathfrak L \bigl( \mathcal B ( \mathcal H ) \bigr) \bigr) \subseteq {\mathfrak L}_q \bigl( \mathcal B ( \mathcal H ) \bigr)\, $ which is the kernel of the (multiplicative) function representation 
$$ r :\> C^* \bigl( \mathfrak L \bigl( \mathcal B ( \mathcal H ) \bigr) \bigr)\> \largerightarrow\> l_{\infty } \bigl( \mathcal P \bigl( \mathcal B ( \mathcal H ) \bigr) \bigr) $$
for the pure states of $\, \mathcal B ( \mathcal H )\, $ onto a corresponding order ideal 
$\, {\mathfrak J}^{\mu ,\upsilon }_0 \subseteq C^* \bigl( \mathfrak L ( I ) \bigr)\, $ with quotient 
$\, {\mathfrak L}^{\mu , \upsilon } ( I ) = C^* \bigl( \mathfrak L ( I ) \bigr)\, /\, {\mathfrak J}^{\mu , \upsilon }_0\, $. Put 
$\, {\mathfrak L}^{\mu , \upsilon }_0 ( I )\, =\, {\mathfrak L}_1 \bigl( {\mathfrak L}^{\mu , \upsilon } ( I ) \bigr)\, $. 
If $\, \underline{\mathcal C}\, ,\, \underline{\mathcal D}\in \mathfrak L ( I )\, $ are positive basic elements
which agree modulo $\, {\mathfrak J}^{\mu , \upsilon }_0\, $ let $\, \underline{\mathcal C}^{\mu }\, ,\, \underline{\mathcal D}^{\mu }\, $ denote their images in 
$\, \mathfrak L ( R )\, $. Then 
$$ \bigl( \underline{\mathcal C}^{\mu } \bigr)^2\> \equiv\> \underline{\mathcal C}^{\mu }\cdot \underline{\mathcal C}^{\mu }\> =\> \mathcal M \bigl( \underline{\mathcal C}\cdot \underline{\mathcal C} \bigr)\> ,\quad \bigl( \underline{\mathcal D}^{\mu } \bigr)^2\> \equiv\> 
\underline{\mathcal D}^{\mu }\cdot \underline{\mathcal D}^{\mu }\> =\> 
\mathcal M \bigl( \underline{\mathcal D}\cdot\underline{\mathcal D} \bigr)\qquad \mod {\mathfrak J}^r_0 $$
and 
$$ \Upsilon \bigl( \bigl( \underline{\mathcal C}^{\mu } \bigr)^2 \bigr)\> =\> 
\underline{\mathcal C}^2\> ,\quad \Upsilon \bigl( \bigl( \underline{\mathcal D}^{\mu } \bigr)^2 \bigr)\> =\> \underline{\mathcal D}^2\> .  $$
The relation $\, \Upsilon \bigl( \bigl( \underline{\mathcal C}^{\mu }\bigr)^2 \bigr) \geq 
\underline{\mathcal C}^2\, $ follows from the Schwarz inequality applied to $\,\upsilon\, $ and the reverse inequality follows from $\, \bigl( \underline{\mathcal C}^{\mu } \bigr)^2 \leq \mathcal M \bigl( \underline{\mathcal C}^2 \bigr)\, $. Then 
$$ \underline{\mathcal C}^2\> \equiv\> \underline{\mathcal C}\cdot \underline{\mathcal C}\> ,\quad \underline{\mathcal D}^2\> \equiv\> \underline{\mathcal D}\cdot \underline{\mathcal D}\qquad\mod {\mathfrak J}^{\mu , \upsilon }_0 $$
so that since the quotient map modulo the order ideal $\, {\mathfrak J}^{\mu , \upsilon }_0\, $ is necessarily a $*$-homomorphism the result 
$$ \underline{\mathcal C}^2\>\equiv\>\underline{\mathcal D}^2\qquad\mod {\mathfrak J}^{\mu , \upsilon }_0 $$
follows, i.e. the basic squaring operation is well defined on the quotient and agrees with the $C^*$-square of the corresponding elements. The argument in case of the operation $\, (13 )\, $ is much the same. Then the quadratic clouds of all types plus the Lie bracket associated with the Jordan squaring operation must be contained in $\, {\mathfrak J}^{\mu ,\upsilon }_0\, $ for any given triple of positive basic elements and $\, {\mathfrak J}_0^r\cap \mathfrak L ( I )\subseteq {\mathfrak J}^{\mu , \upsilon }_0\, $. 
In particular the image of $\, \mathfrak L ( I )\, $ is a $C^*$-subalgebra of $\, {\mathfrak L}^{\mu , \upsilon }_0 ( I )\, $. Define
$$ {\mathfrak J}^{\mu }_0\> =\> \bigcap_{ \upsilon }\, {\mathfrak J}^{\mu , \upsilon }_0  $$
where $\, \{ \upsilon \}\, $ ranges over all completely positive retractions for $\, \mu\, $
and check that $\, {\mathfrak L}^{\mu }_0 ( I )\, $ resp. $\, {\mathfrak L}_0 ( I )\, $ as defined in the Theorem has the required properties by diagonal embedding into the direct product $\, {\mathfrak L}^{\mu  }_0 ( I ) \hookrightarrow \prod_{( \mu , \upsilon )}\, {\mathfrak L}^{\mu , \upsilon }_0 ( I )\, $. 
For a general multiplicative injective representation $\, \mu = \iota\circ \pi\, $ let 
$\, {\mathfrak J}^{\mu }_0\, $  denote the preimage of the corresponding ideal 
$\, {\mathfrak J}_0^{\iota } \subseteq \mathfrak L ( J )\, $ under the canonical surjection 
$\, \mathfrak L ( I ) \twoheadrightarrow \mathfrak L ( J )\, $. Then 
$$ {\mathfrak L}^{\mu }_0 ( I )\> =\> {\mathfrak L}_0^{\iota } ( J ) \> . $$
Also from the argument above the canonical multiplicative surjection 
$$ C^* \bigl( \mathfrak L ( I ) \bigr)\> \twoheadrightarrow\> {\mathfrak L}^{\iota } ( J ) $$ 
extends by the method used above  (nonuniquely) to  a  
$*$-homomorphism 
$$  {\mathfrak L}_q ( I )\> \largerightarrow {\mathfrak L}_0^{\mu } ( I )\> .  $$
From the definition of $\, {\mathfrak L}^{\mu }_0 ( I )\, $ it is immediately clear that if 
$\, \mu\, $ is multiplicative then the basic squaring operation $\, ( 1 )\,  $ drops to the quotient 
$\, {\mathfrak L}^{\mu } ( I )\, $ where it coincides with the corresponding $C^*$-operation since 
the quotient map $\, \Pi : \mathfrak L ( I ) \twoheadrightarrow \mathfrak L ( J )\, $ sends 
a basic square $\, \underline{\mathcal C}^2\, $ to $\, \Pi ( \underline{\mathcal C} )^2\, $. 
Then check that if $\, \mathcal M\, $ denotes the set of equivalence classes of multiplicative injective representations of $\, I\, $ the approximate quotient $\, {\mathfrak L}_0^m ( I )\, $ of 
$\, {\mathfrak L}_q ( I )\, $ as defined in the Theorem has the required properties by diagonal embedding into the direct product
$$ {\mathfrak L}_0^m ( I ) \hookrightarrow \prod_{\mu\in\mathcal M}\, {\mathfrak L}_0^{\mu } ( I ) \> . $$
By considering the case of a faithful injective representation one sees  from uniqueness of $\, \pi\, $
that it factors over $\, {\pi }^m_0\, $\qed
\par\bigskip\bigskip\bigskip\noindent
{\bf 2.\quad Injective envelopes and enveloping von Neumann algebras.}
\par\bigskip\noindent 
The following can be seen as a preamble to the whole text which we have chosen to present here because of its close relation with the results of Theorem 2. 
None of its results is really new, c.f. \cite{A}, Lemma 8, only the proof using transfinite induction seems to be. Also the assumption that the order injection $\, \iota\, $ of the Theorem is multiplicative restricted to $\, A\, $ is unnecessary by the Remark following Lemma 8 in \cite{A} that given a positive linear unital order injection $\, \iota : A \hookrightarrow B\, $ of commutative $C^*$-algebras any positive linear retraction $\, B \twoheadrightarrow A\, $ for $\, \iota\, $ is necessarily a $*$-homomorphism restricted to $\, C^* ( \iota ( A ) )\subseteq B\, $. One should also remark that injectivity of a monotone complete abelian $C^*$-algebra with respect to linear lattice maps is a direct consequence of injectivity with respect to $*$-homomorphisms of abelian $C^*$-algebras since any linear sublattice of an abelian $C^*$-algebra is normdense in its $C^*$-envelope, c.f. \cite{Pe2}, Lemma 4.3.2. Note also that from Corollary 2 of section 3 every $C^*$-quotient of an injective $C^*$-algebra (in particular any quotient of $\, l_{\infty } ( Z )\, $) is sequentially monotone complete.
\par\bigskip\noindent
{\bf Theorem A.}\quad Let $\, A\, $ be a commutative $C^*$-algebra with injective envelope $\, I = I ( A )\, $. Given a unital order isomorphic embedding $\, \iota :\> I\> \hookrightarrow\> B\, $ into a commutative $C^*$-algebra which is multiplicative restricted to $\, A\, $ there exists a positive retraction $\, \rho :\> B\> \twoheadrightarrow\> I\, $ which is a $*$-homomorphism, in particular any commutative injective $C^*$-algebra is a $C^*$-quotient of $\, l_{\infty } ( Z )\, $ for some set $\, Z\, $. Any commutative injective $C^*$-algebra $\, I\, $ is injective in the three categories of positive linear maps of function systems, $*$-homomorphisms of abelian $C^*$-algebras and lattice maps of linear function lattices (which are order isomorphic to a linear sublattice of a commutative $C^*$-algebra). 
\par\bigskip\bigskip\bigskip\noindent
{\it Proof.}\quad Choose a well ordered basis 
$\, \{ c_{\lambda } {\}}_{\lambda\in \Lambda }\, $ consisting of selfadjoint elements of norm one, assuming that the subset $\, \{ c_{\lambda } \} \cap \iota ( A )\, $ is a basis for $\, \iota ( A )\, $ 
and exhausts the leading halfopen interval of all indices $\, 1 \leq \lambda < {\lambda }_0\, $. To save notation one puts $\, {\lambda }_0 = 0\, $ disregarding all indices $\, \lambda < {\lambda }_0\, $ since for these the map $\, \rho \, $ is canonically defined. Let $\, {\mathfrak Y}_{\lambda }\, $ denote the $*$-linear subspace generated by the set $\, \{ c_{\kappa }\,\vert\, \kappa \leq \lambda \}\, $ and 
$\, {\mathfrak Y}_{< \lambda }\, $ the linear span of $\, \{ c_{\kappa }\,\vert\, \kappa < \lambda \}\, $. Inductively define $\, {\rho }_{\lambda } : {\mathfrak Y}_{\lambda } \rightarrow I\, $ by linear extension of  
$$ {\rho }_{\lambda } ( c_{\lambda } )\> =\> \inf\, \bigl\{ {\rho }_{< \lambda } ( y )\bigm\vert y\in {\mathfrak Y}_{< \lambda }\, ,\, y \geq c_{\lambda } \bigr\}\> ,  $$
assuming by induction that $\, {\rho }_{< \lambda }\, $ is well defined and positive.
One checks positivity of $\, {\rho }_{\lambda }\, $. Let 
$\, x = b + \gamma\, c_{\lambda } \geq 0\, $ with $\, b\in {\mathfrak Y}_{< \lambda }\, $ be given. First assume that $\, \gamma > 0\, $. Then 
$$ {\rho }_{\lambda } ( x )\> =\> \inf\, \bigl\{ {\rho }_{< \lambda } ( y )\bigm\vert y \in {\mathfrak Y}_{< \lambda }\, ,\, y \geq x \bigr\}\> \geq\> 0 \> . $$
On the other hand if $\, \gamma < 0\, $ then 
$$ {\rho }_{\lambda } ( x )\> =\> \sup\, \bigl\{ {\rho }_{< \lambda } ( y )\bigm\vert y\in {\mathfrak Y}_{< \lambda }\, ,\, y \leq x \bigr\}\> \geq\> 0 $$
since $\, 0\, $ is contained in $\, {\mathfrak Y}_{< \lambda }\, $ and $\, x \geq 0\, $. Proceeding by induction this results in a positive $*$-linear map $\, \rho : B \rightarrow I\, $. 
Then one needs to check that $\, \rho\, $ is also multiplicative, i.e. a $*$-homomorphism, which is equivalent to validity of the Schwarz identity 
$$ \rho ( x^2 )\> =\> \rho ( x )^2 $$
for selfadjoint elements $\, x\, $. Since the Schwarz inequality $\, \rho ( x^2 ) \geq \rho ( x )^2\, $ holds by virtue of positivity one only needs to check the inverse Schwarz inequality 
$\, \rho ( x^2 ) \leq \rho ( x )^2\, $. This can be done inductively on each fixed subspace $\, {\mathfrak Y}_{\lambda }\, $ assuming the Schwarz identity for all elements in $\, {\mathfrak Y}_{< \lambda }\, $. Let $\, b\in {\mathfrak Y}_{< \lambda } + {\mathbb R}_+\, c_{\lambda } \geq 0\, $ be given. Then 
$$ \rho ( b^2 )\> \leq\> \inf\, \bigl\{ \rho ( a^2 )\bigm\vert a\in {\mathfrak Y}_{< \lambda }\, ,\, a \geq b \bigr\}\> =\> \inf\, \bigl\{ \rho ( a )^2\bigm\vert a\in {\mathfrak Y}_{< \lambda }\, ,\, a \geq b \bigr\} $$
$$\quad =\> \left( \inf\, \bigl\{ \rho ( a )\bigm\vert a \in {\mathfrak Y}_{< \lambda }\, ,\, a \geq b \bigr\} \right)^2\> =\> \rho ( b )^2\> . $$
The first equality in the second line follows by normality of the squaring operation in (monotone complete) commutative $C^*$-algebras. Any selfadjoint element of $\, {\mathfrak Y}_{\lambda }\, $ can be written as a difference of elements as above. If $\, b\, ,\, c \geq 0\, $ are two such elements then
$$ 2\,\rho \bigl( b c \bigr)\> =\> \rho \bigl( ( b + c )^2 \bigr)\> -\> \rho \bigl( b^2 \bigr)\> -\> \rho \bigl( c^2 \bigr)\> =\> 2\, \rho \bigl( b \bigr)\, \rho \bigl( c \bigr) $$ 
whence 
$$ \rho \bigl( ( b - c )^2 \bigr)\> =\> \rho \bigl( b^2 \bigr)\> -\> 2\, \rho \bigl( b c \bigr)\> +\> \rho \bigl( c^2 \bigr)\> =\> \rho \bigl( b - c \bigr)^2\>  $$
proving the Schwarz identity for arbitrary selfadjoint elements in $\, {\mathfrak Y}_{\lambda }\, $.
The induction starts by noting that the Schwarz identity trivially holds for elements of $\, A\, $. Therefore the construction yields a well defined $*$-homomorphism  which necessarily is a retraction since it extends the identity map of $\, A\, $. Since any injective $C^*$-algebra is its own injective envelope and admits a faithful $*$-representation in some $\, l_{\infty } ( Z )\, $ this also proves that any injective $C^*$-algebra can be represented as a $C^*$-quotient of some $\, l_{\infty } ( Z )\, $. 
Then basically the same argument shows that given an inclusion $\, A \subseteq B\, $ of commutative $C^*$-algebras and a $*$-homomorphism $\, \rho : A \longrightarrow I\, $ into a commutative injective $C^*$-algebra there exists a positive linear extension $\, \overline\rho : B \longrightarrow I\, $ of $\, \rho\, $ which satisfies the Schwarz identity $\, \overline\rho ( x^2 ) = \overline\rho ( x )^2\, $ for each selfadjoint element $\, x\in B^{sa}\, $, i.e. $\,\overline\rho\, $ is a $*$-homomorphism. In particular every commutative injective $C^*$-algebra $\, A\, $ is a complete function lattice, i.e. the lattice operations are normal. Namely using a unital $*$-homomorphic embedding $\, j : A \hookrightarrow l_{\infty } ( Z )\, $ there exists an extremal retraction  $\, r:  l_{\infty } ( Z ) \twoheadrightarrow A\, $ which is a $*$-homomorphism, in particular a lattice map so that 
$$ \sup_{\mu }\, \bigl ( x_{\mu } \wedge y \bigr)\> =\> r \left( \sup_{\mu } \bigl( j ( x_{\mu } )\wedge j ( y ) \bigr) \right)\> =\> r \left( \bigl( \sup_{\mu } j ( x_{\mu } ) \bigr) \wedge j ( y ) \right)\> =\> \bigl( \sup_{\mu }\, x_{\mu } \bigr) \wedge y \> . $$
Now the above argument can be used to show that $\, I\, $ is injective also for lattice maps. Suppose given a sublattice inclusion of linear function lattices 
$\, A \subseteq B\, $ (both sublattices of the selfadjoint part of some commutative $C^*$-algebra). One proceeds as above choosing a well ordered real basis $\, \{ c_{\lambda } {\}}_{\lambda\in \Lambda }\, $ of $\, B\, $ such that the leading halfopen interval $\, \{ c_{\lambda }\,\vert\, \lambda < {\lambda }_0 \}\, $ is a basis for $\, A\, $ and extends the given lattice map $\, \rho : A \longrightarrow I\, $ to a positive linear map $\, \overline\rho : B \longrightarrow I\, $ in the manner above. Then one needs to show that this extremal extension is a lattice map. One proceeds by induction. Assume that for given fixed index $\, \lambda\in \Lambda\, $ the identity 
$$ \overline\rho ( a\vee y )\> =\> \overline\rho ( a )\vee \overline\rho ( y ) $$ 
holds for all elements 
$\, a\in {\mathfrak Y}_{< \lambda }\, ,\, y\in B\, $. Let 
$\, x = b + \gamma\, c_{\lambda }\, $ be given with $\, \gamma\in \mathbb R\, $.  One has  
$$ \overline\rho ( x\vee y )\> =\> \overline\rho \bigl( \bigl( b\vee ( y - \gamma c_{\lambda } ) \bigr) + \gamma\, c_{\lambda } \bigr)\> =\> \overline\rho \bigl( b\vee ( y - \gamma\, c_{\lambda } ) \bigr)\> +\> \gamma\, \overline\rho \bigl( c_{\lambda } \bigr) $$
$$\quad =\> 
\overline\rho \bigl( b )\vee \overline\rho \bigl( y - \gamma\, c_{\lambda } \bigr)\> +\> \gamma\, \overline\rho \bigl( c_{\lambda } \bigr)\> =\> \overline\rho \bigl( x \bigr)\vee \overline\rho \bigl( y \bigr) $$   
by linearity of $\, \overline\rho\, $. Thus by induction the relation will hold if it holds for all $\, x\in A\, ,\, y\in B\, $. Assume by induction that for some fixed index $\, \mu\in\Lambda\, $ the relation holds whenever $\, y\in {\mathfrak Y}_{< \mu }\, $. Put $\, y = c + \gamma\, c_{\mu }\, $ with 
$\, \gamma > 0\, $.
By monotonicity one has 
$\, \overline\rho ( x\vee y )\>\geq\> \overline\rho ( x )\vee \overline\rho ( y )\, $ so we only need to prove the reverse inequality. Then
$$ \overline\rho ( x\vee y )\>\leq\> \inf\, \bigl\{ \overline\rho ( x\vee d )\bigm\vert d\geq y\, ,\, d\in {\mathfrak Y}_{< \mu } \bigr\}\> =\> \inf\, \bigl\{ \overline\rho ( x )\vee \overline\rho ( d )\bigm\vert 
d\geq y\, ,\, d\in {\mathfrak Y}_{< \mu } \bigr\} $$
$$\quad =\> \overline\rho ( x )\vee
\inf\, \bigl\{ \overline\rho ( d )\bigm\vert d\geq y\, ,\, d\in {\mathfrak Y}_{< \mu } \bigr\}\> =\> 
\overline\rho ( x ) \vee \overline\rho ( y )\> .  $$
Now suppose that $\, \gamma < 0\, $. One has 
$$ \overline\rho \bigl( ( x - y )\vee 0 \bigr)\> =\> \overline\rho \bigl( x - y \bigr) \vee 0 $$
from the argument above so that by linearity of $\, \overline\rho\, $ one gets 
$$ \overline\rho \bigl( x \vee y \bigr)\> =\> \overline\rho \bigl( ( x - y )\vee 0 \bigr)\> +\> \overline\rho ( y )\> =\> \left( \bigl( \overline\rho ( x )\> -\> \overline\rho ( y ) \bigr) \vee 0 \right)\> +\> \overline\rho ( y )\> =\> \overline\rho ( x ) \vee \overline\rho ( y )\> . $$
Therefore this relation holds by induction for all $\, x\, ,\, y\in B\, $. The relation 
$$ \overline\rho ( x \wedge y )\> =\> \overline\rho ( x ) \wedge \overline\rho ( y ) $$
follows by the symmetry $\, x \mapsto - x\, ,\, y \mapsto - y\, $, exchanging 
$\, \wedge\, $ for $\, \vee\, $\qed
\par\bigskip\noindent 
If $\, \mathfrak X\subseteq \mathcal B ( \mathcal H )\, $ is an operator system its injective envelope is denoted $\, I ( \mathfrak X )\, $. A positive linear map 
$\, \sigma : \mathcal B ( \mathcal H ) \rightarrow  {\mathfrak L}_1 \bigl( I ( \mathfrak X ) \bigr)\, $ will be called a {\it retraction} (for $\,\iota\, $) if given a completely isometric unital embedding 
$\, \iota : I ( \mathfrak X ) \hookrightarrow \mathcal B ( \mathcal H )\, $ one has 
$\, ( \sigma \circ \iota ) ( x ) = x\, $. The retraction $\, \sigma\, $ is called {\it proper} if 
$\,\sigma \bigl( \mathcal B ( \mathcal H ) \bigr) = I ( \mathfrak X )\subseteq {\mathfrak L}_1 \bigl( I ( \mathfrak X ) \bigr)\, $. For a given $C^*$-algebra $\, A\, $ recall from \cite{Wr} the definition of its regular $\sigma $-completion denoted $\, {\sigma }^r ( A )\, $ which is a minimal sequentially monotone complete $C^*$-algebra containing $\, A\, $ and completely isometric to the sequentially monotone completion of $\, A\, $ in $\, {\mathfrak L}_1 ( A )\, $. By definition $\, {\sigma }^r ( A )\, $ is the quotient of the Baire envelope $\, \mathfrak B ( A ) \subseteq A^{**}\, $ of $\, A\, $, i.e. the sequential monotone closure of $\, A\, $ in $\, A^{**}\, $, modulo a sequentially monotone complete $*$-ideal. 
If $\, A\, $ is separable then $\, {\sigma }^r ( A )\, $ is monotone complete. Recall some definitions of \cite{Wr2}: $\, \mathcal M\subseteq A^{**}\, $ is the subspace of elements such that the set $\, \bigl\{ \rho\in \mathcal P ( A )\bigm\vert \rho ( x ) \neq 0 \bigr\}\, $ is a meagre subset in the relative topology, $\,\mathcal N \subseteq \mathcal M\, $ is the linear span of the positive elements of $\,\mathcal M\, $. Putting $\, \mathcal I = \mathcal N \cap \mathfrak B ( A ) = \mathcal M \cap \mathfrak B ( A )\, $ one defines $\, {\sigma }^r ( A )\, $ to be the $C^*$-quotient $\, \mathfrak B ( A ) \bigm/ \mathcal I\, $. Related notions are the weak (resp. strong) Baire envelope $\, {\mathfrak B}^w ( A )\, $ 
and $\, {\mathfrak B}^s ( A )\, $ respectively which are defined to be the smallest subspaces of $\, A^{**}\, $ containing $\, A\, $ and closed under taking weak (resp. strong) limits of bounded sequences 
both of which are (sequentially monotone complete) sub-$C^*$-algebras of $\, A^{**}\, $ (see \cite{Pe1}, chap. 4.5). Another related notion we use below is the 
{\it ultraweak Baire envelope} $\, {\mathfrak B}^{\omega } ( A )\, $. If given any bounded sequence of elements $\, \{ x_n \}\subseteq A^{**}\, $ and any free ultrafilter $\, \omega\, $ on $\, \mathbb N\, $ one obtains the $\, \omega$-limit $\, x^{\omega } = \lim_{\omega }\, \{ x_n \}\in A^{**}\, $ viewing selfadjoint elements of $\  A^{**}\, $ as bounded affine functions on the state space $\, \mathcal S ( A ) \subseteq A^*\, $ and taking the pointwise limit of the functions $\, \{ f_n \}\, $ corresponding to the $\, \{ x_n \}\, $ along $\,\omega\, $. Then we say that a subspace (or subset) $\, V \subseteq A^{**}\, $ is {\it ultraweakly sequentially closed} if it is closed under taking  arbitrary ultrafilter limits of bounded sequences. Clearly this implies that $\, V\, $ is weakly sequentially closed, in particular 
$\, \mathfrak B ( A ) \subseteq {\mathfrak B}^s ( A ) \subseteq {\mathfrak B}^w ( A ) \subseteq {\mathfrak B}^{\omega } ( A )\, $. A straightforward transposition of the corresponding results proved in \cite{Wr}, \cite{Wr2} 
gives that $\, \mathcal M\, $ and hence also $\, \mathcal N\, $ are ultraweakly sequentially closed subspaces (compare with Lemma 1 of \cite{Wr2}). We will however slightly modify the definition of $\, \mathcal M\, $ and $\, \mathcal N\, $ respectively suited for our purposes. Let $\, X = \overline{\mathcal P ( A )} \subseteq \mathcal S ( A )\, $ denote the closure of the subset of pure states in the state space of $\, A\, $ yielding a function representation of $\, A \subseteq C ( X )\, $ by continuous functions on $\, X\, $ and a corresponding normal function representation $\, p_X : A^{**} \rightarrow C ( X )^{**}\subseteq B ( Y )\, $ as bounded functions on $\, Y = \mathcal S ( C ( X ) )\, $ by evaluation.  
Moreover $\, C ( X )\, $ is easily seen to be generated as an abelian $C^*$-algebra by the image of $\, A\, $. Also consider the abelian sub-$C^*$-algebra $\, C^*_{\pi } ( A ) \subseteq {\mathfrak L}_1 ( A )\, $ generated by the image of $\, A\, $ in its positive injective envelope. The identity map of $\, A\, $ extends to a positive linear map $\, {\mathfrak L}_1 ( A ) \hookrightarrow {\mathfrak L}_1 ( C ( X ) )\, $ which admits a positive linear retraction $\, {\mathfrak L}_1 ( C ( X ) ) \twoheadrightarrow {\mathfrak L}_1 ( A )\, $. This map is necessarily a $*$-homomorphism restricted to the sub-$C^*$-algebra of $\, {\mathfrak L}_1 ( C ( X ) )\, $ generated by the image of $\, {\mathfrak L}_1 ( A )\, $ (cf.  \cite{Ch-E}, Theorem 4.1) containing $\, C ( X )\, $ leading to a $*$-homomorphic surjection $\, C ( X ) \twoheadrightarrow C^*_{\pi } ( A )\, $ corresponding to a continuous embedding $\, 
spec ( C^*_{\pi } ( A ) ) \subseteq X\, $. Since both spaces are compact closures of $\, \mathcal P ( A )\, $ the latter being a uniquely determined subset in both $\, spec ( C^*_{\pi } ( A ) )\, $ and $\, X\, $ one finds from continuity of the embedding that $\, spec ( C^*_{\pi } ( A ) ) \hookrightarrow X\, $ is onto, whence $\, C^*_{\pi } ( A ) \simeq C ( X )\, $.  Restriction to $\, {\mathfrak B}^{\omega } ( A )\, $ yields a sequentially normal representation $\, p_X : {\mathfrak B}^{\omega } ( A ) \rightarrow {\mathfrak B}^{\omega } ( X )\, $ with $\, {\mathfrak B}^{\omega } ( X )\, $ the ultraweak sequential closure of $\, C ( X )\, $ in $\, B ( Y )\, $. Let 
$\, {\mathcal M}_X\, $ be the subspace of functions in $\, B ( Y )\, $ vanishing outside a meagre subset of $\, X\, $. Define $\, {\mathcal I}^{\omega } = p_X^{-1} \bigl( {\mathcal M}_X \cap p_X \bigl( {\mathfrak B}^{\omega } ( A ) \bigr) \bigr)\, $ and put $\, {\sigma }^r_{\omega } ( A ) :=  {\mathfrak B}^{\omega } ( A ) \bigm/ {\mathcal I}^{\omega }\, $. It will be shown below that $\, {\mathcal I}^{\omega }\subseteq {\mathfrak B}^{\omega } ( A )\, $ is a (sequentially monotone complete) $*$-ideal so that the quotient map is a (sequentially normal) $*$-homomorphism. We call $\, {\sigma }^r_{\omega } ( A )\, $ the {\it ultraweak regular envelope} of $\, A\, $ where the term regular in this context is used to indicate that each selfadjoint element $\, x\in {\sigma }^r_{\omega } ( A )\, $ is the supremum of all elements $\, \{ b\in A^{sa}\,\vert\, b\leq x \}\, $ and the infimum of all elements $\, \{ a\in A^{sa}\,\vert\, a \geq x \}\, $ (see below). It may be that $\, {\mathfrak B}^{\omega } ( A ) = A^{**}\, $ in certain cases. Obviously any ultraweak sequential limit point $\, x_{\omega }\, $ of a bounded sequence 
$\, \{ x_n \}\subseteq V \subseteq A^{**}\, $ is a weak limit point of the sequence so the question arises whether iterated adjunction of such limit points gives the weak closure of $\, V\, $ which necessitates firstly by taking $\, V\, $ to be the linear span of a given sequence $\, \{ x_n \}\, $ that the weak closure of $\, \{ x_n \}\, $ is equal to its ultraweak sequential closure, and secondly that the weak closure of a given subspace is equal to its extended weak sequential closure the latter being defined by adjoining  all weak limit points of any bounded subsequence in $\, V\, $. If $\, V\, $ is separable the second problem has an affirmative solution, since the unit ball of $\, V\, $ contains a weakly dense subsequence of the unit ball of the weak closure of $\, V\, $. 
\par\noindent
An operator system $\, J\, $ is called 
{\it weakly injective} iff given any unital completely isometric inclusion of operator systems 
$\, \mathfrak X \subseteq \mathfrak Y\, $ together with a completely positive map 
$\, \varphi : \mathfrak X \rightarrow J\, $ there exists a positive linear map
$\, \widehat\varphi : \mathfrak Y \rightarrow J\, $ extending $\, \varphi\, $. Any weakly injective operator system admits a structure as a Jordan $*$-algebra (= $JC$-algebra, c.f. \cite{E-St}) .
If $\, A\, $ is a $C^*$-algebra then a weakly injective $C^*$-algebra $\, J ( A )\, $ containing $\, A\, $ as a subalgebra is called a weak injective envelope of $\, A\, $ iff any monotonous extension 
$\, \iota : J ( A ) \rightarrow J ( A )\, $ of the identity map of $\, A\, $ is the identity map. From definition it is obvious that the weak injective envelope if it exists is uniquely determined up to Jordan isomorphism.
\par\bigskip\noindent
{\bf Theorem 2.}\quad (i) If $\,\mathfrak X\, $ is an operator system contained in an abelian $C^*$-algebra then $\, I ( \mathfrak X ) = {\mathfrak L}_1 ( \mathfrak X )\, $ and each selfadjoint element 
$\, x\in I ( \mathfrak X )^{sa}\, $ is the least upper bound of the subset 
$\, \{ a_{\lambda }\in {\mathfrak X}^{sa}\,\vert\, a_{\lambda } \leq x \}\, $, and the greatest lower bound of the subset $\, \{ a_{\mu }\in {\mathfrak X}^{sa}\,\vert\, x\leq a_{\mu } \} \, $.  
In particular any monotone complete abelian $C^*$-algebra is injective (this of course is well known, cf.  Theorem 4.3.6 of \cite{A-K}). 
\par\smallskip\noindent
(ii) Let $\, \lambda : A\rightarrow \mathcal B ( \mathcal H )\, $ be a faithful unital $*$-representation of the unital $C^*$-algebra $\, A\, $ with strong closure given by the injective von Neumann algebra $\, R \, $ acting on the separable Hilbert space $\,\mathcal H\, $. 
Let $\, \iota : I ( A ) \hookrightarrow R\, $ be any completely isometric embedding extending the identity map of $\, A\, $ (which exists by injectivity of $\, R\, $). Then there exists canonical $*$-ideal $\, J \vartriangleleft R\, $ which is contained in the kernel of a certain positive linear retraction $\, \sigma :\, R\,\twoheadrightarrow\, {\mathfrak L}_1 \bigl( I ( A ) \bigr)\, $ for $\,\iota\, $ and is trivial only if $\, R\, $ is (completely isometric to) the injective envelope of $\, A\, $. 
\par\smallskip\noindent
(iii) For any $C^*$-algebra $\, A\, $ its ultraweak regular envelope $\, {\sigma }^r_{\omega } ( A )\, $
admits the structure of a sequentially monotone complete $C^*$-algebra and a sequentially monotone complete subspace of $\, {\mathfrak L}_1 ( A )\, $ containing $\, A\, $. Moreover $\, {\sigma }^r_{\omega } ( A )\, $ inherits a {\it ultraweak sequential topology} such that the unit ball of $\, {\sigma }^r_{\omega } ( A )\, $ is sequentially compact for this topology in the sense that every bounded sequence in $\, {\sigma }^r_{\omega } ( A )\, $ has an ultraweak sequential limit point (so there exists a subnet converging to this point although the subnet may not be a subsequence). Equivalently, every countable open covering of a closed bounded (= sequentially compact) subset $\, C\subseteq {\sigma }^r_{\omega } ( A )\, $ has a finite subcovering. For any bounded sequence $\, \{ x_n \} \subseteq {\sigma }^r_{\omega } ( A )\, $ and any ultrafilter $\, \omega\, $ on $\, \mathbb N\, $ the $\, \omega $-limit $\, x_{\omega } = \lim_{\omega }\, \{ x_n \}\, $ is well defined in $\, {\sigma }^r_{\omega } ( A )\, $ and linear in the sense that given two bounded sequences $\, \{ x_n \}\, ,\, \{ y_n \}\, $ one has $\, \lim_{\omega } \{ x_n \} + \lim_{\omega }\, \{ y_n \} = \lim_{\omega }\, \{ x_n + y_n \}\, $. Each $\omega $-limit $\, x_{\omega }\, $ is an ultraweak sequential limit point of the sequence and a subset $\, V\, $ is closed in the ultraweak sequential topology iff it contains each $\omega $-limit of any bounded sequence in $\, V\, $. In particular this applies to an arbitrary injective abelian $C^*$-algebra $\, I\, $ which is equal to its own ultraweak regular envelope. If $\, A\, $ is separable then $\, {\sigma }^r_{\omega } ( A )\, $ is isometric to a weakly injective subspace $\, J ( A ) \subseteq I ( A )\, $ containing $\, A\, $ and is a weak injective envelope of $\, A\, $. In particular there is a positive projection 
$\, \Phi : I ( A ) \rightarrow I ( A )\, $ with range $\, J ( A )\, $. 
\par\smallskip\noindent
(iv) Any von Neumann algebra $\, R\subseteq \mathcal B ( \mathcal H )\, $ is weakly injective, i.e. there exists a positive projection $\, {\Phi }_R : \mathcal B ( \mathcal H ) \rightarrow \mathcal B ( \mathcal H )\, $ with range $\, R\, $.
\par\bigskip\noindent
{\it Proof.}\quad We begin with part (i) concerning an operator subsystem of an abelian $C^*$-algebra. The assumption implies that also 
$\, I ( \mathfrak X )\, $ is abelian since the $C^*$-product of $\, I ( \mathfrak X )\, $ can be defined recurring to the $C^*$-product in some injective abelian $C^*$-algebra $\, A\supseteq \mathfrak X\, $ on extending the inclusion of $\, \mathfrak X\, $ to a unital isometric inclusion $\, I ( \mathfrak X ) \subseteq A\, $ from injectivity, by the formula 
$$ x\cdot y\> =\> \Phi \bigl( x\, y \bigr)\> ,\quad x\, ,\, y\in I \bigl( \mathfrak X \bigr) $$
where $\, \Phi :\, A \rightarrow A\, $ is a positive projection with range equal to $\, I ( \mathfrak X )\, $
(compare with the proof of Theorem 6.1.3 of \cite{E-R}). Let 
$\, x\in I ( \mathfrak X )^{sa}\, $ be given and 
$\, \{ a_{\lambda }\,\vert\, a_{\lambda }\in \mathfrak X\, ,\, a_{\lambda } \leq x \}\, $ be the subset of elements in $\, {\mathfrak X}^{sa}\, $ which are smaller or equal than $\, x\, $. Let $\, \overline x\, $ be the least upper bound of this set in $\, I ( \mathfrak X )\, $ which exists by monotone completeness of the injective envelope (Theorem 6.1.3 of \cite{E-R}). Then 
$\,\overline x \leq x\, $. Consider the subspaces 
$\, A_x = \mathfrak X + \mathbb C\, x \subseteq I ( \mathfrak X )\, $ and 
$\, A_{\overline x} = \mathfrak X + \mathbb C\, \overline x\subseteq I ( \mathfrak X )\, $ and define a map 
$\, \nu : A_x \rightarrow A_{\overline x}\, $ extending the identity map of $\, \mathfrak X\, $ in the obvious way by sending $\, x\, $ to $\,\overline x\, $. We claim that $\,\nu\, $ is positive (and hence completely contractive since unital with $\, I ( \mathfrak X )\, $ abelian). To see this let a positive element in 
$\, A_x\, $ be given which can be written as 
$$\, y\, =\,  a\, +\, \gamma \, x \geq 0 $$
with $\, \gamma \in \mathbb R\, $ and $\, a\in {\mathfrak X}^{sa}\, $.
Suppose that $\, \gamma  < 0\, $. Then since $\,\overline x \leq x\, $ one has 
$\, \nu ( y ) \geq y \geq 0\, $. We may therefore assume $\, \gamma > 0\, $. Then 
$\, \nu ( y )\, $ is equal to the least upper bound of the set $\, \{ \gamma a_{\lambda } + a\,\vert\, 
a_{\lambda }\in \mathfrak X\, ,\, a_{\lambda } \leq x \}\, $ which equals the least upper bound of the set 
$\, \{ b_{\lambda }\in\mathfrak X\,\vert\, b_{\lambda } \leq a + \gamma x \}\, $, hence 
$\, \nu ( y ) \geq 0\, $ as desired. Extending $\,\nu\, $ to a completely positive map of 
$\, I ( \mathfrak X )\, $ into 
$\, I ( \mathfrak X )\, $ and using rigidity gives that $\, x = \overline x\, $. The case of $\, x\, $ being equal to the greatest lower bound of elements in $\, {\mathfrak X}^{sa}\, $ which are larger follows by symmetry. This proves the special Up/Down-property of $\, I ( \mathfrak X )\, $. Now it is easy to see that the monotone complete abelian $C^*$-algebras $\, {\mathfrak L}_1 ( \mathfrak X )\, $ and $\, I ( \mathfrak X )\, $ are naturally (completely) order isomorphic.
If $\, A\, $ is a monotone complete 
abelian $C^*$-algebra then by the foregoing argument each element $\, x\in I ( A )^{sa}\, $ is the least upper bound of all elements $\, \{ a\in A\,\vert\, a\leq x \}\, $. But this set also has a least upper bound 
$\, \overline x\, $ in $\, A\, $ with $\, \overline x \geq x\, $ whereas the set $\, \{ b\in A^{sa}\,\vert\, b\geq x \}\, $ has a greatest lower bound $\, \underline x\, $ in $\, A\, $, so that $\, \overline x \leq \underline x \leq x\leq 
\overline x\, $ and equality follows in each instance, i.e. $\, A = I ( A )\, $ so $\, A\, $ must be injective proving (i).
\par\smallskip\noindent
Let $\, \lambda : A \rightarrow \mathcal B ( \mathcal H )\, $ be a faithful representation of $\, A\, $ as in part (ii) of the Theorem with strong closure given by the injective von Neumann algebra $\, R = \lambda ( A )''\, $. To save notation put $\, I = I ( A )\, $. 
Any completely positive extension $\, \iota : I \hookrightarrow R \subseteq \mathcal B ( \mathcal H )\, $ of $\, \lambda\, $ then is a separable relatively transitive injective representation of $\, I\, $ so that from Proposition 3 the minimal monotonous extension $\, \underline\pi : {\mathfrak L}_q^{\nu } \bigl( \mathcal B ( \mathcal H ) \bigr) \twoheadrightarrow {\mathfrak L}_1 ( I )\, $ of the identity map of $\, I\, $ is boundedly basically increasing normal on the image of $\, {\mathfrak L}^{\nu } ( R )\, $.
Let $\, \mathcal J\subseteq R\, $ be the operator subsystem generated by the affine cone $\, {\mathcal J}^-\, $ of elements which are infima of monotone decreasing nets $\, ( x_{\lambda } {)}_{\lambda }\searrow z\, $ with $\, \{ x_{\lambda } \}\subseteq \iota ( I )\, $. If $\, x\in {\mathcal J}^+\, $ is the supremum of the monotone increasing net $\, ( y_{\lambda } )_{\lambda } \nearrow x\, $ with $\, y_{\lambda }\in \iota ( I )\, $ and $\, \overline x\in \iota ( I )\, $ is the supremum of the same net in $\, \iota ( I )\, $ then the positive element 
$\, z_x = \overline x - x\in {\mathcal J}^-\, $ is in the kernel of $\, \upsilon\, $ for every monotonous retraction $\, \upsilon : \mathfrak L \bigl( \mathcal B ( \mathcal H ) \bigr) \twoheadrightarrow {\mathfrak L}_1 ( I )\, $ for $\,\iota\, $. 
The completely positive linear map $\, \iota\, $ is necessarily an $A$-module map, i.e. 
$\, \iota ( a x b ) = a\, \iota ( x )\, b\, $ for $\, a\, ,\, b\in A\, ,\, x\in I\, $, since $\, \iota\, $ is multiplicative restricted to $\, A\, $. This fact can be checked from considering some $*$-homomorphic Stinespring dilation of $\, \iota\, $ and assuming $\, a\, ,\, b\, $ to be unitary elements since an arbitrary element is a linear combination of unitaries. Therefore
$\, z_{a x a} = a z_x a\, $ so that 
$\,  \upsilon ( a\, z_x\, a ) = 0\, $ for every $\, a\in A^{sa}\, $. 
Let $\, {\mathcal J}_A^+ \subseteq {\mathcal J}^+\, $ denote the subset of elements which are suprema of increasing nets of elements in $\, A^{sa}\, $.
Any element $\, c\in R^{sa}\, $ is the infimum of a monotone decreasing sequence $\, ( b_n {)}_n \searrow c\, $ with $\, b_n\in {\mathcal J}_A^+\, $ the limit of a monotone increasing sequence 
$\, ( a_{n m} {)}_m \nearrow b_n\, $ of elements $\, \{ a_{n m} \} \subseteq A^{sa}\, $. 
Given $\, z_x\, $ as above there exists since $\, \mathcal H\, $ is separable a subsequence 
$\, \{ b_{n_k} {\}}_k\, $ and an assignment $\, n_k \mapsto m_k = m ( n_k )\, $ such that the net 
$\, \bigl\{ a_{n_k m}\, z_x\, a_{n_k m}\bigm\vert m \geq m_k \bigr\}\, $ converges strongly to 
$\, c\, z_x\, c\, $ where $\, ( n_k , m ) \geq ( n_l , r )\, $ if either $\, k \geq l\, $ or else $\, k = l\, $ and $\, m \geq r\, $. By selecting corresponding subsequences of the $\, ( b_n {)}_n\, $ and $\, ( a_{n m } {)}_m\, $ we may assume that $\, \{ a_{n m}\, z_x\, a_{n m} \} \, $ converges strongly to $\, c\, z_x\, c\, $. To save notation we will abbreviate the index set $\, \{ (n , m) \}\, $ by $\, \{ \nu \}\, $ and 
$\, z_{a_{\nu } x a_{\nu }}\, $ by $\, z_{\nu }\, $. Choose any monotonous extension
$\, {\pi }_q^{\nu } : {\mathfrak L}_q ( R ) \rightarrow {\mathfrak L}_q^{\nu } \bigl( \mathcal B ( \mathcal H ) \bigr)\, $ of the canonical map $\, {\pi }^{\nu } :\, \mathfrak L ( R ) \twoheadrightarrow {\mathfrak L}^{\nu } ( R ) \hookrightarrow {\mathfrak L}^{\nu } \bigl( \mathcal B ( \mathcal H ) \bigr)\, $.
Then the image of $\, c\, z_x\, c\, $ in $\, {\mathfrak L}^{\nu } \bigl( \mathcal B ( \mathcal H ) \bigr)\, $ is given by the image of the elements 
$$ \liminf_{\nu }\, z_{a_{\nu }\, x\, a_{\nu }}\> =\> 
\sup_{\nu }\, \inf_{\mu\geq \nu}\, \bigl\{ a_{\mu } ( \overline x - x ) a_{\mu } \bigr\}\> \leq\> 
\inf_{\nu }\, \sup_{\mu \geq \nu} \bigl\{ a_{\mu } ( \overline x - x ) a_{\mu } \bigr\}\> =\> \limsup_{\nu }\, z_{a_{\nu }\, x\, a_{\nu }}  $$
of $\, {\mathfrak L}_q \bigl( R \bigr)\, $.
That these images are the same is seen from the relations
$$ {\pi }_q^{\nu } \bigl( \liminf_{\nu } \{ z_{\nu} \} \bigr)\> \geq\> \sup_{\nu } \, 
\inf_{\mu \geq \nu }\, \bigl\{ {\pi }_q^{\nu } \bigl( \{ z_{\mu} \} \bigr) \bigr\}\> =\> 
\inf_{\nu}\, \sup_{\mu \geq \nu }\, \bigl\{ {\pi }_q^{\nu } \bigl( \{ z_{\mu } \} \bigr) \bigr\} $$
$$\quad \geq\> {\pi }_q^{\nu } \bigl( \limsup_{\nu }\, \{ z_{\nu } \} \bigr)\>\geq\> {\pi }_q^{\nu } \bigl( \liminf_{\nu }\, \{ z_{\nu} \} \bigr) $$
which follow from monotonicity and the fact that $\, {\pi }^{\nu }\, $ is decreasing normal on basic elements and increasing normal on antibasic elements. 
Therefore in $\, {\mathfrak L}^{\nu }_q \bigl( \mathcal B ( \mathcal H ) \bigr)\, $ one has 
$$ \inf\, \{ c\, z_x\,c \}\> =\> \sup_{\mu }\, {\pi }^{\nu } \bigl( \underline{\mathcal A}_{\mu } \bigr) $$
where $\, \underline{\mathcal A}_{\mu } = \inf_{\nu\geq\mu}\, \overline{\bigl\{ z_{\nu } \bigr\}}^{w^*}\in {\mathfrak L}^{\nu } ( R )\, $ is the $w^*$-closed basic element corresponding to the subset 
$\, \bigl\{ z_{\nu } \bigm\vert \nu \geq \mu \bigr\}\, $. Clearly the monotone increasing net 
$\, \bigl( \underline{\mathcal A}_{\mu } {\bigr)}_{\mu }\nearrow \inf\, \bigl\{ c\, z_x\, c \bigr\}\, $ is uniformly boundedly generated. It is easy to see that the minimal monotonous retraction $\, \underline\pi\, $
as above satisfies $\, \underline\pi \bigl( {\pi }^{\nu } \bigl( \underline{\mathcal A}_{\mu } \bigr) \bigr)\> =\> 0\, $ for every $\, \mu\, $ whence also $\, \underline\pi \bigl( c\, z_x\, c \bigr) = 0\, $ follows for every $\, c\in R^{sa}\, $ and every $\, z_x\in {\mathcal J}^-\, $ as above.
Recall the following general scheme from \cite{HSch}, $\S\, 6\, $: given a positive concave (= monotonous) map $\, c : L \rightarrow M\, $  of complete linear function lattices which is increasing normal for a specified family $\, \mathcal F\, $ of monotone increasing nets in $\, L_+\, $ invariant under affine translation, i.e. if $\, c\, $ is increasing normal for the net $\, ( a_{\lambda } {)}_{\lambda } \nearrow a\, $ then also for every translated net $\, ( a_{\lambda } + b {)}_{\lambda } \nearrow a + b\, $ with $\, b\in L_+\, $, then there is a procedure to bend up the map $\, c\, $ step by step to an affine (linear) map using transfinite induction by stepwise affinization of $\, c\, $ for addition of positive scalar multiples of a fixed positive element with respect to a chosen well order on the set $\, L_{+ , 1} = \bigl\{ a\in L_+\,\vert\, \Vert a\Vert = 1 \bigr\}\, $. Namely if $\, c_{<\omega }\, $ denotes the concave map which is the result of the affinization of $\, c\, $ for all elements $\, \{ a_{\kappa }\,\vert \kappa < \omega \}\, $ which is assumed to exist with the property of being positive concave and increasing normal for $\, \mathcal F\, $ by induction assumption then the affinization of $\, c_{<\omega }\, $ with respect to the ray $\, {\mathbb R}_+\, a_{\omega }\, $ is defined by the formula 
$$ c_{\omega } ( a )\> =\> \sup_{t\to\infty }\, c_{<\omega } ( a + t a_{\omega } )\> -\> t c_{<\omega } ( a_{\omega } ) \> . $$
In order that the supremum exists it is convenient to assume that $\, c_{<\omega}\, $ is affine already for addition of positive scalar multiples of the unit element $\, {\bf 1}_L\, $ and satisfies $\, c_{<\omega } ( {\bf 1}_L ) = {\bf 1}_M\, $. Then it is easily checked from monotonicity that the expressions to the right remain within a bounded subset. One checks the following facts: $\, c_{\omega }\, $ is positive concave and increasing normal for $\,\mathcal F\, $. Moreover generally $\, c_{\omega } ( a ) \geq c_{<\omega } ( a )\, $ and $\, c_{\omega } ( a ) = c_{<\omega } ( a )\, $ whenever $\, c_{<\omega }\, $ is already affine for addition of $\, {\mathbb R}_+ a\, $ so that changes are only made for elements which do not yet belong to the affine domain of $\, c_{<\omega }\, $. One proceeds by transfinite induction since the supremum of a monotone increasing net of positive concave and $\mathcal F $-increasing normal maps has the same properties. Thus having completed the affinization procedure for all elements in $\, L_+\, $ one arrives at a positive affine map $\, \overline c : L_+ \rightarrow M_+\, ,\, \overline c \geq c\, $ uniquely extending to a positive linear map $\, L \rightarrow M\, $ which is increasing normal for $\,\mathcal F\, $. Applying this general scheme in the present context the map $\, c\, $ is given by $\, \underline\pi\, $ and the family $\, \mathcal F\, $ is given by translates $\, \bigl\{ \underline{\mathcal A}_{\lambda } + \mathfrak B {\bigr\}}_{\lambda }\, $ of sets of uniformly boundedly generated positive basic elements {\it in the image of $\, \mathfrak L ( {\mathcal J}^- )\, $} where the translates $\, \mathfrak B\, $ are also assumed to be of this form, i.e. $\, \mathfrak B = \sup_{\mu } 
\bigl\{ {\pi }^{\nu } \bigl( \underline{\mathcal B}_{\mu } \bigr) \bigr\}\, $ for some set of uniformly boundedly generated positive basic elements $\, \bigl\{ \underline{\mathcal B}_{\mu } \bigr\}\subseteq \mathfrak L ( {\mathcal J}^- ) \subseteq \mathfrak L ( R )\, $. This is a bit more restrictive than in the general scheme not allowing for arbitrary translates but sufficient since we may start the affinization process with these elements choosing a corresponding well order on the set $\, \bigl\{ \mathfrak P\in {\mathfrak L}_q ( \mathfrak V )\bigm\vert \Vert\mathfrak P\Vert = 1 \bigr\}\, $. Then having completed the affinization process with respect to all elements which are suprema of such type the images of these elements will remain constant under any further affinization, hence the maps will remain increasing normal for $\,\mathcal F\, $. Since $\, \underline\pi \, $ is uniquely determined on (and linear for addition of) elements from $\, \mathfrak L ( {\mathcal J}^- )\, $ by Proposition 3 the restriction to sets of positive basic elements has no impact for the general result since any set of uniformly boundedly generated basic elements can be translated to a set of positive uniformly boundedly generated basic elements upon adding a scalar multiple of the unit element.
Note that $\, {\pi }_c\, $ is linear with respect to addition of elements $\, \mathcal B\in\mathfrak L( {\mathcal J}^- )\, $ since 
$$ \underline\pi  \bigl( \mathfrak A - \mathcal B \bigr)\> \leq\> \underline\pi \bigl( \mathfrak A \bigr)\> -\> {\pi }^-  \bigl( \mathcal B \bigr)\> =\> \underline\pi \bigl( \mathfrak A \bigr)\> +\> {\pi }^- \bigl( - \mathcal B \bigr)\> \leq\> \underline\pi \bigl( \mathfrak A - \mathcal B \bigr)  $$
on using concavity of $\, \underline\pi\, $ and linearity of $\, {\pi }^- : \mathfrak L ( {\mathcal J}^- ) \twoheadrightarrow {\mathfrak L}_1 ( I )\, $. This implies that the values of 
the images of any supremum over a set of uniformly boundedly generated basic elements in the image of $\, \mathfrak L ( {\mathcal J}^- )\, $ must remain fixed under the affinization process since the values of the single basic elements remain fixed. Then the linear limit map $\, \sigma\, $ of the affinization process coincides with $\, \underline\pi\, $ on all such elements. In particular the positive linear retraction $\, \sigma \, $ satisfies $\, \sigma ( c\, z_x\, c ) = 0\, $ for every $\, c\in R^{sa}\, $ and $\, z_x\in {\mathcal J}^-\, $ as above. Then being positive and linear $\, \sigma\, $ satisfies $\, \sigma ( a\, z_x\, b ) = 0\, $ for every 
$\, a\, ,\, b\in R^{sa}\, $, thus the $C^*$-ideal $\, J\vartriangleleft R\, $ generated by all elements 
$\, \bigl\{ z_x\bigm\vert x\in {\mathcal J}^+ \bigr\}\, $ is contained in the kernel of $\, \sigma\, $. If this ideal is trivial then $\, z_x = 0\, $ for any $\, x\in {\mathcal J}^+\, $ implying $\, \mathcal J = \iota ( I )\, $ and a forteriori $\, R = \mathcal J = \iota ( I )\, $ which proves (ii).
\par\noindent
Let $\, A\, $ be any unital $C^*$-algebra with ultraweak Baire envelope $\, {\mathfrak B}^{\omega } ( A )\subseteq A^{**}\, $ and ultraweak regular envelope $\, {\sigma }^r_{\omega } ( A ) = {\mathfrak B}^{\omega } ( A )\bigm/ {\mathcal I}^{\omega }\, $ where $\, {\mathcal I}^{\omega } = p_X^{-1} \bigl( {\mathcal M}_X \cap p_X \bigl( {\mathfrak B}^{\omega } ( A ) \bigr) \bigr)\, $ as above. We wish to show that $\, {\mathfrak B}^{\omega } ( A )\, $ is a sub-$C^*$-algebra of $\, A^{**}\, $ with twosided ideal $\, {\mathcal I}^{\omega }\vartriangleleft {\mathfrak B}^{\omega } ( A )\, $ so that the quotient map is in fact a surjective (sequentially normal) $*$-homomorphism. Note that for any given ultrafilter, any bounded sequence of elements $\, \{ x_n \} \subseteq {\mathfrak B}^{\omega } ( A )\, $ and any countable collection of states $\, \{ {\phi }_{\kappa } \}\subseteq \mathcal S ( A )\, $ it is possible to find a subsequence $\, \{ x_{n_k} {\}}_k\, $ with $\, \{ n_k {\}}_k \subseteq \omega\, $ such that each of the sequences $\, \{ {\phi }_{\kappa } ( x_{n_k} ) {\}}_k\, $ converges to the value of $\, {\phi }_{\kappa } ( x_{\omega } )\, $. From this one finds that multiplication from the left or the right with a fixed element is sequentially ultraweakly continuous, i.e. given 
$\, x_{\omega } = \lim_{\omega }\, \{ x_n \}\, $ with $\, x_n\, ,\, y\in A^{**}\, $ one has 
$$\lim_{\omega }\, \{ y\, x_n \} = y\, x_{\omega }\> ,\quad \lim_{\omega }\, \{ x_n\, y \}\> =\> x_{\omega }\, y \> . $$
In fact for $\, \phi ( x ) = \langle x\, \xi\, ,\, \xi \rangle\, $ with $\, \xi\in \mathcal H\, $, the Hilbert space of the universal representation of $\, A\, $, one computes
$$ \phi \bigl(  (y\, x_n)_{\omega } \bigr)\> =\> \lim_{\omega }\, \langle y\, x_n\,\xi\, ,\, \xi \rangle\> =\> 
\lim_{\omega }\, \langle x_n\,\xi\, ,\, y^*\, \xi \rangle\> =\> \lim_{\omega }\, \langle x_n\, \xi\, ,\, \eta \rangle  $$
where $\, \eta = y\, \xi\, $ so that choosing a subsequence $\, \{ n_k {\}}_k\subseteq \omega\, $ with $\, \{ x_{n_k} {\}}_k\, $ converging on the states corresponding to the onedimensional subspaces $\, \{ \mathbb C\, \xi\, ,\, \mathbb C\, \eta\, ,\, \mathbb C\, \bigl( \xi \pm \eta \bigr)\, ,\, \mathbb C\, \bigl( \xi\,  \pm i\, \eta \bigr) \}\, $ the polarization identity gives that $\, \phi ( y\, x_{\omega } )\, =\, 
\lim_{\omega }\, \phi ( y x_n )\, $ whence the result. Then $\, {\mathfrak B}^{\omega } ( A )\, $ is a sub-$C^*$-algebra of $\, A^{**}\, $. It is proved in \cite{Wr} that for each positive (or invertible) element $\, a\in A\, $ and $\, y\in\mathcal M\, $ one has $\, a\, y\, a^*\in\mathcal M\, $, the proof given there is only for positive elements $\, y\in {\mathcal M}_+\, $ but extends without difficulty to arbitrary elements, compare the proofs of Lemma 4 of \cite {Wr2} and Proposition 2.1 of \cite{Wr}. 
In fact for any positive invertible element $\, c\in A\, $ the map $\, {\tau }_c : \mathcal S ( A ) \rightarrow \mathcal S ( A )\, ,\, {\tau }_c ( \rho ) ( x ) = {\rho } ( c^2 )^{-1}\, {\rho } ( c\, x\, c )\, $ is an affine homeomorphism with inverse $\, {\tau }_{c^{-1}}\, $ so that it maps $\, X = \overline{\mathcal P ( A )}\, $ homeomorphically onto itself proving that for 
$\, y\in p_X^{-1} \bigl( {\mathcal M}_X \cap p_X \bigl( A^{**} \bigr) \bigr)\, $ and any positive invertible $\, c\in A\, $ also $\, p_X \bigl( c\, y\, c \bigr)\in {\mathcal M}_X \cap p_X \bigl( A^{**} \bigr)\, $. From this one obtains the relation 
$\, a\, y + y\, a \in {\mathcal I}^{\omega }\, $ whenever $\, a\in A\, ,\, y\in {\mathcal I}^{\omega }\, $. By induction this property passes to arbitrary elements of 
$\, {\mathfrak B}^{\omega } ( A )\, $, i.e. $\, x\, y + y\, x \in {\mathcal I}^{\omega }\, $ for any $\, x\in {\mathfrak B}^{\omega } ( A )\, ,\, y\in {\mathcal I}^{\omega }\, $. Assume by induction that this relation holds for all elements of a given subspace 
$\,  V \subseteq {\mathfrak B}^{\omega } ( A )\, $ and let $\, \{ x_n \} \subseteq V\, $ be a bounded sequence converging to $\, x_{\omega }\, $ along $\, \omega\, $. Then the union 
$$ \bigcup_{n\in\mathbb N}\, \bigl\{ \rho\in X \bigm\vert \rho ( x_n\, y\, + y\, x_n ) \neq 0 \bigr\} $$
is meagre and for any point $\, \phi\, $ of the complement of this meagre subset one has 
$$ \phi \bigl( x_{\omega }\, y\, + y\, x_{\omega } \bigr)\> =\> \phi \bigl( x_{\omega }\, y \bigr)\> +\> \phi \bigl(  
y\, x_{\omega } \bigr)\> =\> 
\lim_{\omega }\, \phi \bigl( x_n\, y \bigr)\> +\> \lim_{\omega }\, \phi \bigl( y\, x_n \bigr) $$
$$\quad =\> 
\lim_{\omega }\, \phi \bigl( x_n\, y + y\, x_n \bigr)\> =\> 0\> . $$
By induction one finds that $\, x\, y + y\, x \in {\mathcal I}^{\omega }\, $ for any $\, x\in {\mathfrak B}^{\omega } ( A )\, $ and any $\, y\in {\mathcal I}^{\omega }\, $. Then if $\, y\in {\mathcal I}^{\omega }\, $ one obtains $\, y^2\in {\mathcal I}^{\omega }\, $ whence if $\, p_{\pm }\in {\mathfrak B}^{\omega } ( A )\, $ denotes the support projection of $\, y_{\pm }\, $ (the support projection of any element is contained in $\, {\mathfrak B}^{\omega } ( A )\, $ by sequential monotone completeness, compare with section 3 or \cite{KaPe}) $\, p_{\pm }\, $ commutes with $\, y^2\, $ so that $\, y_{\pm }^2 = p_{\pm }\, y^2 = y^2\, p_{\pm }\in {\mathcal I}^{\omega }\, $. Finally $\, {\mathcal I}^{\omega }_+\, $ is obviously closed under taking squareroots whence $\, y_{\pm }\in {\mathcal I}^{\omega }\, $. Thus $\, {\mathcal I}^{\omega }\, $ is positively generated and is an order ideal.
Next we check on the regularity condition of $\, {\sigma }^r_{\omega } ( A )\, $. We first show that for any  element $\, x\in {\mathfrak B}^{\omega } ( A )\, $ corresponding to the bounded function 
$\, f_x\, $ on $\, X\, $ there exists an upper semicontinuous function $\, u\, $ which is near $\, f_x\, $ in the sense that the subset of points 
$\, \bigl\{ \rho\in X\bigm\vert u ( \rho ) \neq f_x ( \rho ) \bigr\}\, $ is meagre. If this is the case there is also a lower semicontinuous function $\, \breve u\, $ which is near $\, f_x\, $ on letting
$$ \breve u ( \rho )\> =\> \sup\, \bigl\{ g ( \rho ) \bigm\vert g\in C ( X )\, ,\, g \leq u \bigr\} $$
(see \cite{Wr}, p. 300). From Corollary 4 of \cite{Wr3} this is true for any Borel measurable function. Let us call a function {\it almost Borel measurable} if it is near to a Borel measurable function in the sense above, which then implies that it is near to a complemented (upper or lower) semicontinuous function. Let $\, \mathfrak B ( Y )\subseteq B ( Y )\, $ denote the subspace of bounded Borel functions and $\, {\mathfrak B}^{\omega } ( X )\subseteq B ( Y )\, $ the ultraweak sequential closure of $\, C ( X )\, $ in $\, B ( Y )\, $. The subspace
$\, {\mathcal I}_X = {\mathcal M}_X \cap \mathfrak B ( Y )\, $ is positively generated and is a twosided ideal of (the abelian $C^*$-algebra) $\, \mathfrak B ( Y )\, $. From \cite{Dx} the quotient 
$\, \mathfrak B ( Y ) \bigm/ {\mathcal I}_X\, $ is isomorphic to the injective envelope 
$\, I ( C ( X ) ) \simeq {\mathfrak L}_1 ( C ( X ) ) \simeq {\mathfrak L}_1 ( A )\, $. 
If any element in $\, {\mathfrak B}^{\omega } ( X )\, $ is congruent to an element of $\, \mathfrak B ( Y )\, $ modulo $\, {\mathcal M}_X\, $ then also any element of $\, {\mathfrak B}^{\omega } ( A )\, $ is congruent to a bounded Borel function modulo 
$\, {\mathcal M}_X\, $. Thus it is sufficient to consider the commutative setting where $\, A = C ( X )\, $ and 
$\, {\mathfrak B}^{\omega } ( A ) = {\mathfrak B}^{\omega } ( X ) $. To prove the statement one proceeds by induction: assume given a subspace $\, V \subseteq {\mathfrak B}^{\omega } ( X )\, $ such that each function in $\, f\in V\, $ is near a complemented concave upper semicontinuous function $\, \Hat f\, $.
Trivially this condition is satisfied for $\, f\in C ( X )\, $. Let $\, \{ f_n \} \subseteq V\, $ be a bounded sequence. Given an ultrafilter $\, \omega\, $ on $\, \mathbb N\, $ the $\omega $-limit of the sequence is given by the function
$$ f_{\omega }\> =\> \inf_{S\in \omega }\, f_S\> =\> \inf_{S\in \omega }\, \sup_k\, \bigl\{ f_{n_k}\bigm\vert S = \{ n_k {\}}_k \bigr\} $$
which putting $\, {\breve f}_S = \sup_k\, \bigl\{ \breve{\Hat f}_{n_k}\bigm\vert S = \{ n_k {\}}_k \bigr\}\, ,\, 
{\Hat f}_S = \sup_k\, \bigl\{ {\Hat f}_{n_k}\bigm\vert S = \{ n_k {\}}_k \bigr\}\, $ is near both functions 
$$ {\breve f}_{\omega }\> =\> \lim_{\omega }\, \bigl\{ \breve{\Hat f}_n \bigr\}\> =\>   
\inf_{S\in \omega }\, \sup_k\, \bigl\{ \breve{\Hat f}_{n_k} \bigm\vert S = \{ n_k {\}}_k \bigr\}\> =\> 
\inf_{S\in \omega }\,\, {\breve f}_S\> , $$
$$ {\Hat f}_{\omega }\> =\> \lim_{\omega }\, \bigl\{ {\Hat f}_n \bigr\}\> =\> \inf_{S\in \omega }\, \sup_k\, \bigl\{ {\Hat f}_{n_k}\bigm\vert S = \{ n_k {\}}_k \bigr\}\> =\> \inf_{S\in \omega }\, {\Hat f}_S  $$
due to the fact that taking $\omega $-limits is linear and $\, {\mathcal M}_X\, $ is ultraweakly sequentially closed. Given any bounded function $\, g\, $ put 
$$ g^{\wedge }\> =\> \inf\, \bigl\{ a\in C ( X ) \bigm\vert a \geq g \bigr\}\> ,\quad g^{\vee }\> =\> \sup\, \bigr\{ b\in C ( X ) \bigm\vert b \leq g \bigr\} \> . $$
Thus we may assume that each function $\, f_n = \breve{\Hat f}_n\, $ is complemented convex lower semicontinuous
whence also each function $\, f_S = {\breve f}_S\, $ is convex lower semicontinuous and is near the complemented concave upper semicontinuous function 
$$  ( f_S )^{\wedge }\> =\> \inf\, \bigl\{ a\in C ( X )\bigm\vert a \geq f_S \bigr\}\> =\> \bigl( {\Hat f}_S \bigr)^{\wedge }\geq {\Hat f}_S\>  . $$
We claim that $\, f_{\omega }\, $ is near the concave upper semicontinuous function 
$$ u\> =\> \inf_{S\in \omega }\, \left( \left( {\Hat f}_S \right)^{\vee } \right)^{\wedge }\> =\> \inf_{S\in \omega }\, \left( f_S \right)^{\wedge } \> .  $$
Note that $\, u\, $ being the infimum of concave upper semicontinuous functions is again concave upper semicontinuous. Being convex lower (resp. concave upper) semicontinuous the functions $\, \{ f_n \}\, $ (resp. 
$\, \{ {\Hat f}_n \}\, $) define a sequence of uniformly bounded antibasic elements 
$\, \bigl\{ \overline{\mathcal B}_n \bigr\}\, ,\, \overline{\mathcal B}_n = \sup\, \bigl\{ d\in C ( X )\bigm\vert d ( \rho ) \leq f_n ( \rho )\, ,\, \forall\, \rho\in Y \bigr\}  \, $ (resp. of uniformly bounded basic elements $\, \bigl\{ \underline{\mathcal A}_n \bigr\}\, ,\, \underline{\mathcal A}_n = \inf\, \bigl\{ c \in C ( X )\bigm\vert c ( \rho ) \geq f_n ( \rho )\, ,\, \forall\, \rho\in Y \bigr\}\, $) of $\, \mathfrak L ( C ( X ) )\, $. Note that basic elements $\, \underline{\mathcal A}\in \mathfrak L ( C ( X ) )\, $ which are infima of maximal representatives with respect to evaluation on the subset of pure states are automatically bounded in a commutative setting since $\, C ( X )\, $ is a lattice and evaluation on a pure state is a lattice map so that if two different elements $\, c\, ,\, d\in C ( X )\, $ dominate $\, \underline{\mathcal A}\, $ then also $\, c \wedge d\, $ will dominate $\, \underline{\mathcal A}\, $. Also if $\, u\, $ is any concave upper semicontinuous function then it is near the concave upper semicontinuous $\, \overline u\, $ representing the basic element $\, \underline{\mathcal A}_{\overline u}\, $ which is the infimum of maximal representative set $\, {\mathcal A}_{\overline u}\, $ such that $\, r ( \underline{\mathcal A}_{\overline u } ) = u {\vert }_X\, $. In particular complemented concave upper semicontinuous functions are always of the form $\, s ( \underline{\mathcal A}_{\overline u} )\, $, hence bounded.
Also note that from Theorem 1 $\, B ( X )\, $ is equal to the minimal monotone completion of $\, r \bigl( \mathfrak L ( C ( X ) ) \bigr)\, $, i.e. $\, B ( X ) \simeq 
{\mathfrak L}_1 \bigl( r \bigl( \mathfrak L ( C ( X ) ) \bigr) \bigr)\, $.
Then 
$$ \inf_{S\in \omega }\, {\Hat f}_S\> =\> \inf_{S\in\omega }\, \sup_k\, {\Hat f}_{n_k}\> = {\Hat f}_{\omega }\> =\> \sup_{S\in \omega }\, \inf_k\, {\Hat f}_{n_k}\> =\> \sup_{S\in \omega }\, {\Hat g}_S\> , $$
$$ \inf_{S\in \omega }\, f_S\> =\> \inf_{S\in \omega }\, \sup_k\, f_{n_k}\> =\> f_{\omega }\> =\> 
\sup_{S\in \omega }\, \inf_k\, f_{n_k}\> =\> \sup_{S\in \omega }\, g_S\> . $$
Put
$$\quad f_{\omega }^c\> =\> \inf_{S\in \omega }\, ( f_S )^{\wedge }\> =\> \inf_{S\in \omega }\, ( {\Hat f}_S )^{\wedge }\> , \quad f_{\omega }^{cc}\> =\> \inf_{S\in \omega }\, \bigl( ( {\Hat f}_S )^{\wedge } \bigr)^{\vee }\> =\> \inf_{S\in \omega }\, \bigl( ( f_S )^{\wedge } \bigr)^{\vee }\> , $$
$$\quad ( f_{\omega } )_c\> =\> \sup_{S\in \omega }\, ( g_S )^{\vee }\> =\> \sup_{S\in \omega }\, ( {\Hat g}_S )^{\vee }\> ,\quad ( f_{\omega } )_{cc}\> =\> 
\sup_{S\in \omega }\, \bigl( ( {\Hat g}_S )^{\vee } \bigr)^{\wedge }\> =\> \sup_{S\in \omega }\, \bigl( ( g_S )^{\vee } \bigr)^{\wedge } $$
where $\, S = \{ n_k {\}}_k\, $ is understood, all of which are elements of $\, B ( Y )\, $. The canonical map $\, \pi : \mathfrak L ( C ( X ) ) \twoheadrightarrow {\mathfrak L}_1 ( C ( X ) )\, $ factors over the injective function representation 
$\, s : \mathfrak L ( C ( X ) ) \twoheadrightarrow B ( Y )\, $. Then the maximal convex (resp. minimal concave) extension 
$$ {\overline\pi }^c\, ,\, {\underline\pi }_c :\> B ( Y )\> \twoheadrightarrow\> {\mathfrak L}_1 ( C ( X ) )\> ,\quad {\overline\pi }^c ( f )\> =\> \inf\, \bigl\{ \pi ( c )\bigm\vert c\in C ( X )\, ,\,  c\geq f \bigr\} $$
of $\, \pi\, $ is boundedly antibasically decreasing normal (resp. boundedly basically increasing normal)  from Proposition 3, i.e. 
$$ {\overline\pi }^c \bigl( \inf_{\lambda }\, s \bigl( \overline{\mathcal B}_{\lambda } \bigr) \bigr) =\, \inf_{\lambda }\, \pi \bigl( \overline{\mathcal B}_{\lambda }\, \bigr) $$
since $\, {\mathfrak L}^{\nu } \bigl( B ( X ) \bigr) \subseteq B ( Z )\, $  by the function representation 
$\, s^{\nu }\, $ where $\, Z = \overline{co \bigl( X )}\subseteq Y\, $ denotes the closed convex hull of $\, X\, $ in $\, Y\, $ and the maximal monotonous extension $\, \overline\pi : {\mathfrak L}^{\nu } \bigl( B ( X ) \bigr) \twoheadrightarrow {\mathfrak L}_1 \bigl( C ( X ) \bigr)\, $ is boundedly antibasically decreasing normal from Proposition 3 so there exists a monotonous extension $\, B ( Y ) \twoheadrightarrow {\mathfrak L}_1 ( C ( X ) )\, $ of the identity of $\, C ( X )\, $ factoring over the (normal) evaluation map to $\, B ( Z )\, $ which is boundedly antibasically decreasing normal for the image of $\, \mathfrak L ( C ( X ) )\, $. But then the maximal monotonous extension also must have this property.
We first show that 
$\, {\overline\pi }^c \bigl( {\Hat f}_{\omega } \bigr)\> =\> {\overline\pi }^c \bigl( f_{\omega } \bigr)\, $. One has 
$$ {\pi }^c \bigl( {\Hat f}_{\omega } \bigr)\> =\> \pi \bigl( f_{\omega }^c \bigr)\> \leq\> 
\inf_{S\in \omega }\, \pi \bigl( ( {\Hat f}_S )^{\wedge } \bigr)\> =\> \inf_{S\in \omega }\, \pi \bigl( f_S \bigr)\> =\> {\overline\pi }^c \bigl( f_{\omega } \bigr) \> . $$
where the last equality follows from Proposition 3 by the fact that the antibasic elements 
$\, \bigl\{ \overline{\mathcal B}_S\bigm\vert \overline{\mathcal B}_S = \sup\, \bigl\{ \overline{\mathcal B}_{n_k}\,\vert\, S = \{ n_k \} \bigr\} \bigr\} \, $ are uniformly bounded. The reverse inequality follows from $\, {\Hat f}_{\omega } \geq f_{\omega }\, $. Consider the order ideal 
$\, \mathfrak N \subseteq B ( Y )\, $ linearly generated by positive elements in the kernel of the convex map $\, {\overline\pi }^c\, $ which subspace is contained in the kernel of any positive linear extension 
$\, {\overline\pi }_q : B ( Y ) \twoheadrightarrow {\mathfrak L}_1 ( C ( X ) )\, $ of the identity map of $\, C ( X )\, $ by maximality of $\, {\overline\pi }^c\, $ among all monotonous extensions. 
We claim that $\, \mathfrak N\, $ contains the element $\, f_{\omega }^c - f_{\omega }\, $. Since $\, {\pi }^c\, $ is linear with respect to addition of elements in $\, \mathfrak L ( C ( X ) )\, $ from Proposition 1 and 
$\, f_{\omega }^{cc}\, $ is the infimum of the net $\, \{ \bigl( ( f_S )^{\wedge } \bigr)^{\vee } {\}}_{S\in \omega }\, $ of images of complemented uniformly bounded antibasic elements it follows from Proposition 3 that 
$\, {\underline\pi }_c \bigl( f_{\omega }^{cc} \bigr) = {\overline\pi }^c \bigl( f_{\omega }^{cc} \bigr) = \inf_{S\in \omega }\, \pi \bigl( \bigl( ( f_S )^{\wedge } \bigr)^{\vee } \bigr) = \inf_{S\in \omega }\, \pi \bigl( f_S \bigr)\, $ whence
$$ {\overline\pi }^c \left( f_{\omega }^c\, -\, f_{\omega }^{cc} \right)\> =\> 
\pi \bigl( f_{\omega }^c \bigr)\> -\> {\underline\pi }_c \bigl( f_{\omega }^{cc} \bigr)\> =\> \pi \bigl( f_{\omega }^c \bigr)\> -\> {\overline\pi }^c \bigl( f_{\omega }^{cc} \bigr)\> =\> 0  $$
since
$$ \pi \bigl( f_{\omega }^c \bigr)\> =\> \inf_{S\in \omega }\, \pi \bigl( f_S \bigr) $$
from the argument above. Thus $\, f_{\omega } - f_{\omega }^{cc}\, $ is contained in $\, \mathfrak N\, $. One computes
$$ {\overline\pi }^c \left( f_{\omega }^{cc}\, -\, f_{\omega } \right)\> =\> 
{\overline\pi }^c \left( \inf_{S\in \omega }\, \bigl( ( f_S )^{\wedge } \bigr)^{\vee }\, -\, \sup_{S\in \omega }\, g_S \right)\> \leq\> {\overline\pi }^c \left( \inf_{S\in \omega }\, \bigl( ( f_S )^{\wedge } \bigr)^{\vee } \right)\> -\> {\overline\pi }_c \left( \sup_{S\in \omega }\, g_S \right) $$
$$\> =\> {\overline\pi }^c \left( \inf_{S\in \omega }\, \bigl( ( f_S )^{\wedge }\bigr)^{\vee } \right)\> -\> {\overline\pi }_c \left( 
\sup_{S\in \omega }\, {\Hat g}_S \right)\> =\> \inf_{S\in \omega }\, \pi  \left( f_S \right)\> -\> \sup_{S\in \omega }\, \pi \left( {\Hat g}_S \right) $$
$$\>\leq\> \inf_{S\in\omega }\, \pi \bigl( f_S\, -\, {\Hat g}_S \bigr)\> =:\> \Delta \> .  $$ 
The last inequality in the first line and the first equality in the second line follow from convexity of $\, {\overline\pi }^c\, $ plus the identity 
$$ {\overline\pi }_c \bigl( {\Hat f}_{\omega } \bigr)\> =\> {\overline\pi }_c \bigl( f_{\omega } \bigr) $$ 
proved analogous to the relation $\, {\overline\pi }^c \bigl( {\Hat f}_{\omega } \bigr) = {\overline\pi }^c \bigl( f_{\omega } \bigr)\, $ as above. The second equation in the second line again follows from Proposition 3. We want to show that 
$\, \Delta \leq 0\, $. Then it is sufficient to show that the element 
$\, \inf_S\, \bigl( f_S\, -\, {\Hat g}_S \bigr)\, $ is negative.
For any given point $\, \rho\in Y\, $ and any $\, S_0\in \omega\, $ one has
$$ \inf_{S\in \omega }\, \bigl\{ f_S  ( \rho )\, -\, {\Hat g}_S ( \rho ) \bigr\}\> =\> \inf_{S {\succ }_{\omega } S_0}\, \bigl\{ f_S ( \rho )\, -\, {\Hat g}_S ( \rho ) \bigr\}\> $$ 
where the notation $\, S {\succ }_{\omega } S_0\, $ is shorthand for $\, S \subset S_0\, ,\, S\in \omega\, $,
on noting that each such subnet is cofinal. But for each state $\, \rho\, $ and given $\, \epsilon > 0\, $ there exists a subsequence 
$\, S_{\rho , \epsilon }\in \omega\, $ with the property that $\, \vert f_{\omega } ( \rho ) - f_{n_k} ( \rho ) \vert \leq \epsilon\, $ and $\, \vert {\Hat f}_{\omega } ( \rho ) - {\Hat f}_{n_k} ( \rho ) \vert \leq \epsilon\, $ for all $\, n_k\in S_{\rho , \epsilon }\, $ so that since $\, {\Hat f}_{\omega } \geq f_{\omega }\, $ one gets 
$$ \inf_{S {\succ }_{\omega } S_{\rho , \epsilon }}\, \bigl\{ f_S ( \rho )\, -\, {\Hat g}_S ( \rho ) \bigr\}\> \leq 2\, \epsilon $$
and since $\, \epsilon > 0\, $ is arbitrary 
$$ \inf_{S\in \omega }\, \bigl\{ f_S ( \rho )\, -\, {\Hat g}_S ( \rho ) \bigr\}\> \leq\> 0 $$
whence also $\, \Delta\> \leq\> 0\, $ and $\, {\overline\pi }^c \bigl( f_{\omega }^c - f_{\omega } \bigr) = 0\, $ follows from positivity of $\, f_{\omega }^c - f_{\omega }\, $ and convexity of $\, {\overline\pi }^c\, $ since 
$$ {\overline\pi }^c \bigl( f_{\omega }^c\, -\, f_{\omega } \bigr)\> \leq\> {\overline\pi }^c \bigl( f_{\omega }^c\, -\, f_{\omega }^{cc} \bigr)\> +\> {\overline\pi }^c\bigl( f_{\omega }^{cc}\, -\, f_{\omega } \bigr)\>\leq\> 0\> . $$  Therefore $\, f_{\omega }^c - f_{\omega }\in \mathfrak N\, $ as desired. Similarly 
$\, {\Hat f}_{\omega }\, -\, ( f_{\omega } )_c \in \mathfrak N\, $ implying 
$\, f_{\omega }^c - ( f_{\omega } )_c\in \mathfrak N\, $. Since both functions $\, f_{\omega }^c\, $ and $\, ( f_{\omega } )_c\, $ are Borel and in the image of $\, \mathfrak L ( C ( X ) )\, $ so that their image in $\, {\mathfrak L}_1 ( C ( X ) )\, $ is uniquely determined this necessitates $\, f_{\omega }^c - ( f_{\omega } )_c\in {\mathcal M}_X\, $ and since $\, {\mathcal M}_X\, $ is an order ideal one gets $\, f_{\omega }^c - f_{\omega }\in {\mathcal M}_X\, $ so both $\, f_{\omega }\, $ and $\, {\Hat f}_{\omega }\, $ are near the upper semicontinuous function 
$\, f_{\omega }^c\, $ as desired completing the induction step.
The proof that $\, {\sigma }^r_{\omega } ( A ) = {\mathfrak B}^{\omega } ( A ) \bigm/ {\mathcal I}^{\omega }\, $ is regular now is completely analogous to the corresponding proof for the regular $\sigma $-completion $\, {\sigma }^r ( A )\, $ as in \cite{Wr}, Theorem 2.3. Therefore $\, {\sigma }^r_{\omega } ( A )\, $ embeds as a  natural sequentially monotone complete subspace of $\, {\mathfrak L}_1 ( A )\, $ and $\, {\mathcal I}^{\omega }\, $ being positively generated and the kernel of a positive map from a $C^*$-algebra into another $C^*$-algebra is necessarily a twosided ideal so that $\, {\sigma }^r_{\omega }\, $ inherits the structure of a (sequentially monotone complete) $C^*$-algebra containing $\, A\, $. Since the ideal $\, {\mathcal I}^{\omega } \vartriangleleft {\mathfrak B}^{\omega } ( A )\, $ is ultraweakly sequentially closed (and hence also weakly sequentially closed) $\, {\sigma }^r_{\omega } ( A )\, $ inherits both topologies from $\, {\mathfrak B}^{\omega } ( A )\, $. By the same argument the $\omega $-limit $\, x_{\omega }\, $ of any bounded sequence $\, \{ x_n \} \subseteq {\sigma }^r_{\omega } ( A )\, $ is well defined in $\, {\sigma }^r_{\omega } ( A )\, $ and is an ultraweak sequential limit point of the sequence whence the unit ball of $\, {\sigma }^r_{\omega } ( A )\, $ is seen to be sequentially compact for the ultraweak sequential topology, moreover any ultrafilter $\, \omega\, $ defines a positive linear 
map $\, \widehat\omega : {\prod }_{\mathbb N}\, {\sigma }^r_{\omega } ( A ) \rightarrow {\sigma }^r_{\omega } ( A )\, $. 
\par\noindent
Assume now that $\, A\, $ is separable and let $\, I ( A )\, $ denote the injective envelope of $\, A\, $. Since $\, A \subseteq {\sigma }^r_{\omega } ( A )\, $ is an injective $*$-homomorphism it is completely positive so that the identity map of $\, A\, $ extends to a completely positive linear map 
$\, \Psi : {\sigma }^r_{\omega } ( A ) \rightarrow I ( A )\, $ and from injectivity of $\, {\mathfrak L}_1 ( A )\, $ there is a (completely) positive linear map $\, \Gamma : I ( A ) \rightarrow {\mathfrak L}_1 ( A )\, $ extending the identity map of $\, A\, $. The composition 
$$ \Gamma\circ\Psi :\> {\sigma }^r_{\omega } ( A )\> \largerightarrow\> {\mathfrak L}_1 ( A ) $$
is equal to the natural inclusion $\, {\sigma }^r_{\omega } ( A ) \subseteq {\mathfrak L}_1 ( A )\, $ by rigidity. Therefore $\, \Psi\, $ is a unital isometric embedding of $\, {\sigma }^r_{\omega } ( A )\, $ into $\, I ( A )\, $. 
Put $\, J ( A ) = \Psi ( {\sigma }^r_{\omega } ( A ) )\, $. We claim that there exists a positive projection 
$\, \Phi : I ( A ) \rightarrow I ( A )\, $ with range $\, J ( A )\, $. On utilizing the map $\, \Psi\, $ the assertion is equivalent to the existence of a positive linear map $\, I ( A ) \twoheadrightarrow {\sigma }^r_{\omega } ( A )\, $ extending the identity map  of $\, A\, $ which from the completely positive embedding $\, I ( A ) \hookrightarrow R\, $ where $\, R = \lambda ( A )''\, $ is the enveloping von Neumann algebra of $\, A\, $ in some separable supertransitive representation will follow from the existence of a corresponding positive linear map $\, \gamma : R \twoheadrightarrow {\sigma }^r_{\omega } ( A )\, $ over the identity map of $\, A\, $. From \cite{A2} each state of $\, A\, $ is in the $w^*$-closure of a sequence of mutually orthogonal vector states of $\, \mathcal B \bigl( {\oplus }_n\, {\mathcal H}_n \bigr)\, $ with each $\, {\mathcal H}_n\, $ corresponding to an irreducible factor of $\, \lambda\, $. Since $\, A\, $ is contained in $\, R\, $ it follows that each state of $\, A\, $  can be $w^*$-approximated by a sequence of convex combinations  of vector states in $\, \cup_n\, \mathcal B ( {\mathcal H}_n )_*\, $.  
Let $\, \delta X \subseteq \mathcal P ( A )\, $ denote the subset of pure states obtained as restrictions of vector states from $\, \cup_n\, \mathcal B ( {\mathcal H}_n )_*\, $ and $\, X = co \bigl( \delta X \bigr)\, $ its convex hull which is $w^*$-dense in $\, \mathcal S ( A )\, $. Each element $\, x\in R^{sa}\, $ determines a bounded (affine) function $\, f_x\, $ on $\, X\, $ by evaluation. Functions corresponding to elements in $\, A\, $ extend uniquely to continuous affine functions on $\, \mathcal S ( A )\, $. If $\, P_n\in R\, $ denotes the support projection of the states $\, \{ p_1\, ,\,\cdots\, ,\, p_n \}\, $ viewed as elements of $\, R_*\, $ let $\, {\mathcal E}^n = P_n\, \mathcal H\subseteq \mathcal H\, $ be the corresponding finitedimensional subspace. Let $\, Y\subseteq A^*\, $ be the subspace generated by $\, X\, $ with $\, Y^n = {\lambda }^* \bigl( ( Ad\, P_n )_* \bigl( \mathcal B ( {\mathcal E}^n )_* \bigr) \bigr)\subseteq Y\, $ the finitedimensional subspace generated by the vector states subordinate to $\, P_n\, $ which is completely isometric to a sum of matrix space duals and put $\, X^n = Y^n \cap \mathcal S ( A )\, $. Then $\, X^n\, $ is $w^*$-closed since $\, Y^n\, $ is finitedimensional. For $\, v\in Y\, $ the assignment $\, v \mapsto v\circ Ad\, P_n\, $ where 
$\, \bigl( v\circ Ad\, P_n \bigr) \bigl( a \bigr) = v \bigl( P_n\, \lambda ( a )\, P_n \bigr)\, $ 
is well defined and completely positive, so by injectivity of $\, Y^n\, $ the map extends to a completely positive map on $\, A^*\, $ denoted $\, Ad\, P_n\, $ by abuse of notation. 
By Proposition 2 there exist for every $\, \epsilon > 0\, $  completely positive $w^*$-continuous maps 
$$ {\rho }^n_{\epsilon } :\> A^*\>\largerightarrow\> Y^n $$
such that $\, {\rho }^n_{\epsilon }\, $ is almost a retraction for the natural $w^*$-continuous embedding $\, {\Lambda }^n : Y^n \hookrightarrow A^*\, $, i.e. putting $\, {\Phi }^n_{\epsilon } = {\Lambda }^n\circ {\rho }^n_{\epsilon }\, $ and $\, {\Psi }^n_{\epsilon } = {\rho }^n_{\epsilon }\circ {\Lambda }^n\, $ one has $\, \Vert ( {\Phi }^n_{\epsilon } )^2 - {\Phi }^n_{\epsilon } \Vert < \epsilon\, $ and $\, \Vert {\Psi }^n_{\epsilon } - id\Vert < \epsilon\, $. Then also the map $\, 1 - {\Phi }_{\epsilon }^n\, $ is $w^*$-continuous and $\, \ker\, {\Phi }^n_{\epsilon }\, $ as well as $\, \bigl( 1 - {\Phi }^n_{\epsilon } \bigr) \bigl( A^* \bigr)\, $ are $w^*$-closed from the Krein-Smulian-Theorem. 
Since the intersection of $\, \ker\, {\Phi }_{\epsilon }^n\, $ with $\, Y^n\, $ is trivial for $\, \epsilon < 1, $ also the intersection of the $w^*$-closure $\, {\Sigma }^n\, $ of $\, {\Gamma }^n : = \ker\, {\Phi }_{\epsilon }^n \cap ( 1 - Ad\, P_n ) ( Y )\, $ with $\, Y^n\, $ is trivial and its intersection with $\, Y\, $ has finite codimension in $\, Y\, $. A finite dimensional (hence $w^*$-closed) complement $\, Z^n\, $ can be chosen in the subspace $\, \bigl( 1 - Ad\, P_n \bigr) ( Y )\, $. Moreover the complement $\, Z^n\, $ can be chosen in order that $\, Y^n + Z^n\, $ is positively generated.
Then one obtains a decomposition of $\, A^*\, $ as 
$$ A^*\> =\> Y^n\> +\> Z^n\> +\> {\Sigma }^n $$
and a decomposition
$$  {\widetilde Y}^n\> :=\> Y^n\> +\> Z^n\> +\> {\Gamma }^n $$  
(if $\, v = y + z + w\, $ with $\, y\in Y^n\, ,\, z\in Z^n\, ,\, w\in {\Gamma }^n\, $ then $\, \Vert y\Vert \leq 
\Vert v\Vert\, $ and there exists a constant $\, C > 0\, $ with $\, \Vert z\Vert\, ,\, \Vert w\Vert \leq C\, \Vert v\Vert\, $ so that $\, Y^n + Z^n + {\Sigma }^n\, $ is $w^*$-closed and coincides with $\, A^*\, $ and the $w^*$-closure of $\, {\widetilde Y}^n\, $ from the Krein-Smulian Theorem). To see this note that since 
$\, Z^n \cap \ker\, {\Phi }^n_{\epsilon } = \{ 0 \}\, $ and $\, Z^n\, $ being finitedimensional there exists a constant $\, d > 0\, $ with $\, \Vert {\Phi }^n_{\epsilon } ( z ) \Vert \geq d\, \Vert z\Vert\, $ so assume that the decomposition $\, Z^n + {\Gamma }^n\, $ is unbounded, i.e. there exists a sequence 
$\, \bigl\{ x_k \bigm\vert \Vert x_k\Vert = 1\, ,\, x_k = z_k + w_k \bigr\}\subseteq Z^n + {\Gamma }^n\, $ with $\, w_k\in {\Gamma }^n\, ,\, z_k\in Z^n\, ,\; \Vert z_k\Vert \geq k\, $. Then $\, \Vert {\Phi }^n_{\epsilon } ( x_k )\Vert = \Vert {\Phi }^n_{\epsilon } ( z_k )\Vert \geq d\, k \geq n\, $ for $\, k \geq n\, d^{-1}\, $ contradicting the fact that $\, \Vert {\Phi }^n_{\epsilon }\Vert \leq n\, $.  
Given $\, x\in R^{sa}\, $ define an element $\, f^n_x\in C_{aff } \bigl( \mathcal S ( A ) \bigr) \simeq A^{sa}\, $ by restriction of the linear function $\, {\widehat f}^n_x : A^* \rightarrow \mathbb C\, $ where $\, {\widehat f}^n_x\, $ is given by linear extension of 
$$ f^n_x ( p )\> =\> f_x ( p )\> ,\quad p\in X^n\> ,\qquad {\widehat f}^n_x ( z )\> =\> 0\> ,\quad z\in Z^n\> +\> {\Sigma }^n\> .  $$ 
That $\, f^n_x\, $ is continuous is easily checked from the fact that $\, Z^n + {\Sigma }^n\, $ is $w^*$-closed. To prove that the assignment $\, x\mapsto f^n_x\, $ is contractive it suffices by continuity to consider the $w^*$-dense subspace $\, {\widetilde Y}^n\, $. If $\, v = y + z + w \in {\widetilde Y}^n\, $ with $\, y\in Y^n\, ,\, z\in Z^n\, ,\, w\in \bigl( 1 - Ad\, P_n \bigr) ( {\Gamma }^n )\, $ and $\, \Vert v\Vert = 1\, $ then $\, \Vert y\Vert \leq 1\, $ since $\, Ad\, P_n\, $ is contractive on $\, Y\, $. Therefore $\, \Vert f^n_x \Vert \leq \Vert f_x\Vert \leq \Vert x\Vert\, $. 
It is clear that the maps $\, ( f^n_x )\, $ converge pointwise to $\, f_x\, $ on $\, Y\, $. 
Choose a free ultrafilter $\, \omega\, $ on $\,\mathbb N\, $ so that the sequence $\, \{ f^n_x \}\, $ converges along $\,\omega\, $ to a bounded affine almost Borel function 
$$ ( f^n_x )\>\buildrel \omega\over\largerightarrow\> f^{\omega }_x $$
determining an element $\, \widehat x\, $ of the ultraweak sequential Baire envelope $\, {\mathfrak B}^{\omega } ( A )\, $, the ultraweak sequential closure of $\, A\, $ in $\, A^{**}\, $. The application $\, x \mapsto \widehat x\, $ is clearly linear and contractive. The composite map $\, \gamma : R^{sa} \rightarrow {\sigma }^r_{\omega } ( A )^{sa}\, ,\, \gamma ( x ) = \upsilon ( \widehat x )\, $ is contractive and linear. We need to show that for $\, x\in A^{sa}\, $ this map restricts to the identity map with respect to the canonical identifications, i.e. $\, \upsilon ( \widehat x ) = \upsilon ( x )\, $ or $\, ( \pi\circ {\Upsilon }^{\omega } ) \bigl( \{ \widehat x \} \bigr) = ( \pi\circ {\Upsilon }^{\omega } ) \bigl( \{ x \} \bigr)\, $ in $\, {\mathfrak L}_1 \bigl( {\sigma }^r_{\omega } ( A ) \bigr) \simeq {\mathfrak L}_1 \bigl( A \bigr)\, $. 
Let $\, {\mathfrak L}^{\nu } \bigl( {\mathfrak B}^{\omega } ( A ) \bigr) = s^{\nu } \bigl( \mathfrak L \bigl( {\mathfrak B}^{\omega } ( A ) \bigr) \bigr) \subseteq 
l_{\infty } \bigl( \mathcal S ( A ) \bigr)\, $ denote the image of $\, \mathfrak L \bigl( {\mathfrak B}^{\omega } ( A ) \bigr)\, $ under the function representation on the space of sequentially normal states which identifies with the state space of $\, A\, $. Then it is plain to see that in the ${\mathfrak L}_1$-completion $\, {\mathfrak L}_q^{\nu } \bigl( {\mathfrak B}^{\omega } ( A ) \bigr)\, $ one has the estimates
$$ \liminf_n\, \bigl\{ x_n \bigr\}\> =\> \sup_n\, \inf_{m\geq n}\, \bigl\{ x_m \bigr\}\> \leq\> \bigl\{ \widehat x \bigr\}\> \leq\> \inf_n\, \sup_{m\geq n}\, \bigl\{ x_m \bigr\}\> =\> \limsup_n\, \bigl\{ x_n \bigr\} $$
the elements $\, \{ x_n \} \subseteq A^{sa}\, $ corresponding to the functions $\, \{ f^n_x \}\subseteq C_{aff} ( \mathcal S ( A ) )\, $. On the other hand the functorial map $\, \mathfrak L ( A ) \rightarrow \mathfrak L \bigl( {\mathfrak B}^{\omega } ( A ) \bigr) \twoheadrightarrow {\mathfrak L}^{\omega } \bigl( {\mathfrak B}^{\omega } ( A ) \bigr)\, $ is clearly injective. Then there is a monotonous  map 
$\, uc : {\mathfrak L}^{\nu }_q \bigl( {\mathfrak B}^{\omega } ( A ) \bigr) \rightarrow \underline{\mathfrak L} ( A )\, $ by taking upper complements in $\, A\, $, i.e. $\, uc \bigl( \mathfrak A \bigr) = \inf\, \bigl\{ a\in A \bigm\vert a \geq \mathfrak A \bigr\}\, $. Then $\, uc\, $ is a retraction for the natural embedding 
$\, \underline{\mathfrak L} ( A ) \subseteq {\mathfrak L}^{\nu } \bigl( {\mathfrak B}^{\omega } ( A ) \bigr)\, $
such that the image $\, uc \bigl( s^{\nu } ( \mathcal A ) \bigr)\, $ of an element $\, \mathcal A\in \mathfrak L \bigl( {\mathfrak B}^{\omega } ( A ) \bigr)\, $ coincides with the image $\, uc \bigl( \mathcal A \bigr)\, $ of the corresponding map $\, uc : \mathfrak L ( {\mathfrak B}^{\omega } ( A ) \bigr) \rightarrow \underline{\mathfrak L} ( A )\, $ of any preimage $\, \mathcal A\in \mathfrak L \bigl( {\mathfrak B}^{\omega } ( A ) \bigr)\, $. Let $\, \upsilon : {\mathfrak B}^{\omega } ( A ) \rightarrow {\sigma }^r_{\omega } ( A )\, $ denote the canonical sequentially normal  and ultraweakly sequentially continuous $*$-homomorphism and $\, \Upsilon  :\, \mathfrak L \bigl( {\mathfrak B}^{\omega } ( A ) \bigr)\rightarrow \mathfrak L \bigl( {\sigma }^r_{\omega } ( A ) \bigr)\, $ the functorial lattice extension of $\, \upsilon\, $. Then one obtains 
$$ \upsilon ( \widehat x )\> \leq\> uc \bigl( \{ \widehat x \} \bigr)\> \leq\> uc \bigl( \sup_{m\geq n}\, \bigl\{ x_m \bigr\} \bigr) $$
for all $\, n\in \mathbb N\, $ which relation is to be understood to hold  in $\, \mathfrak L \bigl( {\sigma }^r_{\omega } ( A ) \bigr)\, $ with respect to the functorial normal embedding 
$\, {\mathfrak L}_q ( A ) \subseteq {\mathfrak L}_q \bigl( {\sigma }^r_{\omega } \bigr)\, $. To see this note that since $\, uc ( \{ \widehat x \} ) \geq \{ \widehat x \}\, $ in $\, \mathfrak L \bigl( {\mathfrak B}^{\omega } ( A ) \bigr)\, $ and $\, \Upsilon \bigl( uc \bigl( \{ \widehat x \} \bigr) \bigr) = uc \bigl( \{ \widehat x \} \bigr)\, $
with $\, \Upsilon \bigl( \{ \widehat x \} \bigr) = \{ \upsilon ( \widehat x ) \}\, $ the result follows from monotonicity of $\, \Upsilon\, $. Any monotonous map 
$$ \mathfrak L \bigl( {\sigma }^r_{\omega } ( A ) \bigr)\> \largerightarrow\> {\mathfrak L}_1 \bigl( A \bigr)\> \simeq\> {\mathfrak L}_1 \bigl( {\sigma }^r_{\omega } ( A ) \bigr) $$ 
extending the identity map of $\, A\, $ is necessarily also the identity map for $\, {\sigma }^r_{\omega } ( A )\, $ whence it is unique and equal to the canonical surjection $\, \pi\, $ so that in particular the image of $\, \inf\, \{ \upsilon ( \widehat x ) \}\, $ in $\, {\mathfrak L}_1 ( A )\, $ is uniquely determined independent of the chosen monotonous extension of the identity of $\, A\, $. 
Then there is a commutative diagram 
$$ \vbox{\halign{ #&#&#\cr 
\hfil $\mathfrak L \bigl( {\sigma }^r_{\omega } ( A ) \bigr)$\hfil & \hfil $\buildrel \overline r\, ,\, \underline r\over\largerightarrow$\hfil &\hfil ${\mathfrak L}_q ( A )$\hfil \cr
\hfil $\pi\, \Bigm\downarrow  $\hfil  &&\hfil ${\pi  }^c \Bigm\downarrow {\pi }_c $\hfil \cr
\hfil ${\mathfrak L}_1 \bigl( {\sigma }^r_{\omega } ( A ) \bigr) $\hfil &\hfil $\buildrel\sim\over\largerightarrow $\hfil &\hfil $ {\mathfrak L}_1 \bigl( A \bigr)$\cr }}\>   $$
where $\, \underline r\, ,\, \overline r\, $ denote the basic and antibasic restriction maps for the functorial inclusion 
$\, {\mathfrak L}_q ( A ) \hookrightarrow {\mathfrak L}_q \bigl( {\sigma }^r_{\omega } ( A ) \bigr)\, $. Since $\, \pi : \mathfrak L ( A ) \rightarrow {\mathfrak L}_1 ( A )\, $ factors over the function representation on the (image in $\, \mathcal S ( A )\, $ of the) normal state space $\, {\mathcal S}^{\nu } ( R )\, $ denoted $\, {\pi }^{\nu } : \mathfrak L ( A ) \rightarrow l_{\infty } \bigl( {\mathcal S}^{\nu } ( R ) \bigr)\, $, i.e. $\, \pi = \overline\pi \circ {\pi }^{\nu }\, $ where $\, \overline\pi : l_{\infty } \bigl( {\mathcal S}^{\nu } ( R ) \bigr) \rightarrow {\mathfrak L}_1 ( A )\, $ is any monotonous extension of the identity map of $\, A\, $, it follows that 
$$ \pi \bigl( \{ \upsilon ( \widehat x ) \} \bigr)\>\leq \inf_n\, \pi \bigl( \sup_{m\geq n}\, \bigl\{ x_n \bigr\} \bigr) \> =\> \inf_n\, \overline\pi \bigl( {\pi }^{\nu } \bigl( \sup_{m\geq n} \bigl\{ x_m \bigr\} \bigr) \bigr)\> .  $$
Clearly the subset $\, \bigl\{ {\pi }^{\nu } \bigl( \sup_{m\geq n}\, \bigl\{ x_m \bigr\} \bigr) {\bigr\}}_n\, $ viewed as elements of $\, {\mathfrak L}^{\nu } ( R )\, $ form a uniformly boundedly generated and monotone decreasing sequence of antibasic elements so that taking 
$\, \overline\pi : {\mathfrak L}^{\nu }_q ( R ) \rightarrow {\mathfrak L}_1 ( A )\, $ to be the maximal convex monotonous extension of the identity of $\, A\, $ one gets 
$$ \inf_n\, \Bigl\{ ( \overline\pi\circ {\pi }^{\nu } ) \bigl( \sup_{m\geq n}\, \bigl\{ x_m \bigr\} \bigr) \Bigr\}\> =\> 
\overline\pi \Bigl( \inf_n\, \bigl\{ {\pi }^{\nu } \bigl( \sup_{m\geq n}\, \bigl\{ x_m \bigr\} \bigr) \bigr\} \Bigr) $$
from Proposition 3, $\, R\, $ being a direct product of type I factors. But the function 
$\, \inf_n\, \bigl\{ s^{\nu } \bigl( \sup_{m\geq n} \bigl\{ x_m \bigr\} \bigr) \bigr\}\, $ is clearly equal to 
$\, s^{\nu } \bigl( \{ x \} \bigr)\, $ since the functions $\, \{ s^{\nu } ( x_n ) \}\, $ converge pointwise to $\, s^{\nu } ( x )\, $ whence the onesided estimate 
$\, \upsilon ( \widehat x ) \leq x\, $ is established. The reverse inequality follows by a symmetric argument we leave to the reader. In particular 
$\, \gamma ( {\bf 1} ) = {\bf 1}\, $ showing that $\, \gamma\, $ is positive since contractive and unital. Thus $\, J ( A ) \simeq {\sigma }^r_{\omega } ( A )\, $ is a weak injective envelope for $\, A\, $. Being equal to the range of a positive projection of a monotone complete $C^*$-algebra $\, {\sigma }^r_{\omega } ( A )\, $ is monotone complete for separable $\, A\, $ (compare the proof of Theorem 6.1.3 of \cite{E-R}). This completes the proof of (iii).
\par\noindent
Let $\, R\, $ be a von Neumann algebra with predual $\, R_*\, $.  The proof that $\, R\, $ is weakly injective is similar to the argument for $\ {\sigma }^r ( A )\, $ but even simpler in that it does not require use of the Baire envelope, nor of Proposition 2, nor the use of the Jordan lattice completion. Consider any faithful representation $\, \gamma : R \rightarrow \mathcal B ( \mathcal H )\, $ which is a sum of irreducible representations $\, \gamma = \oplus {\gamma }_{\nu }\, $ with $\, {\gamma }_{\nu } : R \rightarrow \mathcal B ( {\mathcal H}_{\nu } )\, $ irreducible for each $\, \nu \, $ such that each normal state of $\, R\, $ can be approximated in norm by convex combinations of vector states corresponding to vectors in $\, \cup_{\nu } {\mathcal H}_{\nu }\, $. We will exhibit the existence of a positive retraction $\, \rho  : \gamma ( R )'' \rightarrow R\, $ for the inclusion $\, \gamma ( R ) \subseteq \gamma ( R )''\, $ from which the result follows by injectivity of $\, \gamma ( R )''\, $.
There exists an increasing net of finitedimensional projections 
$\, \{ P_{\lambda } {\}}_{\lambda\in \Lambda } \subseteq \gamma ( R )''\, $ with corresponding finitedimensional subspaces 
$\, \bigl\{ {\mathcal E}^{\lambda } \}\subseteq \mathcal H\, ,\, {\mathcal E}^{\lambda } = P_{\lambda }\, \mathcal H\, $ ordered by inclusion such that each normal state of $\, R\, $ can be approximated in norm by convex combinations of vector states corresponding to vectors in $\, \cup_{\lambda } {\mathcal E}^{\lambda }\, $. Putting $\, Y^{\lambda } = {\gamma }^* \bigl( \bigl( Ad\, P_{\lambda } \bigr)_* \bigl( \mathcal B ( {\mathcal E}^{\lambda } )_* \bigr)\, $ and $\, Y = {\gamma }^* \bigl( \mathcal B ( \mathcal H )_* \bigr)\, $ the assignment $\, Y \rightarrow Y^{\lambda }\, ,\, \sigma \mapsto \sigma\circ Ad\, P_{\lambda }\, $ with 
$\, ( \sigma\circ Ad\, P_{\lambda } ) \bigl( y \bigr) = \sigma \bigl( P_{\lambda }\, a\, P_{\lambda } \bigr)\, $ is well defined (completely) positive linear and contractive for each $\, \sigma\in Y\, $ hence extends by injectivity of $\, Y^{\lambda }\, $ the latter which is completely order isomorphic to a sum of matrix duals to a (completely) positive linear map 
$\, R^* \rightarrow Y^{\lambda }\, $ denoted $\, Ad\, P_{\lambda }\, $ by abuse of notation which is completely bounded by $\, \dim\, P_{\lambda }\, $. To see this extend the completely positive embedding $\, Y \subseteq R^*\, $ to a completely positive embedding $\, Y^+ \subseteq ( R^* )^+\, $ of the corresponding unitizations as in the proof of Proposition 2 endowed with the corresponding matrix orders by which the unitizations are completely order isomorphic with operator systems, restricting to the given matrix orders on the subspaces 
$\, Y \subseteq R^*\, $. The unitization of the completely positive retraction $\, Ad\, P_{\lambda } : Y \twoheadrightarrow Y^{\lambda }\, $ composed with the natural completely positive retraction $\, ( Y^{\lambda } )^+ \twoheadrightarrow Y^{\lambda }\, $ (compare with the proof of Proposition 2) then extends by injectivity of $\, Y^{\lambda }\, $ to a completely positive retraction $\, Ad\, P_{\lambda } : R^* \subseteq ( R^* )^+ \twoheadrightarrow Y^{\lambda }\, $. Then $\, Y\, $ (resp. $\, R^*\, $) admits a decomposition as
$$ Y\> =\> Y^{\lambda }\> +\> \bigl( 1 - Ad\, P_{\lambda } \bigr) \bigl( Y \bigr)\> ,\qquad R^*\> =\> Y^{\lambda }\> + \bigl( 1 - Ad\, P_{\lambda } \bigr) \bigl( R^* \bigr)\> .  $$ 
Let $\, X = {\gamma }^* \bigl( co \bigl( {\mathcal P}^{\nu } \bigl( \lambda ( R )'' \bigr) \bigr) \bigr)\, $ be the image of the convex hull of the vector states corresponding to vectors in $\, \cup_{\lambda } {\mathcal E}^{\lambda }\, $  and $\, X^{\lambda } = X \cap Y^{\lambda }\, $ and note that the map $\, {\gamma }^* :  co \bigl( {\mathcal P}^{\nu } \bigl( \lambda ( R )'' \bigr) \bigr) \rightarrow  X \subseteq R^*\, $ is one-to-one. Given $\, x\in \bigl( \gamma ( R )'' \bigr)^{sa}\, $ one obtains an affine function $\, f_x : X \rightarrow \mathbb R\, $ by evaluation of $\, x\, $ on the corresponding preimages of $\, X\, $. Define a net of affine normcontinuous functions $\, \{ f^{\lambda }_x {\}}_{\lambda\in\Lambda }\, ,\, f^{\lambda }_x : \mathcal S ( R ) \rightarrow \mathbb R\, $ by restriction of the linear function $\, {\widehat f}^{\lambda }_x : R^* \rightarrow \mathbb R\, $ where $\, {\widehat f}^{\lambda }_x\, $ is defined by linear extension of 
$$ f^{\lambda }_x ( p )\> =\> f_x ( p )\> ,\quad p\in X^{\lambda }\> ,\qquad {\widehat f}^{\lambda }_x ( z )\> =\> 0\> ,\quad z\in \bigl( 1 - Ad\, P_{\lambda } \bigr) \bigl( R^* \bigr) \> . $$
Evaluation on the normal states of $\, R\, $ yields affine functions 
$\, \{ g^{\lambda }_x {\}}_{\lambda\in\Lambda }\, ,\, g^{\lambda }_x : {\mathcal S}^{\nu } ( R ) \rightarrow \mathbb R\, $, which being affine are automatically normcontinuous (resp. weakly continuous) thus defining an element $\, x_{\lambda }\in R^{sa}\, $ (compare with \cite{Pe1} , Theorem 3.10.3). 
Moreover the assignment 
$\, x\mapsto g^{\lambda }_x\, $ is positive which follows from the fact that $\, Ad\, P_{\lambda }\, $ is positive. Choosing any free ultrafilter $\, \omega\, $ for the net $\, \Lambda\, $, i.e. an ultrafilter $\,\omega\, $ on $\, \Lambda\, $ containing the subset $\, {\Lambda }_{\mu } = \{ \lambda\in\Lambda \,\vert\, \lambda \geq \mu \} \subseteq \Lambda\, $ for every $\, \mu\in\Lambda\, $.
Then the net $\, \{ g^{\lambda }_x \}\, $ converges pointwise along $\, \omega\, $ to a bounded affine function 
$$ \bigl( g^{\lambda }_x \bigr)\> \buildrel\omega\over\largerightarrow\> g^{\omega }_x $$
corresponding to a unique element $\, \widehat x\in R^{sa}\, $. Define $\, \rho ( x ) = \widehat x\, $.
From the argument above the map $\, \rho : \lambda ( R )'' \rightarrow R\, $ is linear and positive. Then we need to show that it is the identity map for $\, x\in R\, $ but this is fairly obvious from the fact that the net $\, \{ g^{\lambda }_x \}\, $ converges pointwise to the (canonical) function representation by evaluation for $\, x\in R^{sa}\, $\qed  
\par\bigskip\noindent
{\it Remark.}\quad If one could construct a positive map $\, I ( M ) \twoheadrightarrow {\sigma }^r ( M )\, $  for general monotone complete $C^*$-algebras then each monotone complete $C^*$-algebra $\, M\, $ would be weakly injective since the inclusion $\, M \subseteq {\mathfrak L}_1 ( M )\, $ is obviously normal whence $\, M = {\sigma }^r ( M )\, $. This however seems to be unlikely since in case of von Neumann algebras the argument given above makes essental use of the fact that $\, M\, $ admits sufficiently many normal functionals.
\par\bigskip\noindent
Let $\, M\, $ be a von Neumann algebra. A state $\, s\, $ of $\, {\mathfrak L}_q^{\nu } ( M )\, $ is called {\it boundedly basically increasing normal} iff given any monotone increasing net of uniformly bounded $w^*$-closed basic elements $\, \bigl( \underline{\mathcal A}_{\lambda } {\bigr)}_{\lambda } \nearrow \mathfrak A\, $ converging up to $\, \mathfrak A\in {\mathfrak L}_q^{\nu } ( M )\, $ the net $\, \bigl( s ( \underline{\mathcal A}_{\lambda } ) {\bigr)}_{\lambda }\, $ converges to $\, s ( \mathfrak A )\, $. By symmetry this implies that given any monotone decreasing net of uniformly bounded $w^*$-closed antibasic elements $\, \bigl( \overline{\mathcal B}_{\lambda } \bigr)_{\lambda } \searrow \mathfrak B\, $ converging down to $\, \mathfrak B\in {\mathfrak L}_q^{\nu } ( M )\, $ the monotone decreasing net $\, \bigl( s \bigl( \overline{\mathcal B}_{\lambda } \bigr) \bigr)_{\lambda }\, $ converges to $\, s \bigl( \mathfrak B \bigr)\, $. 
\par\bigskip\noindent 
{\bf Corollary.}\quad Let $\, M\, $ be a von Neumann algebra. Then $\, {\mathfrak L}^{\nu }_q ( M )\, $ admits a separating family of boundedly basically increasing normal states.
\par\bigskip\noindent
{\it Proof.}\quad  Let $\, M\, $ be given admitting a natural isometric function representation as bounded affine functions on the set of normal states $\, X = {\mathcal S}^{\nu } ( M )\, $ by evaluation. The natural inclusion $\, {\Lambda }^{\nu } : {\mathfrak L}^{\nu } ( M ) \hookrightarrow {\mathfrak L}^{\nu } ( N )\, $ with $\, N = B ( X )\, $ an abelian von Neumann algebra admits a positive linear extension 
$$ {\Lambda }^{\nu }_q :\> {\mathfrak L}_q^{\nu } ( M )\> \largerightarrow\> {\mathfrak L}_q^{\nu } ( N )\> . $$
Composing $\, {\Lambda }^{\nu }_q\, $ with the minimal (concave) monotonous extension 
$$ {\pi }^{\nu }_c :\> {\mathfrak L}_q^{\nu } ( N )\> \largerightarrow\> {\mathfrak L}_1 ( N )\> \simeq\> N $$
of the identity map of $\, N\, $ as in Proposition 1 yields a concave monotonous extension 
$\, s^{\nu }_c\, $ of the natural function representation $\, s^{\nu } :\> {\mathfrak L}^{\nu } ( M )\,\rightarrow\, B ( X )\, $ to $\, {\mathfrak L}_q^{\nu } ( M )\, $. We claim that $\, {\Lambda }^{\nu }_q\, $ is boundedly basically increasing normal so that since $\, {\pi }^{\nu }_c\, $ is boundedly basically increasing normal from Proposition 3 also $\, s^{\nu }_c\, $ is boundedly basically increasing normal. Note that $\, {\pi }^{\nu }_c\, $ is linear with respect to addition of elements in $\, {\mathfrak L}^{\nu } ( N )\, $ since 
$$ {\pi }^{\nu }_c \bigl( \mathfrak A - \mathcal B \bigr)\> \leq\> {\pi }^{\nu }_c \bigl( \mathfrak A \bigr)\> -\> {\pi }^{\nu } \bigl( \mathcal B \bigr)\> =\> {\pi }^{\nu }_c \bigl( \mathfrak A \bigr)\> +\> {\pi }^{\nu } \bigl( - \mathcal B \bigr)\> \leq\> {\pi }^{\nu }_c \bigl( \mathfrak A - \mathcal B \bigr)  $$
on using concavity of $\, {\pi }^{\nu }_c\, $ and linearity of $\, {\pi }^{\nu }\, $.
Let $\, Y = {\mathcal S}^{\nu } \bigl( B ( X ) \bigr)\, $ denote the space of normal states of the abelian von Neumann algebra $\, B ( X )\, $.
If $\, \bigl( \underline{\mathcal E}_{\lambda } {\bigr)}_{\lambda } \nearrow \mathfrak E\, $ is a monotone increasing net of uniformly bounded $w^*$-closed basic elements converging up to 
$\, \mathfrak E \in {\mathfrak L}^{\nu }_q \bigl( B ( X ) \bigr)\, $ then the pointwise supremum $\, \underline e\, $ of the image net $\, \bigl\{ s^{\nu } \bigl( \underline{\mathcal E}_{\lambda } \bigr) \bigr\} \subseteq B ( Y )\, $
is concave and Lipschitz normcontinuous being the supremum of a net of uniformly Lipschitz normcontinuous concave functions and since any Lipschitz normcontinuous concave function on a simplex $\, Y\, $ can be pointwise approximated by affine functions (corresponding to arbitrary bounded functions on the boundary $\, X = \delta Y\, $) one concludes that this supremum represents a $w^*$-closed basic element $\, \underline{\mathcal E}\in {\mathfrak L}^{\nu } \bigl( B ( X ) \bigr)\, $ which therefore must coincide with the supremum $\, \mathfrak E = \sup_{\lambda } \underline{\mathcal E}_{\lambda }\, $. To see this note that any point of $\, \eta\in Y\, $ can be approximated in norm by finite convex combinations of points in $\, X\, $. Choosing such a finite convex combination $\, \rho = \sum_{k = 1}^n\, {\lambda }_k\, x_k\, $ satisfying 
$\, \Vert \rho - \eta \Vert \leq \epsilon\, $ for given $\, \epsilon > 0\, $ one finds an element 
$\, w_{\eta }^{\rho }\in B ( X )\, $ with $\, \vert \rho ( w_{\eta }^{\rho } )\, -\, \eta ( w^{\rho }_{\eta } ) \vert \leq \epsilon\, $ and 
$\, \sigma ( w_{\eta }^{\rho } )\geq \underline e ( \sigma )\, ,\> \forall \sigma\in {\Delta }_{\rho } = co \bigl( x_1\, ,\cdots\, ,\, x_n \bigr)\, ,\, \rho  ( w_{\eta }^{\rho } ) = \underline e ( \rho )\, $. Moreover the net $\, ( w_{\eta }^{\rho } )\, $ can be chosen uniformly bounded since $\, \underline e\, $ is Lipschitz normcontinuous where $\, \rho \prec \sigma\, $ is defined by the relations $\, {\Delta }_{\rho } \subseteq {\Delta }_{\sigma }\, $ and 
$\, \epsilon ( \sigma ) \leq \epsilon ( \rho )\, $ for $\, \epsilon ( \rho ) = \Vert \rho - \eta\Vert\, $. Then any $w^*$-limit point $\, w_{\eta }\in B ( X )\, $ of the net $\, ( w_{\eta }^{\rho } )\, $ will have the required properties. Thus in the context of $\, {\mathfrak L}^{\nu }_q \bigl( B ( X ) \bigr)\, $ the supremum of any net of uniformly bounded $w^*$-closed basic elements is basic. This in turn can be used to show that $\, \underline{\Lambda }^{\nu }\, $ is boundedly basically increasing normal. Let 
$\, \bigl( \underline{\mathcal A}_{\lambda } {\bigr)}_{\lambda } \nearrow \mathfrak A\, $ be as above and let $\, \underline{\mathfrak A}\in {\mathfrak L}_q ( M )\, $ denote the supremum of the net 
$\, \bigl( \underline{\mathcal A}_{\lambda } {\bigr)}_{\lambda }\, $ viewed as basic elements of $\, \mathfrak L ( M )\, $ under the canonical inclusion $\, {\mathfrak L}^{\nu } ( M ) \hookrightarrow \mathfrak L ( M )\, $. Let $\, \underline{\mathcal A}^{\Lambda }\, $ denote the supremum of the image net $\, \bigl( {\Lambda }^{\nu } \bigl( \underline{\mathcal A}_{\lambda } \bigr) {\bigr)}_{\lambda }\, $ in $\, {\mathfrak L}_q^{\nu } ( N )\, $ which from the argument above is a $w^*$-closed basic element. Then $\, \underline{\mathcal A}^{\Lambda } \geq \Lambda \bigl( \underline{\mathcal A}_{\lambda } \bigr)\, $ for each $\, \lambda\, $ viewed as basic elements of $\, \mathfrak L ( N )\, $  since $\, \Lambda \bigl( \underline{\mathcal A}_{\lambda } \bigr)\, $ being the image of a bounded $w^*$-basic element is $w^*$-closed. Then
clearly $\, \underline{\mathcal A}^{\Lambda } = \overline p \bigl( \Lambda \bigl( \underline{\mathfrak A} \bigr) \bigr)\, $ where $\, \overline p\, $ is any monotonous extension of the canonical projection $\, p : \mathfrak L ( N ) \twoheadrightarrow {\mathfrak L}^{\nu } ( N )\, $ since $\, \Lambda\, $ is normal. Taking $\, \overline p\, $ to be the maximal monotonous extension of $\, p\, $, i.e. 
$$ \overline p \bigl( \mathfrak E \bigr)\> =\> \inf\, \bigl\{ p \bigl( \overline{\mathcal E} \bigr) \bigm\vert \overline{\mathcal E} \in \overline{\mathfrak L} ( N )\, ,\, \overline{\mathcal E} \geq \mathfrak E \bigr\} $$
one finds that $\, \underline{\mathcal A}^{\Lambda }\leq ( {\Lambda }^{\nu }_q\circ \overline p ) \bigl( \underline{\mathfrak A} \bigr) \leq( \overline p\circ \Lambda ) \bigl( \underline{\mathfrak A} \bigr) = \underline{\mathcal A}^{\Lambda }\, $ since 
$$ \bigl( \overline p\circ \Lambda \bigr) \bigl( \underline{\mathfrak A} \bigr)\> =\> \inf\, \bigl\{ ( p\circ \Lambda ) \bigl( \overline{\mathcal B} \bigr)\bigm\vert \overline{\mathcal B} \geq \underline{\mathfrak A} \bigr\} $$
so the result follows from 
$\, \overline p \bigl( \underline{\mathfrak A} \bigr) = \mathfrak A\, $.  One also  gets that 
$\, s_c^{\nu }\, $ which is linear with respect to addition of $w^*$-closed basic elements is also linear with respect to addition of $\, \mathfrak A\, $, i.e. 
$$ s_c^{\nu } \bigl( \mathfrak A + \mathfrak B \bigr)\> =\> {\pi }_c^{\nu } \bigl( \underline{\mathcal A}^{\Lambda }\> +\> {\Lambda }_q^{\nu } \bigl( \mathfrak B \bigr) \bigr)\> =\> s_c^{\nu } \bigl( \mathfrak A \bigr)\> +\> s_c^{\nu } \bigl( \mathfrak B \bigr) $$
for an arbitrary element $\, \mathfrak B\in {\mathfrak L}_q^{\nu } ( M )\, $. 
Evaluation at a point $\, \rho \in X\, $ yields a unital concave boundedly basically increasing normal map 
$$ m_{\rho } :\> {\mathfrak L}_q^{\nu } ( M )\> \largerightarrow\> \mathbb R $$
extending the natural positive linear evaluation map 
$$ s_{\rho } :\> {\mathfrak L}^{\nu } ( M )\> \largerightarrow\> \mathbb R $$
which is loosely referred to as a "Minkowski functional" in the sequel despite the fact that it is suplinear and positively homogenous instead of being sublinear. 
By linearity with respect to addition of suprema of nets of uniformly bounded $w^*$-closed basic elements $\, m_{\rho } = e_{\rho }\circ {\pi }_c\circ {\Lambda }^{\nu }_q\, $ is also increasing normal on each monotone increasing net of the form $\, \bigl( {\mathfrak C}_{\lambda } {)}_{\lambda } \nearrow \mathfrak C \, ,\, {\mathfrak C}_{\lambda } = \underline{\mathcal A}_{\lambda } + \mathfrak B\, ,\, \mathfrak C = \mathfrak A + \mathfrak B\, $ if $\, \bigl( \underline{\mathcal A}_{\lambda } \bigr)_{\lambda } \nearrow \mathfrak A\, $ is a monotone increasing net of uniformly bounded $w^*$-closed basic elements converging up to the element $\, \mathfrak A\, $. 
Given a Minkowski functional $\, m : V \rightarrow \mathbb R\, $ defined on a real ordered Banach space $\, V\, $ which is increasing normal restricted to a certain family $\, \mathcal F\, $ of monotone increasing nets in $\, V\, $ ($ m\, $ is said to be $\mathcal F$-normal in this case) 
assume that the family $\, \mathcal F\, $ remains stable under addition by fixed elements of $\, V\, $, i.e. $\, ( x_{\lambda } )_{\lambda }\in \mathcal F\, ,\, y\in V \Rightarrow ( x_{\lambda } + y )_{\lambda } \in \mathcal F\, $. Then 
there exists a family of $\mathcal F$-normal functionals 
$\, \{ {\rho }_{\kappa } \}\, $ with $\, {\rho }_{\kappa } \geq m\, ,\, \forall\> \kappa\, $ by the method as applied in the proof of part (ii) of Theorem 2 above. 
In our setting where $\, m = s^{\nu }_c ( \rho )\, $ one finds since $\, s_c^{\nu }\, $ is linear for addition of elements in $\, {\mathfrak L}^{\nu } ( M )\, $ that the process above applied to $\, m\, $ does not affect the images of elements in $\, {\mathfrak L}^{\nu } ( M )\, $ so that $\, m \bigl( \mathcal A \bigr) = s_{\rho } \bigl( \mathcal A \bigr) = {\rho }_{\kappa } \bigl( \mathcal A \bigr)\, $ for any boundedly basically increasing normal functional $\, {\rho }_{\kappa }\, $ obtained by applying the procedure as above to $\, m\, $ with respect to some given well order on the positive elements of norm $1$ in $\, {\mathfrak L}_q^{\nu } ( M )\, $. Therefore the resulting functionals $\, \{ {\rho }_{\kappa } \}\, $ are all positive linear extensions of the canonical map $\  s_{\rho } : {\mathfrak L}^{\nu } ( M ) \rightarrow \mathbb R\, $. Since the function representation $\, s^{\nu }\, $ is faithful on $\, {\mathfrak L}^{\nu } ( M )\, $ and choosing one boundedly basically increasing normal functional extension $\, \overline\rho\, $ for each $\, s_{\rho }\, $ one concludes that the set of functionals $\, \bigl\{ \overline\rho\bigm\vert \rho\in X \bigr\}\, $ is separating for $\, {\mathfrak L}_q^{\nu } ( M )\, $\qed
\par\bigskip\noindent
We also introduce the {\it (weak) normal topology} (as opposed to the strong normal topology considered in \cite{KaPe}) of a function system $\, \mathfrak X\, $ which is weaker than (but may coincide in certain cases with) the norm topology and coincides with the usual $w^*$-topology for a $W^*$-algebra $\, M\, $ with respect to its predual. Namely a bounded net $\, ( x_{\lambda } {)}_{\lambda }\, $ in $\, \mathfrak X\, $ is normally convergent to the element $\, x\in\mathfrak X\, $ iff 
$$ \limsup_{\lambda }\, x_{\lambda }\> =\> x\> =\> \liminf_{\lambda }\, x_{\lambda } $$
which in case that $\, \mathfrak X\, $ is not a complete function lattice is to be understood as a relation in $\, {\mathfrak L}_1 ( \mathfrak X )\, $. It is easy to see that a normconvergent net has this property as well as a $w^*$-convergent net in the abelian $W^*$-algebra $\, B ( X )\, $ of bounded functions on the set $\, X \, $. Then a map $\, \varphi : \mathfrak X \rightarrow \mathfrak Y\, $ of function systems is normally continuous iff given any normally convergent net $\, ( x_{\lambda } {)}_{\lambda }\to x\, $ in $\, \mathfrak X\, $ the image net $\, ( \varphi ( x_{\lambda } ) {)}_{\lambda }\, $ is normally convergent with limit $\, \varphi ( x )\, $. Any normally continuous map is normal in the usual sense (compare with the definition below) but a normal map need not be altogether normally continuous.
Given an arbitrary $W^*$-algebra $\, M\, $ there exists a $w^*$-continuous order isomorphic embedding $\, M \hookrightarrow B ( X )\, $ with $\, X = {\mathcal S}^{\nu } ( M )\, $ the set of normal states of $\, M\, $ which embedding extends to a positive linear order isomorphic embedding 
$\, {\iota }: {\mathfrak L}_1 ( M ) \hookrightarrow {\mathfrak L}^{\nu }_q \bigl( B ( X ) \bigr)\, $. Suppose given a normally convergent net $\, ( x_{\lambda } {)}_{\lambda } {\buildrel\nu\over\longrightarrow} x\, $ in $\, M\, $ and choose any $w^*$-convergent subnet $\, \bigl( y_{\mu } {\bigr)}_{\mu } {\buildrel w^*\over\longrightarrow} y\, $. Then also $\, \bigl( y_{\mu } {\bigr)}_{\mu } {\buildrel\nu\over\longrightarrow } x\, $ converges normally to $\, x\, $.
Consider the $w^*$-closed basic elements $\, {\underline C}_{\mu } = \inf_{\kappa\geq\mu }\, \bigl\{ y_{\kappa } \bigr\}\in {\mathfrak L}^{\nu } \bigl( M \bigr)\subseteq {\mathfrak L}^{\nu } \bigl( B ( X ) \bigr)\, $ and the $w^*$-closed antibasic elements $\, {\overline B}_{\mu } = \sup_{\kappa \geq\mu }\, \bigl\{ y_{\kappa } \bigr\}\in {\mathfrak L}^{\nu } ( M )\subseteq {\mathfrak L}^{^\nu } \bigl( B ( X ) \bigr)\, $ so that $\, \inf\, \{ y \} = \sup_{\mu }\, {\underline C}_{\mu } = \inf_{\mu }\, {\overline B}_{\mu }\in {\mathfrak L}^{\nu }_q ( M )\, $. Since $\, \iota \bigl( \sup_{\kappa \geq \mu }\, \bigl\{ y_{\kappa } \bigr\} \bigr) \geq {\overline B}_{\mu }\, $ and $\, \iota \bigl( \inf_{\kappa \geq \mu }\, \bigl\{ y_{\kappa }\bigr\} \bigr) \leq {\underline C}_{\mu }\, $ any positive linear retraction $\, \rho : {\mathfrak L}^{\nu } \bigl( B ( X ) \bigr) \twoheadrightarrow {\mathfrak L}_1 ( M )\, $ for $\, \iota\, $ must satisfy 
$\, \rho \bigl( {\overline B}_{\mu } \bigr) = \sup_{\kappa \geq \mu }\, \bigl\{ y_{\kappa } \bigr\}\, $ and $\, 
\rho \bigl( {\underline C}_{\mu } \bigr) = \inf_{\kappa\geq\mu }\, \bigl\{ y_{\kappa } \bigr\} \, $. Therefore 
$\, y = \rho \bigl( \{ y \} \bigr) = \rho \bigl( \inf_{\mu }\, {\overline B}_{\mu } \bigr) \leq x\, $ and 
$\, y = \rho \bigl( \{ y \} \bigr) = \rho \bigl( \sup_{\mu }\, {\underline C}_{\mu } \bigr) \geq x\, $ proving that $\, y = x\, $, i.e. the net $\, ( x_{\lambda } {)}_{\lambda }\, $ is $w^*$-convergent with limit $\, x\, $. 
Suppose conversely given a $w^*$-convergent net $\, ( x_{\lambda } {)}_{\lambda } {\buildrel w^*\over\longrightarrow } x\, $ in $\, M\, $ and consider the corresponding uniformly boundedly generated  $w^*$-closed basic and antibasic elements $\, {\underline C}_{\lambda } = \inf_{\mu\geq \lambda }\, \bigl\{ x_{\mu } \bigr\}\, ,\, {\overline B}_{\lambda } = \sup_{\mu\geq\lambda }\, \bigl\{ x_{\mu } \bigr\}\, $ in $\, {\mathfrak L}^{\nu } \bigl( B ( X ) \bigr)\, $. Then the minimal monotonous retraction 
$\, \underline\pi : {\mathfrak L}_q^{\nu } \bigl( B ( X ) \bigr) \twoheadrightarrow {\mathfrak L}_1 ( M )\, $ over the identity of $\, M\, $ is boundedly basically increasing normal and the maximal monotonous retraction $\, \overline\pi : {\mathfrak L}_q^{\nu } \bigl( B ( X ) \bigr) \twoheadrightarrow {\mathfrak L}_1 ( M )\, $ is boundedly antibasically decreasing normal, so that 
$\, x = \underline\pi \bigl( \{ x \} \bigr) = \underline\pi \bigl( \sup_{\lambda }\, {\underline C}_{\lambda } \bigr) = \sup_{\lambda }\, \underline\pi \bigl( {\underline C}_{\lambda } \bigr) \leq \liminf_{\lambda }\, x_{\lambda }\, $ whereas $\, x = \overline\pi \bigl( \inf_{\lambda }\, {\overline B}_{\lambda } \bigr) = \inf_{\lambda }\, \overline\pi \bigl( {\overline B}_{\lambda } \bigr)\geq \limsup_{\lambda } x_{\lambda }\, $ whence the relation 
$$ \limsup_{\lambda }\, x_{\lambda }\> =\> x\> =\> \liminf_{\lambda }\, x_{\lambda } $$
in $\, {\mathfrak L}_1 ( M )\, $ follows and $\, ( x_{\lambda } {)}_{\lambda } \, $ is normally convergent to $\, x\, $ proving that the $w^*$-topology of $\, M\, $ coincides with the normal topology (on bounded subsets).
\par\bigskip\bigskip\bigskip\bigskip\noindent
{\bf 3.\quad Projection lattices and $\mathcal P$-algebras.}
\par\bigskip\noindent
A $C^*$-algebra $\, A\, $ will be called a {\it $\Sigma $-algebra} iff $\, A\, $ is monotone complete and in addition the monotone closure of each unital abelian sub-$C^*$-algebra 
is contained in an abelian subalgebra (i.e. each maximal abelian subalgebra is relatively monotone complete).  In the meanwhile the author learned that the second condition is superfluous since automatically true in any monotone complete $C^*$-algebra which moreover has the general property that the monotone closure of any subalgebra $\, D\subseteq A\, $ is contained in its double commutant $\, D'' \subseteq A\, $, cf. \cite{KaPe}.  $\, A\, $ is a complete $\Sigma $-algebra iff $\, M_n ( A )\, $ is a $\Sigma $-algebra for every $\, n\in\mathbb N\, $. Every von Neumann algebra is a complete $\Sigma $-algebra, and also every injective $C^*$-algebra. It should be hard to find a monotone complete $C^*$-algebra $\, A\, $ without $\, M_n ( A )\, $ also being monotone complete, on the other hand the author knows of no argument ensuring that this implication is automatic. Natural candidates are the weak injective envelopes of section 2.
In \cite{Kap} Kaplansky coined the term $AW^*$-algebra for a $C^*$-algebra $\, A\, $ such that every maximal abelian subalgebra is monotone complete (so that $\, A\, $ contains "sufficiently many projections"). Since the prefix $W^*$ seems to indicate the existence of normal functionals at least on the level of maximal abelian subalgebras we will not use this nomenclatura here but instead talk of $A\Sigma $-algebras (i.e. $A\Sigma $-algebra $\, =\, AW^*$-algebra in the sense of Kaplansky). 
Of course an $A\Sigma $-algebra $\, A\, $ with a separating family of completely additive states (so that in particular any maximal abelian subalgebra is a $W^*$-algebra) is automatically a $W^*$-algebra itself, compare \cite{Pe1}, Theorem 3.9.4. We also occasionally consider sequentially monotone complete $C^*$-algebras called $ {\Sigma }_{\omega }$-algebras. 
A $C^*$-algebra $\, A\, $ is a {\it $\mathcal P$-algebra} if any positive element $\, x\in A_+\, $ can be approximated in norm by finite positive linear combinations of pairwise commuting positive projections {\it and if} for any projection $\, p\in \mathcal P ( A )\, $ the algebra $\, p\, A\, p\, $ is again a $\mathcal P$-algebra. $\, A\, $ is called a {\it ${\mathcal P}_{\sigma }$-algebra} if every positive element $\, x\in A_+\, $ can be approximated in norm from below by finite positive linear combinations of  ("almost spectral")  pairwise commuting projections also commuting with $\, x\, $. We adopt the general terminology to call a projection contained in the double commutant $\, \bigl\{ x \bigr\}''\, $ of a selfadjoint element $\, x\, $ a {\it spectral projection of $\, x\, $}, and a projection contained in the commutant $\, \bigl\{ x \bigr\}'\, $ an {\it almost spectral projection of $\, x\, $}. Of course the former notion depends on the enveloping algebra whereas the latter does not. Accordingly we define a {\it spectral $\mathcal P$-algebra} or ${\mathcal P}_s$-algebra to be a $\mathcal P$-algebra such that any positive element can be approximated in norm from below by finite positive combinations of spectral projections of $\, x\, $. Any ${\Sigma }_{\omega }$-algebra is a spectral $\mathcal P$-algebra.
It is shown below that any separable $C^*$-algebra admits an embedding into a separable $\mathcal P$-algebra but not necessarily into a separable ${\mathcal P}_{\sigma }$-algebra which is why we have settled for the more general definition. One should bear in mind however that many $\mathcal P$-algebras are in fact ${\mathcal P}_{\sigma }$-algebras, examples include any $*$-ideal or $C^*$-quotient of a ${\Sigma }_{\omega }$-algebra. $\, A\, $ is called  a {\it $\mathcal P$-lattice algebra} if it is a $\mathcal P$-algebra and its subset of positive projections is a lattice for the induced order. A ${\mathcal P}^+$-algebra is a $\mathcal P$-algebra $\, A\, $ such that any element $\, x\in A\, $ admits a polar decomposition 
$\, x = v\, \vert x\vert = \vert x^*\vert\, v\, $ with $\, v\in A\, $ a partial isometry.  An ${\Sigma }^+$-algebra is a $\Sigma $-algebra which is a ${\mathcal P}^+$-algebra. A monotone increasing net of projections $\, \bigl\{ p_{\lambda } {\bigr\}}_{\lambda } \in \mathcal P ( A )\, $ of a $C^*$-algebra $\, A\, $ is called a {\it projection unit} iff given $\, \epsilon > 0\, $ there exists for every $\, x\in A_+\, $ an element $\, p_{\lambda }\, $ with $\, \Vert p_{\lambda }\, x - x \Vert \leq \epsilon\, $. Obviously any $\mathcal P$-algebra admits a projection unit. If $\, \mathcal J \vartriangleleft A\, $ is a $C^*$-ideal and $\, B\subseteq A\, $ is any given subalgebra then the projection unit $\, \bigl\{ p_{\lambda } \bigr\} \subseteq \mathcal P ( \mathcal J )\, $ will be called {\it quasicentral for $\, B\, $} iff 
$\, \lim_{\lambda\to\infty }\, \bigm\Vert p_{\lambda }\, x\, -\, x\, p_{\lambda } \bigm\Vert = 0\, $ for any $\, x\in B\, $. If $\, A\, $ is a $\mathcal P$-algebra define $\, {\mathcal P}_f ( A ) \subseteq A\, $ to be the dense subset consisting of finite linear combinations of pairwise commuting projections.
A (possibly discontinuous) positively homogenous map $\, \sigma : {\mathcal P}_f ( A ) \rightarrow {\mathcal P}_f ( B )\, $ of unital $\mathcal P$-algebras which is monotonous (and hence continuous) on the intersection of each abelian subalgebra $\, C\subseteq A\, $ with $\, {\mathcal P}_f ( A )\, $ will be called a $\mathcal P$-map iff it maps positive projections to positive projections and if $\, \sigma ( \alpha {\bf 1} + x ) = \alpha \sigma ( {\bf 1} ) + \sigma ( x )\, $ for each element $\, x\in {\mathcal P}_f ( A )\, $ and $\, \alpha \in \mathbb C\, $. We will mostly consider $\mathcal P$-maps which send orthogonal projections $\, p\bot q\, $ to orthogonal projections. However to have the opportunity of considering more general $\mathcal P$-maps we call such maps 
{\it orthogonal $\mathcal P$-maps} or {\it ${\mathcal P}^{o}$-maps}.
A $\mathcal P$-map which maps each pair of complementary projections to a pair of complementary projections will be called {\it complemented} or a {\it ${\mathcal P}^c$-map}. Any 
${\mathcal P}^c$-map is orthogonal, but the converse need not be true. The dual notion is a {\it coorthogonal $\mathcal P$-map} or {\it ${\mathcal P}^{co}$-map} which sends each pair of coorthogonal projections $\, p\top q\, $, i.e. $\, [ p\, ,\, q ] = 0\, ,\, p\vee q = {\bf 1}\, $ to a pair of coorthogonal projections. A complemented $\mathcal P$-map is both orthogonal and coorthogonal. A
$\mathcal P$-map 
which maps each pair of commuting projections to a pair of commuting projections will be called a 
${\mathcal P}^{a}$-map. A ${\mathcal P}^{a}$-map 
such that $\, \sigma ( p \wedge q ) = \sigma ( p ) \wedge \sigma ( q )\, $ for each pair of commuting projections $\, [ p\, ,\, q ] = 0\, $ will be called a {\it  ${\mathcal P}$-$A$-wedge-map or ${\mathcal P}^{a\wedge}$-map}, and a {\it  
$\mathcal P$-$A$-lattice-map} if in addition it sends complementary projections to complementary projections. Correspondingly one may define 
{\it ${\mathcal P}^{a\vee}$-maps, and ${\mathcal P}^{\wedge }$-maps resp. ${\mathcal P}^{\vee }$-maps} of $\mathcal P$-lattice algebras which instead of only considering the operations $\,\wedge , \vee\, $ for two commuting projections (which are defined in any $C^*$-algebra) respect these operations for general pairs of projections. 
The following compiles these notions into a general concept. 
A {\it $*$-decorated $\mathcal P$-map or ${\mathcal P}^*$-map} of $\mathcal P$-lattice algebras is a $\mathcal P$-map with one or more of the additional properties $\, * \in \{ o , co , c , a , a\wedge , a\vee , \wedge , \vee \}\, $. 
A decorated $\mathcal P$-map respecting all three lattice operations (and hence has all other decorations as well) will be called a {\it $\mathcal P$-lattice map}, i.e. a $\mathcal P$-lattice map satisfies 
$$ \sigma ( p\vee q )\> =\> \sigma ( p ) \vee \sigma ( q )\> ,\quad \sigma ( p \wedge q )\> =\> \sigma ( p ) \wedge \sigma ( q )\> ,\quad \sigma ( {\bf 1} - p )\> =\> {\bf 1}\> -\> \sigma ( p )  $$ 
for any pair of positive projections $\, p\, ,\, q\in A\, $. One would like to be dealing only with $\mathcal P$-lattice maps since these have the strongest properties but then these maps occur very sparsely in nature, not even $*$-homomorphisms of $\mathcal P$-lattice algebras are guaranteed to be $\mathcal P$-lattice maps. On the other hand 
$\mathcal P$-A-lattice maps still have extremely strong implications (the restriction to any abelian $\mathcal P$-subalgebra is a $*$-homomorphism), and include the notion of $*$-homomorphisms. Note however that even a $\mathcal P$-lattice map need not be overall monotonous or continuous. Simple counterexamples can be constructed by 
discontinuous bijective ${\mathcal P}^c$-maps $\, \sigma : M_2 ( \mathbb C ) \rightarrow M_2 ( \mathbb C )\, $ on realizing that any such map is a $\mathcal P$-lattice map.
\par\noindent
An (abstract) lattice $\, \Lambda\, $ is said to be of {\it finite type} if each totally ordered subset in a finitely generated sublattice is finite, and a {\it function lattice} iff it satisfies the identities 
$\, ( x \wedge y ) \vee z = ( x \vee z )\wedge ( y \vee z )\, $ and $\, ( x \vee y ) \wedge z = ( x \wedge z ) \vee ( y \wedge z )\, $ for any three elements $\, x\, ,\, y\, ,\, z\in\Lambda\, $. A function lattice is clearly of finite type, another example is the lattice of projections of any finitedimensional $C^*$-algebra, which in fact is of {\it finite depth}, i.e. the size of each totally ordered subset is uniformly bounded by some constant $\, N\, $.  A complete lattice will be called a {\it complete function lattice} if it is a function lattice and in addition 
for each monotone increasing net $\, ( x_{\lambda } ) \nearrow x\, $ converging up to an element $\, x\, $ and any $\, y\, $ one has the identities $\, x \wedge y = \sup_{\lambda }\, ( x_{\lambda } \wedge y )\, $, resp. for any monotone decreasing net $\, ( w_{\mu } ) \searrow w\, $ and any $\, z\, $ one has $\, w\vee z = \inf_{\mu }\, ( w_{\mu }\vee z )\, $. Also recall from \cite{Kap} that a {\it modular lattice} is a lattice satisfying $\, ( e\vee f ) \wedge g = e \vee ( f \wedge g )\, $ whenever $\, e \leq g\, $.
A {\it complemented lattice} is a lattice $\, \Lambda\, $ containing a unique maximal element $\, {\bf 1}\, $ and a unique minimal element $\, 0\, $ which is equipped with an involutary map 
$$ x\>\mapsto\> x^c\> =\> {\bf 1}\> -\> x\> ,\quad ( x^c )^c = x  $$
such that 
$$ x \leq y \iff y^c\leq x^c\> ,\quad ( x \wedge y )^c\> =\> ( x^c\vee y^c )^c  \> . $$
Then $\, x\, $ and $\, y\, $ are {\it orthogonal} iff $\, x\leq y^c\, $. For orthogonal elements 
$\, x\bot y\, $ we also write the lattice operations in additive/ subtractive form 
$\, x \vee y = x + y\, $ and  $\,  y^c\wedge x^c = y^c - x = x^c - y\, $ respectively. Then one assumes the following restricted commutation property: 
$\, ( x + y ) \wedge ( x + z )\, =\, x + ( y \wedge z )\, $. 
A complemented lattice $\, \Lambda\, $ will be called {\it signed} iff given a disjoint decomposition 
$\, \Lambda = {\Lambda }_- \cup {\Lambda }_+\, $ (a signature) such that $\, e\in {\Lambda }_- \iff e^c\in {\Lambda }_+\, $ and $\, e \leq f\, ,\, f\in {\Lambda }_-\, $ implies $\, e\in {\Lambda }_-\, $. A signature will be called {\it polar} iff the following condition holds: given a pair of orthogonal elements $\, e\,\bot\, f\, $ with $\, e + f\in {\Lambda }_+\, $ then either $\, e\in {\Lambda }_+\, $ or else $\, f\in {\Lambda }_+\, $. Correspondingly a lattice $\,\Lambda\, $ is called {\it polar} iff for each pair of nonorthogonal elements $\, c\, ,\, d\in \Lambda\, ,\, c\nleq d^c\, $ there exists a polar decomposition
$\, \Lambda = {\Lambda }_- \cup {\Lambda }_+\, $ with $\, c\, ,\, d\in {\Lambda }_+\, $. 
Any (unital complemented) sublattice $\, {\Lambda }_0 \subseteq \Lambda\, $ of a polar lattice is again polar (by restricting any suitable polar signature from $\, \Lambda\, $ to $\, {\Lambda }_0\, $). 
In fact the notion only requires knowledge of the $A$-lattice structure of $\, \Lambda\, $ (complements and orthogonal sums of elements). Any projection lattice of a commutative $C^*$-algebra is polar since a polar decomposition corresponds to the restriction of a multiplicative functional in this case and a commutative $C^*$-algebra has a separating family of such functionals. Another case of a polar projection lattice is $\, \mathcal P \bigl( M_2 ( \mathbb C ) \bigr)\, $. In fact any signature is polar in this case. Apart from these examples no matrix algebra of higher (or infinite) dimension admits a polar signature. If such a signature would exist on $\, \mathcal P \bigl( M_3 ( \mathbb R ) \bigr)\, $ say, one finds that the corresponding region of the real projective plane consisting of minimal projections of positive signature must be simply connected and convex. Moreover for each point of the boundary the maximal length of the intersection of some geodesic containing the boundary point with this region must be the same and equal to $\, \pi / 2\, $.The 
only solutions for this setting are a disc of diameter $\, \pi / 2\, $, the equilateral geodesic triangle of side length $\, \pi / 2\, $ and more generally any uneven regular geodesic polygon of diameter $\, \pi / 2\, $. In each case there remain orthogonal systems lying completely outside of the region in question. Therefore a polar signature cannot exist for $\, \mathcal P \bigl( M_k ( \mathbb C ) \bigr)\, $ for any $\, k \geq 3\, $. Thus considering full projection lattices of $C^*$-algebras the notion of polar signature is confined more or less to commutative algebras. The situation might be quite different however if one considers (countable) sublattices of projections (for a countable sublattice the connectedness property makes no sense).
A complemented monotonous map of signed lattices $\, r : \Lambda \rightarrow {\Lambda }'\, $  will be called {\it signed} iff $\, r ( {\Lambda }_{\pm } ) \subseteq {\Lambda }_{\pm }'\, $. Specifying a signature on a lattice $\, \Lambda\, $ is equivalent to specifying a complemented monotonous map 
$\, \tau : \Lambda \twoheadrightarrow \{ 0 , {\bf 1} \}\, $ into the trivial lattice by letting 
$\, {\Lambda }_- = {\tau }^{-1} ( 0 )\, ,\, {\Lambda }_+ = {\tau }^{-1} ( {\bf 1} )\, $. Correspondingly, specifying a polar signature on $\, \Lambda\, $ is equivalent to giving an $A$-lattice map 
$\, \tau : \Lambda \twoheadrightarrow \{ 0 , {\bf 1} \}\, $ into the trivial lattice. A {\it presignature} is a disjoint decomposition $\, \Lambda = {\Lambda }_- \cup {\Lambda }_+\, $ such that $\, e \leq f\, ,\, e\in {\Lambda }_+\, $ implies $\, e^c\in {\Lambda }_-\, ,\, f\in {\Lambda }_+\, $. A presignature is equivalent to giving an orthogonal monotonous map into the trivial lattice $\, \{ 0 , {\bf 1} \}\, $.
An {\it $A$-lattice} is a partially ordered set $\, \Gamma\, $ containing a unique maximal element 
$\, {\bf 1}\, $ and a unique minimal element $\, 0\, $ which is equipped with an order reversing involution $\, e \mapsto e^c = {\bf 1} - e\, ,\, ( e^c )^c = e\, $ and an addition $\, ( e , f ) \mapsto e + f\, $ of orthogonal elements $\, e\,\bot\, f \iff e \leq f^c\, $ such that $\, e + e^c = {\bf 1}\, $ and $\, e + f\in \Gamma\, $ is the smallest element  larger than both $\, e\, $ and $\, f\, $, then there is also defined a subtraction operation $\, ( e , f ) \mapsto f - e\, $ for ordered pairs $\, e\leq f\, $ with $\, ( f - e ) + e = f\, $. 
The subset of projections of any $C^*$-algebra is a natural example of an $A$-lattice. 
A faithful $A$-lattice representation is an order isomorphic map $\, \iota : \Gamma \hookrightarrow \Lambda\, $ of an $A$-lattice $\, \Gamma\, $ into a complemented lattice $\, \Lambda\, $ respecting complements and orthogonal sums such that 
$\, \iota ( c )\,\bot\, \iota ( d ) \Rightarrow c\,\bot\, d\, $. An $A$-sublattice $\, \Gamma \subseteq \Lambda\, $ of a complemented lattice will be called {\it proper} iff $\, e\wedge f\in \Gamma\, $ for any pair of commuting elements 
$\, e\, ,\, f\in \Gamma\, ,\, [ e\, ,\, f ] = 0\, $. 
A $\mathcal P$-lattice algebra $\, A\, $ will be called of {\it finite type} iff its subset of projections 
$\, \mathcal P ( A )\, $ is a (complemented) lattice of finite type and a {\it modular $\mathcal P$-lattice algebra} iff its subset of projections is a modular lattice. 
\par\noindent
The bulk of results given below also works if one considers $A\Sigma $-algebras and their $C^*$-quotients instead of $\Sigma $-algebras. The result that the subset of projections in an $A\Sigma$-algebra forms a complete lattice is well known, however the proof given here is a bit different from the one given in \cite{Kap} and has as consequence that 
for any normal inclusion of $\Sigma $-algebras (${\Sigma }_{\omega }$-algebras) one gets a corresponding normal (resp. sequentially normal) inclusion of projection lattices, which also easily follows from \cite{KaPe} but does not seem to have been made explicit in any paper known to the author. It also seems that 
the $A\Sigma $-algebras encountered in nature are mostly monotone complete, so that $\Sigma $-algebras constitute a natural and interesting category. In fact most of the results with few exceptions even extend to the category of $A\Sigma_{ \omega }$-algebras which are defined to be 
${\mathcal P}_s$-algebras characterized by the condition that every monotone increasing sequence $\, x_1\leq\cdots\leq x_n\leq x_{n + 1} \leq\cdots\, $ of mutually commuting selfadjoint elements of $\, A\, $ has a least upper bound {\it within the double commutant $\, \{ x_n \}''\subseteq A\, $} of the sequence (note that it is not required that the increasing sequence has a least upper bound in $\, A\, $). This condition implies that the subset of projections forms a sequentially monotone complete lattice by an easy adaptation of the proof of Theorem P below and that every maximal abelian subalgebra is (absolutely) sequentially monotone complete on replacing the terms monotone complete and normal by sequentially monotone complete and sequentially normal at the corresponding places. 
\par\bigskip\noindent
{\bf Theorem P.}\quad If $\, A \subseteq B\, $ is a separable subalgebra of the ${\Sigma }_{\omega }$-algebra $\, B\, $ then there exists a separable $\mathcal P$-subalgebra $\, C\subseteq B\, $ containing $\, A\, $. If $\, B\subseteq A\, $ is an abelian $\mathcal P$-subalgebra of the $\mathcal P$-algebra $\, A\, $ and $\, \mathcal J\vartriangleleft A\, $ is a closed $*$-ideal then $\, \mathcal J\, $ contains a projection unit quasicentral for $\, B\, $. 
Any $C^*$-quotient and any $C^*$-ideal of a $\mathcal P$-algebra (resp. ${\mathcal P}_{\sigma }$-algebra)  $\, A\, $ is a $\mathcal P$-algebra (resp. ${\mathcal P}_{\sigma }$-algebra). Any $C^*$-ideal of a spectral $\mathcal P$-algebra is a spectral $\mathcal P$-algebra.
If $\, A\, $ is a $\Sigma  $-algebra (resp. ${\Sigma }_{\omega }$-algebra) its subset of positive projections $\, \mathcal P ( A )\, $ forms a complete (resp. sequentially monotone complete) lattice, in particular $\, A\, $ is a $\mathcal P$-lattice algebra. In case that $\, q : A \twoheadrightarrow B\, $ is a normal $*$-epimorphism $\, B\, $ is a $\Sigma $-algebra and $\, q\, $ admits a natural normal (nonunital) $*$-homomorphic cross section, there exists a central projection $\, p_q\in \mathcal Z ( A )\, $ with 
$\, \ker\, q = ( {\bf 1} - p_q )\, A\, $ and $\, B \simeq p_q\, A\, $. If $\, q\, $ is sequentially normal whence $\, B\, $ is sequentially monotone complete, then $\, B\, $ is a $\mathcal P$-lattice algebra and $\, q\, $ is a $\mathcal P$-lattice map. 
Any $\mathcal P$-map $\, \sigma : {\mathcal P}_f ( A ) \rightarrow {\mathcal P}_f ( B )\, $ of a unital $\mathcal P$-algebra $\, A\, $ into a $\Sigma $-algebra $\, B\, $ extends to a positively homogenous map $\, \sigma : A \rightarrow B\, $ with the properties that $\, \sigma\, $ is monotonous and continuous  restricted to each abelian subalgebra and satisfies $\, \sigma \bigl( \alpha {\bf 1} + x \bigr) = \alpha\, \sigma \bigl( {\bf 1} \bigr) + \sigma \bigl( x \bigr)\, $ for any $\, \alpha\in \mathbb R\, $. Moreover if $\, A\, $ is a ${\mathcal P}_{\sigma }$-algebra then the extension is unique.
If $\, B\, $ is any $C^*$-quotient of the 
$\Sigma $-algebra $\, A\, $  such that $\, B\, $ is a $\mathcal P$-lattice algebra with quotient map $\, q\, $ a $\mathcal P$-lattice map (for example if $\, A\, $ is abelian), and $\, \mathcal P\subseteq \mathcal P ( B )\, $ is a complemented countable sublattice there exists a ${\mathcal P}^c$-map $\, \sigma : {\mathcal P}_f ( B ) \rightarrow A\, $ which is a  cross section to the quotient map restricted to $\, {\mathcal P}_f = \bigl\{ x = \sum_{k = 1}^n\, {\alpha }_k\, p_k \in {\mathcal P}_f ( B ) \bigm\vert \{ p_k {\}}_k \subseteq \mathcal P\, ,\, \bigl[ p_k\, ,\, p_l \bigr] = 0 \bigr\}\, $. If $\, \mathcal P \subseteq \mathcal P ( B_0 )\, $ is a normdense sublattice for some separable abelian $\mathcal P$-subalgebra $\, B_0 \subseteq B\, $ then $\, \sigma\, $ extends to $\, B_0\, $ and can be chosen a $\mathcal P$-$A$-lattice map (= $\mathcal P$-lattice map) restricted to $\, \mathcal P ( B_0 )\, $ extending to a $*$-homomorphism $\, \sigma : B_0 \rightarrow A\, $ by continuity. 
In particular if $\, B = q ( A )\, $ is any $C^*$-quotient of a ${\Sigma }^+ $-algebra $\, A\, $ such that  $\, B\, $ admits an order isomorphic representation as bounded operators on a separable Hilbert space
then $\, B\, $ is sequentially monotone complete and any unitary in $\, B\, $ lifts to a unitary in $\, A\, $, and in case that any maximal abelian subalgebra of $\, B\, $ is a $\mathcal P$-subalgebra every normal element of $\, B\, $ lifts to a normal element of $\, A\, $. Alternatively, if $\, A\, $ is of finite type $\,\sigma\, $ can be chosen to be the restriction of a ${\mathcal P}^{\wedge , o}$-map on $\, {\mathcal P}_f\, $ and a ${\mathcal P}^{o}$-map on $\, {\mathcal P}_f ( B )\, $, or (by duality) a ${\mathcal P}^{\vee , co}$-map restricted to $\, {\mathcal P}_f\, $, and a ${\mathcal P}^{co}$-map on $\, {\mathcal P}_f ( B )\, $. In case that $\, q\, $ is sequentially normal $\sigma \, $ can be chosen a (nonunital) continuous
$\mathcal P$-lattice cross section restricted to $\, {\mathcal P}_f\, $ and a ${\mathcal P}^c$-map on $\, {\mathcal P}_f ( B )\, $ (with respect to the unit $\, s ( {\bf 1} )\in \mathcal P ( A )\, $). 
A $\mathcal P$-lattice map (or even a $\mathcal P$-$A$-lattice map) if restricted to any abelian $\mathcal P$-subalgebra extends a $*$-homomorphism on this subalgebra (on any abelian subalgebra if every maximal abelian subalgebra is a $\mathcal P$-algebra).
Any $\mathcal P$-map 
$\,\sigma : {\mathcal P}_f ( A ) \rightarrow {\mathcal P}_f ( B )\, $ of $\mathcal P$-algebras is uniquely determined by its (monotonous) restriction $\, s : \mathcal P ( A ) \rightarrow \mathcal P ( B )\, $ to the subset of positive projections. Further it is monotonous (and hence continuous) restricted to the intersection with any abelian $\mathcal P$-subalgebra and satisfies the identity $\, \sigma ( x^2 ) = \sigma ( x )^2\, $ for any $\, x\geq 0\, $. Thus any convex combination of $\mathcal P$-maps satisfies the Schwarz inequality for positive elements in $\, {\mathcal P}_f ( A )\, $ by operator convexity of the squaring operation. 
Any monotone complete complemented lattice $\,\Lambda\, $ has the following extension property with respect to ${\mathcal P}^c$-maps (monotonous maps commuting with taking complements): given an inclusion of $A$-lattices 
$\, \Gamma \subseteq {\Gamma }' \, $ and a
${\mathcal P}^c$-map $\, s : \Gamma \rightarrow \Lambda\, $ there exists a ${\mathcal P}^c$-extension $\, s' : {\Gamma }' \rightarrow \Lambda\, $ with $\, s' {\vert }_{\Gamma } = s\, $.
Any injective $C^*$-algebra and any von Neumann algebra is a complete ${\Sigma }^+$-algebra. In particular any normal quotient of an injective $C^*$-algebra is injective.
\par\bigskip\noindent
{\it Proof.}\quad It immediately follows from the definition of a $\mathcal P$-algebra (resp. ${\mathcal P}_{\sigma }$-algebra)  that this property drops to quotients and ideals. Similarly it is easy to see that the property of spectral $\mathcal P$-algebra drops to any $*$-ideal but not necessarily to arbitrary quotients of such an algebra.  Let $\, A\subseteq B\, $ be a separable $C^*$-subalgebra of the ${\Sigma }_{\omega }$-algebra $\, B\, $ and $\, \bigl\{ p_n {\bigr\}}_n \subseteq \mathcal P ( A )\, $ be a normdense sequence of projections. For each $\, n\, $ choose a normdense subsequence $\, \bigl\{ x_{n , k} {\bigr\}}_k\subseteq \bigl\{ x\in ( p_n\, A\, p_n )_+ \bigm\vert \Vert x \Vert = 1 \bigr\}\, $. For each $\, ( n , k )\in {\mathbb N}^2\, $ choose a countable subset $\, {\mathcal P}_{n , k} \subseteq \mathcal P ( B )\, $ of projections in $\, B\, $ each smaller than $\, p_n\, $ such that $\, x_{n , k}\, $ can be approximated in norm by finite positive linear combinations of pairwise commuting projections in $\, {\mathcal P}_{n , k}\, $. Then consider the separable subalgebra $\, A \subseteq A^1 \subseteq B\, $ generated by $\, A\, $ and $\, \bigcup_{n , k} {\mathcal P}_{n , k}\, $. Repeating this process infinitely many times leads to an ascending sequence of separable subalgebras $\, \bigl\{ A^n {\bigr\}}_n\, $ of $\, B\, $ containing $\, A\, $. Put $\, C = \overline{\bigcup_n A^n}\, $. We claim that $\, C\, $ is a separable $\mathcal P$-algebra containing $\, A\, $.  Obviously any element $\, x\in C_+\, $ can be approximated in norm by positive linear combinations of pairwise commuting projections in $\, C\, $. Let $\, p\in \mathcal P ( C )\, $ be any projection and $\, x = p\, x\, p \geq 0\, $ be a positive element of $\, p\, C\, p\, $.  Then from continuous functional calculus there exists for every $\, \epsilon > 0\, $ an index $\, m\, $ and a projection $\, p^m\in A^m\, $ with $\, \Vert p^m - p \Vert \leq \epsilon\, $ together with a finite positive linear combination of pairwise commuting projections in $\, A^{m + 1}\, $ each smaller than $\, p^m\, $ approximating $\, p^m\, x\, p^m\, $ in norm up to the order of $\, \epsilon\, $, and again from continuous functional calculus a finite subset of pairwise commuting projections smaller than $\, p\, $ each of which is $\sqrt{\epsilon } $-close in norm to one of the former projections hence approximating $\, x\, $ in norm up to the order of $\, 2\, \sqrt{\epsilon}\, $ proving that $\, C\, $ is a $\mathcal P$-algebra.  
Now let $\, B \subseteq A\, $ be an abelian $\mathcal P$-subalgebra of the $\mathcal P$-algebra $\, A\, $ and $\, \mathcal J \vartriangleleft  A\, $ a $*$-ideal. Let $\, \bigl\{ u_{\lambda } {\bigr\}}_{\lambda\in \Lambda }\subseteq {\mathcal J }_+\, $ be a quasicentral approximate unit for $\, \mathcal J\, $. For each finite subset $\, S =  \bigl\{ x_1\, ,\,\cdots\, ,\, x_n \bigr\} \cup \bigl\{ Q_1\, ,\, \cdots\, ,\, Q_m\bigm\vert i \neq j \Rightarrow Q_i \bot Q_j\, ,\, \sum_{i = 1}^m Q_i = {\bf 1}  \bigr\} \subseteq \bigl\{ x\in {\mathcal J}_+ \bigm\vert \Vert x\Vert = 1 \bigr\} \cup \mathcal P ( B )\, $
and $\, \epsilon > 0\, $ choose $\, u_{{\lambda }_{S , \epsilon }}\, $ with $\, \Vert u_{{\lambda }_{S , \epsilon }}\, x_k - x_k \Vert \leq\epsilon \, $ and 
$\, \Vert u_{{\lambda }_{S , \epsilon }}\, Q_i\, -\, Q_i\, u_{{\lambda }_{S , \epsilon }} \Vert \leq {\epsilon\over m}\, $. Then replace $\, u_{{\lambda }_{S , \epsilon }}\, $ by a contractive element $\, v_{S , \epsilon }\, $ of $\, \mathcal P \bigl( \mathcal J \bigr)_f\, $ commuting with the projections $\, \bigl\{ Q_i {\bigr\}}_{i = 1}^m\, $ and approximating the element $\, \sum_{i = 1}^m\, Q_i\, u_{{\lambda }_{S , \epsilon }}\, Q_i\, $ in norm up to the order of $\, \epsilon\, $ and check that $\, \bigl\{ v_{S , \epsilon } {\bigr\}}_{( S , \epsilon )}\, $ is again a quasicentral approximate unit for $\, \mathcal J\, $. Then if $\, v_{S , \epsilon } = \sum_{k = 1}^l\, {\alpha }_k\, p^{S , \epsilon }_k\, $ for some monotone decreasing string of projections $\, \bigl\{ p^{S , \epsilon }_k \bigm\vert p^{S , \epsilon }_k \geq p^{S , \epsilon }_{k + 1} \bigr\} \subseteq \mathcal P ( \mathcal J )\, $ in $\, \mathcal J\, $ consider the projection unit $\, \bigl\{ p_1^{S , \epsilon } {\bigr\}}_{( S , \epsilon )}\, $ and check that it is quasicentral for $\, B\, $.
\par\noindent
Assume next that $\, A\, $ is a $\Sigma $-algebra (resp. ${\Sigma }_{\omega }$-algebra). We first show that the subset of its positive projections forms a complete (resp. sequentially monotone complete)  lattice. Given two projections 
$\, p , q\in A\, $ consider the two monotone decreasing sequences of alternating products 
$$ ( x_n )_n \searrow x\> ,\quad x_n\> =\> \underbrace{( p\, q )\, ( p\ q)\,\cdots\, ( p\, q )}_{n- times}\, p\> ,\quad ( y_n )_n \searrow y\> ,\quad y_n\> =\> \underbrace{( q\, p )\, ( q\, p )\,\cdots\, ( q\, p )}_{n-times}\, q \> . $$
Clearly, $\, \sqrt x\leq \inf_n\, \sqrt{ x_{2n}}\, =\, \inf_n\, x_n\, =\, x\, $ by monotonicity of the squareroot whence $\, x\, $ (and similarly $\, y\, $) is a projection. Then one checks the relations 
$$ x\, y\, x\> =\> \inf_n\, ( x\, y_n\, x )\> =\> \inf_n\, ( x\, x_{n+1}\, x )\> =\> x\> ,\quad 
y\, x\, y\> =\> \inf_n\, ( y\, x_n\, y )\> =\> \inf_n\, ( y\, y_{n+1}\, y )\> =\> y  $$
since both $\, Ad\, x\, $ and $\, Ad\, y\, $ are normal from which $\, x = y\, $ follows. 
(That $\, Ad\, x\, $ is normal for positive $\, x\, $ in a monotone complete $C^*$-algebra is easily seen if $\, x\, $ is invertible, however in the general case one can approximate $\, x\, $ in norm by invertible positive elements of the form $\, x_{\epsilon } = x + \epsilon {\bf 1}\, $ whence the result). Then $\, x\leq p\, ,\, x\leq q\, $ and we claim that $\, x\, $ is the largest positive element with this property whence $\, x = p\wedge q\, $. Given $\, 0\leq z\leq p , q\, $ one has $\, p\,z\, = \, z\, =\, q\, z\, $ then also $\, x_n\, z\, =\, z\, =\, z\, x_n\, $ for every $\, n\, $ so that 
$\, z\, x\, z = z^2\, $ and $\, z\, ( 1 - x ) = 0 = ( 1 - x )\, z\, $ follows. 
Thus the projections form a lattice putting $\, p\vee q = 1 - ( ( 1 - p )\wedge ( 1 - q ) )\, $. Then suppose given a monotone decreasing net of projections $\, ( p_{\lambda } )_{\lambda } \searrow p\, $. One needs to show that $\, p = p^2\iff \sqrt p = p\, $. By monotonicity of the squareroot one has 
$$ \sqrt p\>\leq\> \inf_{\lambda }\, \sqrt{ p_{\lambda }}\> =\> p $$
and the reverse relation is trivial proving $\, p = p^2\, $. Then the projections form a (sequentially) complete lattice. 
\par\noindent
From the second property of a $\Sigma $-algebra the monotone completion of each 
unital sub-$C^*$-algebra generated by a single positive element $\, x\geq 0\, $ is an abelian subalgebra of the double commutant of $\, x\, $ hence injective so that since any commutative injective $C^*$-algebra is a $C^*$-quotient of some 
$\, l_{\infty } ( Z )\, $ by Theorem A of section 2 any positive element can be approximated in norm from above and from below by finite positive linear combinations of pairwise commuting (spectral) projections. The same is true in case of an $A\Sigma $-algebra or a $\, {\Sigma }_{\omega }$-algebra. On any abelian subalgebra the $\wedge $- and $\vee $-operations defined for projections generalize to the usual lattice operations for arbitrary selfadjoint elements in a commutative setting. Given any positive element 
$\, x\geq 0\, $ in $\, A\, $ there is a minimal projection $\, p_x\, $ with $\, p_x x = x p_x = x\, $ called the support projection. Indeed, if $\, p\, $ and $\, q\, $ are two projections with this property one readily checks that also $\, ( p\wedge q ) x = x\, $, hence the projections $\, \{ p_{\lambda , x} \}\, $ with this property form a monotone decreasing net in $\, \mathcal P ( A )\, $ the subset of all positive projections in $\, A\, $ and the infimum $\, p_x\, $ satisfies 
$$ x\, ( {\bf 1}\> -\> p_x )\, x\> =\> \sup_{\lambda }\, x\, ( {\bf 1}\> -\> p_{\lambda , x} )\, x\> =\> 0 $$
whence $\, p_x\, x = x\, $ follows. The support projection may also be defined in a ${\Sigma }_{\omega }$-algebra by $\, p_x = \sup_n\, \root n\of{ x / \Vert x \Vert}\, $ for which one easily checks the same properties. By duality one defines the {\it cosupport projection $\, r_x\, $} of $\, x\, $ to be the complement of the support projection of $\, z = \Vert x \Vert {\bf 1} - x\, $, i.e. $\, r_x = {\bf 1} - p_z\, $.
Let $\, q : A \twoheadrightarrow B\, $ be a surjective normal $*$-homomorphism. Then the subset of contractive positive elements contained in $\,\ker\, q\, $ has a supremum $\, r_q\, $ since for each $\, x\in \ker\, q\, ,\, 0\leq x\leq {\bf 1}\, $ also $\, p_x\in \ker\, q\, $ then by normality $\, r_q = \sup\, \bigl\{ p_x\bigm\vert x\in\ker q\, ,\, x\geq 0 \bigr\}\, $ lies in $\, \ker\, q\, $ and is the maximal projection with this property. Therefore the kernel of $\, q\, $ is contained in $\, r_q\, A\, r_q\, $ while on the other hand the kernel of $\, q\, $ contains $\, r_q\, A\, +\,  A\, r_q\, $ showing that $\, r_q\, $ is central. Putting $\, p_q = ( {\bf 1} - r_q )\, $ one gets $\, \ker\, q = ( {\bf 1} - p_q )\, A\, $ and 
$\, B \simeq p_q\, A\, $ giving a natural normal $*$-homomorphic lift to the quotient map. 
Now let $\, q : A \twoheadrightarrow B\, $ be a surjective sequentially normal $*$-homomorphism. 
We need to show that $\, q\, $ is a lattice map, i.e. for any pair of positive projections $\, e\, ,\, f\,\in A\, $ the element $\, q ( e\wedge f ) (=: q ( e )\wedge q ( f ))\, $ is the largest positive element in $\, B\, $ dominated by $\, q ( e )\, $ and $\, q ( f )\, $. Since $\, B\, $ is sequentially monotone complete the formulas above yield that $\, B\, $ is a $\mathcal P$-lattice algebra and that 
$\, q\, $ is a lattice map. 
\par\noindent
Now suppose that $\, q : A \twoheadrightarrow B\, $ is any surjective $*$-homomorphism which is a $\mathcal P$-lattice map with $\, B\, $ a $\mathcal P$-lattice algebra (e.g. for commutative algebras any $*$-homomorphism is a lattice map) and $\, \mathcal P\subseteq \mathcal P ( B )\, $ is a countable sublattice. As a warm up for the more complicated task of constructing $\mathcal P$-maps with fancy decorations we begin by constructing a ${\mathcal P}^c$-map so that certain general features of the construction are more apparent and the reader can compare and recourse to this much simpler construction on examining the others. The first task is to construct a monotonous complemented map 
$$ s :\> \mathcal P ( B )\> \longrightarrow\> \mathcal P ( A )   $$
which is a cross section restricted to  $\, \mathcal P\, $. 
Using a transfinite enumeration of the set of positive projections 
$\, \mathcal P ( B ) = \{ e_{\gamma }\,\vert\, \gamma\in\Gamma \}\, $ with $\, \Gamma\, $ a well ordered set such that the complementary projection of $\, e_{\gamma }\, $ is either the direct successor or the direct precessor of $\, e_{\gamma }\, $ starting with the pair $\, ( 0 , {\bf 1} )\, $ then followed by the projections in $\, \mathcal P\, $ one constructs a monotonous map of the positive projections in $\, B\, $ into the set $\, \mathcal P ( A )\, $ of positive projections in $\, A\, $ by transfinite induction. Put $\, s ( 0 ) = 0\, ,\, s ( {\bf 1} ) = {\bf 1}\, $. Assume by induction that for some given ordinal $\, {\gamma }_0\, $ with direct successor $\, {\gamma }_0 + 1\, $ corresponding to the complementary projection 
$\, {\bf 1} - e_{{\gamma }_0}\, $ one has constructed a coherent lift $\, s\, $ of the subset 
$\, \bigl\{ e_{\gamma } {\bigr\}}_{\gamma < {\gamma }_0}\, $ into the set of projections of $\, A\, $ meaning that for all $\,\gamma\, ,\, \kappa < {\gamma }_0\, $ the relations
$$ s ( {\bf 1}\, -\, e_{\gamma } )\> =\> {\bf 1}\> -\> s ( e_{\gamma } )\> ,\quad s ( e_{\gamma } )\>\leq\> s ( e_{\kappa } ) $$
are satisfied whenever $\, e_{\gamma } < e_{\kappa }\, $. 
Put 
$$ \underline{\mathcal E}_{{\gamma }_0}\> =\> \bigcup_{\gamma < {\gamma }_0} \bigl\{ s ( e_{\gamma } )\bigm\vert e_{\gamma } < e_{{\gamma }_0} \bigr\}\> ,\quad \overline {\mathcal E}_{{\gamma }_0}\> =\> \bigcup_{\gamma < {\gamma }_0} \bigl\{ s ( e_{\kappa } )\bigm\vert e_{{\gamma }_0} < e_{\kappa } \bigr\} \> , $$
$$ {\underline e}_{{\gamma }_0}\> =\> \bigvee_{e\in \underline{\mathcal E}_{{\gamma }_0}}\, \bigl\{ e \bigr\}\> ,\quad {\overline e}_{{\gamma }_0}\> =\> \bigwedge_{f\in \overline{\mathcal E}_{{\gamma }_0}}\, \bigl\{ f \bigr\}  $$
and choosing an arbitrary projection $\, e^s_{{\gamma }_0}\, $ with $\, q ( e^s_{{\gamma }_0} ) = e_{{\gamma }_0} \, $ define 
$$  s ( e_{{\gamma }_0} )\> =\> \left( \left( {\underline e}_{{\gamma }_0} \right)  \vee e^s_{{\gamma }_0} \right) \wedge \left( {\overline e}_{{\gamma }_0} \right)\> ,\quad 
s ( {\bf 1} - e_{{\gamma }_0} )\> =\> {\bf 1}\> -\> s ( e_{{\gamma }_0} ) \> .  $$
One proceeds inductively to construct a monotonous map $\, s : \mathcal P ( B ) \rightarrow \mathcal P ( A )\, $. If $\, e\bot f\, $ are orthogonal projections then $\, e\leq {\bf 1} - f\, $ whence $\, s ( e )\leq s ( {\bf 1} - f ) = {\bf 1} - s ( f )\, $, therefore 
$\, s ( e )\bot s ( f )\, $ are orthogonal. Consider the normdense subset $\, {\mathcal P}_f ( B )\subseteq B^{sa}\, $ of the real subspace of selfadjoint elements consisting of finite inear combinations of pairwise commuting projections. If $\, x\in {\mathcal P}_f ( B )\, ,\, x \geq 0\, $ there is a unique decomposition 
$$ x\> =\> \sum_{k=1}^n\, {\alpha }_k\, p_k\> ,\quad 0 < p_1 <\cdots < p_n \leq {\bf 1}\> ,\quad {\alpha }_k\in {\mathbb R}_+\, ,\, 1\leq k \leq n\>  , $$
whence the $\, \{ p_k \}\, $ are pairwise commuting projections (contained in some abelian subalgebra). This follows since the corresponding orthogonal decomposition gives the spectral projections and eigenvalues of $\, x\, $ which if they exist are unique in any $C^*$-algebra. Using this decomposition extend $\, s\, $ to a lift 
$\, \sigma : {\mathcal P}_f ( B ) \rightarrow {\mathcal P}_f ( A )\, $ by the formula 
$$ \sigma ( x )\> =\> \sum_{k=1}^n\, {\alpha }_k s ( p_k ) \> . $$
One needs to check that $\,\sigma\, $ is monotonous on the intersection of any abelian subalgebra $\, C\subseteq B\, $ with  $\, {\mathcal P}_f ( B )_+\, $. Let 
$$ x\> =\> \sum_{k=1}^n\, {\alpha }_k\, p_k\>\leq\> \sum_{l=1}^m\, {\beta }_l\, q_l\> =\> y $$
denote the canonical decompositions of $\, x\leq y\, $ assuming that $\, [ p_k\, ,\, q_l ] = 0\, $ for all $\, k\, ,\, l\, $. By subtracting a suitable scalar multiple of $\, {\bf 1}\, $ from both sides
one may assume that the support projection $\, p_n\, $ of $\, x\, $ is strictly smaller than $\, {\bf 1}\, $. Then 
$\, \min \{ {\alpha }_n , {\beta }_m \}\, p_n\,\leq\, \min \{ {\alpha }_n , {\beta }_m \}\, q_m\, $. 
Subtracting $\, \min \{ {\alpha }_n , {\beta }_m \} \, s ( p_n )\, $ from the images of both sides the remaining positive term involving the projection $\, s ( q_m )\, -\, s ( p_n )\, $ on the right side may be dropped from consideration since it is orthogonal to the support projection on the left and doesn't affect the validity of the inequality relation (this argument only works in the abelian case). Then one proceeds by induction to prove that $\, \sigma ( x ) \leq \sigma ( y )\, $.  If $\, B\, $ is a spectral $\mathcal P$-algebra so that $\, x\, $ can be approximated from below by elements in $\, \mathcal P_f ( B ) \cap \{ x \} ''\, $ these approximations $\, \bigl\{ x_{\lambda } {\bigr\}}_{\lambda\in \Lambda }\, $ are all contained in a maximal abelian $\mathcal P$-subalgebra of $\, B\, $ and can be seen to form a monotone increasing net of elements converging up to $\, x\, $ so that defining $\, \sigma ( x ) = \sup_{\lambda } \sigma ( x_{\lambda } )\, $ gives an extension of $\, \sigma\, $ to all of $\, B\, $ which is seen to be monotonous and hence continuous restricted to each abelian subalgebra, where continuity is a consequence of monotonicity plus the property $\, \sigma \bigl( \alpha\, {\bf 1} + x \bigr) = \alpha\, \sigma \bigl( {\bf 1} \bigr) + \sigma \bigl( x \bigr)\, $. Nonetheless $\, \sigma\, $ may fail to be monotonous or even continuous in general. 
Positive homogeneity is more than obvious from the definition. Even if $\, B\, $ is not a spectral $\mathcal P$-algebra it may be embedded into a spectral $\mathcal P$-algebra, e.g. by considering some faithful $*$-representation into $\, \mathcal B ( \mathcal H )\, $ which is a spectral $\mathcal P$-algebra, so speaking generally if given any $\mathcal P$-map of $\, B \subseteq \mathcal B ( \mathcal H )\, $ into a $\Sigma $-algebra $\, A\, $ it may be extended to a $\mathcal P$-map $\, \overline s : {\mathcal P}_f ( \mathcal B ( \mathcal H ) ) \rightarrow {\mathcal P}_f ( A )\, $ which then extends to a map $\, \overline\sigma : \mathcal B ( \mathcal H ) \rightarrow A\, $ with the property of being monotonous and continuous restricted to any abelian subalgebra, in particular one finds that the value of $\, \sigma ( x )\, $ for some given positive element $\, x\in B_+\, $ must be larger than any of the approximations $\, \sigma ( x_{\lambda } )\, ,\, x_{\lambda }\in {\mathcal P}_f ( B )\, ,\, x_{\lambda } \leq x\, ,\, \bigl[ x\, ,\, x_{\lambda } \bigr[ = 0\, $ in case that $\, B\, $ is a ${\mathcal P}_{\sigma }$-algebra, and smaller than any of its approximations from above, since for each approximation from below up to the order $\, \epsilon > 0\, $ in norm there exists an approximation up to the same order from above in $\, \mathcal B ( \mathcal H )\, $ by a spectral element commuting with the smaller element while on the other hand these elements being close in norm must also have images which are close in norm. Similarly the value of $\, \sigma ( x )\, $ is necessarily smaller than any of its approximations from above in $\, {\mathcal P}_f ( B )\, $. To see that this uniquely determines the value of $\, \sigma ( x )\, $ note that for any element $\, x_{\lambda }\in {\mathcal P}_f ( B )\, ,\, x_{\lambda } \leq x\, ,\, \Vert x - x_{\lambda } \Vert \leq \epsilon\, $ commuting with $\, x\, $ there exists the element $\, {\overline x}_{\lambda } = x_{\lambda } + \epsilon\, {\bf 1}\, $ satisfies $\, {\overline x}_{\lambda } \geq x\, ,\, \Vert {\overline x}_{\lambda } - x \Vert \leq \epsilon\, $ and as $\, \epsilon > 0\, $ is arbitrary the result follows. That $\, \sigma\, $ can be extended at all with the property of being continuous and monotonous on abelian subalgebras and that at least in case of a ${\mathcal P}_{\sigma }$-algebra the extension is unique is quite stunning.
Returning to the setting as above it follows from construction that $\, \sigma\, $ is a cross section for $\, q\, $ restricted to $\, {\mathcal P}_f\, $.
Suppose that $\, B\, $ is abelian. To construct a $\mathcal P$-$A$-lattice map $\, \sigma  : B \rightarrow A\, $ one begins with an enumeration $\, \{ e_{\gamma } \}\, $ of the set $\, \mathcal P ( B )\, $ as above. 
For simplicity we write $\, \gamma\wedge\kappa\, $ resp. $\, \gamma\vee\kappa\, $ to denote the index corresponding to the projection $\, e_{\gamma }\wedge e_{\kappa }\, $ and $\, e_{\gamma }\vee e_{\kappa }\, $ etc., then assume by induction that for some given finite ordinal 
$\, {\gamma }_0\, $ one already has constructed a coherent lift of the uncomplemented sublattice 
$\, {\Lambda }^{ {\gamma }_0}\, $, i.e. $\, {\Lambda }_{{\gamma }_0}\, $ is closed under the lattice operations $\, \wedge , \vee\, $ but not necessarily under $\, c\, $,  generated by the set of all projections $\, \{ e_{\gamma } {\}}_{\gamma < {\gamma }_0}\, $ into the set of projections of $\, A\, $ meaning that for all 
$\,\gamma\, ,\, \kappa\in {\Lambda }^{{\gamma }_0}\, $ the relations
$$ s ( e_{\gamma}\wedge e_{\kappa } )\> =\> s ( e_{\gamma } ) \wedge s ( e_{\kappa } )\> ,\quad 
s ( e_{\gamma }\vee e_{\kappa } )\> =\> s ( e_{\gamma } )\vee s ( e_{\kappa } )\> ,\quad 
\bigl[ s ( e_{\gamma } )\, ,\, s ( e_{\kappa } ) \bigr]\> =\> 0  $$
are satisfied and in general
$$ e_{\gamma }\>\leq e_{\kappa }\>\Rightarrow\> s ( e_{\gamma } )\>\leq\> s ( e_{\kappa } )\> . $$
Note that $\, {\Lambda }^{{\gamma }_0}\, $ is finite for any finite ordinal $\, {\gamma }_0\, $. These conditions imply that putting
$$ {\underline e}_{\mu}^{{\gamma }_0}\> =\> \bigvee_{\omega } \bigl\{ s ( e_{\omega } )\bigm\vert 
\omega\in {\Lambda }_{{\gamma }_0}\, ,\, e_{\omega } \leq e_{\mu } \bigr\}\> ,\quad 
{\overline e}_{\mu }^{{\gamma }_0}\> =\> \bigwedge_{\rho } \bigl\{ s ( e_{\rho } )\bigm\vert \rho\in {\Lambda }_{{\gamma }_0}\, ,\> e_{\rho }\geq e_{\mu } \bigr\}  $$
for arbitrary $\,\mu\, $ one has 
$$ {\underline e}_{\mu}^{{\gamma }_0} \wedge s ( e_{\gamma } )\> \leq\> {\underline e}_{\mu\wedge \gamma }^{{\gamma }_0}\> ,\quad {\overline e}_{{\gamma }_0}\vee s ( e_{\kappa } )\> \geq\>  {\overline e}_{\mu\vee\kappa }^{{\gamma }_0}   $$
for any indices  $\, \gamma\, ,\,\kappa\in {\Lambda }_{{\gamma }_0}\, $. In case of infinite 
$\, {\gamma }_0\, $ one uses the fact that any monotone complete abelian subalgebra containing all $\, \{ s ( e_{\gamma } )\,\vert\, e_{\gamma }\in {\Lambda }_{{\gamma }_0} \}\, $ is a complete function lattice by the Corollary of Theorem 2 to get this result.
We first treat the case of finite $\, {\gamma }_0\, $ since in this context the sublattice $\, {\Lambda }_{{\gamma }_0 + 1}\, $ generated by $\, {\Lambda }_{{\gamma }_0}\, $ and $\, e_{{\gamma }_0}\, $ is finite and posesses minimal elements not contained in $\, {\Lambda }_{{\gamma }_0}\, $. We may choose one such element denoted $\, e_{{\delta }_0} \leq e_{{\gamma }_0}\, $. Then inductively choose a sequence of elements $\, \{ e_{{\delta }_k} {\}}_{k = 1}^n\subseteq {\Lambda }_{{\gamma }_0 + 1}\, $ such that if
$\, {\Lambda }_{{\gamma }_0}^k\, $ denotes the sublattice generated by $\, {\Lambda }_{{\gamma }_0}\, $ and $\, \{ e_{{\delta }_0}\, ,\,\cdots\, ,\, e_{{\delta }_k} \}\, $ the element 
$\, e_{{\delta }_{k + 1}}\, $ is a minimal element in the complement 
$\, {\Lambda }_{{\gamma }_0 + 1} \backslash {\Lambda }_{{\gamma }_0}^k\, $ and $\, {\Lambda }_{{\gamma }_0}^n = {\Lambda }_{{\gamma }_0 + 1}\, $.
Assume again by induction that one has already constructed a coherent ($A$-lattice) lift 
on the sublattice $\, {\Lambda }_{{\gamma }_0}^{k - 1}\, $ for given $\, k \leq n\, $. 
One proceeds as above to obtain a projection 
$\, z\in A\, $ with
$$ {\underline e}^{{\gamma }_0}_{{\delta }_k}\leq z\leq {\overline e}^{{\gamma }_0}_{{\delta }_k}\> ,  $$
and such that $\, q ( z ) = e_{{\delta }_k}\, $. 
Then we still have to incorporate the correct lattice and commutation relations. Any projection of the form 
$\, e_{{\delta }_k\wedge\gamma}\, $ strictly smaller than $\, e_{{\delta }_k}\, $ is already contained in $\, {\Lambda }_{{\gamma }_0}^{k - 1}\, $. Then $\, s ( e_{{\delta }_k\wedge\gamma } ) = {\underline e}^{{\gamma }_0}_{{\delta }_k}\wedge s ( e_{\gamma } ) \leq z\wedge s ( e_{\gamma } )\, $ and replacing $\, z\, $ by 
$$ z'\> =\> \left[ z\wedge \bigl( {\bf 1} - s ( e_{\gamma } ) \bigr)\right]\> +\> 
s ( e_{{\delta }_k\wedge\gamma } )\>\leq\> z $$
the new element satisfies
$\, z'\geq {\underline e}^{{\gamma }_0}_{{\delta }_k}\, ,\, q ( z' ) = e_{{\delta }_k}\, ,\, [ z'\, ,\, s ( e_{\gamma } ) ]\, =\, 0\, $ and $\, z'\wedge s ( e_{\gamma } ) = s ( e_{{\delta }_k}\wedge e_{\gamma } )\, $.
For the first inequality note that since $\, {\underline e}^{{\gamma }_0}_{{\delta }_k}\, $ commutes with $\, s ( e_{\gamma } )\, $ one has 
$$ {\underline e}^{{\gamma }_0}_{{\delta }_k}\> =\> \left( {\underline e}^{{\gamma }_0}_{{\delta }_k}\wedge \bigl( {\bf 1}\, -\, s ( e_{\gamma } \bigr)\right)\> +\> s ( e_{{\delta }_k\wedge\gamma } )\>\leq\> z'\> . $$
One proceeds in this manner by induction with respect to the induced enumeration to cover all indices $\, \gamma\in {\Lambda }_{{\gamma }_0}^{k - 1}\, $ leading to a (finite) monotone decreasing sequence of elements 
$\, ( z^{( n ) } )\, $. Then put 
$\, \underline z\, =\, \inf_n  z^{( n )}\geq {\underline e}^{{\gamma }_0}_{{\delta }_k}\, $ which is a lift of $\, e_{{\delta }_k}\, $ such that 
$\, \underline z\wedge s ( e_{\gamma } ) = s ( e_{{\delta }_k\wedge\gamma } )\, $ and $\, \bigl[ \underline z\, ,\, s ( e_{\gamma } ) \bigr] = 0\, $ for each $\,\gamma\in {\Lambda }_{{\gamma }_0}^{k - 1}\, $. Then suppose that 
$\, e_{{\delta }_k\vee\kappa }\in {\Lambda }_{{\gamma }_0}^{k - 1}\, $ for some given element 
$\, e_{\kappa }\in {\Lambda }_{{\gamma }_0}^{k - 1}\, $. Replacing $\, \underline z\, $ by 
$$ {\underline z}'\> =\> \left[ \underline z\> +\> \left( s ( e_{{\delta }_k\vee \kappa } )\> -\> \bigl( \underline z\vee s ( e_{\kappa } ) \bigr) \right) \right] \wedge {\overline e}^{{\gamma }_0}_{{\delta }_k} $$ 
the new element satisfies the same relations as $\,\underline z\, $ plus the relation 
$\, {\underline z}' \vee s ( e_{\kappa } ) = s ( e_{{\delta }_k\vee\kappa } )\, $. Proceeding as before by induction with respect to the induced enumeration of indices $\, \{ \kappa\, \vert\, e_{\kappa }\in {\Lambda }_{{\gamma }_0}^{k - 1}\, ,\, e_{{\delta }_k\vee\kappa }\in {\Lambda }_{{\gamma }_0}^{k - 1} \}\, $ one arrives at an element $\, \overline z\, $ such that 
$$ \overline z\wedge s ( e_{\gamma } )\> =\> s ( e_{{\delta }_k\wedge \gamma } )\> ,\quad \overline z\vee s ( e_{\kappa } )\> =\> s ( e_{{\delta }_k\vee\kappa } )\> ,\quad \bigl[ \overline z\, ,\, s ( e_{\mu } ) \bigr]\> =\> 0 $$
for all $\, \gamma , \kappa , \mu\, $ with $\, e_{{\delta }_k\vee\kappa }\in {\mathcal E}_{{\gamma }_0}^{k - 1}\, $. 
If $\, e_{\mu }\, ,\, e_{\kappa }\in {\mathcal E}_{{\gamma }_0}^{k - 1}\, ,\, e_{\mu }\leq e_{{\delta }_k\vee\kappa }\, $ but $\, s ( e_{\mu } ) \nleq \overline z\vee s ( e_{\kappa } )\, $ replace 
$\, \overline z\, $ by
$$ \left[ \overline z\> +\> \bigl( \overline z\vee s ( e_{\kappa } )\vee s ( e_{\mu } )\> -\> \overline z\vee s ( e_{\kappa } ) \bigr) \right] \wedge {\overline e}^{{\gamma }_0}_{{\delta }_k} $$
etc. and repeating this procedure for all pairs $\, \mu , \kappa\, $ as above define $\, s ( e_{{\delta }_k} )\, $ to be the supremum (maximum) of the corresponding monotone increasing (finite) sequence of elements. Check that $\,s ( e_{{\delta }_k} )\, $ satisfies the same relations as before plus the relations 
$$ e_{\mu }\>\leq\> e_{{\delta }_k}\vee e_{\kappa }\quad\Rightarrow\quad s ( e_{\mu } )\>\leq\> s ( e_{{\delta }_k} )\vee s ( e_{\kappa } ) $$
Note that $\, e_{\mu }\leq e_{{\delta }_k}\vee e_{\kappa }\, $ implies 
$\, s ( e_{\mu } )\leq {\overline e}^{{\gamma }_0}_{{\delta }_k}\vee s ( e_{\kappa } )\, $ since we are working in a commutative setting. Defining 
$$ s ( e_{{\delta }_k\vee\kappa } )\> :=\> s ( e_{{\delta }_k} ) \vee s ( e_{\kappa } ) $$
then gives a coherent $A$-lattice lift for $\, {\Lambda }_{{\gamma }_0}^k\, $ completing the induction step for finite $\, {\gamma }_0\, $. We may assume that the union of all uncomplemented sublattices 
$$ {\bigcup }_{\gamma < \infty}\, {\Lambda }_{\gamma } $$ 
is a complemented sublattice. Then the $A$-lattice condition implies that the map $\, s\, $ is complemented restricted to this union. In case of infinite $\, {\gamma }_0\, $ the complement $\, {\Lambda }_{{\gamma }_0 + 1} \backslash {\Lambda }_{{\gamma }_0}\, $ may contain no minimal elements.  Assume given a coherent 
$A$-lattice-map $\, s\, $ on the sublattice $\, {\Lambda }_{{\gamma }_0}\, $ containing all elements $\, \{ e_{\gamma }\,\vert\, \gamma < {\gamma }_0 \}\, $ together with their finite products (in a commutative setting the wedge of two projections is equal to their product) and vees.
Since we don't have to worry about the cross section property we may simply define 
$$ s ( e_{{\gamma }_0}\wedge e_{\gamma } )\> =\> {\underline e}_{{\gamma }_0}\wedge s ( e_{\gamma } )\> ,\quad s ( e_{{\gamma }_0}\vee e_{\kappa } )\> =\> {\underline e}_{{\gamma }_0}\vee s ( e_{\kappa } )  $$ 
where 
$$ {\underline e}_{{\gamma }_0}\> =\> \sup\, \bigl\{ s ( e_{\gamma } )\bigm\vert e_{\gamma }\in {\Lambda }_{{\gamma }_0}\, ,\, e_{\gamma } < e_{{\gamma }_0} \bigr\} $$  
and check that it satisfies all the required relations. This is because $\, {\underline e}_{{\gamma }_0}\, $ commutes with all images of elements in $\, {\Lambda }_{{\gamma }_0}\, $ so that we are working in a complete function lattice whence $\, \bigl( \sup_{\lambda } x_{\lambda } \bigr)\wedge y\, =\, 
\sup_{\lambda }\, ( x_{\lambda }\wedge y )\, $ for any monotone increasing net $\, ( x_{\lambda } )\, $ and an arbitrary element $\, y\, $. 
\par\noindent
Suppose now that $\, A\, $ is a $\Sigma $-algebra  with $\, q : A \twoheadrightarrow B\, $ a surjective $\mathcal P$-lattice homomorphism such that any abelian subalgebra is contained in an abelian $\mathcal P$-subalgebra. If $\, x = v + i\, w\, ,\, \bigl[ v\, ,\, w \bigr] = 0\, $ is a normal element of $\, B\, $, then $\, C^* ( x , x^* )\, $ is contained in an abelian $\mathcal P$-subalgebra $\, C \subseteq B\, $ and there exists from the argument above a $*$-homomorphic lift $\, \sigma : C \rightarrow A\, $ (even without the assumption that $\, q\, $ is an overall $\mathcal P$-lattice map since any $*$-homomorphism is an $A$-lattice map and only commuting projections are involved in the construction) whence $\, \sigma ( x )\, $ is a normal lift for $\, x\, $. Assuming instead of the second condition that $\, A\, $ is a ${\Sigma }^+$-algebra and $\, B\, $ admits a faithful order isomorphic (not necessarily $*$-homomorphic or even linear) representation on a separable Hilbert space implying that any well ordered strictly monotone increasing subnet in $\, B\, $ is at most countable, let 
$\, u = x + i y\in B\, $ be a unitary with selfadjoint elements $\, x\, ,\, y\in B^{sa}\, $. Then $\, x\, , y\, $ lift to selfadjoint elements $\, \overline X\, ,\, \overline Y\in A^{sa}\, $ which admit a polar decomposition of the form $\, \overline X = \overline V\, \vert \overline X\vert\, ,\, \overline Y = \overline W\, \vert \overline Y\vert\, $ such that $\, \overline V = {\overline V}^* = {\overline V}^{-1}\, $ is a selfadjoint unitary commuting with $\, \vert \overline X\vert\, $, and $\, \overline W = {\overline W}^* = {\overline W}^{-1}\, $ is a selfadjoint unitary commuting with $\, \vert \overline Y\vert\, $. Putting $\, v = q ( \overline V )\, ,\, w = q ( \overline W )\, $ one obtains a corresponding polar decomposition of $\, x\, $ and $\, y\, $ in $\, B\, $. Since $\, u\, $ is unitary one gets $\, \vert y \vert = \sqrt{ {\bf 1} - \vert x {\vert }^2}\, $ and hence
$\, [ v\, ,\, \vert y \vert ] = 0 = \bigl[ w\, ,\, \vert x \vert \bigr] = \bigl[ v\, ,\, w \bigr] \, \vert x\vert\, \vert y \vert\, $. Put $\, p = {{\bf 1} - v\over 2}\, ,\, q = {{\bf 1} - w\over 2}\, $ and 
$$ r\> =\> \bigl( p \wedge q \bigr)\> +\> \bigl( p \wedge ( {\bf 1} - q ) \bigr)\> +\> \bigl( ( {\bf 1} - p ) \wedge q \bigr)\> +\> \bigl( ( {\bf 1} - p ) \wedge ( {\bf 1} - q ) \bigr) $$
and check that $\, r\, \vert x\vert\, \vert y \vert = \vert x\vert\, \vert y \vert\, $ whence $\, r > 0\, $ unless $\, \vert x\vert = 0\, $ or $\, \vert y\vert = 0\, $.
This follows since the relation $\, \bigl[ p\, ,\, q \bigr]\, \vert x \vert\, \vert y \vert = 0 = \vert x \vert\, \vert y \vert\, \bigl[ p\, ,\, q \bigr]\, $ implies that 
$\, p\, \vert x \vert\, \vert y \vert\, q = p\, \vert x \vert\, \vert y \vert\, q\, p = q\, p\, \vert x \vert \, \vert y \vert\, q\, $ is smaller than both $\, p\, ,\, q\, $ hence smaller than $\, p \wedge q\, $ and similarly with $\, p\, $ replaced by $\, {\bf 1} - p\, $ or $\, q\, $ replaced by $\, {\bf 1} - q\, $ so that 
$\, \vert x\vert\, \vert y \vert\, $ is seen to be smaller than $\, r\, $. The size condition on $\, B\, $ implies that $\, B\, $ is sequentially monotone complete from part (i) of Theorem 3 below so that in particular $\, \mathcal P ( B )\, $ is a sequentially monotone complete lattice.
One may also show that the subset of projections $\, Ann ( x ) = \bigl\{ e \bigm\vert e\, \vert x\vert = 0 \bigr\} \subseteq \mathcal P ( B )\, $ annihilating $\, x\, $ is a sequentially complete sublattice meaning that $\, e\, ,\, f \in Ann ( x )\, $ implies $\, e \vee f\in Ann ( x )\, $ and given any monotone increasing sequence $\, \{ e_n \} \subseteq Ann ( x )\, $ there exists a projection $\, e\in Ann ( x )\, $ which is larger than any of the projections $\, \{ e_n \}\, $. Indeed, if $\, e\, ,\, f\in Ann ( x )\, $ considering the monotone increasing sequence $\, \bigl\{ {\bf 1} - \bigl( ( {\bf 1} - e )\, ( {\bf 1} - f)\, ( {\bf 1} - e ) \bigr)^n \bigr\}\, $ each of these elements annihilates $\, x\, $ so that these elements, which can be expressed by powers of a single positive element and $\, {\bf 1}\, $, and $\, \vert x \vert\, $ are contained in a common abelian $\mathcal P$-subalgebra of $\, B\, $ for which there exists a $*$-homomorphic lift into $\, A\, $ whence the images of these elements annihilate the image of $\, \vert x \vert\, $ whence also their supremum whose image in $\, B\, $ is larger than 
$\, e \vee f\, $ has the same property, cf. \cite{KaPe}. Similarly given a monotone increasing sequence 
$\, \{ e_n \}\subseteq Ann ( x )\, $ these are contained together with $\, \vert x \vert\, $ in a common abelian $\mathcal P$-subalgebra of $\, B\, $ so again there exists a $*$-homomorphic lift into $\, A\, $ where the supremum of the image projections is seen to annihilate the image of $\, \vert x \vert\, $ proving the claim. Then from the condition that any well ordered strictly monotone increasing subnet in $\, B\, $ is at most countable $\, Ann ( x )\, $ must contain a unique maximal element $\, e\, $ and this element is seen to commute with both $\, p\, $ and $\, q\, $ since $\, e\in Ann ( x )\, $ implies 
$\, v\, e\, v\, ,\, w\, e\, w\in Ann ( x )\, $ so that by the above argument also $\, e \vee ( v\, e\, v )\, $ and $\, e  \vee  ( w\, e\, w )\, $ are contained in $\, Ann ( x )\, $ and must agree with $\, e\, $ by maximality. But it is plain to see that $\, e \vee ( v\, e\, v )\, $ commutes with $\, p\, $ and $\, e \vee ( w\, e\, w )\, $ commutes with $\, q\, $.
Since $\, p\, ,\, r\, ,\, e\, ,\, p\wedge q\, ,\, ( {\bf 1} - p ) \wedge q\, ,\, p \wedge ( {\bf 1} - q )\, ,\, ( {\bf 1} - p ) \wedge ( {\bf 1} - q )\, ,\, \vert x \vert\, ,\, \vert y \vert\, $ each commute with each other and with $\, \vert x \vert\, ,\, \vert y \vert \in C^* ( \vert x \vert , {\bf 1} )\, $ they are all contained in a common separable abelian $\mathcal P$-subalgebra $\, C \subseteq B\, $ and there exists a $*$-homomorphic cross section $\, l : C \rightarrow A\, $. Put $\, P = l ( p )\, ,\, R = l ( r )\, ,\, E = l ( e )\, ,\,\vert X \vert = l ( \vert x \vert )\, ,\, \Vert Y \vert = l ( \vert y \vert )\, $ with $\, X = ( {\bf 1} - 2\,P )\, \vert X \vert\, $. Then in particular one has $\, \vert X {\vert }^2 + \vert Y {\vert }^2 = 1\, ,\, \vert X \vert\, \vert Y \vert \leq R\, $ so that $\, ( {\bf 1} - R )\, \vert X \vert\, \vert Y \vert = 0\, $. Define the projection $\, \overline E\in \mathcal P ( A )\, $ by $\, \overline E = \lim_n \vert Y {\vert }^n \geq E\, $ which gives the projection onto the kernel of $\, \vert X \vert\, $ by the relation $\, \vert X {\vert }^2 + \vert Y {\vert }^2 = {\bf 1}\, $, hence commutes with $\, P\, ,\, R\, ,\, \vert X \vert\, $  and $\, \vert Y \vert\, $. By maximality $\, e = q ( \overline E )\, $ and one can choose a preimage of $\, q\, $ in $\, \mathcal P ( A )\, $ larger than $\, l \bigl( p \wedge q \bigr)\, ,\, l \bigl( ( {\bf 1} - p ) \wedge q \bigr)\, $ and commuting with both $\, R\, ,$ and $\, \overline E\, $ and put $\, Y = ( {\bf 1} - 2\, Q ) \vert Y \vert\, $, finally putting $\, U = X + i\, Y\, $ which is seen to be a preimage of $\, u\, $. One easily checks that $\, U^*\, U = \vert X {\vert }^2 + \vert Y {\vert }^2 = {\bf 1}\, $ so that $\, U\, $ is an isometry implying $\, U\, U^* \leq U^*\, U\, $. On the other hand one checks that $\, U\, U^* = \vert X {\vert }^2 + W\, \vert Y {\vert }^2 \, W\, $ implying $\, W\, \vert Y {\vert }^2\, W \leq \vert Y {\vert }^2\, $ and since $\, W\, $ is a symmetry this implies $\, U\, U^* = {\bf 1} = U^*\, U\, $, thus $\, U\, $ is a unitary lift for $\, u\, $ as claimed.
\par\noindent
We now come to the case of arbitrary (noncommutative) quotient $\, B\, $ (of finite type). The proof uses in parts the special case above but the problem is much more complex in that we have to consider relations of noncommuting projections. One observes the following fact: for any pair of projections $\, e\, ,\, f\in \mathcal P ( B )\, $ there exists a unique minimal projection $\, e^f\, $ larger than $\, e\, $ and commuting with $\, f\, $ and a unique maximal projection $\, e_f\, $ smaller than $\, e\, $ and commuting with $\, f\, $. To see  this put 
$$ e_f\> =\> e\wedge f\> +\> e\wedge \bigl( {\bf 1} - f \bigr)\>\leq\> e\>\leq\> 
\left[ \bigl( {\bf 1} - e \bigr)\wedge f\> +\> \bigl( {\bf 1} - e \bigr)\wedge \bigl( {\bf 1} - f \bigr) \right]^c\> =\> e^f \> . \leqno{(*)}$$
One has $\, [ e^f\, ,\, f ] = 0 = [ e_f\, ,\, f ]\, $ and $\, e^f = e = e_f\, $ if and only if $\, [ e\, ,\, f ] = 0\, $. From monotonicity of the assignment 
$\, e \mapsto e^f\, $ one finds that $\, e^f\, $ is the minimal projection commuting with $\, f\, $ and larger than $\, e\, $. Correspondingly $\, e_f\, $ is the maximal projection smaller than $\, e\, $ commuting with $\, f\, $. Similarly given $\, e\, $ and a multiplet $\, F = ( f_1\, ,\cdots ,\, f_n )\, $ of projections there exists a unique minimal projection $\, e^F\geq e\, $ 
and a unique maximal projection $\, e_F \leq e\, $ commuting with each $\, f_k\, $ for 
$\, k = 1 ,\cdots , n\, $. For example if $\, n = 2\, $ we may consider the monotone increasing sequence 
$$ e\>\leq\> e^{f_1}\>\leq\> ( e^{f_1} )^{f_2}\> \leq\> ( ( e^{f_1} )^{f_2} )^{f_1}\>\leq\> \cdots $$
and check that each projection in this sequence is smaller than any projection larger than $\, e\, $ and commuting with both $\, f_1\, $ and $\, f_2\, $. Since all projections are contained in the  sublattice generated by $\, \{ e , f_1 , f_2 \}\, $ the sequence must become stationary after finitely many steps. Then the corresponding element commutes with both $\, f_1\, $ and $\, f_2\, $ and obviously is the minimal projection with this property that dominates $\, e\, $. By induction this argument extends to any finite set $\, F = ( f_1\, ,\cdots ,\, f_n )\, $ and the case of 
$\, e_F\, $ follows by symmetry since $\, ( e^F )^c = ( e^c )_F\, $. 
Now let $\, q : A \twoheadrightarrow B\, $ be as above and assume that $\, A\, $ is of finite 
type which implies that $\, B\, $ too is of finite type. We want to construct a ${\mathcal P}^{o}$-map $\, \sigma : B \rightarrow A\, $ which is a cross section for $\, q\, $ and a ${\mathcal P}^{\wedge }$-map restricted to a chosen countable sublattice ${\mathcal P}_0\subseteq \mathcal P ( B )\, $ as above. Choose a transfinite enumeration 
$\, \bigl\{ e_{\gamma } \bigr\}\, $ of the elements of $\, \mathcal P ( B )\, $ beginning with a sequential enumeration of the countable sublattice $\, {\mathcal P}_0\, $ whose linear span is dense in some separable subalgebra containing $\, B_0\, $. Put 
$\, {\mathcal E}_{{\gamma }_0} = \bigl\{ e_{\gamma }\bigm\vert \gamma < {\gamma }_0 \bigr\}\, $ and let $\, {\Lambda }_{{\gamma }_0}\, $ denote the sublattice generated by $\, {\mathcal E}_{{\gamma }_0}\, $. We first treat the case of finite $\, {\gamma }_0\, $. Assume by induction that one has already constructed a coherent lift
$\, s = s_{{\gamma }_0} : {\Lambda }_{{\gamma }_0} \rightarrow \mathcal P ( A )\, $ which is a cross section for $\, q\, $ meaning that the following must be satisfied: 
$$ s ( e \wedge f )\> =\> s ( e )\wedge s ( f )\> ,\qquad e\,\bot\, f\quad\Longrightarrow\quad s ( e )\, \bot\, s ( f )  $$
also implying 
$$ e\leq f\quad\Longrightarrow\quad s ( e ) \leq s ( f )\> . $$
We also assume by induction that if $\, d\in {\Lambda }_{{\gamma }_0}\, $ is a central element then 
$\, s ( d )\vee s ( e ) = s ( d\vee e )\, $ for any $\, e\in {\Lambda }_{{\gamma }_0}\, $ and $\, s ( d^c ) = s ( d )^c\, $.
Of course one presets the values $\, s ( 0 ) = 0\, ,\, s ( {\bf 1} ) = {\bf 1}\, $. Note that if $\, \tilde s\, $ is a $\mathcal P$-lattice map and a cross section then $\, \tilde s \bigl( {\Lambda }_{{\gamma }_0} \bigr)\subseteq \mathcal P ( A )\, $ is a sublattice so that for any $\, e_{\gamma }\in {\Lambda }_{{\gamma }_0}\, $ the element $\, \tilde s ( e_{\gamma } )\, $ is minimal (unique)  in the lattice $\, \tilde s \bigl( {\Lambda }_{{\gamma }_0} \bigr)\, $ with the property that it is a preimage of $\, e_{\gamma }\, $, in other words the lattice generated by the images of
$\, {\Lambda }_{{\gamma }_0}\, $ contains no elements which lie in the kernel of $\, q\, $.
Let $\, {\mathcal Z}_{{\gamma }_0 + 1}\subseteq {\Lambda }_{{\gamma }_0 + 1}\, $ denote the center of $\, {\Lambda }_{{\gamma }_0 + 1}\, $ consisting of those projections commuting with any other projection. Since $\, {\Lambda }_{{\gamma }_0 + 1}\, $ is finite so is its subset of minimal central projections $\, \{ c_1\, ,\cdots ,\, c_n \}\, $ which generate the abelian lattice $\, {\mathcal Z}_{{\gamma }_0 + 1}\, $. We may then extend the given lift $\, s\, $ to a coherent lift on 
$\, {\Lambda }_{{\gamma }_0} \cup {\mathcal Z}_{{\gamma }_0 + 1}\, $ which is an $A$-lattice-map restricted to $\, {\mathcal Z}_{{\gamma }_0 + 1}\, $ in the manner above and on noting that each element in $\, {\Lambda }_{{\gamma }_0 + 1}\, $ has a canonical central decomposition as 
$$ e_{\gamma }\> =\> \sum_{k = 1}^n\, c_k\, e_{\gamma } $$
extend this lift to a coherent $\, {\mathcal P}^{\wedge , a}$-lift on the sublattice generated by $\, {\Lambda }_{{\gamma }_0}\, $ and $\, {\mathcal Z}_{{\gamma }_0 + 1}\, $ by defining 
$$ s \left( \sum_{k = 1}^n\, c_k\, e_{\gamma } \right)\> =\> \sum_{k = 1}^n\, s ( c_k )\, s ( e_{\gamma } ) \> . $$
Thus $\, s\, $ is an $A$-lattice-map restricted to $\, {\mathcal Z}_{{\gamma }_0 + 1}\, {\mathcal Z}_{{\gamma }_0}\, $.
Let $\, \bigl\{ d_1\, ,\cdots ,\, d_m \bigr\}\, $ be an enumeration of the minimal central elements of 
$\, {\mathcal Z}_{{\gamma }_0 + 1}\, {\Lambda }_{{\gamma }_0}\, $, i.e. $\, d_l \leq c_{k_l}\, $ for some 
minimal element $\, c_{k_l}\in {\mathcal Z}_{{\gamma }_0 + 1}\, $.
Then choose elements $\, \{ e_{{\delta }_k} \}\subseteq  {\Lambda }_{{\gamma }_0 + 1} \backslash \bigl( {\mathcal Z}_{{\gamma }_0 + 1}\, {\Lambda }_{{\gamma }_0} \bigr)\, $ such that if 
$\, {\Lambda }_{{\gamma }_0}^k\, $ is the sublattice of $\, {\Lambda }_{{\gamma }_0 + 1}\, $ generated by $\, {\mathcal Z}_{{\gamma }_0 + 1}\, {\Lambda }_{{\gamma }_0} \cup 
\bigl\{ e_{{\delta }_1}\, ,\cdots ,\, e_{{\delta }_k} \bigr\}\, $ then $\, e_{{\delta }_k}\, $ is minimal in the complement of $\, {\Lambda }_{{\gamma }_0}^{k - 1}\, $. By abuse of notation let $\, d_k\, $ be the minimal central element of 
$\, {\Lambda }_{{\gamma }_0}^{k - 1}\, $ exceeding $\, e_{{\delta }_k}\, $.
Define 
$$ {\underline e}_{{\delta }_k}\> =\> \sup_{e_{\gamma }\in {\Lambda }_{{\gamma }_0}^{k - 1}}\, \bigl\{ s ( e_{\gamma } )\bigm\vert e_{\gamma } < e_{{\delta }_k} \bigr\}\>  ,\quad {\overline e}_{{\delta }_k}\> =\> \inf_{e_{\kappa }\in {\Lambda }_{{\gamma }_0}^{k - 1}}\, \bigl\{ s ( e_{\kappa } )\bigm\vert e_{\kappa } > e_{{\delta }_k} \bigr\} \> . $$
For each $\, k\, $ we want to extend $\, s\, $ to a coherent lift of the sublattice $\, d_k\, {\Lambda }_{{\gamma }_0}^k\, $ beginning with a lift of $\, e_{{\delta }_k}\, $ leaving the images of 
$\, \bigl( {\bf 1} - d_k \bigr)\, {\Lambda }_{{\gamma }_0}^{k - 1}\, $ unchanged. Choose any preimage 
$\, z_{{\delta }_k}\, $ of $\, e_{{\delta }_k}\, $ with 
$$ {\underline e}_{{\delta }_k}\>\leq z_{{\delta }_k}\>\leq\> {\overline e}_{{\delta }_k}  $$
then considering the sublattice $\, \widetilde{\Lambda }_{{\gamma }_0}^k\, $ generated by $\, s \bigl( d_k\, {\Lambda }_{{\gamma }_0}^{k - 1} \bigr)\, $ and $\, z_{{\delta }_k}\, $ it contains being of finite type a minimal element 
$\, {\underline s} ( d_k ) \, $ which is a preimage of $\, d_k\, $. Obviously being minimal 
$\, {\underline s} ( d_k )\, $ must be central, since $\, {\underline s} ( d_k ) = {\underline s} ( d_k )_{s \bigl( {\Lambda }_{{\gamma }_0}^{k - 1}\bigr) \cup \{ z_{{\beta }_k} \}}\, $. Then 
$\, \eta = s ( d_k ) - {\underline s} ( d_k )\, $ is a central element in the kernel of $\, q\, $.
Define 
$$ {\underline s} ( e_{{\delta }_k} )\> =\> \bigl( {\bf 1}\, -\, \eta \bigr)\, z_{{\delta }_k}\> ,\quad 
{\underline s} ( e_{\gamma } )\> =\> \bigl( {\bf 1}\, -\, \eta \bigr)\, s ( e_{\gamma } ) $$
for $\, e_{\gamma }\in {\Lambda }_{{\gamma }_0}^{k - 1}\, $. One checks that each element in the sublattice $\, \bigl( {\bf 1} - \eta \bigr)\, \widetilde{\Lambda }_{{\gamma }_0}^k\, $ is minimal with the property of being a preimage for its image under $\, q\, $. Therefore the elements of 
$\, \bigl( {\bf 1} - \eta \bigr)\, \widetilde{\Lambda }_{{\gamma }_0}^k\, $ and the lattice operations are in one-to-one correspondence with the elements and lattice operations of $\, d_k\, {\Lambda }_{{\gamma }_0}^k\, $ 
so that $\, \underline s\, $ extends uniquely to a $\mathcal P$-lattice map. This in turn implies that 
the complementary map 
$$ e_{\gamma }\>\mapsto\> \eta\, s ( e_{\gamma } )\> =:\> {\eta }_{\gamma } $$
for $\, e_{\gamma }\in {\Lambda }_{{\gamma }_0}^{k - 1}\, $ is a ${\mathcal P}^{\wedge , o}$-map which extends to a $\mathcal P$-map on $\, d_k\, {\Lambda }_{{\gamma }_0}^k\, $ putting
$$ {\eta }_{\lambda }\> =\> \eta\, {\underline e}_{\lambda }\> ,\quad 
{\underline e}_{\lambda }\> =\> \sup_{e_{\gamma }\in {\Lambda }_{{\gamma }_0}^{k - 1}}\, \bigl\{ 
s ( e_{\gamma } )\bigm\vert e_{\gamma } < e_{\lambda } \bigr\}\> . $$
Monotonicity and orthogonality of this map are obvious. Then we need to check that 
it is a wedge-map. Since we are considering finite lattices the suprema are attained so that 
$\, {\underline e}_{\lambda } = s ( e_{{\gamma }_{\lambda }} )\, $ where 
$\, e_{{\gamma }_{\lambda }} < e_{\lambda }\, $ is the unique maximal element in 
$\, {\Lambda }_{{\gamma }_0}^{k - 1}\, $ which is smaller than $\, e_{\lambda }\, $. Then the equality 
$\, e_{{\gamma }_{\lambda }}\wedge e_{{\gamma }_{\mu }} = e_{{\gamma }_{\lambda\wedge\mu }}\, $ is straightforward and $\, s\, $ being a wedge map on $\, {\Lambda }_{{\gamma }_0}^{k - 1}\, $ by induction assumption the result follows. Define the extension of $\, s\, $ to $\, d_k\, {\Lambda }_{{\gamma }_0}^k\, $ by 
$$ s ( e_{\lambda } )\> =\> {\underline s} ( e_{\lambda } )\> +\> {\eta }_{\lambda }  $$
which again is a ${\mathcal P}^{\wedge , o}$-map. By induction the argument carries on to give a 
$\, {\mathcal P}^{\wedge , o}$-map extension of $\, s\, $ to all of $\, {\Lambda }_{{\gamma }_0 + 1}\, $
which again satisfies the extra assumption that the images of central elements remain central completing the induction step. This process leads to a ${\mathcal P}^{\wedge , o}$-map defined on  
$\, {\mathcal P}_0\, $. It is then easy to extend this to a ${\mathcal P}^{o}$-map on 
$\, \mathcal P ( B )\, $ by inductively defining
$$ s ( e_{{\gamma }_0} )\> =\> \sup\, \bigl\{ s ( e_{\gamma } )\bigm\vert e_{\gamma } < e_{{\gamma }_0}\, ,\, \gamma < {\gamma }_0 \bigr\}  $$
for infinite $\, {\gamma }_0\, $. Monotonicity and orthogonality of this extension are obvious. 
By duality one obtains a $\, {\mathcal P}^{co}$-map $\, s\, $ which is a $\, {\mathcal P}^{\vee , co}$-map restricted to $\, {\mathcal P}_0\, $ putting 
$$ \widetilde s ( e_{\gamma } )\> =\> {\bf 1}\> -\> s ( e_{\gamma }^c )\> . $$  
Next assume that $\, q\, $ is sequentially normal. Applying the above construction consider the sequence of (nonunital) cross sections $\, {\underline s}_{{\gamma }_0} : {\Lambda }_{{\gamma }_0} \rightarrow \mathcal P ( A )\, $ corresponding to the minimal $\mathcal P$-lattice-map part of 
$\, s\, $ restricted to $\, {\Lambda }_{{\gamma }_0}\, $ for finite index $\, {\gamma }_0\, $ as above. Then this leads for fixed $\, {\Lambda }_{{\gamma }_0}\, $ to a 
monotone decreasing sequence of cross sections $\, \bigl\{ {\underline s}^{{\gamma }_0}_{\gamma } \bigr\}\, $ which are the restrictions of $\, {\underline s}_{\gamma }\, $ to $\, {\Lambda }_{{\gamma }_0}\, $. Since each element of $\, {\mathcal P}_0\, $ has a finite index the limit of these maps yields a well defined cross section $\, \underline s : {\mathcal P}_0 \rightarrow \mathcal P ( A )\, $ from the fact that $\, q\, $ is sequentially normal. 
Let us show that it is a $\mathcal P$-lattice map with respect to the unit element $\, \underline s ( \bf 1 )\, $. Since all images of a pair of complementary projections $\, \{ e , {\bf 1} - e \}\, $ under the family of ${\mathcal P}^c$-maps $\, \{ {\underline s}_{\gamma } \}\, $ (with respect to 
the unit element $\, {\underline s}_{\gamma } ( {\bf 1} )\, $) commute with each other it is easy to see that the limit map is again complemented. Then one only needs to show that it respects wedges. Being represented as a limit of a family of monotone decreasing ${\mathcal P}^{\wedge }$-maps this is again obvious. Fixing the unit $\, \underline s ( {\bf 1} )\in \mathcal P ( A )\, $, i.e. replacing 
$\, \mathcal P ( A )\, $ by $\, \mathcal P \bigl( \underline s ( {\bf 1} )\, A\, \underline s ( {\bf 1} ) \bigr)\, $ this map may be extended to a ${\mathcal P}^c$-map on the whole of $\, \mathcal P ( B )\, $ in the manner described above. 
\par\noindent
Then let $\, I\, $ be an injective $C^*$-algebra and 
$\, \lambda : I \rightarrow \mathcal B ( \mathcal H )\, $ a faithful unital $*$-representation. Choose any completely positive retraction $\, r : \mathcal B ( \mathcal H ) \rightarrow I\, $ for 
$\,\lambda\, $. Then $\, r\, $ is an $I$-module map in the sense that 
$\, r ( x \lambda ( y ) ) = r ( x )\, y\, $ and $\, r ( \lambda ( y ) x ) = y r ( x )\, $ (compare Lemma 6.1.2 of \cite{E-R}). If $\, x\in I\, $ is arbitrary then 
$$ x\> =\> r \bigl( \lambda ( x ) \bigr)\> =\> r \bigl( v\, \lambda ( \vert x\vert ) \bigr)\> =\> r ( v )\, \vert x \vert\> =\> \vert x \vert\, r ( v ) $$
for some partial isometry $\, v\in \mathcal B ( \mathcal H )\, $ with $\, v^*\, v \leq \lambda ( p_x )\, $ 
and $\, p_x\, $ the support projection of $\, x^*\, x\, $. Then the support projection of 
$\, r ( v^* )\, r ( v )\, $ is smaller or equal to $\, p_x\, $ from the Schwarz inequality. Then 
$$ x^*\, x\> =\> \vert x \vert\, r ( v^* )\, r ( v )\, \vert x\vert\> =\> r ( v^* )\, r ( v )\, ( x^*\, x ) $$
implying $\, r ( v^* )\, r ( v ) = p_x\, $ and similarly $\, r ( v )\, r ( v^* ) = p_{x^*}\, $ so that 
$\, r ( v )\, $ is a partial isometry giving the polar decomposition of $\, x\, $. Thus $\, I\, $ is an 
${\Sigma }^+$-algebra as claimed. 
\par\noindent
Next suppose given a $\mathcal P$-map of $\mathcal P$-algebras $\, \sigma : {\mathcal P}_f ( A ) \rightarrow {\mathcal P}_f ( B )\, $. Then $\, \sigma \, $ is determined by its restriction $\, s : \mathcal P ( A ) \rightarrow \mathcal P ( B )\, $ to the subset of positive projections. For this 
suppose given $\, x = \sum_{k=1}^n {\alpha }_k\, p_k\in {\mathcal P}_f ( A )\, $ with 
$\, 0\, <\, p_1\, <\,\cdots\, <\, p_n\, \leq\, {\bf 1}\, $ a finite string of increasing projections and positive coefficients
$\, {\alpha }_k > 0\, ,\, 1\leq k\leq n\, $. 
By monotonicity 
$\, \sigma ( x )\, $ must be larger than any of the elements
$\, \bigl\{ \bigl( \sum_{j = k}^n {\alpha }_j \bigr) s ( p_k ) \bigr\}\, $, while on the other hand it must be smaller or equal than each of the elements 
$\, \bigl\{ \bigl( \sum_{j = k}^n {\alpha }_j \bigr) s ( p_k ) + \bigl( \sum_{j = 1}^n {\alpha }_j \bigr) \, \bigl( {\bf 1}\, -\, s ( p_k ) \bigr) \bigr\}\, $. 
These conditions imply that $\, \sigma  ( x )\, $ has the form 
$\, ( {\overline x}_{ij} )_{ij}\, ,\, 1\leq i , j\leq n\, $ with respect to the matrix decomposition of 
$\, A\, $ corresponding to the finite set of orthogonal projections 
$\, \bigl\{ s ( p_k )\, -\, s ( p_{k - 1} )\,\vert\, k = 1\, ,\cdots\, ,\, n \bigr\}\, $ putting $\, p_0 = 0\, $ such that $\, {\overline x}_{ij} = 0\, $ whenever $\, i\neq j\, $, and $\, {\overline x}_{ii} = \bigl( \sum_{j = k}^n\, {\alpha }_i \bigr) s ( p_k )\, -\, s ( p_{k - 1} )\, $ for $\, k = 1\, ,\cdots\, ,\, n\, $. Therefore on the real subspace of selfadjoint elements it is uniquely determined by the values $\, \{ s ( p_k ) \}\, $ and the condition $\, \sigma ( \alpha {\bf 1} + x ) = \alpha \sigma ( {\bf 1} ) + \sigma ( x )\, $ (or else $\, \sigma  ( x ) = \sigma ( x_+ ) - \sigma ( x_- )\, $). Since the spectral projections of $\, x^2\, $ in the canonical decomposition as above coincide with those of $\, x\, $ for a positive element $\, x\in {\mathcal P}_f ( A )\, $ one checks the relation $\, \sigma ( x^2 ) = \sigma ( x )^2\, $ for such elements and the case of general positive elements follows by continuity. If $\, \sigma\, $ is a homogeneized ${\mathcal P}^o$-map one has for an element in $\, x\in {\mathcal P}_f ( A )\, $ that the spectral projections of $\, x\, $ corresponding to positive eigenvalues are orthogonal to the spectral projections corresponding to negative eigenvalues. Therefore 
$$ \sigma ( x^2 )\> =\> \sigma ( x_+^2 + x_-^2 )\> \geq\> \sigma ( x_+^2 )\> +\> \sigma ( x_-^2 )\> =\> \sigma ( x_+ )^2\> +\> \sigma ( x_- )^2\> =\> \sigma ( x )^2 $$
proving the Schwarz inequality for selfadjoint elements in $\, \mathcal P ( A )_f\, $. The general case again follows by continuity. 
This proves the theorem\qed
\par\bigskip\noindent
{\it Remark.}\quad (i) As a Corollary of these results one proves the approximate triple cone property of 
a $*$-homomorphic quotient map $\, q : A \twoheadrightarrow B\, $ of ${\mathcal P}_{\sigma }$-algebras (even when $\, q\, $ is not an overall $\mathcal P$-lattice map) which may be useful in certain instances.
Since $\, q\, $ is a $*$-homomorphism it has the following general approximate double cone property: given elements $\, a < b\, $ in $\, A\, $ and an element $\, q ( a )\leq x\leq q ( b )\, $ in $\, B\, $ there exists 
$\, a\leq y\leq b + \epsilon {\bf 1}\, $ with $\, q ( y ) = x\, $.
This result is fairly obvious in the case where $\, a = 0\, ,\, b = {\bf 1}\, $ in which case one has an exact double cone property, i.e. $\,\epsilon\, $ can be chosen to be zero which is seen by choosing a unital completely positive lift of the separable abelian $C^*$-algebra $\, C^* ( x, {\bf 1} )\, $ for $\, q\, $. The method above which applies in any case to abelian $\mathcal P$-algebras generated by a single selfadjoint element regardless of whether $\, q\, $ is generally a lattice map
(a $*$-homomorphism of abelian algebras is always a lattice map), shows that there even exists a $*$-homomorphic cross section. This proves that the double cone problem for $\, a = 0\, ,\, b = {\bf 1}\, $ has an exact solution for arbitrary $C^*$-algebras. Then also the double cone problem for $\, a = 0\, $ and $\, b = p\, $ a projection has an exact solution which follows by passing to the hereditary subalgebra $\, p\, A\, p\, $. The general case may be transformed to the case $\, a = 0\, ,\, b = {\bf 1}\, $ by first applying a linear transformation sending $\, a\, $ to $\, 0\, $ and then replacing 
$\, x\, $ by $\, \sqrt{q ( b ) + \epsilon\, {\bf 1}}^{-1}\, x\, \sqrt{q ( b ) + \epsilon {\bf 1}}^{-1}\, $ and 
$\, q ( b )\, ,\, b\, $ by $\, \bf 1\, $ respectively for given $\, \epsilon > 0\, $ proving (only) the approximate double cone property for arbitrary $\, a\leq b\, $. If $\, A\, $ is a ${\mathcal P}_{\sigma }$-algebra $\, q\, $ also posesses a corresponding (approximate) triple cone property: for given $\,\epsilon > 0\, $, elements $\, a\, ,\, b\geq 0\, $ in $\, A\, $ and an element
$\, 0\leq x\leq q ( a )\, ,\, q ( b )\, $ there exists an element $\, 0\leq y < a + \epsilon\, {\bf 1}\, ,\, b + \epsilon {\bf 1}\, $ with $\, q ( y ) = x\, $. To see this replace $\, b\, ,\, q ( b )\, $ by $\, {\bf 1}\, $ and 
$\, a\, $ by $\, a' = \sqrt{b + \delta\, {\bf 1}}^{-1}\, a\, \sqrt{b + \delta\, {\bf 1}}^{-1}\, $, also replacing
 $\, q ( a )\, $ and $\, x\, $  by  
$\,\sqrt{q ( b ) + \delta\, {\bf 1}}^{-1}\, q ( a )\, \sqrt{q ( b ) + \delta\, {\bf 1}}^{-1}\, $ and 
$\, x' = \sqrt{q ( b ) + \delta\, {\bf 1}}^{-1}\, x\, \sqrt{q ( b ) + \delta\, {\bf 1}}^{-1}\, $ respectively where $\, \delta = \epsilon / 2\, $.
Then solving the approximate double cone problem for $\, 0\leq a'\, ,\, 0\leq x'\leq q ( a' )\, $ gives an element $\, 0\leq y'\leq a'\, +\, ( \delta / \Vert b\Vert )\, {\bf 1}\, $ with $\, q ( y' ) = x'\, $. Let $\, z\in {\mathcal P}_f ( A )\, $ be an element commuting with $\, y'\, $ and satisfying 
$\, \Vert z - y' \Vert < \epsilon / ( 3 \Vert b\Vert )\, $. If $\, p_z^{\bf 1}\, $ denotes the spectral projection of $\, z\, $ corresponding to all eigenvalues larger or equal than $\, 1 + \epsilon / (3 \Vert b\Vert )\, $ then $\, p_z^{\bf 1}\, $ is in the kernel of $\, q\, $ and $\, \Vert ( {\bf 1} - p_z^{\bf 1} )\, y' \Vert < 1 + \epsilon / (3 \Vert b\Vert )\, $ giving an element 
$\, 0\leq y < a + \epsilon\, {\bf 1}\, ,\, b + \epsilon\, {\bf 1}\, $ with 
$\, q ( y ) = x\, $. By induction one finds that there is an approximate multiple cone property for each $\, n\in\mathbb N\, $ as long as all cones except possibly one are pointing in the same direction, i.e for any $\,\epsilon > 0\, $ and given positive elements $\, 0\leq a_1\, ,\cdots\, ,\, a_n\, $ and $\, 0\leq x\leq q( a_1 )\, ,\cdots\, ,\, q ( a_n )\, $ there exists $\, 0\leq y < a_1 + \epsilon\, {\bf 1}\, ,\cdots\, ,\, a_n + \epsilon\, {\bf 1}\, $ with $\, q ( y ) = x\, $. 
\par\smallskip\noindent
(ii) The proof above that the subset of projections in a $\Sigma $-algebra is a complete lattice easily extends to $A\Sigma$-algebras and more generally to $A{\Sigma }_{\omega }$-algebras (with sequentially complete lattice). 
In a $\mathcal PA{\Sigma }_{\omega }$-algebra (= $A{\Sigma }_{\omega }$-algebra) $\, A\, $ every positive element $\, x\geq 0\, $ has a well defined support projection. For if $\, ( x_k ) \nearrow x\, $ is a monotone increasing normconvergent sequence of positive elements in $\, {\mathcal P}_f ( A )\, $ converging up to $\, x\, $  then each $\, x_k\, $ has a well defined support projection $\, p_k\, $ and it is easily checked that its supremum $\, p_x = \sup_k\, p_k\, $ satisfies $\, p_x\, x\, =\, x\, p_x\, =\, x\, $. On the other hand if 
$\, q\in \mathcal P ( A )\, $ is any other projection with this property then $\, q \geq p_k\, $ for all $\, k\, $ whence $\, q \geq p_x\, $ follows. Alternatively one can define $\, p_x\, $ for $\, 0\leq x\leq {\bf 1}\, $ to be $\, p_x = \sup_n\, \sqrt[n]{x}\, $ on checking that the monotone increasing sequence 
$\, ( \sqrt[n]{x}\, {)}_n\, $ is contained in an abelian $\mathcal P$-subalgebra whence its supremum is well defined in the $\mathcal P$-commutant of the sequence (which is the $\mathcal P$-commutant of $\, \{ x \}\, $) and is contained in the double $\mathcal P$-commutant of the sequence (in $\, \{ x {\}}_{\mathcal P}''\, $). From continuous function calculus this supremum is a projection $\, p\, $ with $\, p \geq p_k\, $ for each $\, k\, $ as above so that $\, p \geq p_x\, $ while on the other hand $\, p_x \geq \sqrt[n]{x}\, $ implies $\, p_x \geq p\, $ hence equality. 
It is also true that for any fixed projection $\, p\in A\, $ and any finite subset $\, G \subseteq \mathcal P ( A )\, $ (in fact any countable subset) of projections the upper and lower commutants $\, p^G\, $ and $\, p_G\, $ of $\, p\, $ with $\, G\, $ are well defined in $\, \mathcal P ( A )\, $. If $\, G = \{ g_1\, ,\cdots ,\, g_n \}\, $ consider the monotone increasing sequence defined recursively by $\, p_0 = p\, ,\, p_{k + 1} = p_k^{n_{k + 1}}\, $ where 
$\, k = l\, n + n_k\, ,\, 1 \leq n_k \leq n\, ,\, l\in \mathbb N\, $.  Then the supremum 
$\, q = \sup_k\, p_k\, $ commutes with each projection $\, g_i\, $ for $\, i = 1 ,\cdots , n\, $. Namely fixing $\, i\, $ and selecting all indices $\, k\, $ with $\, n_k = i\, $ the supremum over this cofinal subsequence $\, \{ p_{i_k} \}_\subseteq \{ p_k \}\, $ is again equal to $\, q\, $ and since each single projection $\, p_{i_k}\, $ commutes with $\, g_i\, $ so does the supremum. It is easy to see that $\, q\, $ is the minimal projection larger than $\, p\, $ with this property whence $\, q = p^G\, $. The argument for $\, p_G\, $ follows by symmetry.
\par\smallskip\noindent
(iii) It should be interesting to find means of constructing $\mathcal P$-maps sending equivalent projections to equivalent projections. Such a map is likely to be (overall) continuous or even monotonous and might have applications to $K$-theory, but this seems to be a much harder task in general.
\par\bigskip\noindent
A {\it projection filter} is a subset $\, \mathcal F \subseteq \mathcal P ( A )\, $ of the set of positive projections in a $\mathcal P$-lattice algebra $\, A\, $ satisfying 
$$ e \leq f\> ,\> e\in \mathcal F\quad\Longrightarrow\quad f\in \mathcal F\> ,\> e^c\notin \mathcal F\> , \qquad e\, ,\, f\in \mathcal F\quad \Longrightarrow\quad e \wedge f\in \mathcal F \> . $$
A projection filter defines a presignature which corresponds to a ${\mathcal P}^{\wedge}$-map $\, s : \mathcal P ( A ) \twoheadrightarrow \mathbb C\, $  putting $\, s^{-1} ( {\bf 1} ) = \mathcal F\, $. A $\mathcal P$-map with range $\,\mathbb C\, $ will also be called a $\mathcal P$-functional.
An {\it ultrafilter} is a projection filter $\, \mathcal F\subseteq \mathcal P ( A )\, $ such that for any 
$\, p\in \mathcal P ( A )\, $ either $\, p\in \mathcal F\, $ or else $\, p^c\in \mathcal F\, $. One easily finds that an ultrafilter defines a polar signature on $\, \mathcal P ( A )\, $, which corresponds to a $\mathcal P$-lattice map $\, {\omega }_{\mathcal F} : A \twoheadrightarrow \mathbb C\, $ from the second condition of a projection filter which is supposed to hold whether or not $\, e\, $ commutes with $\, f\, $. A projection filter $\, \mathcal F\, $ defines a {\it lattice ideal} iff in addition the following condition holds
$$ f\in \mathcal F \> ,\> p\in \mathcal P ( A )\quad\Longrightarrow\quad f_p\> =\> f\wedge p\> +\> f\wedge p^c\>\in\> \mathcal F\> . $$
Define an equivalence relation on 
$\, \mathcal P ( A )\, $ by concatenation of simple equivalences 
$$ f\in \mathcal F\quad\Longrightarrow\quad p\> {\sim }_{\mathcal F}\> p\vee f^c\> ,\quad p\> {\sim}_{\mathcal F}\> p\wedge f\> . $$
The equivalence relation is compatible with taking complements since 
$$ ( p \vee f^c )^c\> =\> p^c \wedge f\> {\sim }_{\mathcal F}\> p^c\quad ,\quad ( p \wedge f )^c\> =\> p^c \vee f^c\> {\sim }_{\mathcal F}\> f^c \> . $$ 
It is also compatible with the wedge-operation since $\, f\in \mathcal F\, $ implies 
$\, ( p \vee f^c - p )^c\in \mathcal F\, $ for arbitrary $\, p\in \mathcal P ( A )\, $ whence
$$ ( p \vee f^c ) \wedge q\> {\sim }_{\mathcal F}\> 
\left[ \bigl( ( p \vee f^c\> -\> p )\> +\> p \bigr) \wedge q \right] \wedge \bigl( p \vee f^c\> -\> p \bigr)^c\> =\> p\wedge q\> , $$
$$ ( p \wedge f ) \wedge q\> =\> ( p \wedge q ) \wedge f\> {\sim }_{\mathcal F}\> p\wedge q\> .  $$
Then for any two equivalent elements $\, p {\sim }_{\mathcal F} p'\, $ there is a finite chain $\, \{ f_1\, ,\, \cdots\, ,\, f_n \} \subseteq \mathcal F\, $ with 
$$ p'\> =\> \bigl(\bigl( \cdots \bigl( p\vee f_1^c \bigr) \wedge f_2 \bigr)\cdots \bigr) \vee f_n^c $$
so that if $\, p \leq q \leq p'\, $ one gets $\, p' \leq p\vee f^c {\sim }_{\mathcal F} p\, $ with $\, f =  f_1\wedge\cdots \wedge f_n\in \mathcal F\, $ whence 
$$ q\> {\sim }_{\mathcal F}\> q\vee f^c\> =\> p \vee f^c\> {\sim }_{\mathcal F}\> p \> . $$
Thus the quotient set modulo the equivalence relation is partially ordered by the order relation induced from $\, \mathcal P ( A )\, $ and is a complemented lattice denoted $\, \mathcal P ( A )_{\mathcal F} = \mathcal P ( A ) \bigm/ {\sim }_{\mathcal F}\, $. A projection filter generating a lattice ideal will be called an {\it ideal filter}. An ultrafilter is automatically an ideal filter since either 
$\, p\in \mathcal F\, $ or else $\, p^c\in \mathcal F\, $ whence $\, f\in \mathcal F\> \Longrightarrow\> f_p\in \mathcal F\, $. Any intersection of ideal filters is an ideal filter. To a projection filter one may also associate a ${\mathcal P}^{\vee }$-map into $\,\mathbb C\, $, interpreting the set 
$\, {\mathcal F}^c = \{ e\in \mathcal P ( A )\,\vert\, e^c\in \mathcal F \}\, $ as the kernel 
$\, s^{-1} ( 0 )\, $ of the corresponding ${\mathcal P}^{\vee }$-functional.  
\par\noindent
A $*$-ideal $\, \mathcal J \vartriangleleft A\, $ of a $\mathcal P$-lattice algebra will be called a {\it ($C^*$-)lattice ideal} if the induced quotient map $\, q : A \twoheadrightarrow A / \mathcal J\, $ is a lattice map. This implies that $\, \mathcal P ( \mathcal J ) \subseteq \mathcal P ( A )\, $ is a lattice ideal (of complemented lattices). However, if $\, \mathcal K ( \mathcal H ) \vartriangleleft \mathcal B ( \mathcal H )\, $ denotes the ideal of compact operators in the algebra of bounded operators on a separable Hilbert space $\, \mathcal H\, $ then $\, \mathcal P \bigl( \mathcal K ( \mathcal H ) \bigr) \subseteq \mathcal P \bigl( \mathcal B ( \mathcal H ) \bigr)\, $ is a lattice ideal without $\, \mathcal K ( \mathcal H ) \vartriangleleft \mathcal B ( \mathcal H )\, $ being a $C^*$-lattice ideal (see below).
A unital abelian $C^*$-algebra $\, D\, $ is {\it separably injective} iff given a unital inclusion of function systems (abelian $C^*$-algebras, linear function lattices) $\, \mathfrak X \subseteq \mathfrak Y\, $ such that $\, \mathfrak X\, $ is separable together with a completely positive linear map ($*$-homomorphism, linear lattice map) 
$\, s : \mathfrak X \rightarrow D\, $ there exists a  completely positive linear ($*$-homomorphic, linear lattice map) extension $\, \overline s : \mathfrak Y \rightarrow D\, $ of $\, s\, $. 
\par\bigskip\noindent
{\bf Corollary 1.}\quad Let $\, A\, $ be an $A\Sigma$-algebra and given a surjective $*$-homomorphism $\, q : A \twoheadrightarrow B\, $  which is a $\mathcal P$-lattice map. Then every separable abelian $\mathcal P$-subalgebra $\, D_{\lambda }\subseteq B\, $ is contained in an injective (abelian) sub-$C^*$-algebra $\, I_{\lambda }\simeq I ( D_{\lambda } )\subseteq B\, $ which admits a $*$-homomorphic cross section 
$\, s_{\lambda } : I_{\lambda } \rightarrow A\, $ to the quotient map. If $\, q\, $ is an arbitrary surjective $*$-homomorphism then any selfadjoint element $\, x\in B^{sa}\, $ is contained in an injective abelian sub-$C^*$-algebra of $\, B\, $ (which admits a $*$-homomorphic lift to $\, A\, $). Any quotient of an injective abelian $C^*$-algebra is separably injective.
The ${\mathcal P}^{\wedge }$- and the ${\mathcal P}^{\vee }$-functionals corresponding to an ideal filter are monotonous (hence continuous). In particular the $\mathcal P$-lattice functional corresponding to an ultrafilter is monotonous. Therefore any $\mathcal P$-lattice map into an abelian $C^*$-algebra is monotonous.
\par\bigskip\noindent
{\it Proof.}\quad Given a surjective $*$-homomorphism $\, q : A \twoheadrightarrow B\, $ with 
$\, A\, $ an $A\Sigma $-algebra assume that $\, B\, $ is a $\mathcal P$-lattice algebra and $\, q\, $ is a $\mathcal P$-lattice map. Let $\, D_{\lambda }\subseteq B\, $ be a separable abelian $\mathcal P$-subalgebra. From the proof of Theorem P there exists a $*$-homomorphic cross section $\, s_{\lambda } : D_{\lambda } \rightarrow C \subseteq A \, $ for the quotient map where $\, C\, $ is some maximal abelian, hence injective sub-$C^*$-algebra of $\, A\, $. Then from the Theorem A of section 2 there exists a $*$-homomorphic extension of $\, s_{\lambda }\, $ to the injective envelope $\, I ( D_{\lambda } )\, $ and from rigidity one gets that the extension is injective. Therefore its composition with the quotient map $\, q\, $ results in an injective $*$-homomorphic extension $\, {\iota }_{\lambda } : I ( D_{\lambda } ) \hookrightarrow B\, $ of $\, D_{\lambda } \subseteq B\, $. Even if $\, q\, $ is not an overall $\mathcal P$-lattice map its restriction to any maximal abelian (= injective) sub-$C^*$-algebra $\, C \subseteq A\, $ is, so that any separable sub-$C^*$-algebra 
$\, D_{\lambda } \subseteq D = q ( C )\, $ is contained in an injective subalgebra isomorphic to 
$\, I ( D_{\lambda } ) \subseteq D\, $. Therefore $\, D\, $ is separably injective. Assume given an ideal filter $\, \mathcal F\subseteq \mathcal P ( A )\, $ on a $\mathcal P$-lattice algebra $\, A\, $. 
Let $\, \overline\sigma\, ,\, \underline\sigma :  A \twoheadrightarrow \mathbb C\, $ be the induced 
${\mathcal P}^{\wedge }$-functional (resp. ${\mathcal P}^{\vee }$-functional) whose kernel restricted to $\, \mathcal P ( A )\, $ is given by $\, \mathcal P ( A ) \backslash \mathcal F = {\overline s}^{-1 } ( 0 )\, $ (resp.
by $\, {\mathcal F}^c = {\underline s}^{-1} ( 0 )\, $). It is sufficient to prove that $\, \overline\sigma ( x ) \leq \overline\sigma ( y )\, $ and $\, \underline\sigma ( x ) \leq \underline\sigma ( y )\, $ for elements $\, x\, ,\, y\in {\mathcal P}_f ( A )\, $ with $\, 0 < x \leq y < {\bf 1}\, $. Write 
$$ x\> =\> \sum_{k = 1}^n\, {\alpha }_k\, p_k\> ,\quad 0\> <\> p_1\>\leq\>\cdots\>\leq\> p_n\> ,\qquad y\> =\> \sum_{l = 1}^m\, {\beta }_l\, q_l \> ,\quad 0\> <\> q_1\>\leq\>\cdots\>\leq\> q_m $$ 
with $\, N = \min \{ n , m \}\, $ and assume by induction that either relation is proved for all pairs 
$\, x' \leq y'\, $ in $\, {\mathcal P}_f ( A )\, $ such that either $\, x'\, $ or $\, y'\, $  is a positive linear combination of at most $\, N - 1\, $ pairwise commuting projections. Then also the other relation follows for such pairs by the symmetry 
$\, \underline\sigma ( x ) = - \overline\sigma ( - x )\, $, or $\, \underline s ( p ) = 
( \overline s ( p^c ) )^c\, $.
Note that always $\, p_n \leq q_m\, $ giving that 
$\, x' \leq \beta\, q \, $ implies $\, \underline\sigma ( x' ) \leq \beta\, \underline s ( q )\, $ and 
$\, \overline\sigma ( x' ) \leq \beta\, \overline s ( q )\, $ whence by the symmetry $\, x \mapsto {\bf 1} - x\, $ also 
$\, \alpha\, p \leq y'\, $ implies $\, \alpha\, \underline s ( p ) \leq \underline\sigma ( y' )\, $ and $\, \alpha \overline s ( p ) \leq \overline\sigma ( y' )\, $ for an induction start. Assume first that $\, n \leq m\, $. Then supposing that $\, \overline\sigma ( x ) > \overline\sigma ( y )\, $ for some pair as above we get from the induction hypothesis that 
$\, \overline s ( p_k ) = {\bf 1}\, $ for each $\, k\, $ so that in particular $\, q_m - p_n \in {\mathcal F}^c\, $. Then $\, q_l\in \mathcal F \iff (q_l)^{q_m - p_n} = \bigl( q_l \vee ( q_m - p_n ) \bigr)\wedge \bigl( q_l \vee ( q_m - p_n )^c \bigr)\in \mathcal F\, $ since $\, q_l {\sim }_{\mathcal F}\, \bigl[ q_l \vee ( q_m - p_n ) \bigr]\, $. Therefore we may replace each projection $\, q_l\, $ by the larger projection $\, (q_l)^{q_m - p_n}\, $ without changing the value of $\, \overline\sigma ( y )\, $. For 
$\, {\alpha }_n \leq {\beta }_m\, $ one gets by subtraction of $\, {\alpha }_n\, p_n\, $ from both sides that $\, x' = \sum_{k = 1}^{n - 1}\, {\alpha }_k\, p_k \leq {\alpha }_n\, ( q_m - p_n )\, +\, 
\sum_{l = 1}^m\, {\beta }_l\, (q_l)^{q_m - p_n}\, +\, ( {\beta }_m - {\alpha }_n )\, q_m\, $ 
whence $\, x' \leq \sum_{l = 1}^{m - 1}\, {\beta }_l\, (q_l)^{q_m - p_n}\, +\, ( {\beta }_m - {\alpha }_n )\, q_m\, = y'\, $ giving 
$\, \overline\sigma ( x' ) \leq \overline\sigma ( y' )\, $ from induction assumption implying the relation $\, \overline\sigma ( x ) \leq \overline\sigma ( y )\, $. On the other hand if 
$\, {\alpha }_n > {\beta }_m\, $ one similarly gets $\, x' = \sum_{k = 1}^{n - 1} {\alpha }_k\, p_k\, +\, ( {\alpha }_n - {\beta }_m )\, p_n\,\leq \sum_{l = 1}^{m - 1}\, {\beta }_l\, (q_l)^{q_m - p_n}\, = y'\, $ implying $\, p_n \leq (q_{m - 1})^{q_m - p_n}\, $ and the process can be repeated until the coefficient of $\, p_n\, $ eventually is decreased to zero. Then the induction assumption applies. 
The case where $\, n > m\, $ is very similar and is left to the reader. If $\, \sigma : A \rightarrow B\, $ is a $\mathcal P$-lattice map with $\, B\, $ abelian then the composition of $\, \sigma\, $ with any multiplicative functional on $\, B\, $ yields a $\mathcal P$-lattice functional on $\, A\, $ which by the argument above is monotonous. Since there is a separating family of such functionals $\, \sigma\, $ is monotonous\qed 
\par\bigskip\noindent
{\it Counterexample.}\quad Consider the Toeplitz extension $\, \mathcal T = C^* ( v , v^* )\, $ of the circle algebra $\, C ( \mathbb T ) = C^* ( u , u^* )\, $ generated by a single isometry $\, v^* v = {\bf 1}\, $. If the Calkin algebra $\, \mathcal Q ( \mathcal H ) = \mathcal B ( \mathcal H ) / \mathcal K ( \mathcal H )\, $ would admit an order isomorphic representation by bounded operators on a separable Hilbert space then the induced surjection
$\, q_0 : \mathcal T \twoheadrightarrow C ( \mathbb T )\, $ would admit a $*$-homomorphic  cross section $\, \sigma : C ( \mathbb T ) \rightarrow \mathcal T\, $, lifting the unitary $\, u\, $ to a unitary in $\, \mathcal T\, $ and any Fredholm operator would have index $\, 0\, $ which clearly is not true.
in fact $\, \mathcal Q ( \mathcal H )\, $ is not even a $\mathcal P$-lattice algebra from Corollary 1.9 of \cite{Vo} together with  \cite{KaPe}). 
On the other hand it is easy to see that  
$\, \mathcal F = \bigl\{ P\in \mathcal B ( \mathcal H )\bigm\vert q ( P ) = {\bf 1} \bigr\}\, $ defines an ideal filter in $\, \mathcal P \bigl( \mathcal B ( \mathcal H ) \bigr)\, $ such that the set of complements consists of all compact projections. One only needs to check that if $\, p\, $ is a compact (hence finitedimensional) projection then $\, q ( P\vee p ) = q ( P )\, $ for any $\, P\in \mathcal P \bigl( \mathcal B ( \mathcal H ) \bigr)\, $ whence by duality $\, q ( P \wedge p^c ) = q ( P )\, $. Thus the 
$A$-lattice surjection $\, \mathcal P \bigl( \mathcal B ( \mathcal H ) \bigr) \twoheadrightarrow \mathcal P \bigl( \mathcal Q ( \mathcal H ) \bigr)\, $ factors over the $\mathcal P$-lattice quotient 
$\, \mathcal P \bigl( \mathcal B ( \mathcal H ) \bigr) \twoheadrightarrow  \mathcal P \bigl( \mathcal B ( \mathcal H ) \bigr)_{\mathcal F}\, $. In the case where $\, q\, $ is a $\mathcal P$-lattice map the lattice quotient obtained from the ideal filter will coincide with the projection lattice of the quotient algebra since if two projections $\, P\, ,\, Q\, $ are identified modulo $\, q\, $ their images also agree with $\, q ( P \vee Q )\, $ whence the projections $\, P \vee Q - P\, $ and $\, P \vee Q - Q\, $ are contained in the kernel and $\, P {\sim }_{\mathcal F} Q\, $ follows. 
\par\bigskip\noindent
Given a set $\, X\, $ let $\, B ( X )\, $ denote the abelian $C^*$-algebra of bounded (complexvalued) functions on $\, X\, $. A {\it $\mu $-measure} $\, m\, $ on $\, X\, $ (also called a premeasure in the literature) is a monotonous set function
$$ m :\> \mathcal P ( X )\> =\> \bigl\{ A\bigm\vert A \subseteq X \bigr\}\> \largerightarrow\> {\mathbb R}_+\cup \{ + \infty \} $$
satisfying
$$ m ( \emptyset )\> =\> 0\> ,  $$
$$ A \subseteq B\> \Longrightarrow\> m ( A )\> \leq\> m ( B ) \> . $$
If $\, m ( X ) = 1\, $ then $\, m\, $ is called a {\it probability $\mu $-measure}, and a {\it finite $\mu $-measure} if $\, m ( X ) < + \infty\, $. To each finite $\mu $-measure $\, m\, $ corresponds a positively homogenous monotonous integral ($\mu $-integral) 
$$ {\int}^m :\> B ( X )^{sa}\> \largerightarrow\> \mathbb R\> ,\quad  f\>\mapsto\> {\int }^m\, f $$
with the following properties
$$ \int^m\, \bigl( \alpha\, f \bigr)\> =\> =\> \alpha\, \left( {\int}^m\, f \right)\> ,\quad \forall\> \alpha\in {\mathbb R}_+\> , $$
$$ \int^m\, \bigl( f\> +\> \alpha\, {\bf 1}_X \bigr)\> =\> \int^m\, f\> +\> \alpha\, m ( X )\> ,\quad \forall\> \alpha\in \mathbb R\> ,\> f\in B ( X )\>  , $$
$$ f \leq g\> \Longrightarrow\> \int^m\, f\> \leq\> \int^m\, g\> . $$ 
To define $\, \int^m\, $ consider the subspace of functions $\, B_f ( X )\, $ which are finite linear combinations of characteristic functions $\, f = {\alpha }_0\, {\bf 1}_X\> +\> \sum_{k = 1}^n\, {\alpha }_k\, {\chi }_{A_k}\, $ where $\, {\alpha }_0\in \mathbb R\> ,\> \bigl\{ {\alpha }_k {\bigr\}}_{k = 1}^n\subset {\mathbb R}_+\, $ and $\, A_1 \supseteq A_2 \supseteq\cdots \supseteq A_n\, $ is a decreasing string of nonempty proper subsets of $\, X\, $ with $\, {\chi }_A\, $ denoting the characteristic function of $\, A\subseteq X\, $. For such $\, f\, $ define 
$$ \int^m\, f\> =\> {\alpha }_0\, m ( X )\> +\> \sum_{k = 1}^n\, {\alpha }_k\, m ( A_k ) $$
and check that this yields a positively homogenous monotonous map on $\, B_f ( X )\, $. 
For general $\, f\in B ( X )\, $ define 
$$ \int^m\, f\> =\> \sup\,\left\{ \int^m g\Bigm\vert g\in B_f ( X )\, ,\, g \leq f \right\}\> =\> \inf\, \left\{ \int^m\, h\Bigm\vert h\in B_f ( X )\, ,\, h \geq f \right\}  $$
and check that this definition yields a well defined monotonous and positively homogenous map satisfying the properties as above and extending the original definition on $\, B_f ( X )\, $ (compare with the definition of $\mathcal P$-maps above, in fact a $\mathcal P$-map can be regarded as a noncommutative projection valued $\mu $-measure). If $\, U\subseteq X\, $ is a subset we adopt the notation 
$$ \int_U^m\, f\> =\> \int^m\, \bigl( f\, {\chi }_U \bigr) \> . $$
In case that $\, m\, $ is actually a (outer) measure, i.e. if 
$$ A \subseteq\bigcup_{k = 1}^{\infty }\, C_k\> \Longrightarrow\> m ( A )\> \leq\>\sum_{k = 1}^{\infty }\, m ( C_k ) $$
and $\, f\, $ is a $m$-measurable function (c.f. \cite{Ro}) the $\mu $-integral of $\, f\, $ coincides with the (linear) integral associated with $\, m\, $, i.e. 
$$ \int^m\, f\> =\> \int\, f\, dm \> . $$
If $\, X\, $ is a compact metric space there is a canonical finite $\mu $-measure defined by 
$$ d ( A )\> :=\> \sup\, \bigl\{ d ( x , y ) \bigm\vert x\, ,\, y \in A \bigr\}\> ,  $$ 
which can be normalized to a probability $\mu $-measure by 
$$ \delta ( A )\> =\> d ( A )\, d ( X )^{-1} \> . $$
There also are other natural $\mu $-measures connected with the metric. For example one can define such a measure by 
$$ r ( A )\> =\> \sup\, \bigl\{ r\bigm\vert U_r ( x )\subseteq A \bigr\}\> ,\qquad \rho ( A )\> =\> r ( A )\, r ( X )^{-1} $$
where $\, U_r ( x ) = \{ y\in X\,\vert\, d ( x , y ) < r \}\, $ denotes the open ball of radius $\, r\, $ around a point $\, x\in X\, $. In this setting any subset $\, A\subseteq X\, $ without interior points will have measure zero. For any monotone increasing sequence $\, \{ A_n\,\vert A_n \subseteq A_m\, ,\, n \leq m  \}\, $ one has $\, \sup_n\, d ( A_n ) = d ( A )\, $ where $\, A = \bigcup_n\, A_n\, $. From this one induces that the integral $\, \int^d\, $ is increasing normal, i.e. if given a monotone increasing sequence of functions $\, ( f_n )_n \nearrow f\, $ converging pointwise to the positive function $\, f\, $  
one has 
$$ \int^d\, f\> =\> \sup_n\, \int^d\, f_n \> . $$
For this it is sufficient to assume that the $\, \{ f_n \}\, $ and $\, f\, $ are positive simple functions bounded by $\, {\bf 1}_X\, $ of the form 
$$ f_n\> =\> \sum_{k = 1}^m\, {1\over m}\, {\chi }_{C_n^k}\> ,\quad f\> =\> \sum_{k = 1}^m\, {1\over m}\, {\chi }_{C^k} $$
since arbitrary functions $\, \{ f_n \}\, ,\, f\, $ can be uniformly approximated in norm by simple functions of the form above and the $\mu$-integral $\, \int^d\, $ is continuous. Then 
$\, C^k = \bigcup_n\, C_n^k\, $ and the result follows from increasing normality $\, d ( C^k ) = \sup_n d ( C_n^k )\, $ of the $\mu$-measure $\, d\, $. In contrast to this one has $\, r ( \cap_n A_n ) = \inf_n r ( A_n )\, $ leading to 
$$ ( g_n )_n \searrow g\quad\Longrightarrow\quad \int^r\, g\> =\> \inf_n\, \int^r\, g_n \> . $$
In fact an increasing normal $\mu $-measure is also obtained from a finite distance function 
$\, \delta : X \times X \rightarrow [ 0 , R ]\, $ assigning a distance $\, \delta ( x , y )\, $ to any pair of different points of $\, X\, $ without assuming the triangle inequality for $\,\delta\, $.
\par\noindent
We then introduce the following notion: given a function system $\, \mathfrak X\, $ consider the space 
$\, {\mathfrak X}^{*\mu}\, $ of differences of positively homogenous and monotonous maps $\, \varphi :\> {\mathfrak X}_+\> \longrightarrow\> {\mathbb R}_+\, $ satisfying
$$ \varphi ( \alpha\, x )\> =\> \alpha\, \varphi ( x )\> ,\quad \forall\> \alpha\in {\mathbb R}_+\> ,\qquad x \leq y\Longrightarrow \varphi ( x )\>\leq \varphi ( y )\> , $$
$$\quad \varphi \bigl( \alpha\, {\bf 1}_{\mathfrak X} + x \bigr)\> =\> \alpha\, \varphi \bigl( {\bf 1}_{\mathfrak X} \bigr)\> +\> \varphi ( x )\> ,\quad \forall\> \alpha\in {\mathbb R}_+ $$
and note that such a map is automatically continuous. The space $\, {\mathfrak X}^{*\mu }\, $ with  norm given by $\, \Vert \varphi \Vert = \sup\, \bigl\{ \vert \varphi ( x ) \vert \bigm\vert  0\leq x \leq 1 \bigr\}\, $
will be called the {\it $\mu $-dual of $\, \mathfrak X\, $}. Any monotonous map $\, \varphi\, $ as above extends uniquely to a (continuous) positively homogenous and monotonous map 
$\, \overline\varphi : \mathfrak X \rightarrow \mathbb R\, $ 
also satisfying
$$ \overline\varphi \bigl( \alpha\, {\bf 1}_{\mathfrak X} + x \bigr)\> =\> \alpha\, \overline\varphi \bigl( {\bf 1}_{\mathfrak X} \bigr)\> +\> \overline\varphi ( x )\> ,\quad \forall\> \alpha\in \mathbb R $$
on putting 
$$ \overline\varphi ( x )\> =\> \varphi \bigl( x + \Vert x\Vert {\bf 1} \bigr)\> -\> \Vert x \Vert\, \varphi ( {\bf 1} )\> . $$
A general monotonous positively homogenous map $\, \phi : \mathfrak X \rightarrow \mathfrak Y\, $ of function systems which is linear for addition of scalar multiples of the identity will be called a {\it $\mu $-map}. The concept of $\mu $-map $\, \phi : \mathfrak M \rightarrow \mathfrak N\, $ makes sense for {\it $\mu $-systems}  where $\, \mathfrak M\, $ and $\, \mathfrak N\, $ respectively are normclosed unital subsets of the positive cone of function systems $\, \mathfrak X\, $ and $\, \mathfrak Y\, $ generated by {\it essential elements} respectively with induced order relation which are closed under addition of positive scalar multiples of the identity and closed under multiplication by positive real scalars, where an element $\, x\in {\mathfrak X}_+\, ,\, \Vert x\Vert = 1\, $ of norm $\, 1\, $ is called {\it essential} iff either $\, x = {\bf 1}_{\mathfrak X}\, $ or else $\, {\alpha }_x = \sup\, \bigl\{ \alpha\in {\mathbb R}_+\bigm\vert \alpha {\bf 1}_{\mathfrak X} \leq x \bigr\} = 0\, $. One easily checks that $\, {\mathfrak X}^{*\mu }\, $ is a normed real ordered space by the addition $\, ( \varphi + \psi ) ( x ) = \varphi ( x ) + \psi ( x )\, $ and order determined by the closed cone $\, {\mathfrak X}^{*\mu }_+\subseteq {\mathfrak X}^{*\mu }\, $ of monotonous (increasing) functions. The subset of (increasing) monotonous elements of norm $\, 1\, $ in $\, {\mathfrak X}^{*\mu }\, $ is denoted $\, {\mathcal S}^{\mu } ( \mathfrak X )\, $ and its members are called $\mu $-states. If instead one considers the ordering determined by the closed subcone of positive definite functions of $\, {\mathfrak X}^{*\mu }\, $, i.e. functions which are positive on $\, {\mathfrak X}_+\, $ then $\, {\mathcal S}^{\mu } \bigl( \mathfrak X \bigr)\, $ is a complete function lattice for this second order with $\, ( \varphi \vee \psi ) ( x ) = \max \bigl\{ \varphi ( x ) , \psi ( x ) \bigr\}\, ,\, ( \varphi \wedge \psi ) ( x ) = \min \bigl\{ \varphi ( x ) , \psi ( x ) \bigr\}\, $ restricted to essential elements $\, x\, $ and 
$\, \varphi ( x ) = \sup_{\lambda }\, {\varphi }_{\lambda } ( x )\, $ (resp. $\, \varphi ( x ) = \inf_{\lambda } {\varphi }_{\lambda } ( x )\, $) if $\, \varphi = \sup_{\lambda }\, {\varphi }_{\lambda }\, $ (resp. $\, \varphi = \inf_{\lambda } {\varphi }_{\lambda }\, $) denotes the supremum (resp. infimum) of a monotone increasing (resp. decreasing)  bounded net in $\, {\mathcal S}^{\mu } \bigl( \mathfrak X \bigr)\, $. To distinguish the different orders we let $\, \mathcal M ( \mathfrak X )\, $ denote the real ordered Banach space which is the completion of the normed vector space $\, {\mathfrak X}^{*\mu }\, $ but order determined by the subcone $\, \mathcal M ( \mathfrak X )_+\subseteq \mathcal M ( \mathfrak X )\, $ of positive definite elements and $\, \mathcal M ( \mathfrak X )_r\subseteq \mathcal M ( \mathfrak X )_+\, $ the convex subset of elements with $\, \gamma \bigl( {\bf 1}_{\mathfrak X} \bigr) = r\, $. Then 
$\, {\mathfrak X}^{*\mu }_{+ , r} = {\mathfrak X}^{*\mu }_+ \cap \mathcal M ( \mathfrak X )_r\, $. The notation $\, \varphi \geq \psi\, $ generally refers to the order by positive definite elements, i.e. means
$\, \varphi ( x ) \geq \psi ( x )\, $ for all $\, x\in \mathfrak X\, $ even when talking about elements of $\, {\mathfrak X}^{*\mu }\, $. Each subset $\, {\mathfrak X}^{*\mu }_{+ , r}\, $ is a complete function lattice for pointwise maxima and minima on essential elements and taking pointwise suprema resp. infima of monotonous increasing (resp. decreasing) nets. The $\mu $-functional $\, {\eta }_r\, $ defined by $\, {\eta }_r ( x ) = {\alpha }_x\, r\, $ if 
$\, {\alpha }_x = \sup\, \bigl\{ \alpha\in {\mathbb R}_+\bigm\vert \alpha\, {\bf 1}_{\mathfrak X} \leq x \bigr\}\, $ is the unique minimal element of $\, {\mathfrak X}^{*\mu }_{+ , r}\, $. Then $\, {\mathfrak X}^{*\mu }_+\, $ is seen to be a complete lattice. Let two $\mu $-functionals $\, \varphi\, ,\, \psi\, $ be given with $\, \varphi ( {\bf 1}_{\mathfrak X} ) = r\, $ and $\, \psi ( {\bf 1}_{\mathfrak X} ) = s\, $ for $\, r \leq s\, $. It is straightforward to see that (affine and homogenous continuation) of the pointwise maximum $\, ( \varphi \vee \psi ) ( x ) = \max\, \bigl\{ \varphi ( x )\, ,\, \psi ( x ) \bigr\}\, $ on essential elements $\, x \in {\mathfrak X}_{ess}\, $ is the unique least upper bound of $\, \{ \varphi\, ,\, \psi \}\, $ in $\, {\mathfrak X}^{*\mu }_+\, $.
To define the wedge operation consider for each $\, t \leq r\, $  the subset of elements $\, C_{\varphi , \psi }^t\, =\, \bigl\{ \kappa\in {\mathfrak X}^{*\mu }_{+ , t}\bigm\vert \kappa \leq \varphi\, ,\, \psi \bigr\}\, $. Then from the argument above if it is not empty each subset $\, C_{\varphi , \psi }^t\, $ has a well defined supremum $\, {\gamma }_{\varphi , \psi }^t\in {\mathfrak X}^{*\mu }_{+ , t}\, $ and it is easy to see that 
$\, {\gamma }_{\varphi , \psi }^t \leq \varphi\, ,\, \psi\, $ for each index $\, t\, $ with 
$\, C_{\varphi , \psi }^t \neq \emptyset\, $. On the other hand since $\, {\eta }_t\, $ is minimal in 
$\, {\mathfrak X}^{*\mu }_{+ , t}\, $ and $\, {\eta }_t \leq {\eta }_r\, $ for $\, t\leq r\, $ the set $\, C_{\varphi , \psi }^t\, $ is nonempty for all values $\, t \leq r\, $ while on the other hand there can be no element $\, \kappa\in\mathcal M ( \mathfrak X )_t\, $ such that $\, \kappa \leq \psi\, $ if $\, t > r\, $. We claim that for $\, t \leq r\, $ every element $\, \kappa\in C_{\varphi , \psi }^t\, $ is majorized by an element of $\, C_{\varphi , \psi }^r\, $. Namely put 
$\, {\kappa }' = \kappa + {\eta }_{r - t}\, $ and check that $\, {\kappa }' \leq \varphi\, ,\, \psi\, $.
Therefore the supremum $\, {\gamma }_{\varphi , \psi }^r =: \varphi \wedge \psi\, $ of $\, C_{\varphi , \psi }^r\, $ is the unique largest element of $\, {\mathfrak X}^{*\mu }_+\, $ majorized by both elements $\, \varphi\, ,\, \psi\, $ but unlike the case of the $\vee $-operation need not be given by the pointwise minimum of $\, \varphi\, $ and $\, \psi\, $ (not even for essential elements $\, x \neq {\bf 1}\, $). The positive definite order admits an order unit given by the norm function (on positive elements) of $\, \mathfrak X\, $, i.e. 
$$ {\bf 1}_{\mu } ( x )\> =\> \Vert x\Vert \> , $$
i.e. $\, {\bf 1}_{\mu }\, $ is the unique maximal $\mu $-state of $\, \mathfrak X\, $.
In particular $\, \mathcal M ( \mathfrak X )\, $ is a function system with unit $\, {\bf 1}_{\mu }\, $. Although $\, {\bf 1}_{\mu }\, $ is contained in the positive cone $\, {\mathfrak X}^{*\mu }_+\subseteq \mathcal M ( \mathfrak X )_+\, $ of monotonous elements it is not an order unit for the corresponding order. Every 
element of $\, \mathcal M ( \mathfrak X )\, $ has a unique minimal decomposition by positive elements, namely if $\, \gamma = \varphi - \psi\, $ is any decomposition by monotonous functionals then putting $\, \overline\varphi = \varphi - \varphi \wedge \psi\, $ and $\, \overline\psi = \psi - \varphi \wedge \psi\, $ these are positive definite functions with $\, \gamma = \overline\varphi - \overline\psi\, $, orthogonal on essential elements, hence the decomposition is necessarily unique and minimal and satisfies the relation $\, \Vert \gamma\Vert = \max\, \bigl\{ \Vert \overline\varphi \Vert\, ,\, \Vert \overline\psi\Vert \bigr\}\, $. By continuity the decomposition extends to general elements of $\, \mathcal M ( \mathfrak X )\, $. From this one infers that $\, \mathcal M ( \mathfrak X )\, $ is a function lattice, where writing $\, {\gamma }_1 = {\varphi }_1 - {\psi }_1\, ,\, {\gamma }_2 = {\varphi }_2 - {\psi }_2\, $ with $\, {\varphi }_i\, ,\, {\psi }_i\in {\mathfrak X}^{*\mu }_+\, $ for $\, i = 1 , 2\, $ the lattice operations are given by the formulas 
$$ {\gamma }_1 \vee {\gamma }_2\> =\> \bigl( {\varphi }_1 + {\varphi }_2 \bigr)\> -\> \bigl[ \bigl( {\varphi }_1 + {\psi }_2 \bigr) \wedge \bigl( {\psi }_1 + {\varphi }_2 \bigr) \bigr]\> ,  $$
$$ {\gamma }_1 \wedge {\gamma }_2\> =\> \bigl[ \bigl( {\varphi }_1 + {\psi }_2 \bigr) \wedge \bigl( {\varphi }_2 + {\psi }_1 \bigr) \bigr]\> -\> \bigl( {\psi }_1 + {\psi }_2 \bigr) $$
showing by normcompleteness that $\, \mathcal M ( \mathfrak X )\, $ isometrically embeds as a closed subspace of bounded continuous functions on the set of essential elements of $\, \mathfrak X\, $ (in fact this embedding is a complete $\vee $-map  and a complete $\wedge $-map restricted to each convex subspace $\, {\mathfrak X}^{*\mu }_{+ , r}\, $ since for these the lattice operations coincide with the usual pointwise lattice operations of bounded functions on essential elements).  On the subset of extended essential elements $\, {\mathfrak X}_{ess} \cup \{ 0 \} = \bigl\{ x\in A_+\bigm\vert \Vert x\Vert = \Vert {\bf 1} - x\Vert = 1\bigr\} \cup \{ 0 , {\bf 1} \}\, $ there is a natural involution $\, x \mapsto x^{o} = {\bf 1} - x\, $ inducing a linear (nonpositive) involution on the abelian $W^*$-algebra of bounded functions $\, B ( {\mathfrak X}_{ess} )\, $ on the set $\, {\mathfrak X}_{ess}\, $ by 
$$ \omega\> \mapsto \widehat\omega\> ,\quad \widehat\omega ( x )\> =\> \omega ( {\bf 1} )\> -\> \omega \bigl( {\bf 1} - x \bigr)\> . $$
This involution is seen to restrict to a (nonmonotonous) involution of  
$\, {\mathfrak X}^{*\mu }_+\subseteq B ( {\mathfrak X}_{ess} )_+\, $. 
We let $\, {\mathfrak X}^{*\mu , o}_+ \subseteq {\mathfrak X}^{*\mu }_+\, $ denote the subset of elements invariant under 
the involution and note that it contains the subset of positive linear functionals $\, {\mathfrak X}^*_+\, $. Also note that $\, \widehat{{\bf 1}}_{\mu } = {\eta }_1\, $ with $\, {\eta }_1\, $ the minimal element of $\, {\mathfrak X}^{*\mu }_{+ , 1}\, $ of norm one and that the involution is generally order reversing on each subset 
$\, {\mathfrak X}^{*\mu }_{+ , r} = \bigl\{ \xi\in {\mathfrak X}^{*\mu }_+\bigm\vert \xi ( {\bf 1} ) = r \bigr\}\, $. A $\mu $-map $\, \sigma :  {\mathfrak X}_+ \rightarrow {\mathfrak Y}_+\, $ will be called a ${\mu }^{o}$-map iff $\, \sigma \bigl( x \bigr) = \widehat\sigma \bigl( x \bigr) = \Vert x\Vert \sigma \bigl( {\bf 1} \bigr) - \sigma \bigl( \Vert x\Vert {\bf 1} - x \bigr)\, $ for any positive element $\, x\in {\mathfrak X}_+\, $.
Every $\, x\in {\mathfrak X}_+\, $ defines a normal positive linear functional of $\, {\mathfrak X}^{*\mu }_+\, $ by duality.
Also note that every element $\, \gamma\,\in {\mathfrak X}^{*\mu }\subseteq \mathcal M ( \mathfrak X )\, $ admits a minimal decomposition by elements  $\, \varphi \in{\mathfrak X}^{*\mu }_{+ , r}\, ,\, \psi\in {\mathfrak X}^{*\mu }_{+ , s}\, $ for given $\, r\, ,\, s\in {\mathbb R}_+\, $ if any, namely given two such decompositions 
$\, \gamma = {\varphi }_1 - {\psi }1- = {\varphi }_2 - {\psi }_2\, $ one also has 
$\, \gamma = {\varphi }_1\wedge {\varphi }_2 - {\psi }_1\wedge {\psi }_2\, $ which follows by duality from the formulas $\, {\gamma }^{o} = {\varphi }_1^{o} - {\psi }_1^{o} = {\varphi }_2^{o} - {\psi }_2^{o} \Longrightarrow {\gamma }^{o} = {\varphi }_1^{o} \vee {\varphi }_2^{o} - {\psi }_1^{o}\vee {\psi }_2^{o}\, $ using the fact that $\,\vee\, $ is the pointwise maximum on essential elements, 
whence by induction and monotone completeness of $\, {\mathfrak X}^{*\mu }_+\, $ one derives
$$ \gamma\> =\> \inf_{\psi\in {\mathfrak X}^{*\mu }_+}\, \Bigl\{ \varphi\bigm\vert \gamma = \varphi - \psi\, ;\, \varphi \in {\mathfrak X}^{*\mu }_{+ , r} \Bigr\}\> -\> \inf_{\varphi\in {\mathfrak X}^{*\mu }_+}\, \Bigl\{ \psi\bigm\vert \gamma = \varphi - \psi\, ;\, \psi\in {\mathfrak X}^{*\mu }_{+ , s} \Bigr\} $$ 
which decomposition is clearly minimal and unique (compare with the corresponding argument for a minimal decomposition by basic elements in $\, \mathfrak L ( \mathfrak X )\, $). 
Equipped with the $w^*$-topology induced by functionals corresponding to elements of $\, {\mathfrak X}_+\, $ (and linear combinations thereof) the intersection of the unit ball with $\, {\mathfrak X}^{*\mu }_+\, $ (in particular $\, {\mathcal S}^{\mu } ( \mathfrak X )\, $) is compact (compare with \cite{Pe2}, Theorem 2.5.2). 
\par\smallskip\noindent
The definition of a $\mu $-functional may be generalized to the concept of a {\it sub-$\mu $-functional} which is supposed to mean a positively homogenous and monotonous map $\, \xi : {\mathfrak X}_+ \rightarrow {\mathbb R}_+\, $ satisfying $\, \xi \bigl( x + \alpha\, {\bf 1} \bigr) \leq \xi \bigl( x \bigr) + \alpha\, \xi \bigl( {\bf 1} \bigr)\, $ (with $\, \xi \bigl( 0 \bigr) = 0\, $). Again the space of sub-$\mu $-functionals denoted $\, {\mathfrak X}^{*s\mu }_+ = \bigcup_r\, {\mathfrak X}^{*s\mu }_{+ , r}\, $ with unique minimal element $\, {\eta }_r\in {\mathfrak X}^{*s\mu }_{+ , r}\, $ for each $\, r\geq 0\, $ as above can be seen to be a complete lattice for the pointwise $\vee $-operation on positive elements and the $\wedge $-operation given by taking the (pointwise) supremum over all sub-$\mu $-functionals smaller than both $\, \varphi\, ,\, \psi\, $ , i.e. 
$$ \varphi \wedge \psi\> =\> \sup\, \Bigl\{ \rho\in {\mathfrak X}^{*s\mu }_+ \Bigm\vert \rho \leq \varphi\, ,\, \rho \leq \psi \Bigr\} $$
noting that the subspace of minorants is not empty. The arguments are analogous to the ones above.
\par\smallskip\noindent
A {\it convex metric space} is a convex space $\, X\, $ equipped with a metric $\, d\, $ such that each segment $\, [\, x\, ,\, y\, ] := \{\, ( 1 - \lambda )\, x\, +\, \lambda\, y\,\vert\, 0 \leq \lambda\leq 1\, \}\subseteq X\, $ is geodesic, and putting $\, \overline{xy}_{\lambda } = \lambda\, x\, +\, ( 1 - \lambda )\, y\, $ one has
$$ d \bigl( \overline{xz}_{\lambda }\, ,\, \overline{yz}_{\lambda } \bigr)\>\leq\> \lambda\,  d ( x , y )  $$
for any three points $\, x\, ,\, y\, ,\, z\in X\, $ and $\, 0\leq \lambda\leq 1\, $.  For $\, z = x\, $ this implies 
$$ d ( x , \overline{xy}_{\lambda } )\> =\> ( 1 - \lambda )\, d ( x , y )\> ,\quad \forall\> x\, ,\, y \in X\, ,\> 0\>\leq\lambda\leq\> 1  $$
from the triangle inequality which then also implies the relation
$$ d\bigl( \overline{xy}_{\lambda }\, ,\, z \bigr)\>\leq\> \lambda\, d ( x , z )\> +\> ( 1 - \lambda ) d ( y , z )\> . $$
Moreover, whenever $\, d ( x , y )\> >\> \max\, \bigl\{ d ( x , z )\, ,\, d ( y , z ) \bigr\}\, $ there exists $\, 0 < \lambda < 1\, $ with
$$ d ( \overline{xy}_{\lambda } , z )\> <\> \min\, \bigl\{ d ( x , z )\, ,\, d ( y , z ) \bigr\} \> . $$
$\, X\, $ is called a {\it strictly convex metric space} iff in addition for any $\, x , y , z\in X\, ,\, x \neq y\, $ and $\, 0 < \lambda < 1 $ one has
$$ d ( \overline{xy}_{\lambda } , z )\> <\> \max\, \bigl\{ d ( x , z )\, ,\, d ( y , z ) \bigr\} \> . $$
\par\noindent
The {\it geometric center} of a compact convex subset $\, C \subseteq X\, $ is defined as follows: For any given $\, r > 0\, $ choose a maximal (finite) collection of points 
$\, F = \{ c_1 , \cdots , c_n \} \subseteq C\, $ such that the corresponding collection of open balls 
$\, \bigl\{\, U_r ( c_k ) \bigm\vert k = 1 , \cdots , n\, \bigr\}\, $ of radius $\, r\, $ are pairwise disjoint. Then consider the element 
$$ c_{F , r}\> =\> \sum_{k = 1}^n\, {1\over n}\, c_k\> \in C \> . $$
The net $\, \{ c_{F , r} \}\, $ where $\, ( F , r ) \geq ( G , s )\, $ iff $\, r \leq s\, $ and 
$\, G \subseteq F\, $ is easily seen to converge to a unique limit point $\, z ( C )\in C\, $ which point is called the geometric center of $\, C\, $. This is because given any two finite sets $\, \bigl\{ x_1\, ,\,\cdots\, ,\, x_n \bigr\}\, ,\, \bigl\{ y_1\, ,\, \cdots\, ,\, y_n \bigr\} \subseteq X\, $ in a convex metric space such that $\, d ( x_k , y_k ) \leq \epsilon\, $ for each $\, k = 1 ,\cdots , n\, $ and any positive partition 
$\, 1 = \sum_{k=1 }^n {\lambda }_k\, ,\, {\lambda }_k\geq 0\, $ one proves by induction on $\, n\, $ that the relation $\, d ( x , y ) \leq \epsilon\, $ holds where $\, x = \sum_{k = 1}^n {\lambda }_k\, x_k\, $ and $\, y = \sum_{k = 1}^n {\lambda }_k\, y_k\, $.
Given an affine isometric embedding of convex metric spaces $\, i : {\Omega }_0 \hookrightarrow \Omega\, $ with $\, {\Omega }_0\, $ compact there exists a canonical retraction $\, r : \Omega \twoheadrightarrow {\Omega }_0\, $ for $\, i\, $ by defining 
$$ r ( x )\> =\>  z \bigl( C_x \bigr)\> , \quad C_x = \bigl\{ y\in {\Omega }_0 \bigm\vert d ( x , y ) = d \bigl( x , {\Omega }_0 \bigr) \bigr\}\subseteq {\Omega }_0  $$
where $\, z ( C_x )\, $ denotes the geometric center of the compact convex subset of points $\, C_x\subseteq {\Omega }_0\, $ minimizing the distance to $\, x\, $. Since $\, {\Omega }_0\, $ is a compact convex metric space the minimum of the distance function is attained and $\, C_x \neq \emptyset\, $ is a nonempty compact convex subset of $\, {\Omega }_0\, $. If the metric is strictly convex $\, C_x\, $ will consist of precisely one point since if there would exist two different such points $\, \{ y_0 , y_1 \}\, $ then $\, C_x\, $ would also contain the whole segment $\, [ y_0 , y_1 ]\, $ and there would exist a point $\, y_{\lambda }\in [ y_0 , y_1 ]\, $ with $\, d ( x , y_{\lambda } ) < d ( x , y_0 )\, $ by strict convexity contradicting minimality of $\, d ( x , y_0 )\, $. Therefore passing to the geometric center of $\, C_x\, $ is redundant in this case. It is then also easy to see (for a strictly convex metric) that the retraction $\, r\, $ is continuous. This result also holds in the case of a convex metric but is a bit more difficult to see, one has to approximate the geometric center of the compact subset minimizing the distance to a given point by the geometric centers of corresponding compact subsets of points whose distance is $\epsilon $-close to the minimal distance, we leave the precise argument to the reader. 
But more is true. If $\, x\in \Omega\, ,\, y\in {\Omega }_0\, $ then 
$\, d ( r ( x ) , y ) \leq d ( x , y )\, $ so that $\, r\, $ is partially contractive (for arbitrary convex metric $\, d\, $), for otherwise $\, d ( r ( x ) , y ) > d ( x , y ) \geq d ( x , r ( x ) )\, $ implying that there exists a point $\, z\in [ y , r ( x ) ]\subseteq {\Omega }_0\, $ with $\, d ( x , z ) < d ( x , r ( x ) )\, $, a contradiction.
\par\noindent
One checks that the Tychonoff cube $\, T = {\prod}_{n\in\mathbb N}\, [ 0 , 1 ]\, $ equipped with the metric 
$\,  d ( x , y ) =  {\sum}_n\, 2^{-n}\, \vert x_n - y_n \vert\, $ is an example of a convex metric space. If instead one considers the topologically equivalent metric 
$$  \widetilde d ( x , y )\> =\>\sup_n\, \sqrt{ \sum_{k = 1}^n 2^{- k}\, \vert x_k - y_k{\vert }^2} $$
$\, T\, $ turns into a strictly convex metric space. One also notes that the convex metric $\, d\, $ above can be pointwise approximated by the field of strictly convex metrics $\, {\widetilde d}_t ( x , y ) = t\, d ( x , y ) + ( 1 - t )\, \widetilde d ( x , y )\, $ for $\, t\to 1\, $.
The disadvantage of $\, d\, $ not being strictly convex is compensated by the fact that it is of product type, i.e. of the form $\, d = \sum_n\, d_n\, $ with each $\, d_n\, $ a convex metric on the quotient $\, [ 0 , 1 ]\, $ corresponding to the $n$-th component. We prefer to use product type convex metrics which are well behaved with respect to lattice operations.
It follows that if $\, \mathfrak X\, $ is a separable function system then its state space $\, \mathcal S ( \mathfrak X )\, $ (resp. $\, {\mathcal S}^{\mu } ( \mathfrak X )\, $) endowed with the relative $w^*$-topology can be given the structure of a compact (strictly) convex metric space. Indeed choosing a normdense subsequence $\, \{ x_n\,\vert\, x_n \geq 0\, ,\, \Vert x_n\Vert = 1\, ,\, \forall\> n\in\mathbb N \}\subseteq {\mathfrak X}^+_1\, $ of positive elements of norm $1$ consider the sequence of continuous affine maps 
$$ {\widehat x}_n :\> \mathcal S ( \mathfrak X )\> \largerightarrow\> [ 0\, ,\, 1 ]\> ,\quad {\widehat x}_n ( \rho )\> =\> \rho ( x_n ) $$ 
which separate the points of $\, \mathcal S ( \mathfrak X )\, $, whence the product map into the Tychonoff cube 
$$ \prod_n\, {\widehat x}_n :\> \mathcal S ( \mathfrak X )\> \largerightarrow\>  {\prod}_{n\in\mathbb N} [ 0 , 1 ]
$$ 
is injective, affine and continuous, and the image of $\, \mathcal S ( \mathfrak X )\, $ with the induced topology is compact and Hausdorff, thus homeomorphic to $\, \mathcal S ( \mathfrak X )\, $ with relative $w^*$-topology by rigidity, so that $\, \mathcal S ( \mathfrak X )\, $ inherits the structure of a compact (strictly) convex metric space (compare with \cite{Pe2}, Proposition 1.6.14), the same for 
$\, {\mathcal S}^{\mu } ( \mathfrak X )\, $. In fact the embedding as above can be extended to a (not necessarily injective) continuous linear map of $\, \mathcal M ( \mathfrak X )\, $ into the extended Tychonoff space of bounded realvalued sequences with metric as above sending positive definite elements to positive definite sequences.
Then the ensuing topology of $\, {\mathfrak X}^{*\mu }\, $ pulled back from the image of this map coincides with the $w^*$-topology on bounded subsets of $\, {\mathfrak X}^{*\mu }_+\, $. 
\par\noindent
If $\, B\, $ is an abelian $\mathcal P$-algebra then a {\it $\mu $-integral map} 
$\, \psi : B \twoheadrightarrow \mathfrak X\, $ into a function system $\, \mathfrak X\, $ is a $\mu $-map which arises by $\mu $-integration of a ($\, {\mathfrak X}_+$-valued) $\mu $-measure $\, {\psi }_{\mathcal P} : \mathcal P ( B ) \rightarrow {\mathfrak X}_+\, $, i.e. a monotonous map of the subset of projections of $\, B\, $ into the positive elements of $\, \mathfrak X\, $ so that if given a monotone decreasing string of projections $\, \{ p_k\,\vert\, p_{k + 1} < p_k \}\, $ the value of 
$\, b = \sum_{k = 1}^n {\alpha }_k p_k\, $ is given by $\, \psi ( b ) = \sum_{k = 1}^n {\alpha }_k {\psi }_{\mathcal P} ( p_k )\, $.
Let $\, \Psi ( B )\subseteq {\mathcal S}^{\mu }_{\Psi} ( B ) \subseteq {\mathcal S}^{\mu } ( B )\, $ denote the subset of $\mathcal P$-functionals of $\, B\, $ and its closed convex hull respectively. Note that any $\mathcal P$-functional $\, \sigma\, $ is an extreme point of $\, {\mathcal S}^{\mu } ( B )\, $ since assuming that $\, \sigma = \lambda\, \varphi + ( 1 - \lambda )\, \psi\, $ for some $\mu $-states $\, \varphi\, ,\, \psi\, $  one  necessarily has $\, \varphi ( q ) = \psi ( q ) = \{ 0 , 1 \}\, $ for any projection $\, q\in B\, $ showing that both $\, \varphi\, $ and $\, \psi \, $ must be $\mathcal P$-functionals and equal to the unique monotonous extension of their restrictions to the subset of projections, i.e. they must both agree with $\, \sigma\, $. Then the same holds for the images of $\mathcal P$-functionals in $\, {\Sigma }_B\, $. Also note that the $w^*$-closure $\, {\mathcal S}^{\mu }_{\mathcal P} ( B )\, $ of $\, {\mathcal S}^{\mu }_{\Psi } ( B )\, $ contains all (linear) states of $\, B\, $ (any state $\, \rho\in \mathcal S ( B )\, $ is in the $w^*$-closure of the convex hull of the multiplicative functionals from the Krein-Milman theorem which are precisely the $\mathcal P$-$A$-lattice-functionals).
\par\noindent
An affine cone $\, \mathfrak C \subseteq \mathfrak X\, $ in a real vector space $\, \mathfrak X\, $ is a subset closed under addition of elements such that $\, r\, x\in \mathfrak C\, $ for each $\, x\in \mathfrak C\, ,\, r\in {\mathbb R}_+\, $. 
A pair of affine subcones $\, \bigl( \mathfrak C\, ,\, \mathfrak D \bigr) \subseteq \mathfrak X\, $ is called {\it exhausting} iff $\, \mathfrak C - \mathfrak D = \mathfrak X\, $ and {\it generating} if $\, \mathfrak X\, $ is the linear hull of $\, \mathfrak C \cup \mathfrak D\, $.  A convex generalized Minkowski functional $\, m : \mathfrak C \rightarrow \mathbb R \cup \{ \infty \}\, $ is a positively homogenous and sublinear map into the upwards extended real line, i.e. $\, m \bigl( r\, x \bigr) = r\, m \bigl( x \bigr)\, ,\, m \bigl( x + y \bigr) \leq m \bigl( x \bigr) + m \bigl( y \bigr)\, $ with the convention $\, 0\cdot\infty = 0\, ,\, r\cdot\infty = \infty\, $ for $\, r > 0\, $. A concave generalized Minkowski functional $\, \nu : \mathfrak D \rightarrow \mathbb R \cup \{ - \infty \}\, $ is a positively homogenous suplinear map into the downwards extended real line, i.e. $\, \nu \bigl( r\, x \bigr) = r \, \nu \bigl( x \bigr)\, ,\, \nu \bigl( x + y \bigr) \geq \nu \bigl( x \bigr) + \nu \bigl( y \bigr)\, $ with the convention $\, o\cdot ( - \infty ) = 0\, ,\, r\cdot ( - \infty ) = - \infty\, $ for $\, r > 0\, $. The following gives a version of the Hahn-Banach theorem suited for our purposes exhibiting its essentially algebraic nature. 
The Hahn-Banach extension theorem is the magic ingredient algebra donates to functional analysis. In particular it should be noted that the theorem works the same if replacing real vector spaces by $\mathbb Q$-vector spaces (or even free $\mathbb Z$-modules, compare with the Extension theorem for $A$-linear $\mathcal P$-measures given below).
\par\bigskip\noindent
{\bf Affine Hahn-Banach extension theorem.}\quad 
Let $\, \bigl( \mathfrak C\, ,\, \mathfrak D \bigr) \subseteq \mathfrak X\, $ be an exhausting pair of affine subcones in a real vector space and assume given a real subspace $\, \mathfrak Y\subseteq \mathfrak X\, $, together with a real linear functional $\, \varphi : \mathfrak Y \rightarrow \mathbb R\, $ and a  convex (generalized) Minkowski functional $\, m : \mathfrak C \rightarrow \mathbb R \cup \{ \infty \}\, $, a concave (generalized) Minkowski functional $\, \nu : \mathfrak D \rightarrow \mathbb R \cup \{ - \infty \}\, $ satisfying $\, \nu {\bigm\vert }_{\mathfrak C \cap \mathfrak D} \leq m {\bigm\vert }_{\mathfrak C \cap \mathfrak D}\, $ plus the condition  
$$ m ( x )\> -\> \nu ( y )\> \geq\> \phi ( x - y ) \leqno{(*)}   \> $$
\par\noindent
whenever $\, x - y\in \mathfrak Y\, $. Then $\,\varphi\, $ admits  an $\mathbb R$-linear extension $\, \varphi : \mathfrak X \rightarrow \mathbb R\, $ with $\, \nu ( y ) \leq \varphi ( y )\, $ for all $\, y\in \mathfrak D\, $ and $\, \varphi ( x ) \leq m ( x )\, $ for all $\, x\in \mathfrak C\, $.
\par\bigskip\noindent 
{\it Proof.}\quad Consider the (generalized) convex Minkowski functional 
$\, \widetilde m : \mathfrak X \rightarrow \mathbb R \cup \{ \infty \}\, $  given by 
$$ \widetilde m \bigl( x \bigr)\> =\> \inf\, \Bigl\{\> m \bigl( y \bigr)\> -\> \nu \bigl( z \bigr) \Bigm\vert\> x\> =\> y \> -\> z \Bigr\}  $$
so that $\, \varphi \bigl( y \bigr) \leq \widetilde m \bigl( y \bigr)\, $ fpr $\, y\in \mathfrak Y\, $ from the $(*)$-condition with $\, \widetilde m \bigl( x \bigr) \leq m \bigl( x \bigr)\, $ for $\, x\in \mathfrak C\, $ and $\, \nu \bigl( y \bigr) \leq - \widetilde m \bigl( - y \bigr)\, $ for $\, y\in \mathfrak D\, $. From \cite{Pe2}, Lemma 2.3.2. there exists an $\mathbb R$-linear extension of $\, \varphi\, $ to all of $\,\mathfrak X\, $ dominated by $\, \widetilde m\, $\qed
\par\bigskip\noindent
If $\, \mathfrak D = \mathfrak C = \mathfrak X\, $ and putting $\, \nu ( x ) = - m ( - x )\, $ the condition $( * )$ is automatically satisfied which is the case usually considered in the books, cf. \cite{Pe2}, Lemma 2.3.2.  Thus both versions are logically equivalent and can be derived from each other. In case of a nonexhausting pair of subcones or generalized Minkowski functionals there might be a slight problem ensuring that $\, \widetilde m \bigl( x \bigr)\, $ as above is well defined. This can be remedied under the mild assumption that there exist universal constants $\, C\, ,\, D \in \mathbb R\, $ with $\, m \bigl( x \bigr) \geq C\, ,\, \nu \bigl( y \bigr) \leq D\, $ which we assume without further mentioning   Note that condition $( * )$ is trivially satisfied for $\, v\, ,\, w\in \mathfrak Y\, $. It is also satisfied in the following cases: either if $\, x - y\in\mathfrak C\, $ with $\, m\, $ $\mathfrak C$-monotonous and $\, \nu\, $ negative, or else if $\, y - x\in \mathfrak D\, $ with $\, - \nu\, $ $\mathfrak D$-monotonous and $\, m\, $ positive.
\par\bigskip\noindent
Let $\, \mathfrak X\, $ be an operator system and $\, \mathcal P \subseteq \mathcal P ( A )\, $ be an $A$-sublattice of the subset of projections of a $C^*$-algebra $\, A\, $. An {\it $\mathfrak X$-valued $\mathcal P$-measure} on $\, \mathcal P\, $ is a monotonous map $\, s : \mathcal P \rightarrow \mathfrak X\, $ with $\, s \bigl( 0 \bigr) = 0\, $, i.e. $\, p \geq q \Longrightarrow s \bigl( p \bigr) \geq s \bigl( q \bigr) \geq 0\, $. A $\mathcal P$-measure $\, s\, $ is {\it $C$-contractive} for some given $\, C \geq 0\, $ iff for any $\, \epsilon > 0\, $ and two given projections $\, p\, ,\, q\in \mathcal P\, $ with $\, \Vert p - q \Vert \leq \epsilon\, $ this implies $\, \Vert s \bigl( p \bigr) - s \bigl( q \bigr) \Vert \leq C\, \epsilon\, $.  If $\, C = s \bigl( {\bf 1} \bigr)\, $ then $\, s\, $ is called {\it selfcontractive}. The space of $C$-contractive realvalued $\mathcal P$-measures is denoted $\, {\mu }^C_{\mathcal P}\, $ and the space of selfcontractive realvalued $\mathcal P$-measures with $\, s \bigl( {\bf 1} \bigr) = C\, $ is denoted $\, {\mu }^C_{\mathcal P , C}\, $
To any $\mathbb R$-valued $\mathcal P$-measure 
$\, s : \mathcal P ( A ) \rightarrow \mathbb R\, $ one can associate a monotonous and positively homogenous functional on $\, A\, $ extending $\, s\, $ in several ways. The minimal such extension $\, {\nu }_{s , 0} : A_+ \rightarrow \mathbb R\, $  is given by
$$ {\nu }_{s , 0} \bigl( x \bigr)\> =\> \sup\, \Bigl\{\> \alpha\, s \bigl( p \bigr) \Bigm\vert \alpha\, p \leq x\> ,\quad \alpha\in {\mathbb R}_+\> ,\quad p\in \mathcal P ( A )\> \Bigr\} \> , $$
another extension $\, {\nu }_s : A_+ \rightarrow \mathbb R\, $  is given by 
$$ {\nu }_s \bigl( x \bigr)\> =\> \sup\, \left\{ \sum_{k = 1}^n\, {\alpha }_k\, s \bigl( p_k \bigr) \Biggm\vert \sum_{k = 1}^n\, {\alpha }_k\, p_k\>\leq\> x\> ,\> \bigl\{ {\alpha }_k {\bigr\}}_{k = 1}^n\subseteq {\mathbb R}_+\, ,\, \bigl\{ p_k {\bigr\}}_{k = 1}^n \subseteq \mathcal P ( A ) \right\}  $$
where the supremum is taken over all finite descending chains of projections $\, p_1 \geq p_2 \geq \cdots \geq p_n\, $ in $\, A\, $. If $\, A\, $ is a $\mathcal P$-algebra then clearly $\, {\nu }_s \bigl( x \bigr) \geq \sigma \bigl( x \bigr)\, $ for any $\, x\in {\mathcal P}_f ( A )\, $ where $\, \sigma\, $ denotes the canonical $\mathcal P$-integral extension of $\, s\, $ to $\, {\mathcal P}_f ( A )\, $.
If $\, s\, $ is $C$-contractive then another $C$-contractive monotonous functional $\, {\mu }_s^C : A_* \rightarrow {\mathbb R}_+\, $ extending $\, s\, $ is defined by
$$ {\mu }_s^C \bigl( x \bigr)\> =\> \inf\, \Bigl\{ \alpha\, s \bigl( r \bigr)\> +\> \beta\, C\, \Bigm\vert\, \alpha\, r\> +\> \beta\, {\bf 1}\> \geq\> x\> ,\quad 0\> \leq\> \alpha\, ,\, \beta\> ,\quad r\in \mathcal P ( A ) \Bigr\} \> . $$
In particular if $\, s\, $ is selfcontractive then $\, {\mu }_s = {\mu }_s^{s( {\bf 1} )}\, $ is a $\mu $-functional, i.e. one has $\, {\mu }_s \bigl( x + t\, {\bf 1} \bigr) = {\mu }_s \bigl( x \bigr) + t\, s \bigl( {\bf 1} \bigr)\, $.
 If $\, \mathcal J \vartriangleleft A \twoheadrightarrow Q\, $ is an exact sequence of $C^*$-algebras any $C$-contractive realvalued $\mathcal P$-measure on $\, A\, $ defines a $C$-contractive $\mathcal P$-measure on $\, \mathcal J\, $, and any realvalued $C$-contractive $\mathcal P$-measure on $\, Q\, $ defines a realvalued $C$-contractive $\mathcal P$-measure on $\, A\, $ which is trivial on (congruent modulo) $\, \mathcal J\, $ by evaluation. To the contrary there exist two complementary affine projections $\, {\pi }_{\mathcal J} :  {\mu }^C_{\mathcal P ( A )} \rightarrow 
 {\mu }^C_{\mathcal P ( A )}\, ,\, {\pi }_Q : {\mu }_{\mathcal P ( A )}^C \rightarrow {\mu }^C_{\mathcal P ( A )}\, $ such that $\, {\pi }_Q \bigl( {\mu }_{\mathcal P ( A )}^C \bigr)\, $ is contained in the image of $\, {\mu }_{\mathcal P ( Q )}^C\, $ satisfying the relations $\, {\pi }_{\mathcal J} \bigl( {\pi }_{\mathcal J} \bigl( s \bigr) \bigr) = {\pi }_{\mathcal J} \bigl( s \bigr)\, ,\, {\pi }_Q \bigl( {\pi }_Q \bigl( s \bigr) \bigr) = {\pi }_Q \bigl( s \bigr)\, $ and if $\, s\, $ is $A$-linear (see below) the additional relations $\, {\pi }_Q \bigl( {\pi }_{\mathcal J} \bigl( s \bigr) \bigr) = 0 = {\pi }_{\mathcal J} \bigl( {\pi }_Q \bigl( s \bigr) \bigr)\, $ and $\, s = {\pi }_{\mathcal J} \bigl( s \bigr) + {\pi }_Q \bigl( s \bigr)\, $. In particular considering selfcontractive $A$-linear $\mathcal P$-measures one has $\, {\pi }_{\mathcal J} \bigl( {\mu }^C_{\mathcal P ( A ) , C} \bigr) \subseteq {\mu }^{C_0}_{\mathcal P ( A ) , C_0}\, $ and $\, {\pi }_Q \bigl( {\mu }^C_{\mathcal P ( A ) , C} \bigr) \subseteq {\mu }^{C_1}_{\mathcal P ( A ) , C_1}\, $ with $\, C_0 + C_1 = C\, $. These maps are defined in the following way: for each projection $\, p\in \mathcal P ( A )\, $ let $\, \bigl\{ u_{\lambda }^p {\bigr\}}_{\lambda }\subseteq \mathcal P ( \mathcal J )\, $ denote a projection unit of $\, \mathcal J\, $ commuting with $\, p\, $. Then for $\, s\in {\mu }^C_{\mathcal P ( A )}\, $ define 
 $$ {\pi }_{\mathcal J} \bigl( s \bigr) \bigl( p \bigr)\> =\> \limsup_{\lambda \to \infty}\, s \bigl( u_{\lambda }^p\, p \bigr)\> ,\quad {\pi }_Q \bigl( s \bigr) \bigl( p \bigr)\> =\> \liminf_{\lambda\to\infty }\, s \bigl( ( {\bf 1} - u_{\lambda }^p )\, p \bigr) \> . $$
 One needs to show that this definition does not depend on the chosen projection unit. Let $\, \bigl\{ v_{\mu }^p {\bigr\}}_{\mu } \subseteq \mathcal P ( \mathcal J )\, $ be another projection unit commuting with $\, p\, $. Then for any fixed index $\, \lambda\, $ one may replace $\, \bigl\{ v_{\mu } {\bigr\}}_{\mu } \, $ by an equivalent projection unit $\, \bigl\{ v'_{\mu } {\bigr\}}_{\mu }\, $ exactly commuting with $\, u_{\lambda }^p\, $ and $\, p\, $ where equivalent is supposed to mean that $\, \Vert v_{\mu }^p - (v'_{\mu })^p \Vert \to 0\, $ for $\, \mu\to \infty\, $. Clearly equivalent projection units yield the same result. Therefore 
 $$ s \bigl( u_{\lambda }^p\ p \bigr)\>= \> \limsup_{\mu\to\infty }\, s \bigl( u_{\lambda }^p\, (v'_{\mu })^p\, p \bigr)\> \leq\> \limsup_{\mu\to\infty }\, s \bigl( v_{\mu }^p\, p \bigr) $$
 implying $\, \limsup_{\lambda\to\infty }\, s \bigl( u_{\lambda }^p\, p \bigr) \leq \limsup_{\mu\to\infty }\, s \bigl( v_{\mu }^p\, p \bigr)\, $ and the reverse relation follows by symmetry.  A similar argument applies for $\, {\pi }_Q\, $. The fact that $\, {\pi }_{\mathcal J} \bigl( s \bigr)\, $ and $\, {\pi }_Q \bigl( s \bigr)\, $ define (monotonous) $\mathcal P$-measures is immediate from definition since if $\, p \leq q\, $ one can find a projection unit in $\, \mathcal J\, $ commuting with both $\, p\, $ and $\, q\, $. The relations $\, {\pi }_{\mathcal J} \bigl( {\pi }_{\mathcal J} \bigl( s \bigr) \bigr) = {\pi }_{\mathcal J} \bigl( s \bigr)\, ,\, {\pi }_Q \bigl( {\pi }_Q \bigl( s \bigr) \bigr) = {\pi }_Q \bigl( s \bigr)\, $ when defined are now obvious. In case of an $A$-linear $\mathcal P$-measure one also has $\, = {\pi }_Q \bigl( {\pi }_{\mathcal J} \bigl( s \bigr) \bigr) = 0 = {\pi }_{\mathcal J} \bigl( {\pi }_Q \bigl( s \bigr) \bigr)\, $ and $\, s = {\pi }_{\mathcal J} \bigl( s \bigr) + {\pi }_Q \bigl( s \bigr)\, $. To see that $\, {\pi }_{\mathcal J} \bigl( s \bigr)\, $ and $\, {\pi }_Q \bigl( s \bigr)\, $ are again $C$-contractive let $\, p\, ,\, q\in \mathcal P ( A )\, $ be given with $\, \Vert p - q \Vert \leq \epsilon\, $. Choose projection units $\, \bigl\{ u_{\lambda }^p {\bigr\}}_{\lambda }\, $ commuting with $\, p\, $ and $\, \bigl\{ v_{\mu }^q {\bigr\}}_{\mu }\, $ commuting with $\, q\, $ respectively, and a quasicentral approximate unit $\, \bigl\{ w_{\nu } {\bigr\}}_{\nu } \subseteq {\mathcal J}_+\, $. We first show that 
 $$ \limsup_{\lambda }\, s \bigl( u_{\lambda }^p\, p \bigr) \> =\> \lim_{\nu }\, {\mu }^C_s \bigl( w_{\nu }\, p\, w_{\nu } \bigr)\> =\> \lim_{\nu }\, {\mu }^C_s \bigl( p\, w_{\nu }\, p \bigr) \> . $$
 Clearly $\, {\mu }^C_s \bigl( w_{\nu }\, p\, w_{\nu } \bigr) \geq \limsup_{\lambda }\, {\mu }^C_s \bigl( w_{\nu }\, u_{\lambda }^p\, p\, w_{\nu } \bigr)\, $ so that 
$$ \lim_{\nu }\, {\mu }^C_s \bigl( w_{\nu }\, p\, w_{\nu } \bigr)\> \geq\> \limsup_{\lambda }\, \lim_{\nu }\, {\mu }^C_s \bigl( w_{\nu }\, u_{\lambda }^p\, p\, w_{\nu } \bigr)\> =\> \limsup_{\lambda }\, s \bigl( u_{\lambda }^p\, p \bigr) \> . $$
On the other hand for fixed index $\, \nu\, $ one can find for each $\, \delta > 0\, $ an index $\, {\lambda }_{\nu , \delta }\, $ with $\, \Vert u_{{\lambda }_{\nu , \delta }}\, w_{\nu }\, u_{{\lambda }_{\nu , \delta }} - w_{\nu } \Vert \leq \delta\, $ so that since $\, {\mu }^C_s\, $ is $C$-contractive hence continuous one has 
$$ {\mu }^C_s \bigl( p\, w_{\nu }\, p \bigr)\> \leq \> {\mu }^C_s \bigl( p\, u_{{\lambda }_{\nu , \delta }}\, w_{\nu }\, u_{{\lambda }_{\nu , \delta }}\, p \bigr)\> +\> \delta\>\leq\> s \bigl( u_{{\lambda }_{\nu , \delta }}\, p \bigr) \> +\> \delta $$
implying the reverse relation 
$$ \lim_{\nu }\, {\mu }^C_s \bigl( w_{\nu }\, p\, w_{\nu } \bigr)\> \leq\> \limsup_{\lambda }\, s \bigl( u_{\lambda }^p\, p \bigr) \> . $$ 
Then since $\, \Vert p - q \Vert \leq \epsilon\, $ implies 
$\, \Vert w_{\nu }\, p\, w_{\nu } - w_{\nu }\, q\, w_{\nu } \Vert \leq \epsilon\, $ for each index $\, \nu\, $ the relation 
$$ \bigm\vert {\pi }_{\mathcal J} \bigl( s \bigr) \bigl( p \bigr)\> -\> {\pi }_{\mathcal J} \bigl( s \bigr) \bigl( q \bigr) \bigm\vert\> =\> \bigm\vert \lim_{\nu }\, {\mu }_s^C \bigl( w_{\nu }\, p\, w_{\nu } \bigr)\> -\> \lim_{\nu }\, {\mu }_s^C \bigl( w_{\nu }\, q\, w_{\nu } \bigr) \bigm\vert\> \leq\> C\, \epsilon $$
follows from the fact that $\, {\mu }_s^C\, $ is $C$-contractive and an analogous argument gives that $\, {\pi }_Q \bigl( s \bigr)\, $ is $C$-contractive..
However in case that $\, s\, $ is $A$-linear and selfcontractive whence $\, {\mu }_s = {\mu }_s^{s ( {\bf 1} )}\, $ is a $\mu $-functional this estimate can be improved. Consider the adjoint $\mu $-functional $\, {\widehat\mu }_s : A \rightarrow \mathbb R\, $ given by 
$$ {\widehat\mu }_s \bigl( x \bigr)\> =\> {\bf 1}_{\mu } \bigl( x \bigr)\, s \bigl( {\bf 1} \bigr)\> -\> {\mu }_s \bigl( {\bf 1}_{\mu } \bigl( x \bigr)\, {\bf 1}\> -\> x  \bigr) $$
$$\quad =\> \sup\, \Bigl\{ \Bigl( {\bf 1}_{\mu } \bigl( x \bigr) - \beta \Bigr)\, s \bigl( {\bf 1} \bigr)\> -\> \alpha\, s \bigl( r \bigr) \Bigm\vert \beta\, {\bf 1}\> +\> \alpha\, r\>\geq\> {\bf 1}_{\mu } \bigl( x \bigr)\, {\bf 1}\> -\> x \> ,\quad \alpha\, ,\, \beta\> \geq\> 0 \Bigr\} $$
$$\quad =\> \sup\, \Bigl\{ \beta\, s \bigl( {\bf 1} \bigr)\> -\> \alpha\, s \bigl( r \bigr) \Bigm\vert \beta\, {\bf 1}\> -\> \alpha\, r\> \leq\> x\> ,\quad \alpha\> \geq\> \beta - k_x\>\geq\> 0 \Bigr\} $$
where $\, {\bf 1}_{\mu }\, $ denotes the maximal $\mu $-state of $\, A\, $ restricting to the normfunction on $\, A_+\, $  and $\, k_x : A \rightarrow \mathbb R\, $ denotes the (minimal) $\mu $-state given by 
$$ k_x\> =\> \sup\, \Bigl\{ k\in \mathbb R \Bigm\vert k\, {\bf 1} \leq x \Bigr\} \> . $$
Since $\, s\, $ is $A$-linear $\, {\widehat\mu }_s\, $ is again a (selfcontractive) monotonous extension of $\, s\, $. As above one shows that  
$$ {\pi }_{\mathcal J} ( s ) \bigl( p \bigr)\> =\> \lim_{\nu\to\infty  }\, {\widehat\mu }_s \bigl( w_{\nu }\, p\, w_{\nu } \bigr)\> ,\quad {\pi }_Q ( s ) \bigl( p \bigr)\> =\> \lim_{\nu\to\infty }\, {\widehat\mu }_s \bigl( ( {\bf 1} - w_{\nu } )\, p\, ( {\bf 1} - w_{\nu } ) \bigr)\> .  $$
Thus 
$$ \lim_{\nu\to\infty }\, {\widehat\mu }_s \bigl( w_{\nu }\, x\, w_{\nu } \bigr)\> =\> \lim_{\nu\to\infty }\, \sup\, \Bigl\{ \beta\, s \bigl( {\bf 1} \bigr)\> -\> \alpha\, s \bigl( r \bigr) \Bigm\vert \beta\, {\bf 1}\> -\> \alpha\, r\> \leq\> w_{\nu }\, x\, w_{\nu } 
 \> ,\> \alpha\> \geq\> \beta\>\geq\> 0 \Bigr\} $$
 $$ =\> \lim_{\nu\to\infty }\, \sup\, \Bigl\{ \beta\, s \bigl( {\bf 1} - r \bigr)\> -\> ( \alpha - \beta )\, s \bigl( r \bigr) \Bigm\vert \beta ( {\bf 1} - r )\> -\> ( \alpha - \beta )\, r\>\leq\> w_{\nu }\, x\, w_{\nu }\> ,\> \alpha\>\geq\> \beta \>\geq\> 0 \Bigr\} $$ 
where the condition $\, \alpha \geq \beta \geq 0\, $ in the last two lines follows from $\, k_{w_{\nu }\, x\, w_{\nu }} \leq 0\, $  plus the fact that only values $\, \beta \geq 0\, $ contribute to the supremum. Thus $\, {\bf 1} - r\, $ must be contained in $\, \mathcal J\, $ and on replacing $\, r\, $ with a smaller projection $\, r'\in \mathcal J\, $ approximately satisfying $\, \beta ( {\bf 1} - r ) - \alpha r'\, \leq w_{\nu }\, x\, w_{\nu }\, $ and using continuity the functional $\, x\mapsto \lim_{\nu } {\widehat\mu }_s \bigl( w_{\nu }\, x\, w_{\nu } \bigr)\, $ is seen to be generally smaller than the $\mu $-functional $\, {\widehat\mu }_{s , \mathcal J} : A \rightarrow \mathbb R\, $ given by 
$$ {\widehat\mu }_{s , \mathcal J} \bigl( x \bigr)\> =\> \sup\, \Bigl\{ \beta\, {\pi }_{\mathcal J} ( s ) \bigl( {\bf 1} \bigr)\> -\> \alpha\, {\pi }_{\mathcal J} ( s ) \bigl( r \bigr) \Bigm\vert \beta\, {\bf 1} - \alpha\, r\>\leq\> x\> ,\quad \alpha\>\geq\> \beta\>\geq\> 0 \Bigr\} $$
which being $C_0$-contractive for $\, C_0 = {\pi }_{\mathcal J} ( s ) \bigl( {\bf 1} \bigr)\, $ and given any pair of projections with $\, \Vert p - q \Vert \leq \epsilon\, $ implies the relations
$$ \lim_{\nu\to\infty }\, {\widehat\mu }_s \bigl( w_{\nu }\,  p \, w_{\nu } \bigr)\> \leq\> \lim_{\nu\to\infty }\, {\widehat\mu }_s \bigl( w_{\nu } ( q + \epsilon\, {\bf 1} ) w_{\nu } \bigr)\> =\> \lim_{\nu\to\infty }\,  \lim_{\mu\to\infty }\, {\widehat\mu }_s \bigl( w_{\nu }\, v_{\mu }^q\, ( q + \epsilon\, {\bf 1} )\, w_{\nu } \bigr) $$
$$ \>\leq\> \lim_{\mu\to\infty }\, {\widehat\mu }_{s , \mathcal J} \bigl( v_{\mu }^q ( q + \epsilon\, {\bf 1} ) \bigr)\> =\> {\pi }_{\mathcal J} ( s ) \bigl( q \bigr)\> +\> \epsilon\, C_0\> .  $$
Here we have used the identity $\, {\widehat\mu }_{s , \mathcal J} \bigl( q \bigr) = {\pi }_{\mathcal J} ( s ) \bigl( q \bigr) = s \bigl( q \bigr)\, $ for any projection $\, q\in \mathcal P ( \mathcal J )\, $. To see this consider the adjoint $\mu $-functional given by
$$ {\mu }_{s , \mathcal J} \bigl( x \bigr)\> =\> \Vert x \Vert\, C_0\> -\> {\widehat\mu }_{s , \mathcal J} \bigl( \Vert x \Vert\, {\bf 1} - x \bigr) $$
$$\quad =\> \inf\, \Bigl\{ \beta\, C_0\>+\> \alpha\, {\pi }_{\mathcal J} ( s ) \bigl( r \bigr) \Bigm\vert \beta\, {\bf 1}\> +\> \alpha\, r \geq x\> ,\quad \alpha\, ,\, \beta\>\geq\> 0 \Bigr\} $$
$$\quad =\> \lim_{\kappa\to \infty }\, \inf\, \Bigl\{ \beta\, s \bigl( v_{\kappa }^r \bigr)\> +\> \alpha\, s \bigl( v_{\kappa }^r\, r \bigr) \Bigm\vert \beta\, {\bf 1}\> +\> \alpha\, r\>\geq\> x \Bigr\} $$ 
where $\, \bigl\{ v_{\kappa }^r {\bigr\}}_{\kappa } \subseteq \mathcal P ( \mathcal J )\, $ is a projection unit commuting with $\, r\, $. Inserting $\, x = {\bf 1} - q\, $ one has that $\, v_{\kappa }^r\, ( {\bf 1} - q )\, v_{\kappa }^r\, $ is increasingly close to a projection 
$\, v_{\kappa }^r - q_{\kappa }\in \mathcal P ( \mathcal J )\, $ with $\, q_{\kappa } \leq v_{\kappa }^r\, $ and $\, v_{\kappa }^r - q_{\kappa } \leq \beta v_{\kappa }^r + \alpha\, v_{\kappa }^r\, r\, $ so that $\, \beta\, s \bigl( v_{\kappa }^r \bigr) + \alpha\, s \bigl( v_{\kappa }^r\, r \bigr) \geq s \bigl( v_{\kappa }^r \bigr) - s \bigl( q_{\kappa } \bigr) - {\delta }_{\kappa }\, $ with $\, \bigl( {\delta }_{\kappa } {\bigr)}_{\kappa } \to 0\, $ and $\, \bigl( s \bigl( q_{\kappa } \bigr) {\bigr)}_{\kappa } \to {\pi }_{\mathcal J} \bigl( {\bf 1} - q \bigr)\, $. From this one gets $\, {\mu }_{s , \mathcal J} \bigl( {\bf 1} - q \bigr) = {\pi }_{\mathcal J} ( s ) \bigl( {\bf 1} - q \bigr)\, $ and since $\, {\pi }_{\mathcal J} ( s )\, $ is $A$-linear by duality $\, {\widehat\mu }_{s , \mathcal J} \bigl( q \bigr) = s \bigl( q \bigr)\, $. The argument for $\, {\pi }_Q ( s )\, $ with respect to $\, C_1 = {\pi }_Q ( s ) \bigl( {\bf 1} \bigr)\, $ is quite similar and we leave it to the reader.
\par\noindent
 A $\mathbb R$-valued $\mathcal P$-measure $\, l\ : \mathcal P \rightarrow \mathbb R\, $ is called {\it $\mathcal P$-linear} iff 
$\, \sum_{i = 1}^n  l \bigl( p_i \bigr) = \sum_{j = 1}^m\, l \bigl( q_j \bigr)\, $ whenever $\, \sum_{i = 1}^n\, p_i = \sum_{j = 1}^m \, q_j\, $ in $\, A\, $, and {\it $A$-linear} iff $\, l \bigl( p + q \bigr) = l \bigl( p \bigr) + l \bigl( q \bigr)\, $ for any pair of orthogonal projections $\, p \bot q\, $. If $\, r\in \mathcal P\, $ then a difference $\, r = \sum_{i = 1}^n\,p_i\, -\ \sum_{j = 1}^m\, q_j\, $ of elements in $\, \mathcal P\, $ is an {\it $A$-linear representation of $\, r\, $} iff there exists a permutation of the summands such that placing certain brackets (and changing the signs according to these brackets) each $+$-sign represents an orthogonal sum of projections and each $-$-sign represents an orthogonal difference of projections. Then one arrives at an equivalence relation $\, \sum_{i = 1}^n\, p_i {\sim}_A \sum_{j = 1}^m\, q_j\, $ whenever there exist a finite subset of projections $\, \bigl\{ s_k {\bigr\}}_{k = 1}^{o}\, $ such that $\, \Bigl( \sum_{i = 1}^n\, p_i + \sum_{k = 1}^{o}\, s_k \Bigr) - \Bigl( \sum_{j = 1}^m\, q_j + \sum_{k = 1}^{o}\, s_k \Bigr)\, $ is an $A$-linear representation of $\, 0\, $. The ($\mathbb Z$)-linear space $\, {\mathcal A \mathcal P}^{\infty , \infty }\, $ generated by $\, \mathcal P\, $ modulo $\, {\sim }_A\, $ is called the $A$-linear envelope of $\, \mathcal P\, $. Moreover define $\, {\mathcal A\mathcal P}^n.\subseteq {\mathcal A\mathcal P}^{\infty , \infty }\, $ to be the subset consisting of elements which can be written as a sum of at most $\, n\, $ projections in $\, \mathcal P\, $ with $\, {\mathcal A\mathcal P}^{n , m} = {\mathcal A\mathcal P}^n - {\mathcal A\mathcal P}^m\, $ and $\, {\mathcal A\mathcal P}^{\infty } = \bigcup_n\, {\mathcal A\mathcal P}^n\, $.
One similarly defines the $\mathcal P$-linear envelope $\, {\mathcal P}{\infty , \infty } \subseteq A\, $ as the $\mathbb Z$-linear subspace generated by $\, \mathcal P\, $. There is no harm in passing to rational coefficients, i.e. considering $\, \mathbb Q \otimes {\mathcal P}^{\infty , \infty } = \mathbb Q\, {\mathcal P}^{\infty , \infty } \subseteq A\, $ but note that $\, \mathbb R \otimes {\mathcal P}^{\infty , \infty }\, $ need not embed into $\, A\, $. In the same way the $A$-linear envelope of $\, \mathcal P\, $ is generally an extension of its image space in $\, A\, $. If $\, l : \mathcal P \rightarrow \mathbb R\, $ extends to a (nonpositive) $\mathcal P$-linear (resp. $A$-linear) realvalued map on the $\mathcal P$-linear (resp. $A$-linear) envelope then $\, l\, $ is called a {\it signed $\mathcal P$-measure} iff $\, l = l_+ - l_-\, $ for a pair of positive realvalued $\mathcal P$-linear (resp. $A$-linear) $\mathcal P$-measures. We will see below that any realvalued $\mathcal P$-linear (resp. $A$-linear) $\mathcal P$-measure is signed.
A pair of  subcones $\, \mathcal C\, ,\, \mathcal D \subseteq {\mathcal P}^{\infty , \infty }\, $ (resp. $\, \mathcal C\, ,\, \mathcal D\subseteq {\mathcal A\mathcal P}^{\infty , \infty }\, $) containing $\, 0\, $ and closed under addition of elements is called {\it exhausting} iff $\, \mathcal C - \mathcal D = {\mathcal P}^{\infty , \infty }\, $ (resp. $\, = {\mathcal A\mathcal P}^{\infty , \infty }\, $). A function $\, m : \mathcal C \rightarrow \mathbb R \cup \{ \infty \}\, $ which is bounded below is {\it $\mathcal P$-convex} iff $\, m \bigl( 0 \bigr) = 0\, $ and $\, m \bigl( x + y \bigr) \leq m \bigl( x \bigr) + m \bigl( y \bigr)\, $. Correspondingly a function $\, \nu : \mathcal D \rightarrow \mathbb R \cup \{ - \infty \}\, $ which is bounded above is {\it $\mathcal P$-concave} iff $\, \nu \bigl( 0 \bigr) = 0\, $ and $\, \nu \bigl( x + y \bigr) \geq \nu \bigl( x \bigr) + \nu \bigl( y \bigr)\, $.
If $\, \mathcal F \subseteq \mathcal P\, $ is an $A$-sublattice define $\, {\mathcal F}^n\, $ (resp. $\, {\mathcal A\mathcal F}^n\, $) to be the subset of elements in $\, {\mathcal P}^n\, $ (resp. in $\, {\mathcal A\mathcal P}^n\, $) which can be written as a difference of sums of elements.in $\, \mathcal F\, $. One similarly defines $\, {\mathcal A\mathcal F}^{\infty }\, ,\, {\mathcal A\mathcal F}^{\infty , \infty }\, $ etc.. 
\par\bigskip\noindent
{\bf Extension theorem for $\mathcal P$-linear (resp. $A$-linear) $\mathcal P$-measures.}
Assume given an $A$-sublattice $\, \mathcal F \subseteq \mathcal P\, $, a (not necessarily positive) realvalued $\mathcal P$-linear (resp. $A$-linear for the relative $A$-linear structure of the $A$-sublattice $\, \mathcal F \subseteq \mathcal P\, $) $\mathcal P$-measure $\, l_0 : \mathcal F \rightarrow \mathbb R\, $ together with an exhausting pair of subcones $\, \mathcal C\, ,\; \mathcal D \subseteq {\mathcal P}^{\infty , \infty }\, $ (resp. $\, {\mathcal A\mathcal P}^{\infty , \infty }\, $), a $\mathcal P$-convex function $\, m : \mathcal C \rightarrow \mathbb R \cup \{ \infty \}\, $ and a $\mathcal P$-concave function $\, \nu : \mathcal D \rightarrow \mathbb R \cup \{ - \infty \}\, $ satisfying 
$$ l_0 \bigl( z \bigr)\> \leq\> m \bigl( x \bigr)\> -\> \nu \bigl( y \bigr) \leqno{(*)} $$
\par\noindent
whenever $\, z = x - y\, $ for $\, x\in \mathcal C\, ,\, y\in \mathcal D\, $ with $\, z\in {\mathcal F}^{\infty , \infty }\, $ (resp. $\, z\in {\mathcal A\mathcal F}^{\infty , \infty }\, $). Then $\, l_0\, $ extends to a $\mathcal P$-linear (resp. $A$-linear) $\mathcal P$-measure $\, l : \mathcal P \rightarrow \mathbb R\, $ satisfying $\, \nu \bigl( y \bigr) \leq l \bigl( y \bigr)\, $ for $\, y\in \mathcal D\, $ and $\, l \bigl( x \bigr) \leq m \bigl( x \bigr)\, $ for $\, x\in \mathcal C\, $.
\par\bigskip\noindent
{\it Proof.}\quad We only write out the proof in case of an $A$-linear $\mathcal P$-measure the argument for $\mathcal P$-linear measures being completely analogous. Consider the $\mathcal P$-convex function $\, \widetilde m : {\mathcal A\mathcal P}^{\infty , \infty } \rightarrow \mathbb R \cup \{ \infty \}\, $ given by 
$$ \widetilde m \bigl( z \bigr)\> =\> \inf\, \Bigl\{\> m \bigl( x \bigr)\> -\> \nu \bigl( y \bigr)\> \Bigm\vert\> z\> =\> x\> -\> y\> ,\quad x\in \mathcal C\> ,\quad y\in \mathcal D\> \Bigr\} $$
(assuming that $\, \widetilde m \bigl( z \bigr) \geq - \infty\, $ for all $\, z\, $) and check that $\, - \widetilde m \bigl( - z \bigr) \leq l_0 \bigl( z \bigr) \leq \widetilde m \bigl( z \bigr)\, $ for $\, z\in {\mathcal A\mathcal F}^{\infty , \infty }\, $ from the $(*)$-condition. Thus we may assume that $\,  \mathcal C = {\mathcal P}^{\infty , \infty }\, ,\, \mathcal D = \{ 0 \}\, $ replacing $\, m\, $ with $\, \widetilde m\, $ and $\, \nu\, $ with $\, 0\, $. The proof is by induction. Choose an element $\, c\in \mathcal P \backslash \mathcal F\, $ and let $\, {\mathcal F}_c \subseteq \mathcal P\, $ denote the $A$-sublattice generated by $\, \mathcal F \cup \{ c \}\, $. If $\, x\in {\mathcal A{\mathcal F}_c}^{\infty , \infty }\, $ there exists $\, n\in\mathbb Z\, $ with $\, x' = x + n\, c\in {\mathcal A\mathcal F}^{\infty , \infty }\, $. If $\, l_c : {\mathcal A{\mathcal F}_c}^{\infty , \infty } \rightarrow \mathbb R\, $ with $\, l_c \bigl( c \bigr) = \alpha\, $ is any extension of $\, l_0\, $ such that $\, l_c \bigl( x \bigr) \leq m \bigl( x \bigr)\, $ for all $\, x\, $ one must have
$$ m\Bigl( l_0 \bigl( x' \bigr)\> -\> m \bigl( x \bigr) \Bigr)\> \leq\> r\, \alpha\>\leq\> n\, \Bigl( m \bigl( y \bigr)\> -\> l_0 \bigl( y' \bigr) \Bigr) $$
for $\, r = n\, m\, $ and $\, n\, ,\, m\in \mathbb N\, $ whenever $\, x'\, ,\, y'\in {\mathcal A\mathcal F}^{\infty , \infty }\, $ with $\, x = x' - n\, c\, ,\, y = y' + m\, c\, $. In particular putting 
$$ A\> =\> \sup\,\left\{\> {l_0 \bigl( x' \bigr)\> -\> m \bigl( x \bigr)\over n} \>\Biggm\vert\> x\> =\> x' - n\, c\> ,\quad x'\in {\mathcal A\mathcal F}^{\infty , \infty }  \>\right\}\> , $$
$$ B\> =\> \inf\, \left\{ \> {m \bigl( y \bigr)\> -\> l_0 \bigl( y' \bigr)\over m} \>\Biggm\vert\> y\> =\> y' + m\, c\> ,\quad y'\in {\mathcal A\mathcal F}^{\infty , \infty } \> \right\} $$ 
one must have $\, A \leq B\, $ which is an easy consequence of $\mathcal P$-convexity of $\, m\, $ together with $A$-linearity of $\, l_0\, $ and the induction assumption. Proceeding step by step using transfinite induction gives the result\qed
\par\bigskip\noindent
The theorem above can be used to show the existence of a Jordan type positive decomposition for any $A$-linear (resp. $\mathcal P$-linear) nonpositive realvalued $\mathcal P$-measure. We only do the $A$-linear case. Assume without loss of generality given a nonpositive $C$-contractive $\mathcal A$-linear $\mathcal P$-measure $\, l : \mathcal P \rightarrow \mathbb R\, $ with $\, l \bigl( {\bf 1} \bigr) = - 1\, $ and putting $\, c := \sup \bigl\{ \vert l \bigl( p \bigr) \vert \bigm\vert p\in \mathcal P \bigr\}\, $ choose a sequence of projections $\, \bigl\{ p_0\, ,\, p_1\, ,\,\cdots \bigr\} \subseteq \mathcal P\, $ with $\, - 1 \geq l \bigl( p_k \bigr) = - c_k \geq  - c\, $ such that $\, \inf_k\, \bigl\{ l \bigl( p_k \bigr) \bigr\} = - c\, $ . Consider the $\mathcal P$-convex function $\, {\bf 1}_{\mu }^{p_0} : {\mathcal A\mathcal P}^{\infty , \infty } \rightarrow \mathbb R\, $ defined by 
$$ {\bf 1}_{\mu }^{p_0} \bigl( x \bigr)\> =\> {\bf 1}_{\mu } \bigl( p_0\, \overline x\, p_0 \bigr) $$
where $\, \overline x\, $ is the image of $\, x\, $ in $\, A\, $ and $\, {\bf 1}_{\mu }\, $ denotes the maximal $\mu $-state on $\, A\, $. Let 
$\, {\mathcal C}_0 = {\mathcal D}_0 = \mathbb N\, ( {\bf 1} - p_0 ) - {\mathcal A\mathcal P}^{\infty }  \subseteq {\mathcal A\mathcal P}^{\infty , \infty }\, $ denote the subcone of elements which are negative sums of projections modulo the $1$-dimensional subspace generated by $\, p_0\, $. Moreover define the $\mathcal P$-convex function 
$\, r_x : {\mathcal A\mathcal P}^{\infty , \infty } \rightarrow \mathbb Z\, $ by 
$$ r_x\> =\> \inf\, \Bigl\{ r\in \mathbb Z \bigm\vert r_x\, {\bf 1} + x \in {\mathcal A\mathcal P}^{\infty } \Bigr\} \> . $$
Define the $\mathcal P$-convex function $\, m_0 : {\mathcal A\mathcal P}^{\infty , \infty } \rightarrow \mathbb R\, $ by 
$$ m_0 \bigl( x \bigr)\> =\> ( c_0 - 1 )\, r_x\> +\> l \bigl( x \bigr) $$
and the $\mathcal P$-concave function $\, {\nu }_0 : {\mathcal A\mathcal P}^{\infty , \infty } \rightarrow \mathbb R\, $ by 
$$ {\nu }_0 \bigl( x \bigr)\> =\> - c_0\,{\bf 1}_{\mu }^{p_0} \wedge l \bigl( x \bigr) $$
so that $\, {\nu }_0 \bigl( x \bigr) \leq m_0 \bigl( x \bigr)\, $ for $\, x\in {\mathcal C}_0\, $ since $\, x\in {\mathcal C}_0\, $ implies $\, r_x \geq 0\, $.
From the extension theorem there exists an $A$-linear extension $\,  l_0\, $ of $\,  l_0  \bigl( {\bf 1} \bigr) = - c_0 = l_ 0 \bigl( p_0 \bigr)\, $.
satisfying the $(*)$-condition
$$ l_0 \bigl( k\, {\bf 1} + l\, p_0 \bigr)\> \leq\> m_0 \bigl( x \bigr) \> -\> {\nu }_0 \bigl( y \bigr) $$
whenever $\, k\, ,\, l\in \mathbb Z\, $ with $\, x - y = k\, {\bf 1} + l\, p_0\, $. 
To see this note that  for $\, v = k\, {\bf 1} + l\, p_0\, $ one has $\, r_v\> =\> - ( k + l )\, $ if $\, l \leq 0\, $ and $\, r_v = - k\, $ for $\, l \geq 0\, $
so that  $\, m_0 \bigl( k\, {\bf 1} + l\, p_0 \bigr) = - ( k + l )\, c_0\, $ if $\, l \geq 0\, $ and $\, m_0 \bigl( k\, {\bf 1} + l\, p_0 \bigr) = - ( k + l )\, c_0 - ( c_0 - 1 )\, l \geq - ( k + l )\, c_0\, $ for $\, l \leq 0\, $ whereas $\, {\nu }_0 \bigl( k\, {\bf 1} + l\, p_0 \bigr) = - ( k + l )\, c_0\, $ for $\, k \geq 0\, $ and $\, {\nu }_0 \bigl( k\, {\bf 1} + l\, p_0 \bigr) = - k - l\, c_0 \leq - ( k + l )\, c_0\, $ for $\, k \leq 0\, $. Consider the $\mathcal P$-convex function $\, {\underline m}_0 : {\mathcal C}_0 \rightarrow \mathbb R\, $ given by 
$$ {\underline m}_0 \bigl( x \bigr)\> =\> \inf\, \Bigl\{ m_0 \bigl( y \bigr)\> -\> {\nu }_0 \bigl( k\, {\bf 1} + l\, p_0 \bigr) \Bigm\vert k\geq 0\> ,\> x = y - \bigl( k\, {\bf 1} + l\, p_0 \bigr)\> ,\> y\in {\mathcal C}_0 \Bigr\} $$  
so that $\, {\nu }_0 \bigl( x \bigr) \leq {\underline m}_0 \bigl( x \bigr) \leq m_0 \bigl( x \bigr)\, $ with $\, {\underline m}_0 \bigl( k\, {\bf 1} + l\, p_0 \bigr) = - ( k + l )\, c_0\, $ if $\, k\, {\bf 1} + l\, p_0 \in {\mathcal C}_0 \iff k \leq l\, $. Next consider the $\mathcal P$-concave function  $\, {\overline\nu }_0 : {\mathcal C}_0 \rightarrow \mathbb R\, $ given by 
$$ {\overline\nu }_0 \bigl( x \bigr)\> =\> \sup\, \Bigl\{ {\nu }_0 \bigl( y \bigr)\> -\> {\underline m}_0 \bigl( k\, {\bf 1} + l\, p_0 \bigr) \Bigm\vert k\>\leq\> l\> ,\> x\> =\> y\> -\> \bigl( k\, {\bf 1} + l\, p_0 \bigr)\> ,\> y\in {\mathcal C}_0  \Bigr\} $$
so that $\, {\nu }_0 \bigl( x \bigr) \leq {\overline\nu }_0 \bigl( x \bigr) \leq {\underline m}_0 \bigl( x \bigr) \leq m_0 \bigl( x \bigr)\, $ with $\, {\overline\nu }_0 \bigl( k\, {\bf 1} + l\, p_0 \bigr) = - ( k + l )\, c_0\, $ if $\, k\, {\bf 1} + l\, p_0 \in {\mathcal C}_0\, $. Then any element $\, z\, $ in $\, \mathbb Z\, {\bf 1} + \mathbb Z\, p_0\, $ can be written as a difference $\, z = v - w\, $ with $\, v\, ,\, w \in {\mathcal C}_0 \cap \Bigl( \mathbb Z\, {\bf 1} + \mathbb Z\, p_0 \Bigr)\, $ and the $(*)$-condition follows from $\, l_0 \bigl( v \bigr) = {\overline\nu }_0 \bigl( v \bigr)\, ,\, l_0 \bigl( w \bigr) = {\underline m}_0 \bigl( w \bigr)\, $ plus convexity of $\, {\underline m}_0\, $ and concavity of $\, {\overline\nu }_0\, $.  We calim that $\, - c \leq l_0 \bigl( p \bigr) \leq l \bigl( p \bigr)\, $ for each $\, p\in \mathcal P\, $. For this note that since $\, p - {\bf 1} \in {\mathcal C}_0\, $ with $\, m_0 \bigl( p - {\bf 1} \bigr) \leq c + l \bigl( p \bigr)\, $ one has $\, l_0 \bigl( p \bigr) \leq l_0 \bigl( {\bf 1} \bigr) + m \bigl( p - {\bf 1} \bigr) = l \bigl( p \bigr)\, $. On the other hand one has $\, p - p_0\in {\mathcal C}_0\, $ with $\, {\nu }_0 \bigl( p - p_0 \bigr) \geq 0 \wedge l \bigl( p - p_0 \bigr) \geq - c + c_0\, $ so that $\, l_0 \bigl( p \bigr) \geq - c\, $ follows. One also has that $\, l_0\, $ is $C_0$-contractive for $\, C_0 = \max\, \bigl\{ C\, ,\, c_0 \bigr\}\, $. Indeed assume that $\, p\, ,\, q\in \mathcal P\, $ with $\, \Vert p - q \Vert \leq \epsilon\, $. Then $\, p - q - p_0\in {\mathcal C}_0\, $ with 
$\, {\nu }_0 \bigl( p - q - p_0 \bigr) \geq \min\, \bigl\{ c_0\, ( 1 - \epsilon )\, ,\, c_0 - C\, \epsilon \bigr\}\, $ so that $\, l_0 \bigl( p \bigr) - l_0 \bigl( q \bigr) \geq - \max\, \bigl\{ c_0 \, ,\, C \bigr\}\, \epsilon\, $ and by symmetry also $\, l_0 \bigl( q \bigr) - l_0 \bigl( p \bigr) \geq - \max\, \bigl\{ c_0\, ,\, C \bigr\}\, \epsilon\, $ giving the result.  One now proceeds by induction sucessively constructing a sequence of $C_n$-contractive $A$-linear $\mathcal P$-measures $\, \bigl\{ l_n  {\bigr\}}_n\, $ such that $\,  - c\leq l_{n + 1} \bigl( p \bigr) \leq l_n \bigl( p \bigr)\, $ for all $\, p\in \mathcal P\, $ with $\, C_{n + 1} = \max\, \bigl\{ c_n\, ,\, C_n \bigr\}\, $ and $\, l_n \bigl( {\bf 1} \bigr) = - c_n\, $. Then putting $\, - l_-  \bigl( p \bigr) = \inf_n\, \bigl\{ l_n \bigl( p \bigr) \bigr\}\, $ defines a negative $A$-linear $\mathcal P$-measure dominated by $\, l\, $ which is $D$-contractive for $\, D = \max \bigl\{ c\, ,\, C \bigr\}\, $ and satisfies $\, - l_- \bigl( {\bf 1} \bigr) = - c\, $ whence putting $\, l_+ = l + l_- \geq 0\, $ one arrives at the desired Jordan type positive decomposition of $\, l\, $.
\par\smallskip\noindent
If $\, \bigl\{ {\omega }_k {\bigr\}}_k\subseteq \mathcal H\, $ is an orthonormal basis of the Hilbert space $\, \mathcal H\, $ consider the linear (nonpositive)  involution $\, o\, $ on $\, \mathcal B ( \mathcal H )\, $ such that $\, x_{ii }^{o} = x_{ii}\, $ is the identity on the diagonal with respect to the chosen basis and changes signs on the off-diagonal, i.e.  and $\, x_{ij}^{o} = - x_{ij}\, $ if $\, i\neq j\, $.  We continue to denote the dual involution on $\, \mathcal B ( \mathcal H )^*\, $ by $\, o\, $.  An involution of this form will be called {\it basic}. 
\par\bigskip\noindent
{\bf Theorem 3.}\quad If $\, Q \simeq A / \mathcal J\, $ is a  $C^*$-quotient of a ${\Sigma }_{\omega }$-algebra $\, A\, $ modulo a twosided ideal $\, \mathcal J\vartriangleleft A\, $ and $\, \{ x_n {\}}_{n\in\mathbb N}\subseteq Q\, $ is a bounded monotone increasing sequence there exists for any countable collection $\, \{ y_l {\}}_{l\in \mathbb N} \subseteq Q\, $ with $\, y_l \geq x_n\, $ for all 
$\, n\, , l\, $ an element $\, z\in Q\, $ satisfying 
$$ x_n\>\leq\> x_{n + 1}\> \leq\> \cdots\> \leq \> z\>\leq\>  y_l\quad \forall\> l\in \mathbb N\> . $$
In particular, if $\, Q\, $ admits an order isomorphic representation on a separable Hilbert space then $\, Q\, $ is sequentially monotone complete.
\par\bigskip\noindent
{\it Proof.}\quad Let $\, A\, $ be a ${\Sigma }_{\omega }$-algebra and $\, Q \simeq A / \mathcal J\, $ a $C^*$-quotient of $\, A\, $.  Then $\, Q\, $ is a $\mathcal P$-algebra. Let a monotone increasing sequence $\, 0\leq\cdots\leq x_k\leq x_{k + 1}\leq\cdots\leq {\bf 1}\, $ of positive elements in $\, Q_1\, $ be given and $\, \{ y_l {\}}_{l\in\mathbb N}\in Q^+\, $ be any countable collection of elements with every $\, y_l\, $ larger than each of the $\, \{ x_k \}\, $. We first show that there exists an element $\, y\in Q^+\, $ with
$$ 0\>\leq\>\cdots\>\leq\> x_k\>\leq\> x_{k + 1}\>\leq\>\cdots\>\leq\> y\>\leq\> \bigl\{ y_l \bigm\vert l\in\mathbb N \bigr\} \> . $$
Let $\, \epsilon > 0\, $ be given and $\, \{ {\overline y}_1 \}\in A^+\, $ be a lift of $\, y_1\, $. From the Remark after Theorem P there exists a lift $\, {\overline x}_1^{\epsilon }\in A^+\, $ of $\, x_1\, $ satisfying $\, 0 \leq {\overline x}_1^{\epsilon } \leq  {\overline y}_1 + \epsilon\, {\bf 1}\, $. Choose a lift $\, {\overline y}_2\in A^+\, $ of $\, y_2\, $ larger than $\, {\overline x}_1^{\epsilon }\, $. Then again since $\, {\overline y}_1 + \epsilon\, {\bf 1}\, ,\, {\overline y}_2 + \epsilon\, {\bf 1}\, $ are invertible there exists a lift $\, {\overline x}_2^{\epsilon }\in A^+\, $ of $\, x_2\, $ satisfying 
$$ {\overline x}_1^{\epsilon }\>\leq\> {\overline x}_2^{\epsilon }\>\leq\> \left\{ {\overline y}_1 + \epsilon\, {\bf 1}\, ,\, {\overline y}_2 + \epsilon\, {\bf 1} \right\}\> . $$
Proceeding inductively one finds a sequence of lifts $\, \{ {\overline y}_l \}\subseteq A^+\, $ of the $\, \{ y_l \}\, $ and a sequence of lifts $\, \{ {\overline x}_k^{\epsilon } \}\subseteq A^+\, $ of the $\, \{ x_k \}\, $ satisfying
$$ 0\>\leq\>\cdots\>\leq\> {\overline x}_k^{\epsilon }\>\leq\> {\overline x}_{k + 1}^{\epsilon }\>\leq\>\cdots\>\leq \left\{ {\overline y}_l + \epsilon\, {\bf 1} \right\} $$
For $\, z^{\epsilon } = \sup_k\, {\overline x}_k^{\epsilon }\, $ and $\, y^{\epsilon } = q ( z^{\epsilon } )\, $ one gets 
$$ 0\>\leq\>\cdots\>\leq\> x_k\>\leq\> x_{k + 1}\>\leq\>\cdots\>\leq\> y^{\epsilon }\>\leq\> \left\{ y_l + \epsilon\, {\bf 1} \right\} \> . $$
Considering the sequence $\, \{ y^{1 / k} {\}}_{k\to\infty }\, $ choose a lift $\, {\overline x}_{k + 1}\in A^+\, $ of $\, x_{k + 1}\, $ satisfying 
$$ 0\>\leq\> {\overline x}_1\>\leq\> \cdots\> \leq\> {\overline x}_k\>\leq\> {\overline x}_{k + 1}\> \leq\> 
\left\{ {\overline y}^{1 / l}\, +\, {1\over l}\, {\bf 1}\bigm\vert l = 1\, ,\cdots\, ,\, k \right\} $$  
where $\, {\overline y}^{1 / l}\in A^+\, $ is a lift of $\, y^{1 / l}\, $ assuming by induction that a compatible lift $\, \{ {\overline x}_l\, ,\, {\overline y}^{1 / l} {\}}_{l = 1}^k\, $ with 
$$ 0\>\leq\>{\overline x}_1\>\leq\> \cdots\>\leq\> {\overline x}_k\>\leq\> \left\{ {\overline y}^{1 / l} + {1\over l}\, {\bf 1}\bigm\vert l = 1\, ,\,\cdots\, ,\, k \right\} $$
has already been chosen. Then choose a lift $\, {\overline y}^{1 / ( k + 1 )}\in A^+\, $ of $\, y^{1 / ( k + 1)}\, $ satisfying 
$$ {\overline x}_{k + 1}\>\leq {\overline y}^{1 / ( k + 1 )} $$
completing the induction step $\, k \to k + 1\, $. Putting $\, z = \sup_k\, {\overline x}_k\, $ and $\, y = q ( z )\, $ it is easy to see that 
$$ 0\>\leq\>\cdots\>\leq\> x_k\>\leq\> x_{k + 1}\>\leq\cdots\>\leq\> y\>\leq \bigl\{ y_l \bigr\} $$
as claimed. To prove the second assertion choose any maximal well ordered monotone decreasing net $\, \bigl\{ y_{\gamma }\bigm\vert y_{\gamma + 1} < y_{\gamma } {\bigr\}}_{\gamma\in\Gamma }\, $ with each element $\, y_{\gamma }\, $ exceeding each element $\, x_k\, $ which from the argument above either posseses a minimal element or else is uncountable. Then the assertion follows from the fact that any well ordered bounded monotone decreasing net in the bounded operators of a separable Hilbert space is at most countable\qed
\par

\end{document}